\documentclass[11pt,a4paper]{report}

\usepackage{amsmath, amssymb, amsthm}
\usepackage{graphicx,color}
\usepackage[all]{xy}
\usepackage[left=1.5in, right=1in, top=1in, bottom=1.4in,
headheight=13.6pt]{geometry} \usepackage[square, comma, numbers,
sort&compress]{natbib} \usepackage[T1]{fontenc}

\usepackage{times}
\usepackage[hang, small, bf]{caption}
\numberwithin{equation}{section}

\usepackage{fancyhdr}
\makeatletter
\renewcommand{\chaptermark}[1]{\markboth{\textsc{\@chapapp}\ \thechapter:\ #1}{}}
\makeatother

\usepackage{sectsty}
\chapterfont{\large\sc\centering}
\chaptertitlefont{\centering}
\subsubsectionfont{\centering}

\newtheorem{theorem}{Theorem}[section] \newtheorem{lemma}[theorem]{Lemma}
\newtheorem{prop}[theorem]{Proposition}
\newtheorem{definition}[theorem]{Definition}
\newtheorem{corollary}[theorem]{Corollary}
\newtheorem{conjecture}[theorem]{Conjecture}

\newtheorem*{thm1}{Theorem \ref{zetafunc}}
\newtheorem*{thm2}{Corollary \ref{qcontz}}
\newtheorem*{thm3}{Corollary \ref{qcontz1}}
\newtheorem*{thm4}{Corollary \ref{yo}}
\newtheorem*{thm5}{Theorem \ref{thint1}}
\newtheorem*{thm6}{Theorem \ref{stronger}}
\newtheorem*{thm7}{Theorem \ref{klzetash21}}
\newtheorem*{thm8}{Theorem \ref{bc1}}

\begin{document}

    \title{
    \huge{\textbf{Nonstandard Mathematics and New Zeta and
L-Functions}}\\[1.2cm]
       \Large{ B.\ Clare} \\[1.2cm]
        \Large{Thesis submitted to The University of Nottingham
\\ for the degree of Doctor of Philosophy} \\ \vspace{1cm}
 \Large{September 2007}}
 \author{} \date{}
 \maketitle

  \newpage
  \chapter*{Abstract}
  \textsc{I define new zeta functions in a nonstandard setting
and examine some of their properties.  I further develop $p$-adic
interpolation in the nonstandard setting and define the concept
of interpolation with respect to two primes.  The final section
of the dissertation examines the work of M.\ J.\ Shai Haran and
makes initial attempts of viewing it from a nonstandard
perspective.}


  \tableofcontents
 \newpage
 \pagenumbering{arabic}

              \chapter{Introduction}       \label{chint}

There are many examples of modern work in number theory being
developed by fusing it with another field of mathematics.
Examples include quantum field theory, with work
by D.\ Broadhurst and D.\ Kreimer, and the more well known field
of geometry, with the work by A.\ Connes (\cite{Co1} and
\cite{Co2}). This enables studying old problems from a new
direction. The other alternatives in tackling number theory
include proving a known theorem in a new way and trying to prove
an old problem directly (which is becoming less popular due to
the difficulty of such problems - Riemann hypothesis, Birch
Swinnerton-Dyer conjecture, $\ldots$).

The approach in this work is via number-theory-fusion with the
area of nonstandard mathematics. The
work splits into four related parts with the underlying theme
being nonstandard mathematics. This chapter and chapter \ref{chrobo} (first
part) are introductory; setting out the work and introducing nonstandard
mathematics from a model theoretic basis. A review of the main results of
Robinson on nonstandard algebraic number theory is provided as the work in
following chapters further develops his work and approaches. Chapters
\ref{chhypr} and \ref{chhypd} (second part) introduce nonstandard versions of
the Riemann and Dedekind zeta functions. They also prove some nonstandard
analytical properties. The third part consists of chapters \ref{chadic} and
\ref{chapp}. In the first of these chapters the ideas of $p$-adic Mahler
interpolation are interpreted in a nonstandard setting using
shadow (standard part) maps (which enables nonstandard objects to
be viewed in a standard setting). These ideas are pursued in the
next chapter by viewing the Morita gamma function and
Kubota-Leopoldt zeta function in a nonstandard way. One
advantage of using nonstandard spaces ($\mathbb{^{\ast}Q}$ for
example) is the property of being able to consider finite primes
in a symmetrical way (for example subquotients of
$\mathbb{^{\ast}Q}$ include $\mathbb{Q}_{p_{1}}$ and
$\mathbb{Q}_{p_{2}}$ resulting in functions on this nonstandard
space being able to be interpreted as as $p_{1}$-adic and
$p_{2}$-adic functions). This enables double interpolation with
respect to two distinct finite rational primes to be defined,
and more generally interpolation with respect to a finite set of
rational primes. An example given is the double interpolation of
the Riemann zeta function. The final chapter (the final part) is
on the work of Shai Haran. A review of some of his work is given
in conjunction with a review of other authors' work which
relates to his. A small section is dedicated to the beginning of
a nonstandard interpretation of his work.

\section{Chapter \ref{chhypr} Overview}

The first application of nonstandard analysis in this work is in
relation to the Riemann zeta function. A nonstandard version of
the Riemann zeta function, the hyper Riemann
zeta function, is defined. In
particular some of the first analytical properties are examined.

The framework developed from the model theory provides a
foundation for most of the work. Further (analytical) tools are
needed and in particular the shadow maps. They exist in a general setting on hyper
topological spaces with a shadow map corresponding to a non trivial
absolute value. 
Essentially they link (part of) an object on a nonstandard space
to a related object on the standard space. On
$\mathbb{^{\ast}R}$ the shadow map of interest is the real
shadow map corresponding to the real prime $\eta$ and the usual
absolute value (definition \ref{realsh}). This map can also be
extended to a shadow map on $\mathbb{^{\ast}C}$ (definition
\ref{complexsh}). Further it is shown how the shadow map can act
on functions (definition \ref{shimage}) which is of most use in
this chapter. Other tools introduced include the hyperfinite
version of sums, products and integrals along with hyper
sequences.

The central object in this chapter is the new function, the \textit{hyper
Riemann zeta function}, $\zeta_{\mathbb{^{\ast}Q}}$. The definition follows from
defining the internal ideals of $\mathbb{^{\ast}Z}$ and their corresponding
norm. These internal ideals are in bijection with the elements of
$\mathbb{^{\ast}N}$ (lemma \ref{L1}). Also in order to define
$\zeta_{\mathbb{^{\ast}Q}}$ the definitions of
$\text{$^{\ast}$}\exp$ (the \textit{hyper exponential function})
and $\text{$^{\ast}$}\log$ (the \textit{hyper logarithm
function}) are needed in order to define the power function
$z^{s}$ for general $z,s\in\mathbb{^{\ast}C}$. The approach
taken with these hyper functions is different to those in the
literature, for example \cite{S-L}, section 2.4. In this work
most objects are defined from first principles for completeness.

Combining the ideas from the previous paragraph gives definition
\ref{hyzeta} \[
\zeta_{\mathbb{^{\ast}Q}}(s)=\sum_{n\in\mathbb{^{\ast}N}\setminus\{0\}}
\frac{1}{n^{s}}.\] The first analytical properties are proved in
the $Q$-topology which is the
finer topology. The rest of section 2.3 develops these results
along with the first parts of section 2.4. Examples include the
region of $Q$-convergence and continuity in the hyper half plane
$^{\ast}\Re{s}>1$. Like the classical function there is also a
product decomposition into "Euler factors" via proposition
\ref{L8} \[ \zeta_{\mathbb{^{\ast}Q}}(s) = \prod_{p}
(1-p^{-s})^{-1}.\] The product above is taken over all
positive primes in $\mathbb{^{\ast}Z}$. $\mathbb{^{\ast}Z}$
is a ring so primes (or
prime ideals) can be defined. Some of the $Q$-analytical
results require some hyper complex analysis which is taken from
the work of Robinson, for example the hyper Cauchy theorem.

The rest of the chapter is spent proving the main result in this
chapter (theorem \ref{zetafunc} and corollaries \ref{qcontz}
and \ref{qcontz1}), the \textit{functional equation} and
\textit{$Q$-analytic continuation} of
$\zeta_{\mathbb{^{\ast}Q}}$. \begin{thm1}
$\zeta_{^{\ast}\mathbb{A}}(s)=\zeta_{^{\ast}\mathbb{A}}(1-s)$.
\end{thm1} Here $\zeta_{\mathbb{^{\ast}A}}$ is the completed
hyper Riemann zeta function, that is the Euler product of all
the local zeta functions including the contribution from the
prime at infinity $\pi^{-(s/2)}\text{} ^{\ast}\Gamma (s/2)$.
\begin{thm2}    $\zeta _{_{^{\ast}\mathbb{A}}}(s)$ can be
Q-analytically continued on $\mathbb{^{\ast}C}$. \end{thm2}
\begin{thm3}  $\zeta _{_{^{\ast}\mathbb{Q}}}(s)$ can be
Q-analytically continued on $\mathbb{^{\ast}C}$, \[ \zeta
_{_{^{\ast}\mathbb{Q}}}(1-s)=\frac{\pi^{1/2-s} \text{}
^{\ast}\Gamma(\frac{s}{2})\zeta
_{_{^{\ast}\mathbb{Q}}}(s)}{^{\ast}\Gamma(\frac{1-s}{2})}, \]
with a pole at s=1 and trivial zeros at
$s=-2\mathbb{^{\ast}Z_{>\text{0}}} $ . \end{thm3}

Contained in this expression and in the proof are nonstandard versions of
classical functions, many of which introduced are new objects.

The most basic is the \textit{hyper gamma function} which is a
(hyperfinite) product for the hyper naturals (the \textit{hyper factorial}), \[
^{\ast}\Gamma(n)=\prod_{1\leq j \leq n} j := n!,\] for
$n\in\mathbb{^{\ast}N}$ with $^{\ast}\Gamma(0)=0! =1$ (definition
\ref{hypgam}). The definition of the hyper gamma function for
general hyper complex numbers follows by interpolation of the
hyper factorial (definition \ref{hygamf}) \[
^{\ast}\Gamma(s)=\int_{\mathbb{^{\ast}R}}^{\ast}\text{}
^{\ast}\exp(-y)y^{s-1}dy.\] The \textit{functional equation} for
this function is given by theorem \ref{L21} \[ ^{\ast}\Gamma(s+1)
=s\text{} ^{\ast}\Gamma(s).\]

The method of proof of the functional equation follows the
classical version involving theta functions. A \textit{hyper
theta function} is defined (via definition \ref{hytheta}) \[
^{\ast}\Theta(s)=\sum_{n\in\mathbb{^{\ast}Z}}\text{} ^{\ast}
\exp(\pi in^{2}s),\] for $^{\ast}\Im (s)>0$. In order to obtain
this, hyper Fourier transforms and hyper Poisson summation are
needed. Results in this area require evaluation of definite hyper
integrals and the interchanging of hyper products and hyper
integrals. As a result the functional equation is given by \[
^{\ast}\Theta(-1/s)=(s/i)^{1/2}\text{}^{\ast}\Theta(s).\] The
functional equation of $\zeta_{\mathbb{^{\ast}Q}}$ follows from
this as $\zeta_{\mathbb{^{\ast}Q}}$ can be expressed as a hyper
integral involving $^{\ast}\Theta$.

The main methods of proof are a combination of some analysis and
a transfer of results from the classical setting for
$\zeta_{\mathbb{Q}}$. It is emphasized that results do not
transfer directly and so some analytical methods have to be used
to bind the proofs together. A common example of proof is when
one has a function which is a sum over $\mathbb{^{\ast}N}$ and
has a classical counterpart which is a sum over $\mathbb{N}$.
Usually only properties of "an approximation" to the hyper
function can be deduced from the classical counterpart which can
then be extended to the hyper function using analysis. From the
basic definitions of the hyper functions it is seen that they
have a similar definition to their classical counterparts and it
is unsurprising that proofs of nonstandard results follow a
similar structure to the classical case.

\section{Chapter \ref{chhypd} Overview}

This chapter can be considered as a further application of the
nonstandard and analytical tools from the previous chapter. Just
as the Dedekind zeta function generalises the Riemann zeta
function so the \textit{hyper Dedekind zeta function} generalises the
hyper Riemann zeta function by being a zeta function for hyper
number fields. The obvious question is why the first chapter is
included since it can be deduced from the work of this chapter.
The main reason is that the first chapter serves to introduce
some nonstandard tools not covered in the
introduction to model theory and techniques of proof used in
later chapters. Also some of the nonstandard versions of
classical functions are used in this chapter (for example the
hyper gamma function) and in other chapters. The above reasons
could have yielded a chapter focussing on nonstandard
tools and some classical functions in a nonstandard setting. With
 some extra work the hyper Riemann zeta function and some of its
properties have been proven.

In order to define the hyper Dedekind zeta function a small amount of
nonstandard analytic number theory is needed to give a solid meaning to hyper
(algebraic) number fields and their properties. An important characteristic of
any hyper number field is the ring of hyperintegers, section 3.2.3. This
naturally leads into ideals which have been studied to an extent by Robinson
(detailed in chapter \ref{chrobo}). Finally combining all the above the hyper
Dedekind zeta function for a given hyper number field $^{\ast}K$ can be defined
(definition \ref{hypdede}) \[ \zeta_{^{\ast}K}(s)=\sum_{\textbf{a}}\frac{1}{(
N(\textbf{a}))^{s}},\] where the sum runs over the set of internal ideals of
the ring of hyperintegers of $^{\ast}K$ and $N$ is the norm of such an ideal.
The first properties in the $Q$-topology are given including convergence and a
product formula.

The second half of the chapter is devoted to proving the main
result of the chapter; the \textit{functional equation}
\textit{$Q$-analytic continuation} of $\zeta_{^{\ast} K}$ to
$\mathbb{^{\ast}C}$ (corollary \ref{yo}). \begin{thm4}
$\zeta_{^{\ast}K}(s)$ has a Q-analytic continuation to
$\mathbb{^{\ast}C}\setminus\{1\}$. It has a simple pole at $s=1$ with
residue \[ \frac{2^{r_{1}}(2\pi)^{r_{2}}}{w|d_{^{\ast}K}|^{1/2}}
h_{^{\ast}K}\text{}^{\ast}R.\] It also satisfies the functional
equation \[\zeta_{^{\ast}K}(1-s)=\text{}^{\ast}A(s)
\zeta_{^{\ast}K}(s).\] Here \[^{\ast}A(s)= 2^{n}\times
(2\pi)^{-ns} |d_{^{\ast}K}|^{1/2 -s}(^{\ast}\cos(\pi
s/2))^{r_{1}+r_{2}}(^{\ast}\sin(\pi s/2))^{r_{2}}
(^{\ast}\Gamma(s))^{n}.\] \end{thm4} The proof follows a similar
approach to that of the chapter \ref{chhypr} but in the higher
dimensional case, for example introducing a \textit{higher
dimensional hyper gamma function} and a \textit{higher
dimensional hyper theta function} (definition \ref{hight}) which
can be defined for each complete lattice ($L$) in the
\textit{hyper Minkowski space} \[ ^{\ast}\Theta_{L}(z) =
\sum_{g\in L}\text{} ^{\ast}\exp(\pi i \langle gz,z\rangle ),\]
where $\langle, \rangle$ is a hermitian scalar defined in
section 3.2.

The functional equation (theorem \ref{highttr} follows from a
Poisson summation type formula, \[^{\ast}\Theta_{L}(-1/z)=
\frac{\sqrt(z/i)}{\operatorname{vol}(L)}\text{}
^{\ast}\Theta_{L'}(z),\] where $L'$ is the dual lattice defined
in theorem \ref{latdual}.

To finish off the proof a relationship between $\zeta_{^{\ast}K}$ and $^{\ast}
\Theta_{L}$ (for some lattice) is needed in order to use the functional
equation of the hyper theta function. This is provided essentially by a hyper
Mellin transform (section 3.1).

These first two chapters do follow, and rely upon, the classical theory
for the reasons given in the previous section. The reason for
including them are that several nonstandard objects are defined
rigorously with the hope they may have uses outside of the work in later
chapters (which include in chapter 5 the hyper Riemann zeta function and hyper
gamma function being $p$-adically interpolated in a nonstandard setting. It
should also be noted that for a concise and non-repetitive presentation some of
the arguments are not detailed in depth, for example the evaluation of definite
hyper integrals. This is because in the first chapter rigorous proofs have been
given and the method is the same each time, though clear referencing is given
to enable the result to be proven in depth if so required.

\section{Chapter \ref{chadic} Overview}
The concept of $p$-adic interpolation is wide in its application and in the
actual method of obtaining a $p$-adic object from a real object. This chapter
takes the concept of $p$-adic interpolation by the method of Mahler and
interprets it in a nonstandard way via $p$-adic shadow maps. There are two main
theorems in this chapter \ref{thint1} and \ref{stronger} with the
latter a strengthening of the former.

\begin{thm5}Let $f:\mathbb{N}\rightarrow
\mathbb{Q}_{p}$ be a uniformly continuous function, with respect
to the $p$-adic metric, on $\mathbb{N}$ and let
$^{\ast}f:\mathbb{^{\ast}N}\rightarrow \mathbb{^{\ast}Q}_{p}$ be
the extension to its hyper function. Then
$\operatorname{sh}_{p}(^{\ast}f):\mathbb{Z}_{p}\rightarrow\mathbb{Q}_{p}$
is the the unique $p$-adic function obtained by Mahler
interpolation. \end{thm5}

\noindent Here $\operatorname{sh}_{p}$ is defined in definition
\ref{shimage}, as it is not defined on all of
$^{\ast}\mathbb{Q}_{p}$ since it contains unlimited elements.

\begin{thm6}Let $f:\mathbb{N}\rightarrow
\mathbb{Q}_{p}$ be a uniformly continuous function, with respect
to the $p$-adic metric, on $\mathbb{N}$. Then there exists a hyper function
$^{\ast}g:\mathbb{^{\ast}N}\rightarrow \mathbb{^{\ast}Q}^{\lim_{p}}$ such that
$\operatorname{sh}_{p}(^{\ast}g):\mathbb{Z}_{p}\rightarrow\mathbb{Q}_{p}$ is
the the unique $p$-adic function obtained by Mahler interpolation.\end{thm6}

Mahler interpolation (theorems \ref{ma1} and \ref{ma2}) can be stated in two
stages. Firstly any uniformly continuous function
$f:\mathbb{N}\rightarrow\mathbb{Q}_{p}$ can be extended to
a uniformly and continuous function $F:\mathbb{Z}_{p} \rightarrow
\mathbb{Q}_{p}$ such that if $x\in\mathbb{N}$ then $F(x)=f(x)$.
The second part shows that any continuous function
$h:\mathbb{Z}_{p}\rightarrow\mathbb{Q}_{p}$ can be written as a
series over $\mathbb{N}$. The first theorem (\ref{thint1}) shows
how this interpolation can be viewed in a nonstandard way.
Indeed, any function $f:\mathbb{N} \rightarrow \mathbb{Q}_{p}$
can be extended to a hyperfunction $^{\ast}f:
\mathbb{^{\ast}N}\rightarrow\mathbb{^{\ast}Q}_{p}$ (using the
tools of the introduction by regarding
$^{\ast}f(^{\ast}n)=(f(n_{1}),f(n_{2}),\ldots )$
($^{\ast}n=(n_{1},n_{2},\ldots)$) modulo the equivalence
relation) with $^{\ast}f(n)=f(n)$ for $n\in\mathbb{N}$. Applying
the $p$-adic shadow map to this function does not give the
original function $f$ because $\operatorname{sh}_{p}(
\mathbb{^{\ast}N})\neq \mathbb{N}$. As shown in chapter 2 this
result is stronger as $\operatorname{sh}_{p}(\mathbb{^{\ast}N})
= \mathbb{Z}_{p}$. So actually a function
$g:\mathbb{Z}_{p}\rightarrow \mathbb{Q}_{p}$ with $g(n)=f(n)$
for $n\in\mathbb{N}$ is obtained. The main work in the theorem
is showing that the resulting function $g$ from the $p$-adic
shadow map is indeed the function which would have been obtained
by Mahler interpolation of $f$.

The proof shows that the hyper function ($^{\ast}f$) can be written as a sum
over $\mathbb{^{\ast}N}$ which is basically the sum over $\mathbb{N}$ for $f$
with added nonstandard terms (proposition \ref{pro1}). The properties of
$^{\ast}f$ can now be partially deduced from the properties of $f$ and some more
analysis, in particular showing the continuity and convergence properties. Some
care is needed when taking the shadow map. Although the properties of
convergence and continuity carry through actually taking the shadow map on a sum
over $\mathbb{^{\ast}N}$ is not trivial unlike taking it over a hyperfinite
set.

The second main theorem uses the fact that $\operatorname{sh}_{p}(\mathbb{
^{\ast}Q}^{\lim_{p}})=\mathbb{Q}_{p}$. So instead of considering Q-continuous
hyper functions with values lying in $\mathbb{^{\ast}Q}_{p}$ consider the subset
of it $\mathbb{^{\ast}Q}^{\lim_{p}}$ and in particular uniformly Q-continuous
hyper functions of the form $g:\mathbb{^{\ast}N}\rightarrow
\mathbb{^{\ast}Q}^{\lim_{p}}\in \mathcal{C}_{p}(\mathbb{^{\ast}N},
\mathbb{^{\ast}Q}^{\lim_{p}})$, which is the set of all $p$-adically
uniformly Q-continuous functions from $\mathbb{^{\ast}N}$ to
$\mathbb{^{\ast}Q}^{\lim}_{p}$. This theorem shows that the shadow map is a
surjective homomorphism from $\mathcal{C}_{p}(\mathbb{^{\ast}N},
\mathbb{^{\ast}Q}^{\lim_{p}})$ to $\mathcal{C}_{p}(\mathbb{Z}_{p},
\mathbb{Q}_{p})$, the space of $p$-adically uniformly continuous functions from
$\mathbb{Z}_{p}$ to $\mathbb{Q}_{p}$. Showing the map is a homomorphism and that
the map lies in $\mathcal{C}_{p}(\mathbb{Z}_{p},
\mathbb{Q}_{p})$ follows from the definition of the $p$-adic shadow map. The
hardest part is to show that the map is surjective. To show this a hyper
sequence of hyper functions ($^{\ast}f_{m}:\mathbb{^{\ast}N}\rightarrow
\mathbb{^{\ast}Q}^{\lim_{p}}$ ($m\in\mathbb{^{\ast}N}$)) is constructed from
the original function. It is eventually shown that there exist hyper functions
in this hyper sequence ($^{\ast}f_{m}, m\in\mathbb{^{\ast}N}\setminus\mathbb{N}$) which
are infinitely close to the original function by placing a hyper norm on
$\mathcal{C}_{p}(\mathbb{^{\ast}N}, \mathbb{^{\ast}Q}^{\lim_{p}})$. This
hyper norm is based on the classical norm acting on functions in
$\mathcal{C}_{p}(\mathbb{Z}_{p}, \mathbb{Q}_{p})$. In particular it is shown
that functions of $\mathcal{C}_{p}(\mathbb{^{\ast}N},
\mathbb{^{\ast}Q}^{\lim_{p}})$ have a monad (hyper functions of
$\mathcal{C}_{p}(\mathbb{^{\ast}N}, \mathbb{^{\ast}Q}^{\lim_{p}})$ which are
infinitely close to a given function under the hyper norm). In each monad there
is only one standard function infinitely close to the original function. This
can also be interpreted as the hyper norm is defined for functions in $\mathcal{C}_{p}(\mathbb{^{\ast}N},
\mathbb{^{\ast}Q}^{\lim_{p}})$ and in $\mathcal{C}_{p}(\mathbb{Z}_{p},
\mathbb{Q}_{p})$ and for each $f\in\mathcal{C}_{p}(\mathbb{Z}_{p},
\mathbb{Q}_{p})$       there exists at least one hyper function in
$\mathcal{C}_{p}(\mathbb{^{\ast}N}, \mathbb{^{\ast}Q}^{\lim_{p}})$ such that
the $p$-adic shadow map of this function is $f$.

A potential advantage of considering Mahler interpolation in this form is that
on one level the notion of $p$-adic spaces almost "disappears" as the
interpolated function is essentially a hyper function
$^{\ast}f:\mathbb{^{\ast}N}\rightarrow \mathbb{^{\ast}Q}$ which has many of the
same properties as functions from $\mathbb{N}$ to $\mathbb{Q}$. In future this
could make work studying Mahler interpolation simpler.

\section{Chapter \ref{chapp} Overview}

The theme of this chapter is $p$-adic interpolation but in a more general
setting than the previous chapter. There are many more methods of $p$-adic
interpolation such as those used to interpolate the Riemann zeta function or the
gamma function. There is no reason why these cannot be viewed from a nonstandard
viewpoint as well.

The gamma function provides an ideal candidate for interpolation and the study
of the Morita gamma function is in some ways easier than Mahler interpolation.
In a lot of literature on nonstandard analysis calculus and certain elements of
analysis are presented and defined entirely in a nonstandard way. This
philosophy can be carried through to $p$-adically interpolating the gamma
function. The aim of the interpolation is to find a uniformly $Q$-continuous
hyper function from $\mathbb{^{\ast}N}$ to $\mathbb{^{\ast}Q}^{\lim_{p}}$ which
interpolates the gamma function. Simply defining the hyper factorial (definition
\ref{hypgam}) does not work due to lack of continuity from the powers of $p$. By
defining a restricted factorial (as in section 5.1) such a hyper function can
be found, for all
$n\in\mathbb{^{\ast}N}$ with $n\geq2$ \[ ^{\ast}\Gamma_{p} :
\mathbb{^{\ast}N}\rightarrow \mathbb{^{\ast}Z}, \qquad
^{\ast}\Gamma_{p}(n)= ( -1)^{n}\prod_{1\leq j < n, p\nmid n}
j,\] with $^{\ast} \Gamma(0) =-\text{} ^{\ast} \Gamma(1) = 1$. It is shown
in section 5.1 that it has a \textit{functional equation}. One can now study
this object as the interpolation of the gamma function. Naturally the link to
the standard world is provided by the $p$-adic shadow map and results in the
Morita gamma function.

The next section searches for a version of the Kubota-Leopoldt zeta function in
the nonstandard world. A key point is that the hyper Riemann zeta function is
used to provide values which combine to give the nonstandard function which
interpolates the Riemann zeta function. By proving continuity properties the
Kubota-Leopoldt function is obtained under the $p$-adic shadow maps. One
observation is the role of the hyper Riemann zeta function. In section 2.4.1 it
was shown that the shadow map onto the complex numbers takes
$\zeta_{\mathbb{^{\ast}Q}}$ to $\zeta_{\mathbb{Q}}$ and similarly in this
section $\zeta_{\mathbb{^{\ast}Q}}$ is mapped to the Kubota-Leopoldt zeta
functions via the $p$-adic shadow maps by
\begin{thm7}
For a fixed $\sigma_{0}\in\{-1,1,3,\ldots p-3\}$ \[\operatorname{sh}_{p}
((^{\ast}f_{\sigma_{0}}(\sigma)) = \zeta_{p,\sigma_{0}+1}( \operatorname{sh}_{p}
(\sigma)).\]
\end{thm7} (Here $\zeta_{p,s_{0}}$ are the branches from the classical
Kubota-Leopoldt function and $^{\ast}f_{\sigma_{0}}$ is the nonstandard function
interpolating the Kubota-Leopoldt function defined in section \ref{621}.)  In
this way $\zeta_{\mathbb{^{\ast}Q}}$ can be seen as a source of these zeta
functions.

Probably one of the most important concepts contained in this work is that of
\textit{double interpolation}, that is interpolating a real object with respect
to two distinct finite primes. Such a concept in the standard world is
seemingly not possible, apart from the trivial case when the values lie in
$\mathbb{Q}$, because such a double interpolated function would have
to take values in $\mathbb{Q}_{p_{1}}$ and $\mathbb{Q}_{p_{2}}$ which are
non-isomorphic. By treating interpolation from a nonstandard perspective a
solution can be found. An interpolating set of hyper values are needed with the
requirement that they are uniformly $Q$-continuous with respect to both prime
valuations. This is not too restrictive and they form a subset of hyper
functions which are interpolated with respect to a single prime. The
nonstandard space of $\mathbb{^{\ast}Q}$ provides an ideal sanctuary for such
hyper functions since the only requirement is that the values lie in
$\mathbb{^{\ast}Q}^{\lim_{p_{1}}}\cap \mathbb{^{\ast}Q}^{\lim_{p_{1}}}$. So a
function for double interpolation is going to be of the form \[
^{\ast}f:\mathbb{^{\ast}N}\rightarrow \mathbb{^{\ast}Q}^{\lim_{p_{1}}}\cap
\mathbb{^{\ast}Q}^{\lim_{p_{2}}}.\] By taking the respective shadow maps
standard functions are obtained. The aim of double interpolation is the same as
that of single interpolation in that by looking at a real problem from a
different perspective new information about the object can be found. As a note
the process of double interpolation extends to interpolation with respect to a
finite set of primes and in a special case to all finite primes.

An explicit example of double interpolation is given of the Riemann zeta
function. The ideal choice of numbers to interpolate would be the set
$\{(1-p^{m})(1-q^{m})\zeta_{\mathbb{^{\ast}Q}}(-m)\}_{m\in\mathbb{^{\ast}N}}$.
This follows from the work on the Kubota-Leopoldt zeta function.
Firstly these need to be made continuous and by use of the
Kummer congruences the double interpolation can take place. The
next part of this looks at the resulting $p$-adic function (one
only needs to consider one of the primes by symmetry) coming
from the shadow map of the double Riemann zeta function. By
extending some results of Katz on $p$-adic measures (in
particular theorem \ref{bc1}) it is shown that a $p$-adic
measure exists (lemma \ref{bc3}) which corresponds to the
$p$-adic function and enables it to be written as an integral
over $\mathbb{Z}_{p}^{\times}$ (lemma \ref{bc5})

\begin{thm8}
Let $a\in\mathbb{N}$ with $a\geq 2$ and $(a,p)=1$. Also let
$r\in\mathbb{N}$, $(r,p)=1$. Then for all $m\in\mathbb{N}$
\[ (1-a^{m+1})r^{m}\zeta_{\mathbb{Q}}(-m)=\left( t\frac{d}{dt}
\right) ^{m}\Psi_{r}(t)\mid_{t=1}.\]
Here $\Psi_{r}(t)=(1-t^{ra})^{-1}\sum_{b=1}^{a}\xi_{r}(br)t^{br}$
and

\begin{align}
\xi_{r}:\mathbb{Z}&\rightarrow \mathbb{Z}\notag, \\
n&\mapsto \left\{ \begin{array}{ll}
0& \mbox{$r\nmid n$,}\\
 1 & \mbox{$r\mid n$, $ra\nmid n$,}\\ 1-a
& \mbox{$r\mid n$, $ra \mid n$.} \end{array} \right. \end{align}
\end{thm8}

Although the Morita gamma function was simple to interpolate
using nonstandard methods the "double gamma" function is the
trivial value 1. The actual problem in determining its
trivial existence is from elementary number theory. The problem
is as follows. Let $M_{n}=\{x:1\leq x\leq n,p\nmid x, q\nmid
x\}$ where $p$ and $q$ are fixed primes. Then does there exist a
$j\in\mathbb{N}\setminus\{0,1\}$ with $p,q\nmid j$ such that $j$
has an inverse in $M_{p^{r}}$ and $M_{q^{s}}$ for all
$r,s\in\mathbb{N}$ (where $M_{n}$ can be considered as a
multiplicative set of elements modulo $n$)?  The proof is given
in theorem \ref{ivanko}.

As stated above there is a special case of interpolation with respect to all
finite primes (which would not be possible for the Riemann zeta function because
it would be the interpolation of the trivial function 1). It is the function
$n\rightarrow n^{s}$ for some fixed $n\equiv 1 \pmod{p}$. This section shows
that there is a hyper function which is continuous with respect to each finite
prime and by taking the shadow map for any prime $p$ the original $n^{s}$
function is obtained.

The final section gives a method of constructing a nonstandard Teichm\"{u}ller
character. It has the feel of an artificial method  (section 5.8.1) as the
hyperfunction is constructed based on a set of properties it is expected to
have. Using this a definition is given for a double Hurwitz zeta function but
none of the properties are explored. Finally a brief mention is
given of the difficulty in constructing double $L$-functions.

\section{Chapter \ref{chshai} Overview}
This final chapter gives a conceptual overview of the work of
Shai Haran. The central aim of this chapter is to translate my
extensive studies of his work into a simplified version.  Much
of his work is contained in his book, Mysteries of the Real
Prime (\cite{Ha1}). In some extended lectures on his work he gave
several years after its publication he described his own book as
"very condense and hard to read" and as a result very few people
understand this and hence his research. The hope is that this
chapter may be useful on its own to help people who may wish to
study his book (perhaps as an overview of his work before
tackling his book) enabling the reader to have clear ideas of
what his work involves. Naturally the aim of any research is to
try and produce some original work. The final section of this
chapter contains some very basic attempts to view some parts of
his work in a nonstandard way.

Shai Haran's work is centred around the dictionary between
arithmetic and geometry. Although this dictionary is very
powerful there are a two main problems with it. The first
problem relates the geometric property of adding the point at
infinity to the affine line which produces projective geometry.
The analogue of $\infty$ in the arithmetic picture is the real
prime ($\eta$). This is not an unfamiliar object for example the
completion of $\mathbb{Q}$ with respect to $\eta$ is
$\mathbb{R}$. What is unfamiliar are the real integers,
corresponding to the $p$-adic integers. What is
$\mathbb{Z}_{\eta}$?

The second problem is the lack of, what he calls, an
arithmetical surface. For example in geometry the product of two affine lines
is a plane yet the corresponding entities in the arithmetic picture is trivial
as the product in the category of commutative rings of
$\mathbb{Z}\otimes\mathbb{Z}$ is simply $\mathbb{Z}$. The hunt is then for
a category where this does not happen.

Much of his work could be defined as ultimately trying to find a new language in
which the dictionary between arithmetic and geometry is more complete. His way
of doing this is to try and view all primes of $\mathbb{Q}$ (finite $p$ and the
real prime $\eta$) on an equal footing. Many of his methods in his book are
examples of this technique.

His book mainly covers work in the direction of the first problem. By taking the
view that the real integers should be in some way comparable to the $p$-adic
integers a $p$-adic approach can be made to the real integers since plenty of
information is known about $\mathbb{Z}_{p}$. This approach is refreshing given
that for the past century real results have been applied in the $p$-adic world
yet the $p$-adic world in some respects is simpler.  His "tool of comparison" is
the $q$(uantum)-world. He shows that this world interpolates between the real
and $p$-adic worlds of several objects. By looking at limits of when
$q\rightarrow 0$ and $q\rightarrow 1$ he obtains $p$-adic objects and real
objects respectively. By considering a certain $q$ version of well studied
$p$-adic Markov chains he gives an interpretation of some aspects of the real
integers by looking at the chain when $q\rightarrow 1$. It should be stressed
that it his work is an interpretation and future work will show the validity of
this work. Also, in one sense quite separate from this, though related to the
$q$-world, is the second half of his book where he works on the Riemann zeta
function which is closely related to the work of Connes
(\cite{Co1} and \cite{Co2}).

The second problem is in its infancy and his work takes quite a different
direction to the work which is already being developed in this area. The main
object is the "field of one element" ($\mathbb{F}$). Most other work examines
objects related to this such as schemes, varieties and zeta functions over this
"field". Haran differs in his approach and tries to develop a new geometric
language based on, what he calls $\mathbb{F}$-rings. These are certainly
related to the "field of one element" but in a less than obvious way. Also in
this part I make a small review of the other, more traditional, work related to
the "field of one element".

\chapter{The Nonstandard Algebraic Number Theory of Robinson} \label{chrobo}

\section{Nonstandard Mathematics}
Throughout the history of mathematics infinitely large
numbers and, more commonly, infinitely small numbers (infinitesimals) have
caused problems in establishing their existence. Infinitesimals
initially appeared in the mathematical work of the Greek atomist
philosopher Democritus when he put forward the question of
whether it is possible to generate a cone by piling up circular
plane surfaces with decreasing diameter and in his argument for
the existence of atoms. (This problem and similar ones related
to the the concept of the continuum. The concepts of
discreteness and continuity are still of interest today for
example in philosophy and in the nature of time.) Archimedes
also used infinitesimals to give results regarding areas and
centres of gravity but did not state them as proofs because he
did not believe in them. Most well known is the work of Leibniz
in the development of calculus and (the differential notation
$dx$) though this was later replaced by the $\epsilon-\delta$
method of the nineteenth century as the transfer principle of
Leibniz (that results in the reals can be extended to the
infinitesimals) could not be rigorously justified.

The concept of infinitely large numbers
has a similar history beginning with the Indians in the ancient
Yajur Veda (c. 1200--900 BC) which states that "if you remove a
part from infinity or add a part to infinity, still what remains
is infinity". It continued with  Greeks (for example the
paradoxes of Zeno) and to the modern work in axiomatic set
theory.

The rigorous logical framework for infinitely large and small numbers was
pioneered by Robinson in the 1960s via nonstandard mathematics, a branch of
mathematical logic - model theory. A basic fact in model theory
is that every infinite mathematical structure has nonstandard
models (this basically means that there are non-isomorphic
structures which satisfy the same elementary properties). The
existence of nonstandard models has been known since the 1920s
from the work of Thoralf Skolem. More interest began in the
fifties but it was not until Robinson applied this model
theoretic machinery to analysis that nonstandard analysis was
founded.

Robinson's original presentation (\cite{Ro1}) was considered by
most mathematicians to be unnecessarily complicated because of
the logical formalism needed. Several other approaches appeared
and there exists at least eight different (simpler) presentations
of the methods of nonstandard analysis. They fall into two categories: semantic
, or model-theoretic, approach and the syntactic approach. The most used
presentation is the via superstructures (which is model theoretic like
Robinson's original work) introduced by Robinson and Zakon (it
was also the first one to be purely set-theoretic in nature).

The syntactic approach was created by Nelson in the 1970s using
less model theory by introducing an axiomatic formulation of
non-standard analysis, called \textit{internal set theory}.

\subsection{Overview of Nonstandard Analysis}\label{onsa}
This section introduces some of the basic notions of
nonstandard analysis without focussing on the detailed and
specific formulation which follows later. In approaching it in
this generality some logical rigour is sacrificed but an
initial feel can be obtained. The fundamental object in defining
nonstandard analysis is a universe.

\begin{definition}
A universe, $\mathbb{U}$, is a non empty collection of
"mathematical objects" that is closed under subsets and closed
under basic mathematical operations. These operations are union
of sets, intersection of sets, set difference, ordered pair,
Cartesian product, powerset and function set. A universe is also
assumed to contain (copies of) $\mathbb{N}$, $\mathbb{Z}$,
$\mathbb{Q}$, $\mathbb{R}$ and $\mathbb{C}$. Further $\mathbb{U}$
is assumed to be transitive ($a\in A\in\mathbb{U}\Rightarrow
a\in\mathbb{U}$). \end{definition}
 One point of interest is the notion of "mathematical
objects". These are taken to include all the objects of
mathematics (numbers, sets, functions, relations, ordered
tuples, Cartesian products, etc. ). In fact sets are actually
enough and all the objects can be formalized in the foundational
framework of Zermelo-Fraenkel axiomatic set theory. (For example
a function $f:A\rightarrow B$ can be identified with the set of
pairs $\{\langle a,b\rangle:b=f(a)\}$ which is a subset of the
Cartesian product $A\times B$. Conversely a function from $A$ to
$B$ can be defined set theoretically.)

The nonstandard universe is obtained via a \textit{star map}.
This is a one-to-one map $\ast:\mathbb{U}\rightarrow\mathbb{V}$
between two universes that maps object $A\in\mathbb{U}$ to its
hyper-extension (or sometimes termed nonstandard extension)
$^{\ast}A\in\mathbb{V}$. Further it is assumed that for all
$n\in\mathbb{N}$, $^{\ast}n=n$ and $^{\ast}\mathbb{N}\neq
\mathbb{N}$.

The star map has a powerful property in that it preserves a large
class of properties of the standard universe. This is the
transfer principle (or often called the Leibniz principle). In a
non-rigorous way let $P(a_{1},\ldots a_{n})$ be a property of
standard objects $a_{1},\ldots a_{n}$ which has a bounded
formalization in a language is true iff it is true about the
corresponding hyper-extensions $^{\ast}a_{1},\ldots \text{}
^{\ast}a_{n}$. This will be made more precise by using
mathematical logic.

Combining the previous two paragraphs leads to the following
definition.
\begin{definition} A model of nonstandard analysis is a triple
$(\ast,\mathbb{U},\mathbb{V})$ where $\ast:\mathbb{U}\rightarrow
\mathbb{V}$ is a star-map satisfying the transfer principle.
\end{definition}

The other fundamental principle of nonstandard analysis is the
saturation property. Saturation has a precise definition is model
theory, which will be detailed later. However in terms of set
theory saturation can be given an elementary formulation as an
intersection property.

\begin{definition}
An internal set is any $x$ with $x\in \text{}^{\ast}A$ for some
standard $A$. An external set is an element of the
nonstandard model that is not internal. \end{definition}

\begin{definition}
Let $X$ be a set with $A=(A_{i})_{i\in I}$ a family of subsets of
$X$.  Then the collection $A$ has the finite intersection
property, if any subcollection $J\subset I$ has non-empty
intersection $\bigcap_{i\in J} A_{i}\neq
\emptyset$.
\end{definition}

\begin{definition}[$\kappa$-Saturation]
Let $\kappa$ be an infinite cardinal.  Then the
$\kappa$-saturation principle states that if $I$ is an index set
with cardinality $|I|<\kappa$ and $(A_{i})_{i\in I}$ is a family
of internal sets of an internal set $A$ having the finite
intersection property, then $\bigcap_{i\in I} A_{i}\neq
\emptyset$.
\end{definition}

Robinson's original presentation contained a weaker form of
saturation where he used concurrent relations.
It was Luxemburg who introduced $\kappa$-saturation as
a fundamental tool in nonstandard analysis and in particular for
the nonstandard study of topological spaces. The final point
relates to the existence of a star map which can be constructed
using ultraproducts.

\subsection{Model Theory}
Mathematical logic, like many areas of mathematics, has various
branches. One branch is the study of mathematical structures by
considering the first order sentences true in these structures
and sets definable by first order formulae: \textit{model
theory}. The two main, but connected, reasons for studying model
theory are finding out more about a mathematical structure using
model theoretic techniques and given theories proving general
theorems about their models. Some recent results from model
theory have been related to number theoretic problems. The two
main results are the Mordell-Lang conjecture \footnote{The
generalised Mordell-Lang conjecture states that
the irreducible components of the Zariski closure of a subset of a group of
finite rank inside a semi-abelian variety are translates of closed algebraic
subgroups.} (for function fields in positive
characteristic) and another proof of the Manin-Mumford
conjecture \footnote{One way of stating the Manin-Mumford
conjecture is that a curve C in its Jacobian variety J can
only contain a finite number of points that are of finite order
in J, unless C = J. }, both were proved by Hrushovski. For a
concise overview of model theory see \cite{Mar} and for a more
detailed introduction the textbook of \cite{C-Kei} provides a
solid grounding in this area of logic.

The fundamental objects are structures and the components of a
language.

\begin{definition}[Language] A language is a collection of
symbols of three types: \begin{itemize} \item  set of function
symbols ($\mathcal{F}$) and a $n_{f}\in\mathbb{N}$ for each
$f\in\mathcal{F}$; \item set of relation symbols ($\mathcal{R}$)
and $n_{R}\in\mathbb{N}$ for each $R\in\mathcal{R}$; \item set of
constant symbols ($\mathcal{C}$).\end{itemize} The numbers
$n_{f}$ and $n_{R}$ are the arities of the function $f$ and
relation $R$ respectively. \end{definition}

A simple example is given by the language of rings
$\mathcal{L}_{r}=\{ +,-,\cdot,0,1 \}$ where $+,-$ and $\cdot $
are binary function symbols and 0 and 1 are constants.

\begin{definition}[$\mathcal{L}$-Structure] A
$\mathcal{L}$-structure $\mathcal{M}$ consists of the following:
\begin{itemize} \item a non empty set $M$ (the universe of
$\mathcal{M}$); \item for each function symbol $f\in\mathcal{F}$
a function $f^{\mathcal{M}}:M^{n_{f}}\rightarrow M$; \item for
each relation symbol $R\in\mathcal{R}$ a set
$R^{\mathcal{M}}\subset M^{n_{R}}$; \item for each
constant symbol $c\in\mathcal{C}$ an element $c^{\mathcal{M}}\in
M$. \end{itemize} $f^{\mathcal{M}},R^{\mathcal{M}}$ and
$c^{\mathcal{M}}$ are the interpretation of the symbols $f$, $R$
and $c$ respectively. \end{definition}

When there is no confusion the superscript $\mathcal{M}$ is
often dropped. Naturally maps can be considered between
$\mathcal{L}$-structures and the ones of interest are those which
preserve the interpretation of $\mathcal{L}$.

\begin{definition} Suppose $\mathcal{M}$ and $\mathcal{N}$ are
$\mathcal{L}$-structures with universes $M$ and $N$ respectively.
Then an $\mathcal{L}$-embedding $\theta :\mathcal{M}\rightarrow
\mathcal{N}$ is a one-to-one map $\theta : M\rightarrow N$ such
that: \begin{itemize} \item for all $f\in\mathcal{F}$ and $a_{1},
\ldots a_{n_{f}}\in M$, $\theta(f^{\mathcal{M}}(a_{1}, \ldots
a_{n_{f}}))=f^{\mathcal{N}}(\theta(a_{1}),\ldots
\theta(a_{n_{f}}) )$; \item for all $R\in\mathcal{R}$ and $a_{1},
\ldots a_{m_{R}}$, $(a_{1},\ldots a_{m_{R}})\in R^{\mathcal{M}}$
iff $(\theta(a_{1}),\ldots \theta(a_{m_{R}}))\in
R^{\mathcal{N}}$; \item for $c\in\mathcal{C}$,
$\theta(c^{\mathcal{M}})=c^{\mathcal{N}}$. \end{itemize} Further
$\mathcal{M}$ is a substructure of $\mathcal{N}$ (or
$\mathcal{N}$ is an extension of $\mathcal{M}$) if $M\subset N$
and the inclusion map is a $\mathcal{L}$-embedding. A bijective
$\mathcal{L}$-embedding is called a $\mathcal{L}$-isomorphism and
in the case $M=N$ it is called a $\mathcal{L}$-automorphism.
Finally the cardinality of $\mathcal{M}$ is $|M|$.
\end{definition}

An example of a structure is $(\mathbb{R},+,0)$ which is
an $\mathcal{L}_{g}$-structure where $\mathcal{L}_{g}=\{+,0\}$
with $+$ is a binary function and $0$ a constant. The structure
$(\mathbb{Z},+,0)$ is a substructure of $(\mathbb{R},+,0)$.

In any language (mathematical or linguistical) one aim is to use
it to use it to convey ideas (often via sentences). In logic the
first step is to create formulae to describe properties of
$\mathcal{L}$-structures. Formulae are strings built using the
symbols of $\mathcal{L}$ and the (assumed disjoint) set of
logical symbols (which consists of logical connectives ($\wedge,
\vee, \neg, \rightarrow$ and $\leftrightarrow$), parentheses
($[,]$), quantifiers ($\forall$ and $\exists$), the equality
symbol ($=$) and variable symbols). Simplistically terms are
expressions obtained from constants and variables by applying
functions.

\begin{definition} The set of $\mathcal{L}$-terms is the smallest
set $\mathcal{T}$ such that \begin{itemize} \item for each
$c\in\mathcal{C}$, $c\in\mathcal{T}$, \item each variable symbol
is an element of $\mathcal{T}$, and \item if $f\in\mathcal{F}$
and $t_{1},\ldots t_{n_{f}}\in\mathcal{T}$ then $f(t_{1}, \ldots
t_{n_{f}})\in\mathcal{T}$.\end{itemize}\end{definition}

A $\mathcal{L}$-term (t) has a unique interpretation in a
$\mathcal{L}$-structure $\mathcal{M}$ as a function $t
^{\mathcal{M}}: M^{m}\rightarrow M$. For a subterm $s$ of a term
$t$ and $\bar{a}=(a_{i_{1}},\ldots, _{i_{m}})\in M$.
$s^{\mathcal{M}}(\bar{a})$ can be defined inductively
\begin{itemize} \item if $s$ is the constant symbol $c$ then
$s^{\mathcal{M}}(\bar{a})=c^{\mathcal{M}}$, \item if $s$ is the
variable $v_{i_{j}}$ then $s^{\mathcal{M}}(\bar{a})=a_{i_{j}}$,
\item if $s$ is the term $f(t_{1},\ldots t_{n_{f}})$ (where $f$
is a function symbol of $\mathcal{L}$ and $t_{i}$ are terms) then
$s^{\mathcal{M}}(\bar{a})=f^{\mathcal{M}}(t^{\mathcal{M}}_{1}
(\bar{a}),\ldots, t^{\mathcal{M}}_{n^{f}}(\bar{a})$.
\end{itemize} Finally the definition of $\mathcal{L}$-formulae
can be given.

\begin{definition} $\mathcal{L}$-formulae are defined via atomic
$\mathcal{L}$-formulae. An object $\phi$ is said to be an atomic
$\mathcal{L}$-formula if $\phi$ is either $t_{1}=t_{2}$ (for
terms $t_{1}$, $t_{2}$) or $R(t_{1},\ldots t_{n_{R}})$ (for
$R\in\mathcal{R}$ and $t_{i}$ terms). Then the set of
$\mathcal{L}$-formulae is the smallest set $\mathcal{W}$
containing the atomic formulae such that \begin{itemize} \item if
$\phi$ is in $\mathcal{W}$ then so is $\neg \phi$, \item if
$\phi$ and $\psi$ are in $\mathcal{W}$ then so are $(\phi \wedge
\psi)$ and $(\phi\vee\psi)$, and \item if $\phi$ is in
$\mathcal{W}$ then so are $\forall v_{i} \phi$ and $\exists v_{i}
\phi$ (where $v_{i}$ is a variable).\end{itemize}
\end{definition}

The set of $\mathcal{L}$-formulae split into two depending on the
variables of the formula. Indeed a variable is said to be bound
in a quantifier if it occurs inside a $\exists v$ or
$\forall v$ quantifier, otherwise it is said to be free. A
formula is called a sentence if it has no free variables
otherwise it is called a predicate.

The next concept is \textit{satisfaction}. That is given an
$\mathcal{L}$-formula ($\phi(\bar{x}), \bar{x}=(x_{1}, \ldots
x_{n})$ are free variables) and a $\mathcal{L}$-structure
$\mathcal{M}$ the notion of $\phi(\bar{a})$ being true in
$\mathcal{M}$, $\bar{a}=(a_{1},\ldots a_{n})$ is an $n$-tuple of
elements in $M$. This is denoted by $\mathcal{M}\models
\phi(\bar{a})$ with the negation denoted by $\mathcal{M}\nvDash
\phi(\bar{a})$.

\begin{definition} The notion of $\mathcal{M}\models
\phi(\bar{a})$ is defined inductively by \begin{itemize} \item If
$\phi$ is $t_{1}=t_{2}$ then $\mathcal{M}\models
\phi(\bar{a})$ if $t^{\mathcal{M}}_{1}( \bar{a})=
t^{\mathcal{M}}_{2}(\bar{a})$; \item If $\phi$ is $R(t_{1},
\ldots t_{n_{R}})$ (a $n_{R}$-ary relation) then
$\mathcal{M}\models \phi(\bar{a})$ if
$(t^{\mathcal{M}}_{1}(\bar{a}), \ldots
t^{\mathcal{M}}_{n_{R}}(\bar{a}))\in\mathcal{R}^{\mathcal{M}}$;
\item If  $\phi$ is $\phi_{1}\wedge \phi_{2}$ then
$\mathcal{M}\models \phi(\bar{a})$ if $\mathcal{M}\models
\phi_{1}(\bar{a})$ \textbf{and} $\mathcal{M}\models
\phi_{2}(\bar{a})$; \item If $\phi$ is $\phi_{1} \vee \phi_{2}$
then $\mathcal{M}\models \phi(\bar{a})$ if $\mathcal{M}\models
\phi_{1}(\bar{a})$ \textbf{or} $\mathcal{M}\models
\phi_{2}(\bar{a})$; \item If $\phi$ is $\neg \phi_{1}$ then
$\mathcal{M}\models \phi(\bar{a})$ if $\mathcal{M}\nvDash
\phi(\bar{a})$; \item If $\phi$ is $\exists x\psi(\bar{v},x)$
(where the free variables of $\psi$ are among $\bar{v},x$) then
$\mathcal{M}\models\phi(\bar{a})$ if there is a $b\in M$ such
that $\mathcal{M}\models\phi(\bar{a},b)$; \item If $\phi$ is
$\forall x \psi(\bar{v},x)$ then $\mathcal{M}\models
\phi(\bar{a})$ if $\mathcal{M}\models \psi(\bar{a},b)$ for all
$b\in M$.\end{itemize}\end{definition}

\begin{definition} A $\mathcal{L}$-theory is a set of sentences
of the language $\mathcal{L}$. A model of a theory $T$ is a
$\mathcal{L}$-structure $\mathcal{M}$ which satisfies all the
sentences of $T$, denoted $\mathcal{M}\models T$. Further a
$\mathcal{L}$-theory, $T$, is \textit{satisfiable} iff there
exists a model of $T$ and it is \textit{consistent} iff a formal
contradiction can be derived from $T$.  \end{definition}

(It can be shown as a corollary from the completeness theorem
below that $T$ is satisfiable iff $T$ is consistent.)

 \begin{definition} Let $\mathcal{M}$ and $\mathcal{N}$ be
$\mathcal{L}$-structures with $M\subset N$. $\mathcal{M}$ is an
elementary substructure of $\mathcal{N}$ (or $\mathcal{N}$ is an
elementary extension of $\mathcal{M}$), denoted $\mathcal{M}
\prec \mathcal{N}$, iff for any formula $\varphi(\bar{x})$ and
tuple $\bar{a}$ from $M$ \[ \mathcal{M} \models \varphi(\bar{a})
\leftrightarrow \mathcal{N} \models \varphi(\bar{a}).\] Further a
map $f:\mathcal{M}\rightarrow \mathcal{N}$ is called an
elementary embedding iff it is an embedding and $f(\mathcal{M})
\prec \mathcal{N}$. \end{definition}

\begin{definition}
Let $\phi$ be an $\mathcal{L}$-sentence and $T$ an
$\mathcal{L}$-theory. A proof of $\phi$ from $T$ is a finite
sequence of $\mathcal{L}$-formulae $\psi_{1},\ldots,\psi_{m}$
such that $\psi_{m}=\phi $ and $\psi_{i}\in T$ or $\psi_{i}$
follows from $\psi_{1},\ldots,\psi_{i-1}$ by a simple logical
rule for each $i$. ($T\vdash \phi$ if there is a proof of $\phi$
from $T$.)
\end{definition}

These definitions lead to some important theorems
including the very important \textit{compactness theorem}, one
of the crucial tools of model theorists.

\begin{theorem}[Completeness Theorem] Let $T$ be an
$\mathcal{L}$-theory and $phi$ an $\mathcal{L}$-sentence, then
$T\models \phi$ if and only if $T\vdash \phi$. Moreover, if $T$
has infinite models then $T$ has a model where the model has
cardinality $\kappa$, for all $\kappa \geq |\mathcal{L}| +
\aleph_{0}$.\end{theorem}

\begin{theorem}[Compactness Theorem]  A $\mathcal{L}$-theory $T$
has a model iff every finite subset of $T$ has a model.
\end{theorem}

A proof of the compactness theorem can be given via ultraproducts
though it is also a consequence of the completeness theorem.

It is often useful to work in a very rich model of a theory, for
example it is often easier to prove things in an algebraically
closed field of infinite transcendence degree and in the context
of this work the nonstandard methods in assuming there are
infinite elements when dealing with the reals. This is made
precise by the use of \textit{types} and the property of
\textit{saturation}.

For an $\mathcal{L}$-structure $\mathcal{A}$ let $\mathcal{L}_{A}
= \mathcal{L}\cup \{c_{a}:a\in A\}$ be the expansion of the
language by adjoining constant symbols. This leads onto the
\textit{method of diagrams} but this is not needed in this
introduction. Instead suppose that for some
$\mathcal{L}$-structure $\mathcal{M}$, $A\subset M$ then let
$\operatorname{Th}_{A}(M)$ be the set of all
$\mathcal{L}_{A}$-sentences, $\varphi$, such that $\mathcal{M}
\models \varphi$.

\begin{definition}[Types] An $n$-type over $A$ is a set of
$\mathcal{L}_{A}$-formulas in free variables  $x_{1}\ldots x_{n}$
that is consistent with $\operatorname{Th}_{A}(\mathcal{M})$. A
complete $n$-type is a maximal $n$-type. Let $S_{n}(A)$ be the
set of complete $n$-types over $A$. \end{definition} A formula
$\phi(x_{1},\ldots x_{n})$ is said to be consistent with a
$\mathcal{L}$-theory $T$ iff there exists a model $\mathcal{U}$
which realises $\phi$ (iff for some $n$-tuple of elements
satisfies $\phi$ in $\mathcal{U}$). A more expansive definition
is to say that a complete $n$-type is a set $q$ of
$\mathcal{L}$-formulae consistent with $\operatorname{Th}_{A}(
\mathcal{M})$ in the free variables $x_{1},\ldots ,x_{n}$ such
that for any $\mathcal{L}$-formula, $\varphi(\bar{x})$, either
$\varphi(\bar{x})\in q$ or $\neg\varphi(\bar{x})\in q$.

\begin{definition} Let $\kappa$ be an infinite cardinal. A
structure $\mathcal{M}$ is $\kappa$-\textit{saturated} if for
every $A\subset M$, $|A|<\kappa$ and $p\in S_{1}(A)$ then $p$ is
realized in $\mathcal{M}$. Induction shows that in this case
every $n$-type over $A$ is also realized in $\mathcal{M}$.
$\mathcal{M}$ is said to be \textit{saturated} if it is
$|M|$-saturated. \end{definition}

Finally a nonstandard model can be defined. \begin{definition} A
nonstandard model of, a $\mathcal{L}$ structure, $\mathcal{M}$ is
a saturated elementary extension of $\mathcal{M}$ which is
usually denoted by \[ ^{\ast}: \mathcal{M}\rightarrow \text{}
^{\ast}\mathcal{M}.\]\end{definition}  As mentioned at the end of
section \ref{onsa} an important point is the existence of a
nonstandard model. One method of existence is via ultraproducts.
In fact this is a basic method of constructing models in general
and originated in the work of Skolem in the 1930s and has been
used extensively since the work of \L os in 1955.

\begin{definition}[Filters and Ultrafilters] Let $I$ be a set. A
\textit{filter} on $I$ is a subset $\mathcal{F}$ of (the power
set of $I$) $\mathcal{P}(I)$ satisfying the following properties:
\begin{enumerate} \item $I\in \mathcal{F}$, $\emptyset\notin
\mathcal{F}$; \item if $U\in\mathcal{F}$ and $U\subset V$ then
$V\in\mathcal{F}$; \item if $U,V\in\mathcal{F}$ then $U\cup V
\in \mathcal{F}$.\end{enumerate} An \textit{ultrafilter} on $I$
is a filter on $I$ which such that for any $U\in\mathcal{P}(I)$
either $U\in\mathcal{F}$ or $I\setminus
U\in\mathcal{F}$.\end{definition}

An ultrafilter ($\mathcal{F}$) on a set $I$ is \textit{principal}
if there is $i\in I$ such that $\{i\}\in \mathcal{F}$ (so
$U\in\mathcal{F} \leftrightarrow i\in U$. An ultrafilter is
\textit{non-principal} if it is not principal. The existence of
ultrafilters is provided by following theorem and corollary

\begin{theorem}[Ultrafilter Theorem] If $E\subset \mathcal{P}(I)$
and $E$ has the finite intersection property then there exists an
ultrafilter $\mathcal{F}$ of $I$ such that $E\subset
\mathcal{F}$. \end{theorem}
(A proof can be found in proposition 4.1.4 of \cite{C-Kei}.)
\begin{corollary} Any proper filter of $I$ can be extended to an
ultrafilter over $I$.\end{corollary}

\begin{definition}[Cartesian Products of
$\mathcal{L}$-Structures] Let $I$ be an index set,
$\mathcal{L}$ a fixed language and $(\mathcal{M}_{i})_{i\in I}$ a
family of $\mathcal{L}$-structures. Then the
$\mathcal{L}$-structure $\mathcal{M}= \prod_{i\in I}
\mathcal{M}_{i}$ is defined as follows: \begin{itemize} \item
universe - the cartesian product of $\mathcal{M}_{i}$s (the set
of sequences $(a_{i})_{i\in I}$ such that $a_{i}\in
\mathcal{M}_{i}$ for each $i\in I$. \item constant symbol - for
each $c$ (constant symbol) of $\mathcal{L}$ define
$c^{\mathcal{M}}=(c^{\mathcal{M}_{i}})_{i\in I}$. \item relation
symbol - for $R$ an $n_{R}$-ary relation symbol define $R^{
\mathcal{M}}= \prod_{i\in I}R^{\mathcal{M}_{i}}$. \item function
symbol - for $f$ a $n_{f}$-ary function symbol and
$((a_{1,i})_{i}, \ldots (a_{n,i})_{i})\in M^{n_{f}}$ then
$f^{\mathcal{M}} ((a_{1,i})_{i}, \ldots (a_{n,i})_{i}) = (
f^{\mathcal{M}_{i}}(a_{1,i},\ldots a_{n,i}))_{i\in I}$.
\end{itemize} \end{definition}

Let $I$ be a set and $I$ a filter on this set. Let $\mathcal{M}$
be a cartesian product of $\mathcal{L}$-structures with
$(\mathcal{M}_{i})_{i\in I}$ the related family of
$\mathcal{L}$-structures. An equivalence relation
($\equiv_{\mathcal{F}}$) can be put on $\mathcal{M}=\prod_{i\in
I} \mathcal{M}_{i}$ by \[ (a_{i})_{i}\equiv_{\mathcal{F}}
(b_{i})_{i} \leftrightarrow  \{i\in I:a_{i}=b_{i}\}\in
\mathcal{F}.\] The equivalence class of the element
$(a_{i})_{i}$ is denoted by $(a_{i})_\mathcal{F}$. This is used
to define another $\mathcal{L}$-structure.
\begin{definition}[Reduced Products of $\mathcal{L}$-Structures]
The \textit{reduced product} of the $\mathcal{M}_{i}$s over
$\mathcal{F}$ is denoted by $\prod_{i\in I}\mathcal{M}_{i}\setminus
\mathcal{F}$. It is the quotient structure defined by:
\begin{itemize} \item universe - the quotient of $\prod_{i\in
I}M_{i}$ by $\equiv_{\mathcal{F}}$. \item the interpretation of
a constant symbol $c$ of $\mathcal{L}$ is
$(c^{\mathcal{M}_{i}})_{\mathcal{F}}$. \item for $R$ an
$n_{R}$-ary relation symbol, $f$ an $n_{f}$-ary function symbol
in $\mathcal{L}$ and $a_{1},\ldots ,a_{n}\in\prod_{i\in
I}M_{i}\setminus\mathcal{F}$ (represented by $(a_{1,i})_{i},\ldots,
(a_{n,i})_{i}\in\prod_{i\in I}M_{i}$) then \[ \prod_{i\in
I}\mathcal{M}\setminus \mathcal{F}\models R(a_{1},\ldots, a_{n})
\leftrightarrow \{i\in I: (a_{1,i},\ldots ,a_{n,i})\in
R^{\mathcal{M}_{i}}\}\in \mathcal{F},\] and \[ f^{\mathcal{M}}
(a_{1},\ldots,a_{n})=(f^{\mathcal{M}_{i}}(a_{1,i},\ldots ,
a_{n,i}))_{\mathcal{F}}.\] \end{itemize}\end{definition}
This quotient structure is well-defined by the properties of
filters. In the special case when $\mathcal{F}$ is an ultrafilter
then $\prod_{i\in I}\mathcal{M}_{i}\setminus\mathcal{F}$ is called the
\textit{ultraproduct} of the $\mathcal{M}_{i}$s with respect to
$\mathcal{F}$.

The key result is \L os' theorem which basically connects what formulae are
satisfied in the ultraproduct and in the original structure. The proof is by
structural induction on the complexity of the formulae.

\begin{theorem}[\L os Theorem] Let $I$ be a set, $\mathcal{F}$ an
ultrafilter on $I$ and $(\mathcal{M}_{i})$ ($i\in I$) a family of
$\mathcal{L}$-structures. Let $\varphi(x_{1},\ldots,x_{n})$ be an
$\mathcal{L}$-formula, and let $a_{1},\ldots,a_{n}\in\prod_{i\in
I} M_{i}\setminus\mathcal{F}$ be represented by $(a_{1,i})_{i}, \ldots
,(a_{n,i})_{i}\in\prod_{i\in I} M_{i}$. Then \[\prod_{i\in I}
\mathcal{M}_{i}\setminus\mathcal{F} \models \varphi(a_{1}, \ldots a_{n})
\leftrightarrow \{ i\in I : \mathcal{M}_{i}\models
\varphi(a_{1,i},\ldots a_{n,i})\}\in \mathcal{F}.\]
\end{theorem}
\begin{corollary}\label{mainre} Let $I$ be a set, $\mathcal{F}$
an ultrafilter on $I$ and $\mathcal{M}$ an
$\mathcal{L}$-structure. Then the natural map
$\mathcal{M}\rightarrow \mathcal{M}^{I}\setminus\mathcal{F}$, $a\mapsto
(a)_{\mathcal{F}}$, is an elementary embedding. (Here
$(a)_{\mathcal{F}}$ is the equivalence class of the sequence
with all terms equal to $a$.) \end{corollary}

\subsection{The Hyperreals}

Putting all the work of the previous section together enables the
hyperreals to be defined. Let $I=\mathbb{N}$ and let
$\mathcal{F}$ be a nonprincipal ultrafilter (such ultrafilters
exist on $\mathbb{N}$ by the axiom of choice, see \cite{Go}
corollary 2.6.2 for a proof). Let $\mathcal{L}=\{ +,-,\leq,0\}$
be the language of abelian ordered groups (with $+$ and $-$
binary function symbols, $\leq$ a binary relation symbol and 0 a
constant) and endow $\mathbb{R}$ with its natural
$\mathcal{L}$-structure to get a $\mathcal{L}$-structure $R$.
Let $\mathcal{M}_{i}=R$ for all $i\in\mathcal{N}$ Let
$\mathcal{R}=\prod_{i\in\mathbb{N}}\mathcal{M}_{i}\setminus\mathcal{F} =
R^{\mathbb{N}}\setminus \mathcal{F}$. Then by the corollary \ref{mainre}
$\mathcal{R}$ is an elementary extension of $R$ and further it
is saturated (see chapter 6 of \cite{C-Kei}) - hence it is a
nonstandard model of $R$, therefore let $\mathcal{^{\ast}R}=
R^{\mathbb{N}}\setminus \mathcal{F}$. The map of the corollary is the
$\ast$-map desired with the transfer principle provided by the
map being an elementary embedding. By applying the transfer
principle (to the corresponding statement for
$\{\mathbb{R},+,-,<\}$) the structure $\{
\mathbb{^{\ast}R},+,-,<\}$ is a complete ordered field. One could
also check this by going through the axioms for a field and
checking they hold for $\mathbb{^{\ast}R}$.

Corollary \ref{mainre} implies that $\mathbb{R}$ can be
considered as embedded in $\mathbb{^{\ast}R}$ but are these the
only elements? (In fact a general nonstandard model
$\mathcal{^{\ast}M}$ does contain elements distinct from those
in the original structure because of saturation.) For the
existence of such elements in $\mathbb{^{\ast}R}$ consider the
sets $X_{n}=\{x\in\mathbb{R}: n<x\}$ for each $n\in\mathbb{N}$.
Then $(X_{n})$ is a countable family of sets satisfying the
finite intersection property. Therefore the intersection of
their image ($\cap \text{} ^{\ast}X_{n}$) is non-empty in
$\mathbb{^{\ast}R}$. Then any such element ($\omega$) of the
intersection satisfies $\omega>\text{}^{\ast}x=x$ for all
$x\in\mathbb{R}$. The set of \textit{infinite} numbers is the set
$\{s\in\mathbb{^{\ast}R}:n<|s| \forall n\in\mathbb{N}\}$. The set
of infinite integers are often called \textit{hyperfinite
integers}. A number $s\in\mathbb{^{\ast}R}$ is said to be
\textit{limited} (or \textit{finite}) if $|s|<n$ for some
$n\in\mathbb{N}$.

A similar construction can be used to show that infinitesimal
elements exist by considering the sets $Y_{n}=\{ y\in \mathbb{R}
: 0<x<1/n\}$ for $n\in\mathbb{N}\setminus\{0\}$. The set of
infinitesimals are denoted by $\mu_{\eta}(0)=\{
x\in\mathbb{^{\ast}R}: |x|<1/n \forall n\in\mathbb{N} \}$.

An equivalence relation ($\simeq_{\eta}$) can be defined for
$x,y\in\mathbb{^{\ast}R}$ by $x\simeq_{\eta}y$ iff
$x-y\in\mu_{\eta}(0)$. This leads to the real shadow map from
limited elements of $\mathbb{^{\ast}R}$ to $\mathbb{R}$. It is
shown in full details in definition \ref{realsh} that if
$x\in\mathbb{^{\ast}R}$ is limited then there exists a unique
$\rho \in\mathbb{R}$ such that $\rho\simeq_{\eta} x$ enabling
the definition of the shadow/standard part of $x$ to be
$\operatorname{sh}_{\eta}(x)=\rho$.

A very important concept is to be able to descend from the
nonstandard model to the standard model. This can be done for a
general topological space with a valuation (for example see
\cite{Mac1}). For the hyperreals there is the real shadow map,
denoted $\operatorname{sh}_{\eta}$, briefly described in the
previous paragraph and developed in full detail in chapter 3.

In the construction of $\mathbb{^{\ast}R}$ via ultrafilters, one
essentially views elements of $\mathbb{^{\ast}N}$ as infinite
sequences of real numbers. This construction can also be used to
define functions and subsets on $\mathbb{^{\ast}R}$. For example
the important concept of an \textit{internal set} can be defined.
Suppose there is a given sequence of subsets
$(A_{n})_{\mathbb{N}}$ of $\mathbb{R}$ then a subset $[A_{n}]$ of
$\mathbb{^{\ast}R}$ can be defined by specifying that for each
$r\in\mathbb{^{\ast}R}$, \[r\in [A_{n}]\leftrightarrow
\{n\in\mathbb{N}: r_{n}\in A_{n}\}\in\mathcal{F},\] where
$(r_{n})$ is the equivalence class of $r$ modulo the ultrafilter
$\mathcal{F}$. This can be shown to be well-defined (\cite{Go},
11.1). A set which is not internal is said to be
\textit{external}. Examples of internal sets include
$\mathbb{^{\ast}N}$, $\mathbb{^{\ast}Q}$, $\mathbb{^{\ast}Z}$ and
$\mathbb{^{\ast}R}$. Full
details of internal sets of $\mathbb{^{\ast}R}$ can be found in
\cite{Go}, chapter 11. Similar constructions using sequences of
functions can be used to create hyper functions, and internal
functions can also be defined in an analogous way to internal
sets.

A problem with $\mathbb{^{\ast}R}$, and more generally for a
nonstandard enlargement of a topological space, is the lack of a
canonical topology. Given a topological space $T$, with topology
$\tau$, and enlargement $^{\ast}T$ there are two main topologies
which can be put on this hyper space - the $S$(tandard)-topology
and the finer $Q$-topology. \begin{itemize} \item The basis of
fundamental neighbourhoods for the $S$-topology are generated by
$^{\ast}U$ (where $U$ runs through the open subsets of $T$).
\item The basis of fundamental neighbourhoods for the
$Q$-topology are $^{\ast}V$ (where $V$ is the set of fundamental
neighbourhoods in $\tau$).\end{itemize} In the case of the
hyperreals the $S$-neighbourhoods are of the form
$((r-\epsilon,r+\epsilon))=\{x\in\mathbb{^{\ast}R}:
\operatorname{sh}_{\eta}|r-x|<\epsilon\}$ where
$r\in\mathbb{^{\ast}R}$ and $\epsilon\in\mathbb{R}_{>0}$. The
$S$-open sets are the union of $S$-neighbourhoods and the
$S$-open sets form the $S$-topology on $\mathbb{^{\ast}R}$. The
$Q$-neighbourhoods are of the form $(s-\delta,s+\delta)$ for
$s\in\mathbb{^{\ast}R}$ and $\delta\in\mathbb{^{\ast}R}_{>0}$.
Further details can be found in chapter 4 of \cite{Ro1}.

All the analytical results in this work are proven in the
$Q$-topology as it is the finer topology.

To summarize the hyperreals can be "explicitly" considered as
certain elements of $\mathbb{R}^{\mathbb{N}}$ or just
$\mathbb{R}$ with added elements which are "infinitely close" to
each element of $\mathbb{R}$ and "infinitely large". Most of the
properties of $\mathbb{R}$ hold in the hyperreals by using the
transfer principle. It should be noted that the construction of
the hyperreals is not unique because of the choice of
ultrafilter. In fact if the \textit{continuum hypothesis} is
assumed it can be shown that all quotients of
$\mathbb{R}^{\mathbb{N}}$ with respect to nonprincipal
ultrafilters on $\mathbb{N}$ are isomorphic as ordered fields
(\cite{Go}, 3.16). A guide to the hyperreals can be found in the
form of \cite{Go}.

A slightly more general construction is instead of considering a
structure with universe $\mathbb{R}$ is to consider a
\textit{superstructure}. Indeed let $X$
be a non empty set of atoms (where atoms are objects that can be
elements of sets but are not themselves and are "empty" with
respect to $\in$). Usually it is assumed that (a
copy of) the natural numbers $\mathbb{N}\subset X$.

\begin{definition} The superstructure over $X$ is defined to be
\[ \mathbb{U}(X)=\bigcup_{n\in\mathbb{N}}X_{n},\] where $X_{0}=X$
and by induction $X_{n+1}=X_{n} \cup \mathcal{P}(X_{n})$
($\mathcal{P}(A)$ is the powerset of $A$). \end{definition}
Note that \[ X=X_{0}\subset X_{1} \subset X_{2} \subset \ldots
.\]
It is a simple exercise to check that the conditions for a
universe are satisfied by the superstructure. Superstructures
enable a mathematical object $Z$ to be investigated by ensuring
that (a copy of) $Z$ is contained in $X$. The set elements of
$\mathbb{U}(X)$ are called the entities of $\mathbb{U}(X)$ and
the individuals of $\mathbb{U}(X)$ are the elements of $X_{0}=X$.

\begin{definition} A superstructure based on a set of atoms $X$
is the set $\mathbb{U}(X)$ together with the notions of equality
and membership on the elements of $\mathbb{U}(X)$: $(
\mathbb{U}(X),\in,=)$. \end{definition}
(The notion of equality is assumed given for individuals and
equality between entities (no atom equals any entity) is when
they have the same elements.)

As a solid example consider the superstructure ($\mathcal{N}$)
based on the natural numbers as the atoms. The basic algebra of $\mathbb{N}$
is part of $\mathcal{N}$. For example addition can be taken to be
the following entity \[ S=\{ (a,b,c):a,b,c\in\mathbb{N} ,
a+b=c\}. \] Further number systems can be obtained as entities of
$\mathcal{N}$ by taking pairs of entities. For example
$\mathbb{Z}$ is formed from ordered pairs of $\mathbb{N}$,
$\mathbb{Q}$ is formed from pairs of $\mathbb{Z}$ and
$\mathbb{R}$ from Dedekind cuts. Virtually anything occurring in
classical analysis is an entity of $\mathcal{N}$.

Using the model theory above it has a nonstandard extension along
with the transfer principle. From a formal view an advantage of
a superstructure is that all the "useful" objects such as
functions, metrics,... are already extended to the nonstandard
setting. The definition of internal is clearer as an object $A$
is internal if $A\in\text{}^{\ast}B$ for some $B\in
\mathbb{U}(X)$. Internal objects are important as from the
definition of the transfer principle; transfer takes place from
standard objects to internal objects. For example induction only
takes place on internal subsets of $\mathbb{^{\ast}N}$. For a
full introduction to nonstandard analysis via superstructures see
\cite{S-L}.

In this work there is some abuse of notation, for example an
element of a nonstandard space $^{\ast}X$ should be written as
$^{\ast}x$ but generally the $^{\ast}$ is dropped without
confusion.

\subsection{Applications and Criticisms}

As a final point it should be mentioned about the uses of
nonstandard analysis and the potential criticism. So far
nonstandard analysis has been successfully applied to many areas
such as probability theory (for example certain products of
infinitely many independent, equally weighted random variables),
mathematical economics (for example the behaviour of large
economies) and mathematical physics.

There was initial expectation that nonstandard analysis might
revolutionize the way mathematicians reasoned with the real
numbers but it never happened. Due to the construction of
nonstandard analysis any nonstandard proof can be reinterpreted
using standard techniques. This does not reduce
nonstandard analysis to a mere redundant method of proof since the
methods can actually be quite powerful not only in simplifying
standard proofs, proving/refuting conjectures but also in giving
precise meaning to many informal notions/concepts which do not
make sense classically (like infinitesimals). Despite these
points there is still skepticism about just how much nonstandard
methods add to mathematics. A (well-)known critic is Alain
Connes, as he mentioned in his famous book Non Commutative
Geometry (\cite{Co3}). The reader is left to make their own
opinions.

(Further it cannot be over looked that the 300 year old problem
of formalizing infinitesimals was solved by Robinson using the
power of 20th century logic in the form of nonstandard analysis.
In \cite{Kei} he shows how integral and differential calculus can
be developed entirely using hyperreal numbers.)

Future progress and applications of nonstandard analysis, and of
model theory, can be found in \cite{Fe1} and \cite{Fe2}.

\section{Nonstandard Number Theory}
One of the first applications of nonstandard mathematics was to
algebraic number theory. Unsurprisingly Robinson produced
much work in this area and there have been important
contributions from Roquette and from MacIntyre. This section
gives a brief summary of some of the work in this area both to
give an application of nonstandard mathematics and to form a
basis for which this thesis is an extension.

The main work of Robinson in this area is contained in the papers
\cite{Ro2}--\cite{Ro5} and the joint paper with Roquette
(\cite{R-R}). The first few papers make use of a certain external
nonstandard ideal in a nonstandard extension of a Dedekind ring
in order to look at properties of ideals in the standard
Dedekind ring.

\subsection{The Quotient Ring $\Delta$}
Let $D$ be a Dedekind domain which possesses at least one proper
ideal (where a proper ideals is any ideal other than $D$ or
the zero ideal). $D$ can be enlarged to give an integral domain
$^{\ast}D$. (To obtain this extension $D$ and $\mathbb{N}$ are
embedded in a structure $M$ and then $M$ is enlarged. From the
work above this is done so as to be consistent in making
references to hyperfinite integers and other nonstandard objects.)
Further let $\Omega$ be the set of ideals of $D$ and this
enlarges to the set of internal ideals of $D$, $^{\ast}\Omega$.
\begin{definition}\label{monadu} Define the monad $\mu$ to be the
intersection of all proper standard ideals in $^{\ast}D$. (An
internal ideal $A\in\text{}^{\ast}\Omega$ is standard if there
exists an ideal $B$ in $D$ such that $A=\text{}^{\ast}B$.)
\end{definition} The following properties of $\mu$ are easy to
establish and full proofs can be found in section 3 of the first
paper in \cite{Ro2}.
\begin{prop}\label{moprop}\begin{enumerate}
\item $\mu$ is an external ideal in $^{\ast}D$. \item The only
element of $D$ in $\mu$ is the zero element. \item For any
ideal, $B$, in $D$ then $a\in D\cap\text{}^{\ast}B$ iff $a\in
B$. \item There exists an internal proper ideal $J$ in
$^{\ast}D$ such that $J\subset \mu$. \end{enumerate} \end{prop}
\begin{definition} \label{delmo} Let $\Delta$ be the quotient
ring of $^{\ast}D$ with respect to $\mu$, $\Delta= \text{}
^{\ast}D/\mu$.\end{definition} The main result from this paper
is the following theorem. \begin{theorem}\label{moth} There
exists a one-to-one multiplicative mapping from the proper
ideals of $D$ into the classes of associated elements of
$\Delta$.\end{theorem}

The second paper of \cite{Ro2} continues the above development.
In particular the factorization laws of internal ideals in
$\Delta$ are examined and the introduction of prime ideals.
There are several theorems in the paper but the most important is
in section 7 regarding the factorization of elements of $D$ in
terms of prime elements of $^{\ast}D$.

\begin{theorem}\label{dprime} Let $a$ be an element of $D$, which
can be regarded as a subset of $\Delta$. Further decompose $(a)$
in $D$ into its prime ideals as $(a)=P^{n_{1}}_{1}\ldots
P_{j}^{n_{j}}$. Then there exists representative primes
$\pi_{1},\ldots \pi_{j}$ of $\Delta$ and a unit, $\epsilon$, such
that the following decomposition is unique
$a=\epsilon\pi_{1}^{n_{1}}\cdots\pi_{j}^{n_{j}}$.\end{theorem}

Details of nonstandard finite factorization can be found in
section 7 of \cite{Ro2}.

The final paper in this series is \cite{Ro3}. This paper relates
$\Delta$ to the theory of $p$-adic numbers and adeles. In the
previous papers the enlargements were constructed by the use of
concurrent relations. In this one ultrapowers are used in a
condensed way, similar to the previous chapter. This construction
is required to prove the following theorem.

\begin{theorem}\label{padicdel} Let $D$ be a countable Dedekind
ring (not a field) such that the quotient rings $D\setminus P$
are finite. (Here $P$ is a non-trivial prime ideal in $D$.) Let
$^{\ast}D$ be a comprehensive enlargement then the ring $\Delta$
is isomorphic to the strong direct sum of $P$-adic completions of
$D$. \end{theorem}

In section 5 of \cite{Ro4} the results regarding $p$-adic
completions and adeles are extended to algebraic extensions. This
is also summarized in section 2 of \cite{Mac1}.
 
The penultimate paper moves on from the work of the properties
of $\Delta$ and to some other topics in algebraic number theory.
These areas move into, for example, infinite Galois theory and
class field theory. This work is not pursued in this thesis.
 
The work in \cite{Ro5} continues with work on ideals. In
particular the theory of entire ideals in an infinite algebraic
extension of $\mathbb{Q}$ (and also some class field theory which
is not mentioned here because this thesis does not develop this
area). In the previous works the quotient ring $\Delta$ has been
studied for Dedekind rings and finite algebraic number fields.
This work looks at the quotient ring in infinite algebraic
extensions of $\mathbb{Q}$.

Indeed let $F$ be an infinite algebraic extension of
$\mathbb{Q}$ with rings of integers $F^{i}$ and $\mathbb{Q}^{i}$
respectively. The nonstandard enlargements can be considered,
with $^{\ast}F$ an enlargement of $F$ and $^{\ast}F^{i}$ an
enlargement of $F^{i}$. Further let $\Phi=\{F_{n}\}$ be a tower
subfields of $F$ such that $F_{0}=\mathbb{Q}$, each $F_{n}$ a
finite extension of $\mathbb{Q}$ and $F=\bigcup_{n}F_{n}$. In the
enlargement $^{\ast}\Phi$ is a mapping from $\mathbb{^{\ast}N}$
into the subfields of $^{\ast}$ with $^{\ast}\Phi=\{H_{n}\}$ and
$H_{n}=F_{n}$ for finite $n$. In particular for any infinite
$n=\omega$ (set $H=H_{\omega}$) $F\subset H \subset \text{}
^{\ast}F$.

\begin{definition} Let $\mu$ be the subset of $^{\ast}F^{i}$
defined as $\mu=\{x:x\in\text{}^{\ast}F^{i} $ and $x$ is
divisible by all non zero standard rational integers $\}$.
\end{definition} This naturally leads to the quotient ring
$\Delta = \text{} ^{\ast}F^{i}/\mu$ and the canonical mapping
$\delta:\text{}^{\ast}F^{i}\rightarrow \Delta$. Since $\mu$ does
not contain any standard elements it follows that $\delta$
injects $F^{i}$ into $\Delta$. Further let $\mu_{H}=\mu \cap H$
then $\Delta_{H}=H^{i}/\mu_{H}$ can be identified with a subring
of $\Delta$. Let $\delta_{H}$ be the restriction of $\delta$ to
$H$ so that it maps $H$ on $\Delta_{H}$ and injects $F^{i}$ into
$\Delta_{H}$. The main theorem (3.2) proved is the following.

\begin{theorem} Let $S_{H}$ be the set of internal ideals in
$H^{i}$ and $a\in F^{i}$ ($a\neq 0$). Then $\delta_{H}(a)$
is invertible in $\Delta_{H}$ iff all prime ideals $P_{j}\in
S_{H}$ which divide the ideal $(a)_{H}$ generated by $a$ in $H$
have norms $NP_{j}$ that are powers of nonstandard primes.
\end{theorem}

The final papers of interest are \cite{R-R} and \cite{Mac1}. The
paper by \cite{R-R} is extremely important in nonstandard number
theory. It gives a new and simplified proof of the finiteness
theorem of Siegel-Mahler theorem concerning Diophantine
equations by the use of
nonstandard methods. The main idea is to relate the algebraic
geometry of a number field $K$ and the nonstandard arithmetic of
a nonstandard enlargement of $K$, $^{\ast}K$. So an algebraic
function field in one variable over a algebraic number field $K$
can be viewed as a field of functions over $K$ and as a subfield
of $^{\ast}K$. The Siegel-Mahler theorem can then be restated as
\begin{theorem}\label{nssm} If $F$ is an algebraic function field
in one variable over some fixed number field $K$ such that
$K\subset F\subset \text{}^{\ast}K$, and if $F$ has genus $g>0$,
then every non constant element $x$ of $F$ admits at least one
nonstandard prime divisor of $^{\ast}K$ in its denominator.
\end{theorem}

The method of proof is based on a transfer
principle which can symmetrically translate arguments of $F$ into
equivalent functional properties of $^{\ast}F$. A key facet of
the proof is that it does not use the Mordell-Weil theorem
\footnote{The Mordell-Weil theorem states that for an abelian
variety A over a number field K, the group A(K) of K-rational
points of A is a finitely-generated abelian group.} (although
there are areas of the proof which are similar to areas of the
Mordell-Weil theorem). This enables the revelation of an
effective bound relative to effective bounds in Roth's theorem
\footnote{Given an algebraic number $\alpha$ and a given
$\epsilon>0$ the inequality $|\alpha - p/q|< q^{-(\epsilon +2)}$
has only a finite set of solutions for coprime integers $p$ and
$q$}, see section 7 in \cite{Mac1}.

\cite{Mac1} builds on the work in \cite{R-R}. The work was
originally given as a talk intended to outline model-theoretic
methods in Diophantine geometry. The main result relates to
another formulation of Weil's theory of distributions via a
"covering theorem" relating geometric and arithmetical ideles.
    
\chapter{The Hyper Riemann Zeta Function}\label{chhypr}

\section{The Shadow Maps of $\mathbb{^{\ast}Q}$}
The shadow (or sometimes termed standard) maps enable a standard entity to
be taken from a nonstandard entity. In this work they are extremely important.
In $\mathbb{^{\ast}Q}$ the shadow maps used in this work correspond directly to
the valuations of $\mathbb{Q}$. The simplest one to develop is the real shadow
map.

\begin{definition}[Real Shadow Map]\label{realsh} Let $|.|_{\eta}$ be the
standard archimedean valuations on
$\mathbb{R}$ extended to a hyper valuation on
$\mathbb{^{\ast}R}$.
\begin{enumerate}
\item Two elements $x,y\in\mathbb{^{\ast}R}$ are said to be
infinitesimally close with respect to $|.|_{\eta}$ (denoted
$x\simeq_{\eta} y$) if $|x-y|_{\eta}$ is an infinitesimal element
of $\mathbb{^{\ast}R}^{+}=\{y\in\mathbb{^{\ast}R}: y\geq 0\}$. Define
$\mathbb{^{\ast}R}^{\inf_{p}}=\{ x\in\mathbb{^{\ast}R}:
x\simeq_{\eta} 0 \}$.
\item Define the monad of an element $x\in\mathbb{^{\ast}R}$ to be the set of
elements which are infinitely close with respect to $|.|_{\eta}$, $\mu_{\eta}
(x)=\{y\in\mathbb{^{\ast}R}: x\simeq_{\eta}y\}=\{x\}+
\mathbb{^{\ast}R}^{\inf_{p}}$.
\item $x\in\mathbb{^{\ast}R}$ is said to be limited if there exists
$a,b\in\mathbb{R}$ such that $a<|x|_{\eta}<b$. Define
$\mathbb{^{\ast}R}^{\lim_{p}}=\{ x\in\mathbb{^{\ast}R}: x \text{ is limited}
\}$. \item For $x\in\mathbb{^{\ast}R}^{\lim_{\eta}}$ define
$\operatorname{sh}_{\eta}(x)$ to be the unique element of $\mathbb{R}$ to which
$x$ is infinitesimally close. \end{enumerate} \end{definition}

This definition and the justification for the final statement can be found in
\cite{Go}. He further goes onto define (in chapter 18) the $p$-adic
shadow maps but in quite an explicit way unlike the concise way for the real
shadow map. Indeed define the natural map

\begin{align}
\theta_{p} :& \mathbb{^{\ast}Z}\rightarrow\mathbb{Z}_{p},\notag\\
&x\mapsto \langle x
\pmod{p} , x\pmod{p^{2}}
,\ldots, x\pmod{p^{n}},\ldots \rangle.\notag
\end{align}
This is a surjective homomorphism with kernel:
\begin{align}
\mathbb{^{\ast}Z}^{\inf_{p}} &=
\{x\in\mathbb{^{\ast}Z}:\theta_{p}(x)=0\},\notag\\
&= \{ p^{N}q: N\text{ is unlimited and } q\in\mathbb{^{\ast}Z} \}
.\notag
\end{align}
This extends to $\mathbb{^{\ast}Q}$ by setting
\[\theta_{p}(\frac{x}{y})=\frac{\theta_{p}(x)}{\theta_{p}(y)},\]
(with $x/y$ in lowest form). For this to be defined it is required
that $\theta_{p}(y)\neq 0$. In fact this map is well defined for
$x\in\mathbb{^{\ast}Q}^{\lim_{p}}=\{
q\in\mathbb{^{\ast}Q}:|q|_{p} \text{ is limited }\}$. This is also
a ring homomorphism and has kernel
$\mathbb{^{\ast}Q}^{\inf_{p}}=\{
x\in\mathbb{^{\ast}Q}:|x|_{p}\simeq 0 \text{ in }
\mathbb{^{\ast}R}\}.$

So I define the $p$-adic shadow maps on $\mathbb{^{\ast}Q}_{p}$ using the
$p$-adic valuation directly.

\begin{definition}[$p$-adic Shadow Map]\label{padicsh}
\begin{itemize} \item Let $x,y\in\mathbb{^{\ast}Q}_{p}$, $x$ is
$p$-adically infinitely close to $y$ (denoted $x\simeq _{p} y$)
if there exists $N\in\mathbb{^{\ast}N}-\mathbb{N}$ such that
$|x-y|_{p}< p^{-N}$.
\item The $p$-adic monad of $x\in\mathbb{^{\ast}Q}_{p}$ is
defined to be $\mu_{p}(x)=\{ y\in\mathbb{^{\ast}Q}_{p}:
x\simeq_{p} y\}$.
\item For $x\in\mathbb{^{\ast}Q}_{p}^{\lim_{p}}=
\{x\in\mathbb{^{\ast}Q}_{p}:|x|_{p} \text{ is limited} \}$,
$\operatorname{sh}_{p}(x)$ is defined to be the unique element of
$\mathbb{Q}_{p}$ which is infinitely close in the $p$-adic
valuation to $x$.
\end{itemize}
\end{definition}

The last statement is proved by the following theorem.

\begin{theorem} \label{pshadow}
Every $x\in\mathbb{^{\ast}Q}_{p}^{\lim_{p}}$ is infinitely close
to exactly one number in $\mathbb{Q}_{p}$.
\end{theorem}
\textbf{Proof:} Let $x\in\mathbb{^{\ast}Q}_{p}^{\lim_{p}}$.
Dealing with uniqueness first. Assume $x$ is infinitely close to
two elements in $\mathbb{Q}_{p}$, $x\simeq_{p} a$ and
$x\simeq_{p} b$ with $a\neq b$ and $a,b\in\mathbb{Q}_{p}$. Then
$|a-b|_{p}=|(a-x)+(x-b)|_{p}\leq \max\{ |a-x|_{p}, |x-b|_{p} \}$.
Both of these are infinitesimal by definition and hence
$a\simeq_{p} b$. This is a contradiction since
$a,b\in\mathbb{Q}_{p}$ and therefore $a=b$.
 
Now for existence. As $x\in\mathbb{^{\ast}Q}_{p}^{\lim_{p}}$
there exists an expansion,
\[ x=a_{n}p^{n}+\ldots +a_{0} + a_{1}p +\ldots +
a_{N}p^{N} +\ldots ,\]
where $n$ is not negative unlimited ($n\notin
\mathbb{^{\ast}(-N)}-\mathbb{(-N)}$), otherwise $x$ would not be
in $\mathbb{^{\ast}Q}_{p}^{\lim_{p}}$. Also $0\leq  a_{r}<p$.
Then let  \[ \bar{x}=\sum_{r\geq n,
r\in\mathbb{N}}a_{r}p^{r}\in\mathbb{Q}_{p},\] this is the
standard part of $x$. Let \[ \hat{x}=\sum_{r\geq n,
r\in\mathbb{^{\ast}N}-\mathbb{N}} a_{r}p^{r},\] This is the
nonstandard part of $x$. Thus, $x=\bar{x}+\hat{x}$. Claim
$x\simeq_{p} \bar{x}$, equivalently $\hat{x}\simeq_{p} 0$.
From basic $p$-adic analysis, if the sum $\sum_{n\in\mathbb{N}}c_{n}$ is
absolutely $p$-adically convergent for some $c_{n}\in\mathbb{Q}_{p}$ then
$|\sum_{n\in\mathbb{N}}c_{n}|_{p}\leq \max_{n\in\mathbb{N}}\{ |c_{n}|_{p} \}.$
By transfer this can be applied to finding $|\bar{x}|_{p}$, \[ |\bar{x}|_{p}\leq
\max_{r\in\mathbb{^{\ast}N}-\mathbb{N}} \{|a_{r}p^{r}|_{p}\}.\]
As $|a_{r}|_{p}<1$ the result follows since $r$ is nonstandard.
\begin{flushright} \textbf{$\Box$} \end{flushright}

\begin{theorem}[Properties of $\operatorname{sh}_{p}$]
\label{props}
Let $x,y\in\mathbb{^{\ast}Q}_{p}^{\lim_{p}}$ then
\begin{enumerate} \item $\operatorname{sh}_{p}(x\pm y)=
\operatorname{sh}_{p}(x) \pm \operatorname{sh}_{p}(y),$
\item $\operatorname{sh}_{p}(x. y)=\operatorname{sh}_{p}(x)
.\operatorname{sh}_{p}(y),$
\item $\operatorname{sh}_{p}(x/
y)=\operatorname{sh}_{p}(x)/ \operatorname{sh}_{p}(y),$ if
$\operatorname{sh}_{p}(y)\neq 0$.
\item
$\operatorname{sh}_{\eta}(|x|_{p})=|\operatorname{sh}_{p}(x)|_{p}.$
\end{enumerate} \end{theorem}

\textbf{Proof:} The three statements are basic exercises in
manipulation of the $p$-adic expansions of $x$ and $y$ then use
the definition of theorem \ref{pshadow}. The final property is
very similar in nature and relies on a proof by cases.

Firstly suppose
$x\in\mathbb{^{\ast}Q}_{p}^{\inf_{p}}$. Then
$x=\sum_{r\geq N, r\in\mathbb{^{\ast}N}-\mathbb{N}} a_{r}p^{r}$
where $N\in\mathbb{^{\ast}N}-\mathbb{N}$ thus $|x|_{p}=p^{-N}$
and $\operatorname{sh}_{\eta}(|x|_{p})=0$. However
$\operatorname{sh}_{p}(x)=0$ and so
$|\operatorname{sh}_{p}(x)|_{p}=0$.

In the other case $x\in\mathbb{^{\ast}Q}_{p}^{\lim_{p}}$ such
that $x=\sum_{r\geq n, r\in\mathbb{^{\ast}N}} a_{r}p^{r}$ with
$n\in\mathbb{N}$. Thus $|x|_{p}=p^{-n}$ and
$\operatorname{sh}_{\eta}(|x|_{p})=p^{-n}$. By the definition of
the $p$-adic shadow map $\operatorname{sh}_{p}(x)=\sum_{r\geq
n,r\in\mathbb{N}}a_{r}p^{r}$ and
$|\operatorname{sh}_{p}(x)|_{p}=p^{-n}$.

\begin{flushright} \textbf{$\Box$} \end{flushright}

\begin{lemma}\label{surj} $\operatorname{sh}_{p}$ is a surjective
map with kernel $\mathbb{^{\ast}Q}_{p}^{\inf_{p}}$. \end{lemma}
\textbf{Proof:} The map is surjective because
$\mathbb{Q}_{p}\subset\mathbb{^{\ast}Q}_{p}^{\lim_{p}}$ and
$\operatorname{sh}_{p}$ acts trivially on $\mathbb{Q}_{p}$.

The kernel is
$\ker(\operatorname{sh}_{p})=\{x\in\mathbb{^{\ast}Q}_{p}^{\lim_{p}}
:\operatorname{sh}_{p}(x)=0\}.$ By the definition of the $p$-adic
shadow map this happens precisely when
$x\in\mathbb{^{\ast}Q}_{p}^{\inf_{p}}=\{x\in\mathbb{^{\ast}Q}_{p}^{\lim_{p}}:
|x|_{p}\simeq_{p} 0\}$.

\begin{flushright} \textbf{$\Box$} \end{flushright}

\begin{corollary}\label{isop}
\[ \mathbb{^{\ast}Q}_{p}^{\lim_{p}} /
\mathbb{^{\ast}Q}_{p}^{\inf_{p}} \cong \mathbb{Q}_{p}.\]
\end{corollary}

In both the infinite and the finite prime cases the associated
shadow maps are surjective ring homomorphisms such that
$\mathbb{^{\ast}Q}^{\lim_{\eta}}_{\eta}/\mathbb{^{\ast}Q}^{\inf
_{\eta}}_{\eta}\cong \mathbb{R}$ and $\mathbb{^{\ast}Q}^{\lim
_{p}}_{p}/\mathbb{^{\ast}Q}^{\inf _{p}}_{p}\cong
\mathbb{Q}_{p}.$

The shadow maps can be extended from acting on just the above spaces to the
functions on them.

\begin{definition}[Shadow Image]\label{shimage}
For a function $^{\ast}h:\text{}^{\ast}X\rightarrow \text{}^{\ast}Y$ (where
$\text{}^{\ast} X$ and $\text{}^{\ast}Y$ are sets upon which the
hyper valuation $|.|$ is defined) the shadow image of $^{\ast}h$
(with respect to the valuation $|.|$) is denoted by
$\operatorname{sh}_{|.|}(^{\ast}h)$. This is a function with
domain consisting of the standard parts
$\operatorname{sh}_{|.|}(x)$ ($x\in \text{}^{\ast}X$) such that
$^{\ast}h(x)$ is infinitely close with respect to the valuation
to a standard element in $\mathbb{R}$. This gives the image
consisting of these $\operatorname{sh}_{|.|} (^{\ast}h(x))$ for
$x$ in the domain. \end{definition}

Further details of the real shadow map acting on functions can be found in
\cite{Go}.

\section{Nonstandard Tools}
In order to develop the analytical theory one needs some basic
notions in nonstandard analysis.

\subsection{Properties of $\mathbb{^{\ast}Z}$}
To define the hyper Riemann zeta function one first needs the idea notion of
ideals for $\mathbb{^{\ast}Z}$.

Let $\Lambda$ and $^{\ast}\Lambda$ be the sets of all prime ideals
in $\mathbb{Z}$ and all internal prime ideals
in $^{\ast}\mathbb{Z}$ respectively. Also, let $\Pi$ and
$^{\ast}\Pi$ be the sets of primes numbers in $\mathbb{Z}$ and
$^{\ast}\mathbb{Z}$ respectively. \begin{lemma} \label{L2} I is
an internal prime ideal in $^{\ast}\mathbb{Z}$ $\iff$
$I=p^{\ast}\mathbb{Z}$ for a unique $p$ $\epsilon ^{\ast}\Pi$.
\end{lemma} \noindent\emph{\textbf{Proof:}}
\begin{equation}
( \forall \text{ I} \in \Lambda ) ( \exists ! \text{ p}\in \Pi ) (
\text{I=p}\mathbb{Z} ),
\end{equation}
\begin{displaymath}
\downarrow \text{$\ast$-transform,}
\end{displaymath}
\begin{equation}
( \forall \text{ I} \in \text{} ^{\ast} \Lambda ) ( \exists
!\text{ n}\in \text{} ^{\ast} \Pi ) ( \text{I=p}^{\ast} \mathbb{Z}
) .
\end{equation}
Conversely,
\begin{equation}
( \forall \text{ p}\in \Pi ) ( \exists ! \text{ I} \in \Lambda ) (
\text{I=p}\mathbb{Z} ),
\end{equation}
\begin{displaymath}
\downarrow \text{$\ast$-transform,}
\end{displaymath}
\begin{equation}
( \forall \text{ n}\in \text{} ^{\ast}\Pi ) ( \exists ! \text{ I}
\in \text{} ^{\ast}\Lambda ) ( \text{I=p}^{\ast}\mathbb{Z} ).
\end{equation}
\begin{flushright} $\Box$ \end{flushright}
For $I\in \Lambda$ define the norm of the ideal to be
N(I)$=[\mathbb{Z} :I]$. This definition can be transferred to
$I\in \text{} ^{\ast}\Lambda$, N(I)$=[\mathbb{^{\ast}Z} :I]\in\mathbb{^{\ast}Z_{\geq\text{0}}}$.

Further properties of $\mathbb{^{\ast}Z}$ can be found in
\cite{Go}, \cite{Ro1} and in the papers by Robinson on
algebraic integers and Dedekind rings in \cite{Ro2}, \cite{Ro3}
and \cite{Ro4}.

\subsection{Hyperfinite Sums, Products, Sequences and Integrals}

Internal functions can be constructed out of internal sets and
the full definition can be found in chapter 12 of \cite{Go}.

In order to define hyperfinite summation the symbol '$\sum$' has
to be defined for nonstandard integers.  This can be found in
chapter 19 of \cite{Go} enabling
$\sum_{n=1}^{M}\text{}^{\ast}f(n)$ to be defined for all internal
$^{\ast}f:\mathbb{^{\ast}N}\rightarrow\mathbb{^{\ast}C}$ and
unlimited $M$.  This enables the following definition:

\begin{definition}
Let $f:\mathbb{^{\ast}N}\rightarrow\mathbb{^{\ast}C}$ be
an internal function then the infinite sum exists,
$\sum_{n\in \mathbb{^{\ast}Z_{>\text{0}}}} f(n)=S\in
\mathbb{^{\ast}C}$, if \[ (\forall \epsilon \in
\mathbb{^{\ast}R_{>\text{0}}})(\exists v\in
\mathbb{^{\ast}Z_{>\text{0}}})(\forall k\in
\mathbb{^{\ast}Z_{>\text{$v$}}})(^{\ast}d(\sum_{\text{0}<n\leq
\text{$k$}}f(n),S)<\epsilon).\]
\end{definition}

In an analogous way the hyperfinite product can also be defined.
Let $f:\mathbb{N}\rightarrow\mathbb{C}$ be a complex function
then a finite product can be defined for $m\in\mathbb{N}$;
$\prod_{n=1}^{m}f(n)=f(1)\ldots f(m)$. $'\prod'$ can be regarded
as a function from finite sequences to the complex numbers -
where the finite sequence is $\{f(n)\}_{n=1,\ldots,m}$.  By
transfer (see chapter 19 of \cite{Go}) $'\prod'$ extends to act
on all internal hyperfinite sequences.  So $'\prod'$ is then
defined for any $M\in\mathbb{^{\ast}N}$ and internal hyper
complex functions.

\begin{definition}
Let $f:\mathbb{^{\ast}N}\rightarrow\mathbb{^{\ast}C}$ be an
internal function then the infinite product exists, $\prod_{n\in
\mathbb{^{\ast}Z_{>\text{0}}}} f(n)=P\in \mathbb{^{\ast}C}$, if
\[ (\forall \epsilon \in \mathbb{^{\ast}R_{>\text{0}}})(\exists
v\in \mathbb{^{\ast}Z_{>\text{0}}})(\forall k\in
\mathbb{^{\ast}Z_{>\text{v}}})(^{\ast}d(\prod_{\text{0}<n\leq
\text{$k$}}f(n),P)<\epsilon).\]  \end{definition}

As already seen above chapter 19 of \cite{Go} enables
hyper finite internal sequences and internal sequences to be
defined.

Indeed let
$^{\ast}f:\mathbb{^{\ast}N}\rightarrow\mathbb{^{\ast}C}$ be an
internal function.  Then a hyper internal sequence is defined to
be $\{z_{m}\}$ where $z_{m}=\text{}^{\ast}f(m)$.

\begin{definition}\label{hyseq}
\begin{enumerate}

\item A hyper internal sequence $\{z_{m}\}$ Q-converges to a limit $S\in\mathbb{^{\ast}C}$ if
$(\forall \epsilon\in\mathbb{^{\ast}R_{>\text{0}}})
(\exists N\in\mathbb{^{\ast}Z_{>\text{0}}})(\forall n\in\mathbb{^{\ast}Z_{\geq\text{N}}})
(\text{}^{\ast}d(z_{n},S)<\epsilon)$.
\item Let $\{z_{m} \}$ be a hyper internal sequence in $\mathbb{^{\ast}C}$ then $\{z_{m}\}$ is a hyper
internal Cauchy sequence in $\mathbb{^{\ast}C}$ if $(\forall \epsilon\in\mathbb{^{\ast}R_{>\text{0}}})
(\exists N\in\mathbb{^{\ast}Z_{>\text{0}}})(\forall m,n\in\mathbb{^{\ast}Z_{\geq\text{N}}})
(\text{}^{\ast}d(z_{n},z_{m})<\epsilon)$.
\end{enumerate}
\end{definition}
\begin{prop}\label{hyseq1}
If a hyper internal sequence $\{z_{m}\}$ Q-converges then it is a hyper Cauchy internal sequence.
\end{prop}
\noindent\emph{\textbf{Proof:}}
Suppose $\{z_{m}\}$ Q-converges to $S\in\mathbb{^{\ast}C}$. By definition of Q-convergence,
\[ (\forall \epsilon /2\in\mathbb{^{\ast}R_{>\text{0}}})
(\exists N\in\mathbb{^{\ast}Z_{>\text{0}}})(\forall m,n\in\mathbb{^{\ast}Z_{\geq\text{N}}})
(\text{}^{\ast}d(z_{n},S)<\epsilon /2 \wedge \text{}^{\ast}d(z_{m},S)<\epsilon /2).\]
Under these conditions and by using the triangle inequality,
\begin{align}
^{\ast}d(z_{n},z_{m}) &\leq ^{\ast}d(z_{n},S) + ^{\ast}d(z_{m},S),\notag\\
&< \epsilon /2 +\epsilon/2 =\epsilon.\notag
\end{align}
\begin{flushright} $\Box$ \end{flushright}
\begin{prop} \label{hyseq2}
If $\{z_{m}\}$ is a hyper internal Cauchy sequence then it Q-converges to a limit in $\mathbb{^{\ast}C}$.
\end{prop}
\noindent\emph{\textbf{Proof:}}
Using the transfer principle and the definition of a hyper internal Cauchy sequence,
\[ (\forall \epsilon\in\mathbb{^{\ast}R_{>\text{0}}})
(\exists N\in\mathbb{^{\ast}Z_{>\text{0}}})(\forall m,n\in\mathbb{^{\ast}Z_{\geq\text{N}}})
(\text{}^{\ast}d(z_{n},z_{m})<\epsilon), \]
\begin{displaymath}
\downarrow \text{$\ast$-transform,}
\end{displaymath}
\[ (\forall \epsilon\in\mathbb{R_{>\text{0}}})
(\exists N\in\mathbb{Z_{>\text{0}}})(\forall m,n\in\mathbb{Z_{\geq\text{N}}})
(d(z_{n},z_{m})<\epsilon) \]
This implies $\{z_{m}\}$ for $m\in\mathbb{^{\ast}Z_{>\text{0}}}$ is a Cauchy sequence in
$\mathbb{C}$. By the standard theorem this converges to a limit, A say, in $\mathbb{C}$. By the transfer principle
\[ (\forall \epsilon\in\mathbb{R_{>\text{0}}})
(\exists N\in\mathbb{Z_{>\text{0}}})(\forall n\in\mathbb{Z_{\geq\text{N}}})
(d(z_{n},A)<\epsilon) \]
\begin{displaymath}
\downarrow \text{$\ast$-transform,}
\end{displaymath}
\[ (\forall \epsilon\in\mathbb{^{\ast}R_{>\text{0}}})
(\exists N\in\mathbb{^{\ast}Z_{>\text{0}}})(\forall n\in\mathbb{^{\ast}Z_{\geq\text{N}}})
(\text{}^{\ast}d(z_{n},A)<\epsilon), \]
\begin{flushright} $\Box$ \end{flushright}

The work on defining hyperfinite integrals can be found beginning
in chapter 5 of \cite{Ro1}.

\subsection{Q-Topology}
Consider a nonstandard set $X$ with a hypermetric $^{\ast}D:X\times X \rightarrow \mathbb{^{\ast}R_{>\text{0}}}.$
\begin{definition}\label{qtop}
\begin{enumerate}
\item A $\underline{Q-ball}$ is the set $^{\ast}B(y,r)=\{ x\in X
:\text{}^{\ast}D(x,y)<r \}$ for $y\in X$ and
$r\in\mathbb{^{\ast}R_{>\text{0}}}$.
\item A set $W\subset X$ is called a $\underline{Q-neighbourhood}$ of $y\in X$ if it contains a ball $^{\ast}B(y,r)$.
\item A set is $\underline{Q-open}$ if it is a Q-neighbourhood of each of its elements.
\item The complement of a Q-open set with respect to $X$ is termed to be $\underline{Q-closed}$ with respect to $X$.
\end{enumerate}
\end{definition}
Using the triangle inequality it follows that every Q-ball is Q-open. Naturally many more properties can be defined
and results developed via the transfer principle from the standard case but these are not relevant for this work. The
only property needed is the analogue of compactness in the nonstandard case.
\begin{definition}\label{qtop1}
\begin{enumerate}
\item Consider an internal collection of Q-open sets.
These form a $\underline{Q-open covering}$ for $X$ if $X$ is
contained in the union of these sets. A
$\underline{Q-subcovering}$ is a subcollection with the same
property. \item A $\underline{hyperfinite}$ $\underline{
Q-covering}$ of $X$ is a Q-open covering of $X$ consisting of a
hyperfinite number of sets. \end{enumerate} \end{definition}
These lead to the definition of compactness in the nonstandard
case. \begin{definition}\label{qtop2} A set $X$ is
$\underline{hypercompact}$ with respect to the Q-topology iff
every Q-open covering of $X$ contains a hyperfinite
Q-subcovering. \end{definition}
 
\subsection{Hyper Exponential and Logarithm Functions}

On $\mathbb{C}$ take the usual metric and when extended to
$^{\ast}\mathbb{C}$ it makes $^{\ast}\mathbb{C}$ into a hyper
metric space where the metric takes values in
$\mathbb{^{\ast}R_{\geq \text{0}}}$.

\begin{definition}\label{hyreal}
Let $s\in\mathbb{^{\ast}C}$ with $s=u+iv$ $(u,v\in
\mathbb{^{\ast}R})$. Define the hyperreal part to be $^{\ast}\Re
(s)=u \in \mathbb{^{\ast}R}$ and the hyper imaginary part to be
$^{\ast}\Im(s)=v\in\mathbb{^{\ast}R}$. \end{definition}

\begin{definition} \label{hycont}
Let $f(z)$ be a hypercomplex function. $f(z)$ is
\emph{Q-continuous} if \[(\forall \epsilon \in\mathbb{^{\ast}R_{>\text{0}}})
(\exists \delta\in\mathbb{^{\ast}R_{>\text{0}}})
(\forall z,w\in
\mathbb{^{\ast}C},^{\ast}d(z,w)<\delta)(^{\ast}d(f(z),f(w))<\epsilon
). \]
\end{definition}
 
\begin{lemma} \label{L6}
The limit function of a uniformly Q-convergent internal sequence
of Q-continuous functions is itself Q-continuous. \end{lemma}

\noindent\emph{\textbf{Proof:}} Suppose that the functions
$f_{n}(x)$ are Q-continuous and uniformly Q-converge to $f(z)$
on a set E. For any $\epsilon \in \text{} ^{\ast}\mathbb{R_{>
\text{0}}}$, $n\in \text{} ^{\ast}\mathbb{N}$ can be found such
that $^{\ast} d(f_{n}(z),f(z))<\epsilon /3$ $\forall z\in $ E.
Let $z_{0}$ be a point in E. As $f_{n}(z)$ is Q-continuous at
$z_{0}$, $\exists \delta \in
\text{}^{\ast}\mathbb{R_{>\text{0}}}$ such that
$^{\ast}d(f_{n}(z),f_{n}(z_{0}))<\epsilon/3$ $\forall z\in E$
with $^{\ast} d(z,z_{0})<\delta$. Under the same conditions on
z,
\begin{equation}
^{\ast}d(f(z),f(z_{0})) \leq ^{\ast}d(f(z),f_{n}(z)) +
^{\ast}d(f_{n}(z),f_{n}(z_{0})) +
^{\ast}d(f_{n}(z_{0}),f(z_{0})) < \epsilon.
\end{equation}
\begin{flushright} $\Box$ \end{flushright}

\begin{definition}\label{hypgam}
Let $N\in\mathbb{^{\ast}Z_{>\text{0}}}$ and define the hyper
factorial inductively by $0!=1$ and $N!=N\times (N-1)!$. So
$N!=\prod_{1\leq n\leq N} n.$ This function is interpolated by the
hyper gamma function to be defined below.
\end{definition}

\begin{definition}\label{qder}
Suppose $f(z)$ is an internal hypercomplex function then
$f'(z)\in\mathbb{^{\ast}C}$ is the \emph{Q-derivative} of $f(z)$
at $z=z_{0}$ if \[ (\forall \epsilon \in
\mathbb{^{\ast}R_{>\text{0}}})(\exists \delta \in
\mathbb{^{\ast}R_{>\text{0}}})(\forall
h\in\mathbb{^{\ast}C},^{\ast}d(h,0)<\delta) (^{\ast}d
(\frac{f(z_{0}+h)-f(z_{0})}{h},{f'(z_{0})})<\epsilon). \]
\end{definition} \begin{definition}\label{qany} Let $B$ be a set
of points in $\mathbb{^{\ast}C}$ $(f(z)$ as above). Then $f(z)$
is \emph{Q-analytic} in $B$ if $f(z)$ is infinitely
Q-differentiable at all points of $B$. \end{definition}

\begin{definition}\label{hypexp}
Let $s\in\mathbb{^{\ast}C}$ and define the hyper exponential as \[
^{\ast}\exp (s)=\text{}^{\ast}e^{s}=\sum_{n\in\mathbb{^{\ast}Z_{>\text{0}}}}
\frac{s^{n}}{n!} .\]
\end{definition}

\begin{lemma} \label{L20}
Properties of $^{\ast}e^{s}$.
\begin{enumerate}

\item $^{\ast}e^{s}$ uniformly Q-converges $\forall
s\in\mathbb{^{\ast}C}$, so $^{\ast}e^{s}:\mathbb{^{\ast}C} \rightarrow \mathbb{^{\ast}C^{\times}}$.

\item $^{\ast}e^{s}$ is Q-analytic $\forall s\in\mathbb{^{\ast}
C}$.

\item $^{\ast}\exp(z)^{'}=^{\ast}\exp(z)$.

\item $^{\ast}e^{s+t}=^{\ast}e^{s}\text{} ^{\ast}e^{t} $
($s,t\in\mathbb{^{\ast}C}$). In particular $^{\ast}e^{0}=1$.

\item $^{\ast}\exp : \mathbb{^{\ast}C}/2\pi i\mathbb{^{\ast}Z} \rightarrow \mathbb{^{\ast}C^{\times}}$ is a group is a group isomorphism.

\end{enumerate}
\end{lemma}

\noindent\emph{\textbf{Proof:}}
\begin{enumerate}
\item Let $N\in\mathbb{^{\ast}Z_{>\text{0}}}$ and define
$^{(N)}e^{s}=\sum_{n=1}^{N}\frac{n^{s}}{n!}$. In the classical
case $^{(N)}e^{s}$ converges uniformly $\forall s\in \mathbb{C}$.
\[(\forall \epsilon
\in\mathbb{R_{>\text{0}}})(\exists N\in\mathbb{Z_{>\text{0}}})(\forall m,n\in \mathbb{Z_{>\text{N}}})
(\forall s\in\mathbb{C})(d(^{(n)}e^{s},^{(m)}e^{s})<\epsilon),\]
\begin{displaymath}
\downarrow \text{$\ast$-transform,}
\end{displaymath}
\[(\forall \epsilon
\in\mathbb{^{\ast}R_{>\text{0}}})(\exists N\in\mathbb{^{\ast}Z_{>\text{0}}})(\forall m,n\in
\mathbb{^{\ast}Z_{>\text{N}}})(\forall
s\in\mathbb{^{\ast}C})(^{\ast}d(^{(n)}e^{s},^{(m)}e^{s})<\epsilon).\]

\item The arguments of \ref{L10}, \ref{L11} and \ref{L12}, see below, will give the results.

\item Since, by 2, $^{\ast}\exp$ is Q-analytic it remains to find the derivative. In
the standard case $\exp(s)^{'}=\exp(s)$.
\[  (\forall \epsilon
\in\mathbb{R_{>\text{0}}})(\exists N\in\mathbb{Z_{>\text{0}}})(\forall m,n\in \mathbb{Z_{>\text{N}}})
(\forall s\in\mathbb{C})(d(^{(n)}\exp(s),^{(m)}\exp(s)^{'})<\epsilon),\]
\begin{displaymath}
\downarrow \text{$\ast$-transform,}
\end{displaymath}
\[  (\forall \epsilon
\in\mathbb{^{\ast}R_{>\text{0}}})(\exists N\in\mathbb{^{\ast}Z_{>\text{0}}})(\forall m,n\in \mathbb{^{\ast}Z_{>\text{N}}})
(\forall s\in\mathbb{^{\ast}C})(^{\ast}d(^{(n)}\exp(s),^{(m)}\exp(s)^{'})<\epsilon).\]

\item By 3: $\frac{d}{dz}(^{\ast}\text{}e^{z})=\text{}^{\ast}e^{z}.$
Using this and the product rule gives for some constant
$c\in\mathbb{^{\ast}C}$: $\frac{d}{dz}(^{\ast}e^{z}
\text{}^{\ast}e^{c-z})= \text{}^{\ast}e^{z} \text{} ^{\ast}e^{c-z}
+\text{}^{\ast}e{z}(- ^{\ast}e^{c-z})=0.$ So $^{\ast}e^{z} \text{}
^{\ast}e^{c-z}=K$,constant ($K\in\mathbb{^{\ast}C}$). Let $z=0$
then $K=\text{}^{\ast}e^{c}$. Now let $z=s, c=s+t$ and the result
follows. In particular $^{\ast}e^{s} \text{} ^{\ast}e^{-s}=1$
which implies $^{\ast}e^{s}$ is never zero.

\item In the standard case $\exp(s)=1 \iff c\in 2\pi i\mathbb{Z}$.
\[ (\forall n\in\mathbb{Z})(\forall N\in\mathbb{Z_{>\text{0}}})(\exists \epsilon\in\mathbb{R_{>\text{0}}})
(d(^{(N)}e^{2\pi in},1)<\epsilon),\]
\begin{displaymath}
\downarrow \text{$\ast$-transform,}
\end{displaymath}
\[ (\forall n\in\mathbb{^{\ast}Z})(\forall N\in\mathbb{^{\ast}Z_{>\text{0}}})(\exists \epsilon\in\mathbb{^{\ast}R_{>\text{0}}})
(^{\ast}d(^{(N)}e^{2\pi in},1)<\epsilon).\]
Conversely,
\[ (\forall s\in\mathbb{C}-2\pi i\mathbb{Z})(\forall N\in\mathbb{Z_{>\text{0}}})
(\exists \epsilon\in\mathbb{R_{>\text{0}}})(d(^{(N)}e^{s},1)>\epsilon),\]
\begin{displaymath}
\downarrow \text{$\ast$-transform,}
\end{displaymath}
\[ (\forall s\in\mathbb{^{\ast}C}-2\pi i\mathbb{^{\ast}Z})(\forall N\in\mathbb{^{\ast}Z_{>\text{0}}})
(\exists \epsilon\in\mathbb{^{\ast}R_{>\text{0}}})(^{\ast}d(^{(N)}e^{s},1)>\epsilon).\]
\end{enumerate}
\begin{flushright} $\Box$ \end{flushright}

\begin{definition}\label{hyplog}
$^{\ast}\log :\mathbb{^{\ast}C}\setminus (-\mathbb{^{\ast}R_{>\text{0}}})\rightarrow \mathbb{^{\ast}C}$. It is
the inverse to $^{\ast}\exp$. $z=\text{}^{\ast}\log(w)$ is a root of the equation
$^{\ast}\exp(z)=w$. This equation has infinitely many
solutions so the hyper logarithm is multivalued. So define the
principal value of the hyper logarithm to be
$^{\ast}\log(z)=\text{}^{\ast}\log|z|+i\text{Arg}(z)$ where
$|z|>0$ and $-\pi<\text{Arg}(z)\leq \pi$.  Also define
$^{(n)}\log(w)=z$ to be a root of the equation
$^{(n)}\exp(z)=w$ $\forall n\in\mathbb{^{\ast}Z_{>\text{0}}}$
provided $-\pi<\text{Arg}(z)\leq \pi$. \end{definition}

\begin{lemma}\label{logprop}
\begin{enumerate}
\item $^{\ast}\log(st)=\text{}^{\ast}\log(s) + \text{}
^{\ast}\log(t)$ iff $-\pi<\text{Arg}(s)+\text{Arg}(t)\leq \pi$.
\item $^{\ast}\log(s)$ is Q-analytic $\forall s\in
\mathbb{^{\ast}C}\setminus (-\mathbb{^{\ast}R_{>\text{0}}})$.
Also in this region, $\frac{d}{ds}\text{}
^{\ast}\log(s)=\frac{1}{s}$. \end{enumerate} \end{lemma}

\noindent\emph{\textbf{Proof:}}
\begin{enumerate}
\item Let $u=\text{}^{\ast}\log(s)$ and $V=\text{}^{\ast}\log(t)$ then $u$ and $v$
satisfy $^{\ast}\exp(u)=s$ and $^{\ast}\exp(v)=t$. By \ref{L20}
(4) $^{\ast}\exp(u)\text{}^{\ast}\exp(v)=\text{}^{\ast}\exp(u+v)=st$ and the
result follows since from the definition \ref{hyplog} and the
proof in the classical case. \item For the first part use the
arguments of \ref{L10}, \ref{L11} and \ref{L12} below with
$^{(n)}\log(s)$. For the second part use the arguments of
\ref{L20}(3). \begin{flushright} $\Box$ \end{flushright}
\end{enumerate}

\begin{definition}\label{hypower}
For $z\in\mathbb{^{\ast}C}\setminus i\mathbb{^{\ast}R}$ and $s\in\mathbb{^{\ast}C}$ define
$z^{s}=\text{}\exp(s\text{}^{\ast}\log(z))$, where the principal
value of the hyper logarithm is used. \end{definition}
 
\begin{lemma}\label{hypowerprop}
Let $z_{1},z_{2}\in \mathbb{^{\ast}C}\setminus i\mathbb{^{\ast}R}$ and $s,t\in\mathbb{^{\ast}C}$.
\begin{enumerate}
\item $(z_{1}z_{2})^{s}=z_{1}^{s}z_{2}^{s}$ is not true in
general.  It does hold when lemma \ref{logprop} 1 holds. \item
$z^{s+t}=z^{s}z^{t}$. \item $\frac{d}{dz} z^{s}=sz^{s-1}$.
\end{enumerate} \end{lemma}

\noindent\emph{\textbf{Proof:}}
\begin{enumerate}
\item
Suppose that $z_{1}$ and $z_{2}$ satisfy the properties of lemma
\ref{logprop} then
\begin{align} (z_{1}z_{2})^{s}
&=\text{}^{\ast}\exp(s\text{}^{\ast}\log(z_{1}z_{2})), \notag\\
&=\text{}^{\ast}\exp( s\text{}^{\ast}\log(z_{1})
+s\text{}^{\ast}\log(z_{2})),\notag\\ &=\text{}^{\ast}\exp(
s\text{}^{\ast}\log(z_{1}))  \text{}^{\ast}\exp(
s\text{}^{\ast}\log(z_{2}))=z_{1}^{s}z_{2}^{s}.\notag
\end{align} \item \begin{align}
z^{s+t}&=\text{}^{\ast}\exp((s+t)\text{}^{\ast}\log(z)),\notag\\
&=\text{}^{\ast}\exp(s\text{}^{\ast}\log(z)
+t\text{}^{\ast}\log(z)),\notag\\ &=z^{s}z^{t}.\notag
\end{align} \item Using the chain rule, \begin{align}
\frac{d}{dz} z^{s}&=
\frac{d}{dz}(^{\ast}\exp(s^{\ast}\log(z)))\notag,\\
&=sz^{s}\frac{d}{dz}(^{\ast}\log(z))\notag,\\ &=sz^{s-1}\notag.
\end{align} \end{enumerate} \begin{flushright} $\Box$
\end{flushright}

\subsection{Hyper Gamma Function}
\begin{definition} \label{hygamf}
Hyper gamma function. For $s\in\mathbb{^{\ast}C}$ and $^{\ast}\Re
(s)\in\text{} ^{\ast}[\epsilon,r], r,\epsilon\in\mathbb{^{\ast}R_{>\text{0}}}$ and $r>\epsilon$ define \[^{\ast}\Gamma (s)=\int_{^{\ast}
\mathbb{R^{+}}}^{\ast} \text{} ^{\ast}e^{-y} y^{s-1} dy .\]
\end{definition}
\begin{prop} \label{L13}
$^{\ast} \Gamma (s)$ is absolutely Q-convergent for
$^{\ast}\Re(s)>0$.
\end{prop}
\noindent\emph{\textbf{Proof:}} In the standard case \[
\Gamma(s)=\int_{\mathbb{R^{+}}} e^{-y} y^{s-1}dy, \] is absolutely
convergent for $\Re (s)\in [\epsilon,r]$ $(r,\epsilon\in\mathbb{R_{>\text{0}}}$ $and$ $r>\epsilon$.
Let \[I_{n,r}(s)=\int_{0}^{r}\text{}
^{(n)}e^{-y} y^{s-1}dy \text{, where }
^{(n)}e^{s}=\sum_{j=1}^{n}\frac{s^{j}}{j!} .\]
\begin{equation}
(\forall s\in\mathbb{C}, \Re(s)\in[\epsilon,r]) (\forall \delta \in \mathbb{R_{>\text{0}}})(\exists m,n\in
\mathbb{Z_{>\text{0}}})(\exists r,t \in
\mathbb{R_{>\text{0}}})(d(I_{n,r}(s),I_{m,t}(s))<\delta),
\end{equation}
\begin{displaymath}
\downarrow \text{$\ast$-transform,}
\end{displaymath}
\begin{equation}
  (\forall s\in\mathbb{^{\ast}C}, ^{\ast}\Re(s)\in[\epsilon,r])
(\forall \delta \in \text{}
^{\ast}\mathbb{R_{>\text{0}}})(\exists m,n\in \text{}
^{\ast}\mathbb{Z_{>\text{0}}})(\exists r,t \in \text{}
^{\ast}\mathbb{R_{>\text{0}}})^{\ast}(d(I_{n,r}(s),I_{m,t}(s))<\delta),
\end{equation}
So $\{I_{n,r}\}$ form a hyper internal Cauchy sequence with limit function
$^{\ast}\Gamma(s)$.
\begin{flushright} $\Box$ \end{flushright}

\begin{lemma} \label{L22}
For $^{\ast}\Re (s)>0$, $\int_{\mathbb{^{\ast}R^{+}}}^{\ast}
\frac{d}{dy}(y^{s} \text{} ^{\ast}e^{-y}) dy=0$.
\end{lemma}
\noindent\emph{\textbf{Proof:}} The classical result for $\Re
(s)>0$ is $\int_{\mathbb{R^{+}}} \frac{d}{dy}(y^{s} \text{}
e^{-y}) dy=0$. Using a similar method as the previous proposition.
Let $I_{n,r}(s)=\int_{0}^{r} (\frac{d}{dy}(y^{s} \text{}
^{(n)}e^{-y}) dy$. Then,
\[
(\forall \epsilon \in \mathbb{R_{>\text{0}}})(\exists n\in
\mathbb{Z_{>\text{0}}})(\exists r \in
\mathbb{R_{>\text{0}}})(d(I_{n,r}(s),0)<\epsilon),
\]
\begin{displaymath}
\downarrow \text{$\ast$-transform,}
\end{displaymath}
\[(\forall \epsilon \in \mathbb{^{\ast}R_{>\text{0}}})(\exists n\in
\mathbb{^{\ast}Z_{>\text{0}}})(\exists r \in
\mathbb{^{\ast}R_{>\text{0}}})(^{\ast}d(I_{n,r}(s),0)<\epsilon).
\]
\begin{flushright} $\Box$ \end{flushright}

\begin{theorem}\label{L21}
Functional equation of $^{\ast}\Gamma (s)$: $^{\ast}\Gamma
(s+1)=s\text{} ^{\ast}\Gamma (s).$
\end{theorem}
\noindent\emph{\textbf{Proof:}}
$^{\ast}\Gamma(s+1)=\int_{\mathbb{^{\ast}R^{+}}}^{\ast} \text{}
^{\ast}e^{-y} y^{s}dy$. Integrating by parts and then using
\ref{L22} gives,
\begin{align}
^{\ast}\Gamma(s+1)&=\int_{\mathbb{^{\ast}R^{+}}}^{\ast}
\frac{d}{dy}(y^{s} \text{} ^{\ast}e^{-y})
dy+s\int_{\mathbb{^{\ast}R^{+}}}^{\ast} \text{}
^{\ast}e^{-y} y^{s}dy,\notag\\
&=s ^{\ast}\Gamma(s) \notag.
\end{align}
\begin{flushright} $\Box$ \end{flushright}

\begin{corollary}\label{gamcor}
$^{\ast}\Gamma(1)=^{\ast}\Gamma(2)=1$ and for
$N\in\mathbb{Z_{>\text{0}}}$ $^{\ast}\Gamma(N)=(N-1)(N-2)\ldots 1
= (N-1)!$.
\end{corollary}
\noindent\emph{\textbf{Proof:}} Using the classical result
$\Gamma(1)=1$ and transfer using the same method as \ref{L22} but
with $I_{n,r}=\int_{0}^{r}\text{} ^{(n)}e^{-y} dy$ gives
$^{\ast}\Gamma(1)=1$. Similarly for $^{\ast}\Gamma(2)$.

Let $N\in\mathbb{Z_{>\text{0}}}$ and using \ref{L21} repeatedly
gives the result.
\begin{flushright} $\Box$ \end{flushright}

\begin{prop} \label{pr1}
$^{\ast}\Gamma(s)$ has simple poles at $\mathbb{^{\ast}Z_{\leq\text{0}}} $ and is non-zero.
\end{prop}
\noindent\emph{\textbf{Proof:}} The standard gamma function,
$\Gamma(s)$ has a simple pole at $s=0$, since the integral diverges.
Using the functional equation the only other poles are simple and are at
$-1,-2,-3,\ldots $. Let \[ I_{n,r}=
\int_{\mathbb{^{\ast}[\text{0,r}]}}^{\ast} \text{} ^{(n)}e^{-y}
dy.\] Then, \[ (\forall \alpha\in\mathbb{R_{>\text{0}}})(\exists
n\in\mathbb{Z_{>\text{0}}})(\exists
r\in\mathbb{R_{>\text{0}}})(d(I_{n,r},0)>\alpha) ,\]
\begin{displaymath}
\downarrow \text{$\ast$-transform,}
\end{displaymath}
\[ (\forall \alpha\in\mathbb{^{\ast}R_{>\text{0}}})(\exists
n\in\mathbb{^{\ast}Z_{>\text{0}}})(\exists
r\in\mathbb{^{\ast}R_{>\text{0}}})(^{\ast}d(I_{\text{n,r}},0)>\alpha)
.\] So there is a pole at $0$. Using \ref{L21}, and as in the
classical case, the other poles are at the negative hyperintegers.

The function being non-zero again follows from the functional
equation and by transfer principle as $\Gamma(s)$ is non zero.
\begin{flushright} $\Box$ \end{flushright}

\section{Hyper Riemann Zeta Function}

\begin{definition}\label{hyzeta}

 Following the definition of the classical Riemann zeta function, naturally define, \[\zeta
_{_{^{\ast} \mathbb{Q}}}(s)=\sum_{\text{proper internal non-zero
ideals I of $\mathbb{^{\ast}Z}$}} \frac{1}{N(I)^{s}}=\sum_{n\in
\mathbb{^{\ast}Z_{> \text{0}}}} \frac{1}{n^{s}}.\]
\end{definition}

Let $\zeta_{_{N}}(s)=\sum_{0<n\leq N} \frac{1}{n^{s}}$ $(N\in
\text{} ^{\ast} \mathbb{Z},s\in \text{} ^{\ast}\mathbb{C})$. From
classical results, $\zeta_{_{\mathbb{Q}}}(s)$ is absolutely
convergent for $\Re (s)>1+\delta$ $(\delta\in\mathbb{R_{>\text{0}}})$ and the functions $\zeta_{_{N}}(s)$
$(N\in \mathbb{Z},s\in \mathbb{C})$ converge uniformly for $\Re
(s)>1+\delta$. The aim of this section is to develop similar results for
$\zeta _{_{^{\ast} \mathbb{Q}}}(s)$.

\begin{lemma} \label{L5}
$\{\zeta_{_{N}}(s)\} $ are Q-continuous functions $\forall n\in
\text{} ^{\ast} \mathbb{N}$, $^{\ast}\Re(s)>1$.
\end{lemma}
\noindent\emph{\textbf{Proof:}}
\begin{equation}
 (\forall n\in \mathbb{Z_{>\text{0}}}) (\forall
\epsilon \in \mathbb{R_{> \text{0}}})(\exists \delta \in
\mathbb{R_{> \text{0}}})(\forall s,t\in \Re (s)>1,\text{
}d(s,t)<\delta)(d(\zeta_{_{N}}(s),\zeta_{_{N}}(t)))< \epsilon ).
\end{equation}
\begin{displaymath}
\downarrow \text{$\ast$-transform,}
\end{displaymath}
\begin{equation}
(\forall n\in \text{}^{\ast}\mathbb{Z_{\text{0}}}) (\forall \epsilon \in
\text{}^{\ast}\mathbb{R_{> \text{0}}})(\exists \delta \in \text{}
^{\ast} \mathbb{R_{> \text{0}}})(\forall s,t\in \text{} ^{\ast}\Re
(s)>1,\text{ }
^{\ast}d(s,t)<\delta)(^{\ast}d(\zeta_{_{N}}(s),\zeta_{_{N}}(t)))<
\epsilon ).
\end{equation}
\begin{flushright} $\Box$ \end{flushright}

\begin{prop} \label{L4}
$\zeta_{_{N}}(s)$ $(N\in \text{}
^{\ast}\mathbb{Z_{>\text{0}}})$Q-converges uniformly to $\zeta
_{_{^{\ast} \mathbb{Q}}}(s)$ for $^{\ast}\Re (s)>1+\delta$, for every $\delta\in\mathbb{^{\ast}R_{>\text{0}}}$ .
\end{prop}
\noindent\emph{\textbf{Proof:}}
\begin{equation}
(\forall \epsilon \in \mathbb{R_{> \text{0}}})(\exists N\in
\mathbb{Z_{>\text{0}}})(\exists \delta\in \mathbb{R_{>\text{0}}}) (\forall
m,n\in \mathbb{N_{\geq \text{N}}})(\forall s\in\mathbb{C}, \Re (s)
>1+\delta)(d(\zeta_{_{m}}(s),\zeta_{_{n}}(s)))< \epsilon ).
\end{equation}
\begin{displaymath}
\downarrow \text{$\ast$-transform,}
\end{displaymath}
\begin{equation}
(\forall \epsilon \in \text{} ^{\ast} \mathbb{R_{>
\text{0}}})(\exists N\in \text{}^{\ast}\mathbb{Z_{>\text{0}}})(\exists
\delta\in \mathbb{^{\ast}R_{>\text{0}}})(\forall m,n\in \text{}
^{\ast} \mathbb{N_{\geq \text{M}}})(\forall s\in \text{}
^{\ast}\mathbb{C},^{\ast} \Re (s)
>1+\delta)(^{\ast}d(\zeta_{_{m}}(s),\zeta_{_{n}}(s)))< \epsilon ).
\end{equation}
The $\{\zeta_{_{N}}(s)\} $ form a hyper internal Cauchy sequence and a limit
function exists. Fixing $n$ and letting $m\rightarrow \infty$
shows that the limit function is $\zeta_{_{^{\ast}\mathbb{Q}}}(s)$. The
uniform Q-convergence also implies absolute Q-convergence of
$\zeta_{_{^{\ast}\mathbb{Q}}}(s)$ for $^{\ast}\Re (s)>1+\delta$.
\begin{flushright} $\Box$ \end{flushright}

\begin{corollary} \label{L7}
$\zeta_{_{^{\ast}\mathbb{Q}}}(s)$ is Q-continuous for $^{\ast} \Re (s)>1$.
\end{corollary}

\begin{prop} \label{L8}
For $^{\ast} \Re (s)>1$,
\[\zeta_{_{^{\ast}\mathbb{Q}}}(s)=\prod_{\text{Non-zero internal prime ideals p
of} ^{\ast}\mathbb{Z}} (1-[^{\ast}\mathbb{Z}:p]^{-s})^{-1}
=\prod_{p\in \text{}^{\ast} \Lambda}(1-p^{-s})^{-1} \] .
\end{prop}
\noindent\emph{\textbf{Proof:}} Let $M\in \text{}
^{\ast}\mathbb{Z}$ and let $\Omega_{M}$ be the set of all primes
$\leq M$. So,
\begin{equation}
 \prod_{p\in \Omega_{M}} (1-p^{-s})^{-1}=\prod_{p\in \Omega_{M}}
 (1+p^{-s}+p^{-2s}+\ldots)=\sum_{m\in \mathcal{M}}
 \frac{1}{m^{s}}.
\end{equation}
In the above sum $\mathcal{M}$ is the set of all $m\in \text{}
^{\ast} \mathbb{Z_{>\text{0}}}$ such that the prime factors of
$m$ are elements of $\Omega_{M}$. Let
$\xi_{M}(s)=\sum_{m\in\mathcal{M}}\frac{1}{m^{s}}$. Classically
$\xi_{M}(s)$ converge absolutely for $\Re (s)>1+\delta$
($\delta\in\mathbb{R_{\text{0}}})$ to
$\zeta_{_{^{\ast}\mathbb{Q}}}(s)$. So by the transfer principle,
\[ (\forall \delta\in\mathbb{R_{>\text{0}}}) (\forall M,N \in
\mathbb{Z_{>\text{0}}})(\exists \epsilon\in
\mathbb{R_{>\text{0}}})(\forall
s\in\mathbb{C},\Re(s)>1+\delta)(d(\xi_{M}(s),\xi_{N}(s))<\epsilon
/2)\] \begin{displaymath} \downarrow \text{$\ast$-transform,}
\end{displaymath} \[  (\forall
\delta\in\mathbb{^{\ast}R_{>\text{0}}}) (\forall M,N \in
\mathbb{^{\ast}Z_{>\text{0}}})(\exists \epsilon\in
\mathbb{^{\ast}R_{>\text{0}}})(\forall
s\in\mathbb{^{\ast}C},^{\ast}\Re(s)>1+\delta)(^{\ast}d(\xi_{M}(s),\xi_{N}(s))<\epsilon
/2 )\] Hence $(\forall M' \in
\mathbb{^{\ast}Z_{>\text{0}}})(\exists \epsilon '\in
\mathbb{^{\ast}R_{>\text{0}}})(\forall
s\in\mathbb{^{\ast}C},^{\ast}\Re(s)>1+\delta)(^{\ast}d(\xi_{M'}(s),\zeta_{_{^{\ast}\mathbb{Q}}}(s))<\epsilon
' /2 )$ and by \ref{L4}, $(\forall M \in
\mathbb{^{\ast}Z_{>\text{0}}})(\exists \epsilon \in
\mathbb{^{\ast}R_{>\text{0}}})(\forall
s\in\mathbb{^{\ast}C},^{\ast}\Re(s)>1+\delta)(^{\ast}d(\zeta_{_{^{\ast}\mathbb{Q}}}(s),\zeta_{M}(s))<\epsilon
 /2 )$. Let $\delta/2=\max \{\epsilon,\epsilon '\}$ and $N=\max
\{M,M'\}$, then $\forall n>N$ and using the triangle inequality,
\begin{align} ^{\ast} d\left( \sum_{n\in
\mathbb{^{\ast}Z_{>\text{0}}}} \frac{1}{n^{s}},\sum_{n\leq
N}\frac{1}{n^{s}}\right) &=^{\ast}
d\left(\xi_{N'}(s),\zeta_{N}(s) \right),\notag\\ &<
^{\ast}d(\xi_{N'}(s),\zeta_{_{^{\ast}\mathbb{Q}}}(s))+d(\zeta_{_{^{\ast}\mathbb{Q}}}(s),\zeta_{N}(s))
,\notag \\ &<\delta /2 +\delta /2, \notag\\ &=\delta.\notag
\end{align}
Hence the result as the two sequences converge to the same limit
function.
\begin{flushright} $\Box$ \end{flushright}

\begin{lemma} \label{L9}
$\zeta_{_{^{\ast}\mathbb{Q}}}(s)$ has no zeros for $^{\ast} \Re (s)>1$.
\end{lemma}
\noindent\emph{\textbf{Proof:}} In the standard case
$\zeta_{_{Q}}(s)$ for $\Re (s)>1$.
\[ (\forall N\in\mathbb{Z_{>\text{0}}}) (\exists
\epsilon \in\mathbb{R_{>\text{0}}})(\forall
s\in\mathbb{C},\Re(s)>1)(d(\xi_{N}(s),0)>\epsilon),
\]
\begin{displaymath}
\downarrow \text{$\ast$-transform,}
\end{displaymath}
\[ (\forall N\in\mathbb{^{\ast}Z_{>\text{0}}})(\exists
\epsilon \in\mathbb{^{\ast}R_{>\text{0}}})(\forall
s\in\mathbb{^{\ast}C},^{\ast}\Re(s)>1)(^{\ast}d(\xi_{N}(s),0)>\epsilon),
\]

\begin{flushright} $\Box$ \end{flushright}

   \subsection{A Hyper Theta Function}
\begin{definition}\label{hytheta}
\textbf{Hyper Theta Function:} For $^{\ast}\Im (s)>0, $ \[ ^{\ast}\Theta(s)=\sum_{n\in
\text{} ^{\ast}\mathbb{Z}} \text{}^{\ast}e^{\pi in^{2}s}=1
+2\sum_{n\in \text{} ^{\ast}\mathbb{Z_{>\text{0}}}}
\text{}^{\ast}e^{\pi in^{2}s}.\]
\end{definition}

Classically the theta function, $\Theta(s)=\sum_{n\in\mathbb{Z}}
e^{\pi in^{2}s}$, is absolutely convergent $\forall
s\in\mathbb{C}$ with $\Im(s)>0$.

\begin{prop}\label{thconv}
$^{\ast}\Theta(s)$ is absolutely Q-convergent $\forall
s\in\mathbb{^{\ast}C}$ with $^{\ast}\Im(s)>0$.
\end{prop}
\noindent\emph{\textbf{Proof:}} Let
$N\in\mathbb{^{\ast}Z_{>\text{0}}}$ and define
$^{(N)}\theta(s)=1+2\sum_{n=1}^{N} \text{}^{\ast}e^{\pi in^{2}s}$
for $^{\ast}\Im(s)>0$.
\[(\forall m,n\in \mathbb{Z_{>\text{0}}})(\exists \epsilon
\in\mathbb{R_{>\text{0}}}) (\forall
s\in\mathbb{C},\Im(s)>0)(d(^{(n)}\Theta(s),^{(m)}\Theta(s))<\epsilon),\]
\begin{displaymath}
\downarrow \text{$\ast$-transform,}
\end{displaymath}
\[(\forall m,n\in
\mathbb{^{\ast}Z_{>\text{0}}})(\exists \epsilon
\in\mathbb{^{\ast}R_{>\text{0}}})(\forall
s\in\mathbb{^{\ast}C},^{\ast}\Im(s)>0)(^{\ast}d(^{(n)}\Theta(s),^{(m)}\Theta(s))<\epsilon),\]
\begin{flushright} $\Box$ \end{flushright}
\subsubsection{Hyper Fourier Transform}
\begin{definition}\label{hyss}
\textbf{Hyper Schwartz Space:} Let $C^{\infty}(\mathbb{^{\ast}R})$
be the space of internal Q-smooth functions of
$\mathbb{^{\ast}R}$ and define the hyper Schwartz space to be:
\[\mathcal{S}=\left\{ f\in C^{\infty} (\mathbb{^{\ast}R}) :
s\in\mathbb{^{\ast}R}, lim_{^{\ast}d(s,0)\rightarrow\infty}
s^{n}\frac{d^{m}f}{ds^{m}}=0 \\ \forall
n,m\in\mathbb{^{\ast}Z_{>\text{0}}} \right\}. \]
\end{definition}

\begin{definition}\label{hyfour}
\textbf{Hyper Fourier Transform:} Let $f\in\mathcal{S}$ and define
its Fourier transform to be,
\begin{equation}
\hat{f}(y)=\int_{^{\ast}\mathbb{R}}^{\ast} f(x) ^{\ast}e^{-2\pi
ixy} dx.
\end{equation}
\end{definition}

In order to examine properties of this transform the idea of hyper
distributions need to be introduced.
\begin{definition}\label{hydd}
\textbf{Hyper Dirac Distribution:}\[
^{\ast}\delta(x)=\lim_{\epsilon\rightarrow 0}
\Delta(x;\epsilon),\] where $x,\epsilon \in\mathbb{^{\ast}R}$ and
$\Delta(x;\epsilon)=0 (^{\ast}d(x,0)<\epsilon),
\frac{1}{2\epsilon} (^{\ast}d(x,0)>\epsilon)$.
\end{definition}
Using this definition the following properties can be derived.
\begin{prop}\label{fourier}
\begin{enumerate}
\item \[\int_{\mathbb{^{\ast}R}}^{\ast} \text{}
^{\ast}\delta(x-x_{0})dx=lim_{\epsilon\rightarrow
0}\int_{^{\ast}[x_{0}-\epsilon,x_{0}+\epsilon]}^{\ast}
\text{}^{\ast}\delta(x-x_{0})dx=1.\]
\item Let $f$ be a Q-continuous hyper real function. Then \[\int_{\mathbb{^{\ast}R}}^{\ast} \text{}
^{\ast}\delta(x-x_{0})f(x)dx=f(x_{0}).\]
\item \[ ^{\ast}\hat{\delta}(y)=\int_{\mathbb{^{\ast}R}}^{\ast} \text{}
^{\ast}\delta(x)\text{} ^{\ast}e^{-2\pi ixy}dx=1. \]
\item \[ ^{\ast}\delta(x)=\int_{\mathbb{^{\ast}R}}^{\ast} \text{}
^{\ast}\delta(x)\text{} ^{\ast}e^{2\pi ixy}dy. \]
\end{enumerate}
\end{prop}

\begin{prop} \label{L50}
\[ \int_{\mathbb{^{\ast}R}}^{\ast}
\int_{\mathbb{^{\ast}R}}^{\ast} f(x) \text{} ^{\ast}e^{-2\pi
iz(x+y)} dz dx
=\int_{\mathbb{^{\ast}R}}^{\ast} dx f(x) \left( dz
\int_{\mathbb{^{\ast}R}}^{\ast} \text{}
^{\ast}e^{-2\pi iz(x+y)} \right).\]
\end{prop}
\noindent\emph{\textbf{Proof:}}
In the standard case this result is true. So let $I(n,r)= \int_{\mathbb{^{\ast}[\text{-r,r}]}}^{\ast}
\int_{\mathbb{^{\ast}[\text{-r,r}]}}^{\ast} f(x) \text{} ^{(n)}e^{-2\pi
iz(x+y)} dz dx$ and let $J(n,r)=  \int_{\mathbb{^{\ast}[\text{-r,r}]}}^{\ast} dx f(x) \left( dz
\int_{\mathbb{^{\ast}[\text{-r,r}]}}^{\ast} \text{}
^{(n)}e^{-2\pi iz(x+y)} \right)$ for $r,n\in\mathbb{^{\ast}Z_{>\text{0}}}$.
 
\[ (\forall \epsilon\in\mathbb{R_{>\text{0}}})(\exists N,R\in\mathbb{Z_{>\text{0}}})
(\forall n,m\in\mathbb{Z_{>\text{N}}})(\forall s,t\in\mathbb{Z_{>\text{R}}})
(d(I(n,s),J(m,t))<\epsilon),\]
\begin{displaymath}
\downarrow \text{$\ast$-transform,}
\end{displaymath}
\[ (\forall \epsilon\in\mathbb{^{\ast}R_{>\text{0}}})(\exists N,R\in\mathbb{^{\ast}Z_{>\text{0}}})
(\forall n,m\in\mathbb{^{\ast}Z_{>\text{N}}})(\forall s,t\in\mathbb{^{\ast}Z_{>\text{R}}})
(^{\ast}d(I(n,s),J(m,t))<\epsilon). \]
\begin{flushright} $\Box$ \end{flushright}

\begin{prop}\label{inv}
Let $f\in\mathcal{S}$ then $\hat{\hat{f}}(y)=f(-y)$.
\end{prop}
\noindent\emph{\textbf{Proof:}} Using the results of the previous
two propositions:
\begin{align}
\hat{\hat{f}}(y) &= \int_{\mathbb{^{\ast}R}}^{\ast}
\int_{\mathbb{^{\ast}R}}^{\ast} f(x) \text{} ^{\ast}e^{-2\pi
iz(x+y)} dz dx, \notag \\
&=\int_{\mathbb{^{\ast}R}}^{\ast} dx f(x) \left( dz
\int_{\mathbb{^{\ast}R}}^{\ast} \text{}
^{\ast}e^{-2\pi iz(x+y)} \right), \notag\\
&=\int_{\mathbb{^{\ast}R}}^{\ast}
dx f(x) \text{} ^{\ast}\delta(x-(-y)), \notag\\
&=f(-y). \notag
\end{align}
\begin{flushright} $\Box$ \end{flushright}

\subsubsection{Hyper Poisson Summation}

Recall Fej\'{e}r's fundamental theorem concerning Fourier series.
Let $f(x)$ be a bounded, measurable and periodic of period 1. Then
the Fourier coefficients of f are given by, \[
c_{n}=\int_{0}^{1}f(x)e^{-2\pi inx}dx,\] for each
$n\in\mathbb{Z}$. The partial sums are defined as
\[S_{N}(x)=\sum_{d(n,0)\leq (N)}c_{n}e^{2\pi inx}. \] When $f(x)$
is continuous and $\sum_{n\in\mathbb{Z}} |c_{n}| <\infty$ then the
function is represented by the absolutely convergent Fourier
series \[ f(x)=\sum_{\mathbb{Z}}c_{n}e^{2\pi inx}.\]

Let $g\in\mathcal{S}$ and periodic of period 1 ($g(x+1)=g(x),
x\in\mathbb{R}$). $f$ is a real valued function and since $g\in\mathcal{S}$, $f$
is a bounded, measurable function of period 1. By the transfer
principle the Fourier coefficients $c_{n}$ and partial sums
$S_{N}$, defined above, can be defined $\forall
n\in\mathbb{^{\ast}Z}$. So, \[ (\forall m,n\in\mathbb{Z})(\exists
\epsilon\in\mathbb{R_{>\text{0}}})(d(S_{N},S_{M})<\epsilon), \]
\begin{displaymath}
\downarrow \text{$\ast$-transform,}
\end{displaymath}
\[ (\forall m,n\in\mathbb{^{\ast}Z})(\exists
\epsilon\in\mathbb{^{\ast}R_{>\text{0}}})(^{\ast}d(S_{N},S_{M})<\epsilon).
\]
Hence $g(x)$ is represented by the absolutely Q-convergent hyper
Fourier series $g(x)=\sum_{n\in\mathbb{^{\ast}Z}}c_{n} \text{}
^{\ast}e^{2\pi inx}. $

\begin{lemma}\label{sumprod}
$ \int_{^{\ast}[\text{0,1}]}^{^{\ast}} \sum_{k\in
\mathbb{^{\ast}Z}} f(x+k) ^{\ast}e^{-2\pi inx}dx
=\sum_{k\in \mathbb{^{\ast}Z}} \int_{^{\ast}[\text{0,1}]}^{^{\ast}} f(x+k) ^{\ast}e^{-2\pi inx}dx.$
\end{lemma}
\noindent\emph{\textbf{Proof:}}
Let \[ I(n,r)= \int_{^{\ast}[\text{1/r,1-1/r}]}^{^{\ast}} \sum_{k\in
\mathbb{^{\ast}Z}} f(x+k) ^{(n)}e^{-2\pi inx}dx\] and \[ J(n,r)=\sum_{k\in \mathbb{^{\ast}Z}} \int_{^{\ast}[\text{1/r,1-1/r}]}^{^{\ast}} f(x+k) ^{(n)}e^{-2\pi inx}dx,\]
for $n,r\in\mathbb{^{\ast}Z_{>\text{0}}}$. Then use the argument of \ref{L50}.
\begin{flushright} $\Box$ \end{flushright}

\begin{lemma} \label{L18}
Hyper Poisson Summation: Let $f\in\mathcal{S}$ then \[\sum_{n\in
\mathbb{^{\ast}Z}} f(n)=\sum_{n\in \mathbb{^{\ast}Z}} \hat{f}(n).\]
\end{lemma}

\noindent\emph{\textbf{Proof:}} Let $g(x)=\sum_{k\in
\mathbb{^{\ast}}} f(k)$ then $g(x)=g(x+1)$ where
$f\in\mathcal{S}$. This periodicity enables $g(x)$ to be written
as $g(x)=\sum_{n\in \mathbb{^{\ast}Z}} a_{n} ^{\ast}e^{2\pi inx}$
where $a_{n}=\int_{^{\ast}[0,1]}^{\ast} g(x) ^{\ast} e^{-2\pi
inx}dx.$ Using the result of the previous lemma.
\begin{align}
a_{n} &= \int_{^{\ast}[0,1]}^{^{\ast}} \sum_{k\in
\mathbb{^{\ast}Z}} f(x+k) ^{\ast}e^{-2\pi inx}dx, \notag\\
&=\sum_{k\in \mathbb{^{\ast}Z}} \int_{^{\ast}[0,1]}^{^{\ast}} f(x+k) ^{\ast}e^{-2\pi inx}dx, \notag\\
&=\hat{f}(n) \notag.
\end{align}
\[ \sum_{n\in \mathbb{^{\ast}Z}} f(n)=g(0)=\sum_{n\in
\mathbb{^{\ast}Z}} a_{n}=\sum_{n\in \mathbb{^{\ast}Z}}
\hat{f}(n).\]
\begin{flushright} $\Box$ \end{flushright}

\begin{lemma} \label{L15}
\[ \int_{^{\ast}\mathbb{R}}^{\ast} \left( \frac{d}{dx}
\text{} ^{\ast}e^{-\pi x^{2}-2\pi ixy} \right) dx =0. \]
\end{lemma}
\noindent\emph{\textbf{Proof:}} Standard result:
$\int_{\mathbb{R}} \frac{d}{dx} (e^{-\pi x^{2}-2\pi ixy}) dx=0)$.
Using a similar method to \ref{L13} let,
\[ I_{n,r}=\int_{[ -r,r]} \left( \frac{d}{dx}
\text{} ^{(n)}e^{-\pi x^{2}-2\pi ixy} \right) dx.\]
\begin{equation}
(\forall \epsilon \in \mathbb{R_{>\text{0}}})(\exists r\in
\mathbb{R_{>\text{0}}})(\exists n\in
\mathbb{N})(d(I_{n,r},0)<\epsilon),
\end{equation}
\begin{displaymath}
\downarrow \text{$\ast$-transform,}
\end{displaymath}
\begin{equation}
(\forall \epsilon \in
\text{}\mathbb{^{\ast}R_{>\text{0}}})(\exists r\in
\text{}^{\ast}\mathbb{R_{>\text{0}}})(\exists n\in
\text{}^{\ast}\mathbb{N})(^{\ast}d(I_{n,r},0)<\epsilon),
\end{equation}
\begin{flushright} $\Box$ \end{flushright}

\begin{lemma}\label{L16}
\[ \int_{^{\ast}\mathbb{R}}^{\ast}
 \text{}^{\ast}e^{-\pi x^{2}} dx=1. \]
\end{lemma}
\noindent\emph{\textbf{Proof:}} In the standard case
$\int_{\mathbb{R}} e^{-\pi x^{2}} dx=1$. Let
$I_{n,r}=\int_{[-r,r]} \text{} ^{\ast}e^{-\pi x^{2}}dx$.
\begin{equation}
(\forall \epsilon \in \mathbb{R_{>\text{0}}})(\exists r\in
\mathbb{R_{>\text{0}}})(\exists n\in
\mathbb{N})(d(I_{n,r},0)<\epsilon),
\end{equation}
\begin{displaymath}
\downarrow \text{$\ast$-transform,}
\end{displaymath}
\begin{equation}
(\forall \epsilon \in
\text{}\mathbb{^{\ast}R_{>\text{0}}})(\exists r\in
\text{}^{\ast}\mathbb{R_{>\text{0}}})(\exists n\in
\text{}^{\ast}\mathbb{N})(^{\ast}d(I_{n,r},1)<\epsilon),
\end{equation}
\begin{flushright} $\Box$ \end{flushright}

\begin{lemma} \label{L17}
Let $h(x)=^{\ast}e^{-\pi x^{2}}$ then $h(x)$ is its own Fourier
transform.
\end{lemma}
\noindent\emph{\textbf{Proof:}}
$\hat{h}(y)=\int_{^{\ast}\mathbb{R}}^{\ast} h(x) ^{\ast}e^{-2\pi
ixy} dx.$ Differentiating,
\begin{align}
\frac{d\hat{h}(y)}{dy}&= -2\pi i\int_{^{\ast}\mathbb{R}}^{\ast}
xh(x) ^{\ast}e^{-2\pi ixy} dx, \notag \\
&= -2\pi i\hat{h}(y), \notag
\end{align}
Integrating by parts and using \ref{L15}. Thus,
$\hat{h}(y)=C^{\ast}e^{-\pi y^{2}}$ for some constant C. Letting
$y=0$ gives $C=1$ (using \ref{L16}) and the result.
\begin{flushright} $\Box$ \end{flushright}

\begin{theorem} \label{th1}
Let $s\in\mathbb{^{\ast}R_{>\text{0}}}$ then,
\[ \sum_{n\in\mathbb{^{\ast}Z}} \text{} ^{\ast}e^{-n^{2}\pi /s} =
s^{1/2}\sum_{n\in\mathbb{^{\ast}Z}} \text{} ^{\ast}e^{-n^{2}\pi
s}.
\]
\end{theorem}
\noindent\emph{\textbf{Proof:}} Using hyper Poisson summation,
\ref{L16} and \ref{L17},
\begin{align}
\sum_{n\in\mathbb{^{\ast}Z}} \text{}^{\ast}e^{-\frac{\pi
n^{2}}{s}} &= \sum_{n\in\mathbb{^{\ast}Z}} \int_{^{\ast}R}^{\ast}
\text{} e^{-2\pi int} \text{}^{\ast}e^{-\frac{\pi
t^{2}}{s}} dt, \notag \\
&= s\sum_{n\in\mathbb{^{\ast}Z}} \int_{^{\ast}R}^{\ast} \text{}
e^{-2\pi
insu-\pi su^{2}} du, \notag \\
&=s\sum_{n\in\mathbb{^{\ast}Z}} \int_{^{\ast}R}^{\ast} \text{}
e^{-\pi
s[(u+in)^{2}+n^{2}]} du, \notag \\
&=s\sum_{n\in\mathbb{^{\ast}Z}} \text{}^{\ast} e^{-\pi sn^{2}}
\int_{^{\ast}R}^{\ast} \text{} e^{-\pi sv^{2}}dv,\notag\\
&=s^{1/2}\sum_{n\in\mathbb{^{\ast}Z}} \text{}^{\ast} e^{-\pi
sn^{2}} \notag.
\end{align}
\begin{flushright} $\Box$ \end{flushright}

\begin{corollary} \label{co1}
Functional equation for the hyper theta function:
\[^{\ast}\Theta(-\frac{1}{z})=(z/i)^{1/2}\text{}^{\ast}\Theta(z).
\]
\end{corollary}
\noindent\emph{\textbf{Proof:}} Let $z=is$ and set $\omega (s)=
\text{} ^{\ast} \Theta (is)$. Then \ref{th1} gives the functional
equation $\omega (1/x)=x^{1/2}\omega (x)$.

\begin{flushright} $\Box$ \end{flushright}

\section{Analytic Properties}

Classically $\zeta _{_{ \mathbb{Q}}}(s)$ is an
analytic function of s with $\Re (s)>1$. It has analytic
continuation to the whole of $\mathbb{C}$ with a simple pole at
$s=1$.

\begin{lemma} \label{L10}
$\zeta _{_{ N}}(s)$ is Q-analytic $\forall N\in \text{} ^{\ast}
\mathbb{N}$ and $^{\ast} \Re (s)>1$.
\end{lemma}
\noindent\emph{\textbf{Proof:}}
\[
(\forall \epsilon \in \mathbb{R_{>\text{0}}})(\exists \delta \in
\mathbb{R_{>\text{0}}})(\forall s\in \mathbb{C},\Re (s)>1
)(\exists h,g\in \mathbb{C})\]
\[(d(\frac{\zeta _{_{
N}}(s+h)-\zeta _{_{ N}}(s)}{h},\frac{\zeta _{_{
N}}(s+g)-\zeta _{_{ N}}(s)}{g})<
\epsilon \wedge d(h,0)< \delta \wedge d(g,0)< \delta),
\]
\begin{displaymath}
\downarrow \text{$\ast$-transform,}
\end{displaymath}
\[
(\forall \epsilon \in \text{} ^{\ast}
\mathbb{R_{>\text{0}}})(\exists \delta \in \text{}
^{\ast}\mathbb{R_{>\text{0}}})(\forall s\in \text{}
^{\ast}\mathbb{C},\Re (s)>1 )(\exists h,g\in \text{}
^{\ast}\mathbb{C} )\]
\[(^{\ast}d(\frac{\zeta
_{_{ N}}(s+h)-\zeta _{_{ N}}(s)}{h},\frac{\zeta _{_{
N}}(s+g)-\zeta _{_{ N}}(s)}{g})<
\epsilon \wedge ^{\ast}d(h,0)< \delta \wedge ^{\ast}d(g,0)< \delta).
\]
\begin{flushright} $\Box$ \end{flushright}
From standard complex analysis a theorem of Weierstrass' enables
the deduction that the limit function of a uniformly convergent
sequence of analytic functions is analytic. The proof of this
relies on Cauchy's theorem.
\begin{theorem} \label{T1}
$f(z)$ is Q-analytic in a region $\Omega \subseteq ^{\ast}
\mathbb{C}$ iff
\begin{equation}
\int_{\gamma}^{\ast} f(z) dz = 0,
\end{equation}
for every cycle, $\gamma$, which is homologous to zero in
$\Omega$.
\end{theorem}
\noindent\emph{\textbf{Proof:}} See section 6.2.3 of \cite{Ro2}.
\begin{flushright} $\Box$ \end{flushright}

\begin{theorem} \label{T2}
\begin{enumerate}
\item Let f be Q-analytic on and inside and on a cycle $\gamma$. Then, if $a$ is inside $\gamma$,
 \[ f(a)=\frac{1}{2\pi i}\int_{\gamma }^{\ast} \frac{f(z)}{z-a}dz .\]
 \item With the same conditions, \[ f^{(n)}(a)=\frac{n!}{2\pi
 i}\int_{\gamma}^{\ast} \frac{f(z)}{(z-a)^{n+1}}dz .\]
 \end{enumerate}
\end{theorem}
\noindent\emph{\textbf{Proof:}} See section 6.2.3 of \cite{Ro1}.
\begin{flushright} $\Box$ \end{flushright}

\begin{theorem}\label{hyqanl}
Suppose $f_{n}(z)$ is analytic in $\Omega_{n} \subseteq
\text{}^{\ast} \mathbb{C}$ and ${f_{n}(z)}$ Q-converges to a limit
function $f(z)$ in a region $\Omega \subseteq \text{}^{\ast}
\mathbb{C}$, Q-uniformly on every hypercompact set (with respect to the
Q-topology) of $\Omega $. Then $f(z)$ is Q-analytic in $\Omega $.
\end{theorem}
\noindent\emph{\textbf{Proof:}} By hyper Cauchy's theorem:
\begin{equation}
f_{n}(z)=\frac{1}{2\pi i} \int_{C}^{\ast} \frac{f_{n}(w)}{w-z}dw,
\end{equation}
where C is a hyper disc $^{\ast} d(w,a)\leq r$ contained in $\Omega$. In the limit
$n\rightarrow \infty$ and by uniform Q-convergence,
\begin{equation}
f(z)=\frac{1}{2\pi i} \int_{C}^{\ast} \frac{f(w)}{w-z}dw,
\end{equation}
and so $f(z)$ is Q-analytic in the disc. Any hypercompact (with respect
to the Q-topology) subset of $\Omega$ can be covered by a hyperfinite
number of such closed discs and therefore convergence is Q-uniform
on every hypercompact (with respect to the Q-topology) subset.
\begin{flushright} $\Box$ \end{flushright}

\begin{corollary} \label{L11}
$\zeta _{_{^{\ast} \mathbb{Q}}}(s)$ is a Q-analytic function for
$^{\ast} \Re (s) >1$.
\end{corollary}

\begin{lemma}\label{L12}
$\zeta_{N}^{\prime }(s)$ Q-converges uniformly to $\zeta
_{_{^{\ast} \mathbb{Q}}}^{\prime}(s)$ for $^{\ast} \Re (s) >1$.
\end{lemma}
\noindent\emph{\textbf{Proof:}} Identical statement as in \ref{L4}
with $\zeta_{N}(s)$ replaced by $\zeta_{N}^{\prime }(s)$.
\begin{flushright} $\Box$ \end{flushright}

\begin{definition}\label{zreal}
\begin{enumerate}
\item $^{\ast}\zeta _{_{\eta}}(s)=\pi
^{-s/2}\text{} ^{\ast}\Gamma(\frac{s}{2})$.
\item $\zeta _{_{^{\ast}\mathbb{A}}}(s)=^{\ast}\zeta _{_{\eta}}(s)
\text{} ^{\ast}\zeta _{_{^{\ast}\mathbb{Q}}}(s)$.
\end{enumerate}
\end{definition}

\begin{prop}\label{zreal1}
\begin{enumerate}
\item $^{\ast}\zeta _{_{\eta}}(s)$ is a Q-analytic function with
poles at $-2\mathbb{^{\ast}Z_{\geq\text{0}}}$.
\item $\zeta _{_{^{\ast}\mathbb{A}}}(s)$ Q-converges for
$^{\ast}\Re(s)>1$.
\end{enumerate}
\end{prop}
\noindent\emph{\textbf{Proof:}}
\begin{enumerate}
\item The only poles of $^{\ast}\zeta _{_{\eta}}(s)$ are those of
$^{\ast}\Gamma(s/2)$.
\item By composition of functions the Q-convergence is dependent
on the Q-convergence of $^{\ast}\zeta _{_{^{\ast}\mathbb{Q}}}(s)$ since
$^{\ast}\zeta _{_{\eta}}(s)$ is Q-convergent on
$^{\ast}\mathbb{C}$ apart from at its poles.
\end{enumerate}
\begin{flushright} $\Box$ \end{flushright}

\begin{lemma}\label{swap}
\[\sum_{n\in \text{}
^{\ast}\mathbb{Z_{>\text{0}}}} \int_{^{\ast}
\mathbb{R^{+}}}^{\ast} \text{} ^{\ast}e^{-\pi n^{2}y} y^{s-1} dy
=\int_{^{\ast} \mathbb{R^{+}}}^{\ast} \left( \sum_{n\in \text{}
^{\ast}\mathbb{Z_{>\text{0}}}} \text{} ^{\ast}e^{-\pi n^{2}y}
\right) y^{s-1} dy.\]
\end{lemma}
\noindent\emph{\textbf{Proof:}}
Use the argument of \ref{L50}.
\begin{flushright} $\Box$ \end{flushright}

\begin{lemma} \label{L14}
\[ \zeta
_{_{^{\ast}\mathbb{A}}}(s)=\frac{1}{2}\int_{^{\ast}\mathbb{R^{+}}}^{\ast}
(^{\ast}\Theta(iy) -1)y^{s/2 -1} dy \text{\; \; \; \; \;     $(\Re (s)> 0)$.}
\]
\end{lemma}
\noindent\emph{\textbf{Proof:}}
\begin{align}
^{\ast}\Gamma(s) &= \int_{^{\ast}\mathbb{R^{+}}}^{\ast} \text{}
^{\ast} e^{-y} y^{s-1} dy, \notag \\
\pi^{-s} \text{}^{\ast}\Gamma (s)\frac{1}{n^{2s}} &=\int_{^{\ast}
\mathbb{R^{+}}}^{\ast} \text{} ^{\ast}e^{-\pi n^{2}y} y^{s-1} dy
\text{\      (by letting y $\mapsto \pi n^{2}y $)}, \notag \\
\pi^{-s} \text{} ^{\ast}\Gamma (s)\zeta _{_{^{\ast}\mathbb{Q}}}(s)
&= \sum_{n\in \text{} ^{\ast}\mathbb{Z_{>\text{0}}}} \int_{^{\ast}
\mathbb{R^{+}}}^{\ast} \text{} ^{\ast}e^{-\pi n^{2}y} y^{s-1} dy \notag, \\
\zeta _{_{^{\ast}\mathbb{A}}}(2s)&=\sum_{n\in \text{}
^{\ast}\mathbb{Z_{>\text{0}}}} \int_{^{\ast}
\mathbb{R^{+}}}^{\ast} \text{} ^{\ast}e^{-\pi n^{2}y} y^{s-1} dy \notag,\\
&=\int_{^{\ast} \mathbb{R^{+}}}^{\ast} \left( \sum_{n\in \text{}
^{\ast}\mathbb{Z_{>\text{0}}}} \text{} ^{\ast}e^{-\pi n^{2}y}
\right) y^{s-1} dy, \text{    ($^{\ast}\Re(s)>0$ )}, \notag \\
&= \frac{1}{2}\int_{^{\ast}\mathbb{R^{+}}}^{\ast}
(^{\ast}\Theta(iy) -1)y^{s -1} dy\notag.
\end{align}
\begin{flushright} $\Box$ \end{flushright}

\begin{lemma} \label{l21}
For $s\neq 0$,
\[\int_{^{\ast}[\text{0,1}]}^{\ast} y^{s/2-1}dy=\frac{2}{s}.\]
Also, for $s\neq 1/2$, \[\int_{^{\ast}\mathbb{R^{+}\setminus
^{\ast}[\text{0,1}]}}^{\ast} y^{-1/2- s/2}dy =1/(s-1/2).\]
\end{lemma}
\noindent\emph{\textbf{Proof:}} Using the methods of \ref{L16} and
the corresponding standard results but with
$I_{n}=\int_{^{\ast}[\text{0,1-1/n}]}^{\ast} y^{s/2-1}dy$
$(n\in\mathbb{Z_{>\text{0}}})$ and
$J_{n}=\int_{^{\ast}[\text{1,n}]} y^{-1/2- s/2}dy
((n\in\mathbb{Z_{>\text{1}}})$ gives the results
\begin{flushright} $\Box$ \end{flushright}

\begin{theorem}\label{zetafunc}
$\zeta_{_{^{\ast}\mathbb{A}}}(s)=\zeta_{_{^{\ast}\mathbb{A}}}(1-s)$.
\end{theorem}
\noindent\emph{\textbf{Proof:}}
\begin{align}
\zeta
_{_{^{\ast}\mathbb{A}}}(s)&=\frac{1}{2}\int_{^{\ast}\mathbb{R^{+}}}^{\ast}
(^{\ast}\Theta(iy) -1)y^{s/2 -1} dy,\notag \\
&=\frac{1}{2}\int_{^{\ast}\mathbb{R^{+}\setminus
^{\ast}[\text{0,1}]}}^{\ast} (^{\ast}\Theta(iy) -1)y^{s/2 -1} dy
+\frac{1}{2}\int_{^{\ast}\mathbb{\text{[0,1]}}}^{\ast}
(^{\ast}\Theta(iy)
-1)y^{s/2 -1} dy, \notag\\
&=\frac{1}{2}\int_{^{\ast}\mathbb{R^{+}\setminus
^{\ast}[\text{0,1}]}}^{\ast} (^{\ast}\Theta(iy) -1)y^{s/2 -1} dy
+\frac{1}{2}\int_{^{\ast}\mathbb{R^{+}\setminus
^{\ast}[\text{0,1}]}}^{\ast} (^{\ast}\Theta(-\frac{1}{iy}))y^{-s/2
-1} dy-\frac{1}{2}\int_{^{\ast}[\text{0,1}]}^{\ast}
y^{s/2-1}dy,\notag \\
&=\frac{1}{2}\int_{^{\ast}\mathbb{R^{+}\setminus
^{\ast}[\text{0,1}]}}^{\ast} (^{\ast}\Theta(iy) -1)y^{s/2 -1} dy
+\frac{1}{2}\int_{^{\ast}\mathbb{R^{+}\setminus
^{\ast}[\text{0,1}]}}^{\ast} (y^{1/2}\text{}
^{\ast}\Theta(iy))y^{-s/2 -1} dy
-\frac{1}{s},\notag \\
&=-\frac{1}{s}+\frac{1}{2s-1}+\frac{1}{2}\int_{^{\ast}\mathbb{R^{+}\setminus
^{\ast}[\text{0,1}]}}^{\ast}
  [^{\ast}\Theta(iy) -1](y^{s/2-1}+y^{-1/2-s/2})dy
 \notag.
\end{align}
This is true for $^{\ast}\Re(s)>1$ however the integral of the
last line converges and absolutely $\forall
s\in\mathbb{^{\ast}C}$. The standard integral in last line
converges absolutely and uniformly $\forall s\in\mathbb{C}$. Using
the transfer principle this last integral does Q-converge
absolutely and uniformly for $s\in\mathbb{^{\ast}C}$. Indeed, for
$n\in\mathbb{Z_{>\text{0}}}$ and $r\in\mathbb{R_{>\text{1}}}$ let
\[ I_{n,r}=-\frac{1}{s}+\frac{1}{2s-1}+\frac{1}{2}\int_{^{\ast}[\text{0,r}]}^{\ast}
  [2\sum_{n\in\mathbb{Z_{>\text{0}}}}\text{} ^{(n)}e^{-\pi
  j^{2}y}](y^{s/2-1}+y^{-1/2-s/2})dy.\] Apply the transfer
  principle as in \ref{L13}.
Also the last integral is invariant under $s\mapsto 1-s$ by
substitution. The result follows.
\begin{flushright} $\Box$ \end{flushright}

\begin{corollary}\label{qcontz}
$\zeta _{_{^{\ast}\mathbb{A}}}(s)$ can be Q-analytically continued
on $\mathbb{^{\ast}C}$.
\end{corollary}
\begin{corollary}\label{qcontz1}
$\zeta _{_{^{\ast}\mathbb{Q}}}(s)$ can be Q-analytically continued
on $\mathbb{^{\ast}C}$ with a pole at s=1 and trivial zeros at
$s=-2\mathbb{^{\ast}Z_{>\text{0}}} $ .
\end{corollary}
\noindent\emph{\textbf{Proof:}} From the definition of $\zeta
_{_{^{\ast}\mathbb{A}}}(s)$, \[ \zeta
_{_{^{\ast}\mathbb{Q}}}(1-s)=\frac{\pi^{1/2-s} \text{}
^{\ast}\Gamma(\frac{s}{2})\zeta
_{_{^{\ast}\mathbb{Q}}}(s)}{^{\ast}\Gamma(\frac{1-s}{2})}. \]
Using \ref{pr1} the only poles of $\zeta
_{_{^{\ast}\mathbb{Q}}}(1-s)$ in the numerator are at
$s=0,-1,-2,\ldots$ and in the denominator $s=1,3,5,\ldots $.
Evaluating at $s=0$ in the equation corresponds to the pole $\zeta
_{_{^{\ast}\mathbb{Q}}}(1)$ since by the transfer principle $\zeta
_{_{\mathbb{Q}}}(s)$ has a pole at $s=1$. For $s=-1,-2,-3,\ldots $
in the equation leads to zeros at $s=-2,-4,-6,\ldots $ since the
left hand side is Q-analytic and non-zero at these values and the
poles of the hyper gamma function cancel with the zeros of the
zeta function. The other values above lead to the same conclusion. These are termed the trivial zeros.
\begin{flushright} $\Box$ \end{flushright}

It has been shown that apart from the trivial zeros above any other zero ($\gamma$)
must lie in the region $0\leq \text{}^{\ast}\Re (\gamma) \leq 1$. Following the standard case,
\begin{conjecture}\label{bigconj}
If $\zeta_{_{^{\ast}\mathbb{Q}}}(s)=0$ and $s$ is not a trivial zero then $^{\ast}\Re (s)=1/2$.
\end{conjecture}

\subsection{The Complex Shadow Map}

The real shadow map can be extended to act on the hyper complex
numbers. Recall the standard norm on $\mathbb{C}$ given for $z\in
\mathbb{C}$ by
$|z|_{\mathbb{C}}=((\Re(z))^{2}+(\Im(z))^{2})^{1/2}$. This then
extends to $\mathbb{^{\ast}C}$ by transfer. Let \[
\mathbb{^{\ast}C}^{\lim}=\{ z\in\mathbb{^{\ast}C} : \exists
r\in\mathbb{R} \text{ such that } |z|_{\mathbb{C}}<r \}.\] Then
$ \mathbb{^{\ast}C}^{\lim} = \mathbb{^{\ast}R}^{\lim} + i
\mathbb{^{\ast}C}^{\lim}$, where $i^{2}=-1$.

\begin{definition}[Complex Shadow Map] \label{complexsh}Let $z\in
\mathbb{^{\ast}C}^{\lim}$ with $x=\text{}^{\ast}\Re(z)$ and
$y=\text{}^{\ast}\Im(z)$.

\begin{align}
\operatorname{sh}_{\mathbb{^{\ast}C}}: \mathbb{^{\ast}C}^{\lim}
&\rightarrow \mathbb{C},\notag\\ x+iy & \mapsto
\operatorname{sh}_{\eta}(x) + i
\operatorname{sh}_{\eta}(y).
\end{align}
\end{definition}

\noindent All the properties of the real shadow map carry
through. As a simple application of this map consider the action
on $\zeta_{\mathbb{^{\ast}Q}}$.

\begin{lemma}\label{cx1}
For all $s\in\mathbb{^{\ast}C}^{\lim}$ with $^{\ast}\Re(s)>1$, \[
\zeta_{\mathbb{^{\ast}Q}}(s)\in\mathbb{^{\ast}C}^{\lim}.\]
\end{lemma}

\noindent\emph{\textbf{Proof:}} Let
$s\in\mathbb{^{\ast}C}^{\lim}$ with $s=x+iy$. Then as
$x\in\mathbb{ ^{\ast}R}^{\lim_{\eta}}$ there exist positive
$a,b\in\mathbb{R}$ such that $a<x<b$. \[
|\zeta_{\mathbb{^{\ast}Q}}(s)|_{\mathbb{^{\ast}C}} \leq
\sum_{n\in\mathbb{^{\ast}N}} n^{-x}\leq
\sum_{n\in\mathbb{^{\ast}N}}
n^{-a}=\zeta_{\mathbb{^{\ast}Q}}(a).\] For a fixed $a$ then $
\zeta_{\mathbb{^{\ast}Q}}(a)$ is a hyperreal series which is an
extension of the real series $\zeta_{\mathbb{Q}}(a)$. Using the
work in \cite{Go} (6.10) one finds that
$\zeta_{N}(a)\simeq_{\eta} \zeta_{\mathbb{Q}}(a)$ for any
$N\in\mathbb{^{\ast}N}/\setminus\mathbb{N}$. Using the work above
there exists an $N\in\mathbb{^{\ast}N}\setminus\mathbb{N}$ such
that $|\zeta_{ \mathbb{^{\ast}Q}}(a)-\zeta_{N}(a)|<\epsilon/2$
and $|\zeta_{N}(a)- \zeta_{ \mathbb{Q}}(a)|<\epsilon/2$ for some
fixed $\epsilon\in\mathbb{^{\ast}R}^{\inf_{\eta}}$. Then by a
simple application of the triangle inequality \begin{align}
|\zeta_{\mathbb{^{\ast}Q}}(a)-\zeta_{\mathbb{Q}}(a)| &\leq
|(\zeta_{ \mathbb{^{\ast}Q}}(a)-\zeta_{N}(a))+(\zeta_{N}(a) -
\zeta_{ \mathbb{Q}}(a))|,\notag\\ &\leq |\zeta_{
\mathbb{^{\ast}Q}}(a)-\zeta_{N}(a)|+|\zeta_{N}(a) - \zeta_{
\mathbb{Q}}(a)|,\notag\\ <\epsilon.\notag \end{align} Therefore
$ \zeta_{\mathbb{^{\ast}Q}}(a)\simeq_{\eta}
\zeta_{\mathbb{Q}}(a)$.

\begin{flushright} $\Box$ \end{flushright}

\begin{corollary}\label{cx2}
For $s\in\mathbb{^{\ast}R}^{\lim_{\eta}}$ and $^{\ast}\Re (s)>1$
then \[ \operatorname{sh}_{\eta}(\zeta_{\mathbb{^{\ast}Q}}(s)) =
\zeta_{\mathbb{Q}}(\operatorname{sh}_{\eta}(s)).\]
\end{corollary}

By using very basic complex analysis a complex sequence can be
considered in terms of its real and imaginary parts and a direct
application of the above leads to the following.

\begin{corollary}\label{cx3}
For $s\in\mathbb{C}$ and $\Re (s)>1$
then  \[
\operatorname{sh}_{\mathbb{C}}(\zeta_{\mathbb{^{\ast}Q}}(s)) =
\zeta_{\mathbb{Q}}(s).\] \end{corollary}

For general $s\in\mathbb{^{\ast}C}^{\lim}$ and $^{\ast}\Re (s)>1$
one uses the properties of the shadow map. Indeed since each term
in $\zeta_{N}(s)$ is limited, \[ \operatorname{sh}_{
\mathbb{C}}(\zeta_{N} (s))=\sum_{n=1}^{N}
\operatorname{sh}_{\mathbb{C}}(n^{-s})=\zeta_{ \mathbb{Q}}
(\operatorname{sh}_{\mathbb{C}}(s)).\]
This last equality follows from the fact that for all
$n\in\mathbb{^{\ast}N}\setminus\mathbb{N}$
$|n^{-s}|_{\mathbb{^{\ast}C}}\simeq_{ \mathbb{^{\ast}C}} 0$. By
the absolute Q-convergence of $\zeta_{\mathbb{^{\ast}Q}}$,
$\zeta_{\mathbb{^{\ast}Q}}(s) \simeq_{\mathbb{^{\ast}C}}
\zeta_{N}(s)$ and so the following is proven.

\begin{theorem}\label{cx4}
For $s\in\mathbb{^{\ast}C}$ and $^{\ast}\Re (s)>1$
then  \[
\operatorname{sh}_{\mathbb{C}}(\zeta_{\mathbb{^{\ast}Q}}(s)) =
\zeta_{\mathbb{Q}}(\operatorname{sh}_{\mathbb{C}}(s)).\]
\end{theorem}
    
\chapter{The Hyper Dedekind Zeta Function}\label{chhypd}

\section{Hyper Mellin Transform}

\begin{definition}
Let $^{\ast}f:\mathbb{^{\ast}R}_{+}^{\times}\rightarrow
\mathbb{^{\ast}C}$ be an internal Q-continuous function with
$^{\ast}f(\infty)=\lim_{y\rightarrow\infty}\text{}^{\ast}f(y)$
existing. Define the set of all such functions to be
$\mathcal{^{\ast}M}$. Then define the hyper Mellin transform to
be the integral
\[^{\ast}L(^{\ast}f,s)=\int_{\mathbb{^{\ast}R}_{+}^{\times}}
^{\ast}
(^{\ast}f(y)-\text{}^{\ast}f(\infty))y^{s}\frac{dy}{y},\]
provided this integral exists.
\end{definition}

\begin{definition}
Suppose $f$ and $g$ are internal hyper complex functions. Define
the notation $f(x)=O(g(x))$ as $|x|\rightarrow \infty$ iff
\[(\exists x_{0}\in\mathbb{^{\ast}R}_{+})(\exists
M\in\mathbb{^{\ast}N})(|f(x)|\leq M|g(x)|\text{ for
$|x|>x_{0}$}).\]
\end{definition}

\begin{theorem}[Mellin]\label{mel} Let $^{\ast}f,\text{}
^{\ast}g\in \mathcal{^{\ast}M}$ with \[ ^{\ast}f(y)=
a_{0}+O(^{\ast}\exp(-cy^{\alpha}),\qquad  ^{\ast}g(y)=
b_{0}+O(^{\ast}\exp(-cy^{\alpha}),\] for $y\rightarrow \infty$
and positive constants $c,\alpha$. Suppose these hyper functions
satisfy \[^{\ast}f(1/y)=Cy^{k}\text{}^{\ast}g(y),\] for some
hyperreal number $k>0$ and some hyper complex number $C\neq 0$.
Then

\begin{enumerate}
\item The integrals $^{\ast}L(^{\ast}f,s)$ and
$^{\ast}L(^{\ast}g,s)$ Q-converge absolutely and uniformly in
$\{s\in \mathbb{^{\ast}C}: ^{\ast}\Re(s)>k\}$. They are therefore
Q-holomorphic on this space and admit Q-holomorphic continuations
to $\mathbb{^{\ast}C}\setminus\{0,k\}$.
\item The hyper Mellin transforms have simple poles at $s=0$ and
$s=k$ with residues $\operatorname{Res}_{s=0} \text{}
^{\ast}L(^{\ast}f, s)= -a_{0}$, $\operatorname{Res}_{s=k}
\text{} ^{\ast}L(^{\ast}f, s)= Cb_{0}$,
$\operatorname{Res}_{s=0} \text{} ^{\ast}L(^{\ast}g, s)= -b_{0}$
and $\operatorname{Res}_{s=k} \text{} ^{\ast}L(^{\ast}g, s)=
C^{-1}a_{0}$. \item They satisfy the functional equation
\[^{\ast}L(^{\ast}f,s)=C\text{}^{\ast}L(^{\ast}g,k-s).\]
\end{enumerate}
\end{theorem}

\noindent \textbf{Proof:}
The Q-convergence follows from the methods used
in the first chapter by looking at partial Mellin transforms \[
^{\ast}L(^{\ast}f,s)_{N} =\int_{^{\ast}(0,N)}
^{\ast}\text{}
(^{\ast}f(y)-\text{}^{\ast}f(\infty))y^{s}\frac{dy}{y},\]
and transfer.

Now let $^{\ast}\Re(s)>k$ then the hyper Mellin transform can be
rewritten by splitting up the interval of integration \[
^{\ast}L(^{\ast}f,s)=\int_{^{\ast}(1,\infty)}^{\ast} (^{\ast}
f(y)-a_{0})y^{s}\frac{dy}{y} + \int_{^{\ast}(0,1]}^{\ast}
(^{\ast} f(y)-a_{0})y^{s}\frac{dy}{y}.\] In the second integral
make the substitution $y\mapsto 1/y$ and use $^{\ast}f(1/y) =
Cy^{k}\text{}^{\ast}g(y)$. \[\int_{^{\ast}(0,1]}^{\ast}
(^{\ast} f(y)-a_{0})y^{s}\frac{dy}{y} = -\frac{a_{0}}{s} +C
\int_{^{\ast}(1,\infty)}^{\ast} (^{\ast}
g(y)-b_{0})y^{k-s-1} dy - \frac{Cb_{0}}{k-s}.\] This Q-converges
absolutely and uniformly for $^{\ast}\Re(s)>k$ by using an
identical method to the above. Hence \[ ^{\ast}L(\text{}^{\ast} f
,s )=-\frac{a_{0}}{s}+\frac{Cb_{0}}{s-k} + \text{} ^{\ast}F(s),\]
where \[ ^{\ast}F(s) = \int_{^{\ast}(1,\infty)}^{\ast} [(^{\ast}
f(y)-a_{0})y^{s}+ C\frac{dy}{y}(^{\ast}
g(y)-b_{0})y^{k-s}]\frac{dy}{y}.\] Then similarly \[
^{\ast}L(\text{}^{\ast} g ,s
)=-\frac{b_{0}}{s}+\frac{C^{-1}a_{0}}{s-k} + \text{}
^{\ast}G(s),\] where \[ ^{\ast}G(s) =
\int_{^{\ast}(1,\infty)}^{\ast} [(^{\ast} g(y)-b_{0})y^{s}+
C^{-1}\frac{dy}{y}(^{\ast} f(y)-a_{0})y^{k-s}]\frac{dy}{y}.\]
These integrals Q-converge absolutely and locally uniformly on
the whole complex plane, so they represent Q-holomorphic
functions. Moreover it is clear that
$^{\ast}F(s)=C\text{}^{\ast}G(k-s)$ and $^{\ast}L(^{\ast}f,s)= C
\text{}^{\ast}L(^{\ast}g,k-s)$.

\begin{flushright} \textbf{$\Box$} \end{flushright}

\section{Nonstandard Algebraic Number Theory}
In order to define the hyper Dedekind zeta function one needs to
develop the notions of algebraic numbers and integers in a
nonstandard setting.

\subsection{Hyper Polynomials and Hyper Algebraic Numbers}
Consider a ring commutative ring $R$ and the set of polynomials
$R[x]$. Both of these sets can be enlarged in an nonstandard
framework to give a hyper commutative ring $^{\ast}R$ and the set
of internal hyperpolynomials $^{\ast}R[x]$ (which is different
from the set of finite polynomials with coefficients from
$^{\ast}R$, $(^{\ast}R)[x]$). By transfer the notion of the
degree of a hyperpolynomial carries through.

\begin{definition}\label{hyalg}
A number $\alpha\in\mathbb{^{\ast}C}$ is a hyper algebraic number
if there exists a $f(x)=a_{N}x^{N}+\ldots
+a_{0} \in\mathbb{^{\ast}Q}[x]$ with $a_{i}$ not all zero and
$f(\alpha)=0$. Further if $\alpha$ is the root of a monic $g
\in\mathbb{^{\ast}Z}[x]$  then it is said to be a hyper algebraic
integer.
\end{definition}

\begin{theorem}\label{minimal}
Let $\alpha$ be a hyper algebraic number. Then there exists a
unique hyperpolynomial $p(x)\in\mathbb{^{\ast}Q}[x]$ which is
monic, irreducible and of smallest degree such that
$p(\alpha)=0$. Moreover, if $f(x)\in\mathbb{^{\ast}Q}[x]$ and
$f(\alpha)=0$ then $p(x)|f(x)$.
\end{theorem}

\noindent \textbf{Proof:}

Let
$\alpha\in\mathcal{^{\ast}A}$ be the set of hyper algebraic
numbers and $^{\ast}S_{\alpha}=\{ g\in\mathbb{^{\ast}Q}[x]:
g(\alpha)=0 \}$. Then by transfer

\[ (\forall \alpha\in\mathcal{A})(\exists p(x)=a_{\deg(p(x))}x^{
\deg(p(x))}+\ldots +a_{0}  \in S_{\alpha})\]\[ (\forall f(x)\in
S_{\alpha}) (\deg (p(x))\leq\deg(f(x)) \wedge a_{\deg(p(x))}=1
\wedge p(x)|f(x)),\]

\begin{displaymath}
\downarrow \text{$\ast$-transform,} \end{displaymath}

\[(\forall
\alpha\in\mathcal{^{\ast}A})(\exists p(x)=a_{\deg(p(x))}x^{
\deg(p(x))}+\ldots +a_{0}  \in ^{\ast}S_{\alpha})
(\forall f(x)\in ^{\ast}S_{\alpha}) \]\[(\deg
(p(x))\leq\deg(f(x))\wedge a_{\deg(p(x))}=1  \wedge
p(x)|f(x)).\]

Suppose this $p(x)$ is not irreducible then it can be written as
a product of two lower degree polynomials in
$\mathbb{^{\ast}Q}[x]$, say $p(x)=g(x)h(x)$. So $p(\alpha)
=g(\alpha)h(\alpha)=0$. By transfer $\mathbb{^{\ast}C}$ is a
hyper integral domain and so either $g(\alpha)=0$ or
$h(\alpha)=0$ which is a contradiction.

For uniqueness suppose there are two such polynomials
$p(x)$ and $q(x)$ satisfying the above then $p(x)|q(x)$ and
$q(x)|p(x)$ which implies $p(x)=q(x)$.

\begin{flushright} \textbf{$\Box$} \end{flushright}

\begin{definition}\label{degree}
Let the degree of $\alpha$ be given by the degree of $p(x)$, the
minimal internal hyperpolynomial of $\alpha$.
\end{definition}

\begin{definition}\label{dimo} A hyper field $^{\ast}K\subset
\mathbb{^{\ast}C}$ is a hyper algebraic number field if its
dimension over $\mathbb{^{\ast}Q}$ is hyper finite. The dimension
of $^{\ast}K$ over $\mathbb{^{\ast}Q}$ is called the degree of
$^{\ast}K$.
\end{definition}

\begin{theorem}\label{single}
Let $\alpha$, $\beta$ be hyper algebraic numbers then there
exists a hyper algebraic number $\gamma$ such that
$\mathbb{^{\ast} Q}(\alpha, \beta)=\mathbb{^{\ast}Q}(\gamma).$
\end{theorem}

\noindent \textbf{Proof:}

Using the transfer principle
\[ (\forall \alpha,\beta\in\mathcal{A})(\exists
\gamma\in\mathcal{A}) (\mathbb{Q}(\alpha,
\beta)=\mathbb{Q}(\gamma) ).\]

\begin{displaymath}
\downarrow \text{$\ast$-transform,} \end{displaymath}

\[ (\forall \alpha,\beta\in\mathcal{^{\ast}A})(\exists
\gamma\in\mathcal{^{\ast}A}) (\mathbb{^{\ast} Q}(\alpha,
\beta)=\mathbb{^{\ast}Q}(\gamma) ).\]

\begin{flushright} \textbf{$\Box$} \end{flushright}

\begin{corollary}\label{single1}
For $n\in\mathbb{^{\ast}N}$ let $\alpha_{1},\ldots, \alpha_{n}$
be a set of hyper algebraic numbers then there exists
$\gamma\in\mathcal{A}$ such that $\mathbb{^{\ast} Q}(\alpha_{1},
\ldots,\alpha_{n})=\mathbb{^{\ast}Q}(\gamma).$
\end{corollary}

For a given hyper algebraic number, $\theta$, let $p(x)$ be its
minimal internal hyperpolynomial of degree $n$. Let
$\theta^{(1)}=\theta$ and $\theta_{(2)},\ldots,\theta^{(n)}$ be
the conjugates of $\theta$. Then
$\mathbb{^{\ast}Q}(\theta^{(i)})$ ($i=2,\ldots,n$) is a
conjugate field to $\mathbb{^{\ast}Q}(\theta)$. Moreover the
maps $\theta\rightarrow\theta^{(i)}$ are embeddings of
$^{\ast}K=\mathbb{^{\ast}Q}(\theta)$ into $\mathbb{^{\ast}C}$.
There are $n$ embeddings which can be split into hyperreal and
hyper complex embeddings depending on whether or not the
conjugate root is hyperreal or hyper complex. The hyper complex
ones then split into pairs because of hyper complex conjugation.
Let $r_{1}$ be the number of hyperreal embeddings and $2r_{2}$
the number of hyper complex embeddings. Then $n=r_{1}=2r_{2}$.

Minkoswki's ideas, on a simplistic level, interpreted an
algebraic number field $K$ over $\mathbb{Q}$ in terms of points
in $n$-dimensional space. The following extends this in the
natural way to the nonstandard algebraic number fields $^{\ast}K$
over $\mathbb{^{\ast}K}$ of hyperfinite degree $n$. There is a
canonical mapping resulting from the $n$ hyper complex embeddings
($\tau$)

\begin{align}
j:^{\ast}K &\rightarrow ^{\ast}K_{^{\ast}C}:= \prod_{\tau}
\mathbb{^{\ast}C},\notag\\
\alpha &\mapsto (\tau(\alpha)).\notag
\end{align}

This $\mathbb{^{\ast}C}$-vector space is equipped with a
hermitian scalar product

\[ \langle x,y \rangle = \sum_{\tau} x_{\tau} \overline{y}
_{\tau},\]
where $\overline{}$ represents hyper complex conjugation. Also
related to the embeddings is a $\mathbb{^{\ast}R}$-vector space
(the hyper Minkowski space $\boldsymbol{^{\ast}R}$) \[
^{\ast}K_{\mathbb{^{\ast}R}}=[ \prod_{\tau} \mathbb{^{\ast}C}
]^{+},\] which consists of points of
$^{\ast}K_{\mathbb{^{\ast}C}}$ which are invariant under the
involution of hyper complex conjugation. These are the points
$(z_{\tau})$ such that $z_{\overline \tau}=\overline z _{\tau}$.
The restriction of the hermitian scalar product from $^{\ast}
K_{\mathbb{^{\ast}C}}$ to $^{\ast}K_{\mathbb{^{\ast}C}}$ yields
a scalar product. There is also the mapping $j:^{\ast}K
\rightarrow ^{\ast}K_{\mathbb{^{\ast}R}}$ since for all
$\alpha\in ^{\ast}K$, $\overline \tau (\alpha)=\overline
\tau(\alpha)$. This mapping can be used to give lattices in the
hyper Minkowski space, see the following section on lattices.

Minkoswki's theory also exists in a multiplicative form and the
$j$ canonical mapping can be restricted to a homomorphism
\[j:^{\ast}K^{\times}\rightarrow
^{\ast}K_{\mathbb{^{\ast}C}}^{\times}=\prod_{\tau}
\mathbb{^{\ast}C}.\] In a similar way there is the space
$^{\ast}K^{\times}_{\mathbb{^{\ast}R}}$.

There is a homomorphism on $^{\ast}K_{\mathbb{^{\ast}C}}
^{\times}$ given by the product of the coordinates
\[\mathcal{N}:^{\ast}K_{\mathbb{^{\ast}K}}^{\times}\rightarrow
\mathbb{^{\ast}C}^{\times}.\] The usual norm defined below on
$^{\ast}K$ is related to this norm by $N_{^{\ast}K\setminus\mathbb{
^{\ast}Q}}(\alpha)=\mathcal{N}(j(\alpha))$. The restriction leads
to the norm $\mathcal{N}$ being defined on the hyper Minkowski
space.

\subsection{Lattices}
Let $m,n\in\mathbb{^{\ast}N}$ ($m\leq n$) and let $\{
e_{1},\ldots, e_{m}\}$ be a set of linearly independent vectors
in $\mathbb{^{\ast}R}^{n}$. The additive subgroup of
$(\mathbb{^{\ast}R}^{n},+)$ generated by this set is called a
lattice of dimension $m$. A lattice is called complete if $m=n$.

A hyper metric can be placed on $\mathbb{^{\ast}R}^{n}$. Indeed
let $x,y\in\mathbb{^{\ast}R}^{n}$ and define $^{\ast}d_{n}(x,y) =
((x_{1}-y_{1})^{2}+\ldots+(x_{n}-y_{n})^{2})^{1/2}$. With this
hyper metric closed balls can be introduced. Let
$x\in\mathbb{^{\ast}R}^{n}$ and
$r\in\mathbb{^{\ast}R}_{\text{$>$0}}$ then $^{\ast}B_{r}[x]=\{
y\in\mathbb{^{\ast}R}^{n}: ^{\ast}d_{n}(x,y)\leq r\}.$ A subset
$X\subset\mathbb{^{\ast}R}^{n}$ is bounded if $X\subset
^{\ast}B_{r}[0]$ for some $r$. A subset of
$\mathbb{^{\ast}R}^{n}$ is discrete if and only if it intersects
every $^{\ast}B_{r}[0]$ in a hyperfinite set. Let $^{\ast}A_{n}$
be the set of additive subgroups of $\mathbb{^{\ast}R}^{n}$, let
$^{\ast}D_{n}$ be the set of discrete subgroups of
$\mathbb{^{\ast}R}^{n}$ and let $^{\ast}L_{n}$ be the set of
lattices of $\mathbb{^{\ast}R}^{n}$.

\begin{theorem}\label{latdisc}An additive subgroup of
$\mathbb{^{\ast}R}^{n}$ is a lattice iff it is discrete.
\end{theorem}

\noindent \textbf{Proof:}
\[(\forall n\in\mathbb{N})(\forall X\in A_{n}\cap
L_{n})(\exists Y\in D_{n}) (X=Y),\]

\begin{displaymath}
\downarrow \text{$\ast$-transform,} \end{displaymath}

\[(\forall n\in\mathbb{^{\ast}N})(\forall X\in ^{\ast}A_{n}\cap
^{\ast}L_{n})(\exists Y\in ^{\ast}D_{n}) (X=Y),\]

The converse is proved in an identical manner.

\begin{flushright} \textbf{$\Box$} \end{flushright}

For each lattice (generated by $\{e_{1},\ldots,e_{n}\}$) a
fundamental domain $T$ can be defined which consists of all
elements $\sum_{r=1}^{n}a_{r}e_{r}$ ($a_{i}\in\mathbb{^{\ast}R}$)
with $0\leq a_{i}<1$.

The next important concept to introduce is the notion of volume.
To generalise a little the spaces being dealt with consider
$^{\ast}V$, a euclidean vector space (a
$\mathbb{^{\ast}R}$-vector space of hyperfinite dimension $n$
with a symmetric, positive definite bilinear form
$\langle,\rangle: ^{\ast}V \times
^{\ast}V\rightarrow\mathbb{^{\ast}R}$). On $^{\ast}V$ a notion of
volume exists with the cube spanned by an orthonormal basis
$e_{1},\ldots,e_{n}$ has volume 1 while the general
parallelpiped, $\phi$ spanned by $n$ linearly independent vectors
$v_{1},\ldots v_{n}$ has volume
\[\operatorname{vol}(\phi)=|\det(A)|,\] where $A$ is the matrix
of the base change. Now let $L$ be a lattice and $T$ be the
associated fundamental domain then \[\operatorname{vol}(L)=
\operatorname{vol}(T).\]

\subsection{Hyper Algebraic Integers}
Let $^{\ast}K$ be a hyper algebraic number field of degree $n$
over $\mathbb{^{\ast}Q}$ and define $\mathcal{O}_{^{\ast}K}$ to
be the set of algebraic hyperintegers. A simple check shows that
this is a ring. Using the transfer principle and an the proof of
\ref{minimal} gives the following lemma.

\begin{lemma}\label{monet}
For $\alpha\in\mathcal{O}_{^{\ast}K}$, its minimal
hyperpolynomial is monic and an element of
$\mathbb{^{\ast}Z}[x].$ \end{lemma}

\begin{lemma}\label{manet}
Let $\alpha\in\mathcal{^{\ast}A}$ then there exists
$m\in\mathbb{^{\ast}Z}$ such that $m\alpha\in \mathcal{O}_{
^{\ast}K}$.
\end{lemma}
 
\noindent \textbf{Proof:}
\[(\forall \alpha\in\boldsymbol{a})(\exists m\in\mathbb{Z})
(m\alpha\in\mathcal{O}_{K}),\]

\begin{displaymath}
\downarrow \text{$\ast$-transform,} \end{displaymath}

\[(\forall \alpha\in\mathcal{^{\ast}A}) (\exists
m\in\mathbb{^{\ast}Z}) (m\alpha\in\mathcal{O}_{^{\ast}K}),\]

\begin{flushright} \textbf{$\Box$} \end{flushright}

Since $^{\ast}K$ is a vector space over $\mathbb{^{\ast}Q}$ there
exists a $\mathbb{^{\ast}Q}$-basis, $\omega_{1},\ldots,
\omega_{n}$ for $^{\ast}K$. Let $\mathcal{B}_{^{\ast}K}$ be the
set of all bases for $^{\ast}K$ over $\mathbb{^{\ast}Q}$.

\begin{definition}\label{ib}
Let $\{\omega_{i}\}\in\mathcal{B}_{^{\ast}K}$ then it is said to
be an integral basis if $w_{i}\in\mathcal{O}_{^{\ast}K}$ for all
$i$ and
$\mathcal{O}_{^{\ast}K}=\mathbb{^{\ast}Z}\omega_{1}+\ldots +
\mathbb{^{\ast}Z}\omega_{n}.$ Let $\mathcal{IB}_{^{\ast}K}$ be
the set of integral bases for $^{\ast}K$ over
$\mathbb{^{\ast}Q}$. \end{definition}

\begin{lemma}\label{exib}
For all $^{\ast}K$ there exists an integral basis.
\end{lemma}

\noindent \textbf{Proof:}
\[(\forall K)(\exists \{\omega_{i}\} \in\mathcal{B}_{K})(
\{\omega_{i}\}\in\mathcal{IB}_{K}),\]

\begin{displaymath}
\downarrow \text{$\ast$-transform,} \end{displaymath}

\[(\forall ^{\ast}K)(\exists \{\omega_{i}\}
\in\mathcal{B}_{^{\ast}K})(
\{\omega_{i}\}\in\mathcal{IB}_{^{\ast}K}),\]

\begin{flushright} \textbf{$\Box$} \end{flushright}

In analogue to the standard case the norm and trace can be
defined. As $^{\ast}K$ is a hyperfinite vector space over
$\mathbb{^{\ast}Q}$ and linear maps can be defined. For any
$\alpha\in ^{\ast}K$ define a map
\begin{align} \Phi_{\alpha}: ^{\ast}K &\rightarrow ^{\ast}K
,\notag\\ x &\mapsto \alpha x.\notag \end{align}
Then the trace can be defined by $\operatorname{Tr}_{^{\ast}K}
(\alpha)= \operatorname{Tr}(\Phi_{\alpha})$ where
$\operatorname{Tr}$ is the usual trace of a linear map. Similarly
the norm can be defined by $N_{^{\ast}K\setminus
\mathbb{^{\ast}Q}}(\alpha)=\det (\Phi_{\alpha})$ where $\det$ is
the
determinant of a linear map. The discriminant of $^{\ast}K$
can also be defined,
\[d_{^{\ast}K}=\det(\omega_{i}^{(j)})^{2},\] where
$\{\omega_{i}\}$ is an integral basis for $^{\ast}K$ and
$\omega_{i}^{j}$ is a conjugate of $\omega_{i}$.

\begin{definition}\label{modn} For a module $m$ with submodule
$N$ define the index of $N$ in $M$ (denoted $[M:N]$) by the
number of elements in $M\setminus N$. \end{definition}

\subsection{Ideals}
Some aspects of ideals in hyper algebraic number fields have been
examined by Robinson (reference). These works have mainly looked
at purely nonstandard ideals. Let
$\Omega_{\mathcal{O}_{^{\ast}K}}$ be the set of proper internal
ideals of $\mathcal{O}_{^{\ast}K}$. The norm of an internal
ideal $\boldsymbol{a}$, denoted by $N(\boldsymbol{a})$, is its index
in $\mathcal{O}_{^{\ast}K}$. Let
$\Omega_{\mathcal{O}_{^{\ast}K},p}$ be the set of internal primes
ideals of $\mathcal{O}_{^{\ast}K}$.

These ideals could be called integral internal ideals as there
are also internal fractional ideals of $\mathcal{O}_{^{\ast}K}$.
These can be defined as an $\mathcal{O}_{^{\ast}K}$-module
contained in $^{\ast}K$, $\mathcal{A}$, such that there exists
$t\in\mathbb{^{\ast}Z}$ with
$t\mathcal{A}\subset\mathcal{O}_{^{\ast}K}$. By taking $t=1$ any
internal integral ideal is necessarily an internal fractional
ideal. Let $\Omega_{\mathcal{O}_{^{\ast}K},F}$ be the set of
internal proper fractional ideals of $\mathcal{O}_{^{\ast}K}$.

\begin{lemma} \label{modinv}
For each internal prime ideal $\wp\in\Omega_{^{\ast}K,p}$ there
is an internal fractional ideal $\wp^{-1}$ such that $\wp
\wp^{-1}=\mathcal{O}_{^{\ast}K}$.
\end{lemma}

\noindent \textbf{Proof:}
\[ (\forall \wp\in\Omega_{K,p}) (\exists
\wp^{-1}\in\Omega_{K,F}) (  \wp
\wp^{-1}=\mathcal{O}_{K}),\]

\begin{displaymath}
\downarrow \text{$\ast$-transform,} \end{displaymath}

\[ (\forall \wp\in\Omega_{^{\ast}K,p}) (\exists
\wp^{-1}\in\Omega_{^{\ast}K,F}) (  \wp
\wp^{-1}=\mathcal{O}_{^{\ast}K}),\]

\begin{flushright} \textbf{$\Box$} \end{flushright}

\begin{lemma}\label{primefact}
For any $\boldsymbol{a}\in\Omega_{^{\ast}K}$, it can be written
uniquely as a product of prime ideals.
\end{lemma}

\noindent \textbf{Proof:}
\[(\forall \boldsymbol{a}\in\Omega_{K})(\exists ! M\in
\mathbb{N}) (\exists ! \wp_{1},\ldots, \wp_{M} \in \Omega_{
K,p}) (\boldsymbol{a}=\wp_{1}\ldots\wp_{M} \wedge N(\boldsymbol{a})=
N(\wp_{1})\ldots N(\wp_{M}) ),\]

\begin{displaymath}
\downarrow \text{$\ast$-transform,} \end{displaymath}

\[(\forall \boldsymbol{a}\in\Omega_{^{\ast}K})(\exists ! M\in
\mathbb{^{\ast}N}) (\exists ! \wp_{1},\ldots, \wp_{M} \in
\Omega_{ ^{\ast}K,p}) (\boldsymbol{a}=\wp_{1}\ldots\wp_{M}) \wedge N(\boldsymbol{a})=
N(\wp_{1})\ldots N(\wp_{M}) ,\]

\begin{flushright} \textbf{$\Box$} \end{flushright}

Using an identical method of proof this result can be extended
to internal fractional ideals.

\begin{lemma}\label{fracfact}
For any $\boldsymbol{a}\in\Omega_{^{\ast}K,F}$, it can be written
uniquely as a quotient of products of prime ideals.
\end{lemma}

Two internal fractional ideals $\mathcal{A}$ and $\mathcal{B}$
can be defined to be equivalent (written $\mathcal{A} \sim
\mathcal{B}$) if there exists $\alpha,\beta\in
\mathcal{O}_{^{\ast}K}$ such that $(\alpha)\mathcal{A} \sim
(\beta)\mathcal{B}$. It can be simply checked that this relation
is an equivalence relation. Let $Cl_{^{\ast}K}$ be the set of
equivalence classes of internal ideals of $^{\ast}K$. A binary
operation can be placed on this set by defining a product of
$\mathcal{I}_{1}$ and $\mathcal{I}_{2}$ in $Cl_{^{\ast}K}$ to be
the equivalence class of $\mathcal{AB}$ where $\mathcal{A}$ and
$\mathcal{B}$ are two representatives of $\mathcal{I}_{1}$ and
$\mathcal{I}_{2}$ respectively. It can be easily checked that
this product is well defined and a group is formed - the
internal ideal class group - with the equivalence class
containing the internal principle ideals as the identity.

\begin{definition}\label{classic}
Denote by $h_{^{\ast}K}$ the cardinality of the internal ideal
class group, the class number.
\end{definition}

\begin{theorem}{hypclass}
For all hyper algebraic number fields $^{\ast}K$, $h_{^{\ast}K}$
is hyperfinite.
\end{theorem}

\noindent \textbf{Proof:}
\[(\forall K)(\exists M\in\mathbb{N}) (
h_{K}<M),\]

\begin{displaymath}
\downarrow \text{$\ast$-transform,} \end{displaymath}

\[(\forall ^{\ast}K)(\exists M\in\mathbb{^{\ast}N}) (
h_{^{\ast}K}<M),\]

\begin{flushright}\textbf{$\Box$}\end{flushright}

Ideals of hyper algebraic integers lead to complete lattices in
the Minkowski space.

\begin{prop}\label{completelattice}
Let $\boldsymbol{a}\in\Omega_{^{\ast}K}$ then $I=j(\boldsymbol)$
is a complete lattice in $\boldsymbol{R}$ with
$\operatorname{vol} (I)=\sqrt (|d_{^{\ast}K}|)
N(\boldsymbol{\alpha})$. \end{prop}

The proof is almost identical to that in \cite{N}, page 31.

\section{Hyper Dedekind Zeta Function and its First Properties}

\begin{definition}\label{hypdede} Let $s\in\mathbb{^{\ast}C}$ and
define the hyper Dedkind zeta function, \[ \zeta_{^{\ast}K}=
\sum_{\boldsymbol{a}\Omega_{\mathcal{^{\ast}K}}}
\frac{1}{(N(a))^{s}}.\] \end{definition}

The partial hyper Dedekind zeta function can be defined for any
$M\in\mathbb{^{\ast}N}$ by \[ \zeta_{N}=\sum_{N(\boldsymbol{a})\leq
M}\frac{1}{(M(\boldsymbol{a}))^{s}}.\]

\begin{lemma}\label{dedeconv}
$\zeta_{M}(s)$ ($M\in\mathbb{^{\ast}N}$) Q-converges uniformly to
$\zeta_{^{\ast}K}(s)$ for $^{\ast}\Re(s)> 1+\delta$ for every
$\delta>0$.
\end{lemma}

\noindent \textbf{Proof:}
\[(\forall \epsilon\in\mathbb{R}, \epsilon>0) (\exists
M\in\mathbb{N}) (\exists \delta\in\mathbb{R}, \delta>0) (\forall
m,n\in\mathbb{N}, m,n\geq M)\]\[ (\forall s\in\mathbb{C},
\Re(s)>1+\delta) ( d(\zeta_{m}(s),\zeta_{n}(s))<\epsilon),\]

\begin{displaymath}
\downarrow \text{$\ast$-transform,} \end{displaymath}

\[(\forall \epsilon\in\mathbb{^{\ast}R}, \epsilon>0) (\exists
M\in\mathbb{^{\ast}N}) (\exists \delta\in\mathbb{^{\ast}R},
\delta>0) (\forall m,n\in\mathbb{^{\ast}N}, m,n\geq M)
\]\[(\forall s\in\mathbb{^{\ast}C}, ^{\ast}\Re(s)>1+\delta)
(^{\ast} d(\zeta_{m}(s),\zeta_{n}(s))<\epsilon),\]

Thus $\{\zeta_{M}(s)\}$ form a hyper internal Cauchy sequence and
a limit function exists and is $\zeta_{^{\ast}K}$.

\begin{flushright}\textbf{$\Box$}\end{flushright}

\begin{lemma}\label{dedeab}
The hyper Dedekind zeta function is absolutely and uniformly
Q-convergent for $^{\ast}\Re(s)\geq 1+\delta$ for every
$\delta>0$. \end{lemma}

\noindent \textbf{Proof:}
The uniform Q-convergence follows from the previous lemma. By
letting $\zeta_{M}^{A}(s)=\sum_{N(\boldsymbol{a})\leq M} \frac{1}{
^{\ast}d((N(\boldsymbol{a}))^{s},0)}$, for all
$M\in\mathbb{^{\ast}N}$ the absolute Q-convergence can be
established by transfer almost identical to the previous lemma.

\begin{flushright} \textbf{$\Box$} \end{flushright}

\begin{lemma}\label{dedecont}
$\zeta_{M}(s)$ is a Q-continuous function for all
$M\in\mathbb{^{\ast}N}$, $^{\ast}\Re(s)\geq1+\delta$
($\delta>0$). \end{lemma}

\begin{corollary}\label{dedecontt}
$\zeta_{^{\ast}K}(s)$ is Q-continuous for
$^{\ast}\Re(s)\geq1+\delta$ for any $\delta>0$.
\end{corollary}

The proof of the lemma is identical to that of lemma \ref{L5} and
the corollary follows from this result and lemma \ref{L6}.

\begin{theorem}\label{deul}
For $^{\ast}\Re(s)\geq1+\delta$ (for any $\delta>0$), \[
\zeta_{^{\ast}K}=\prod_{\wp}(1-(N(\wp))^{-s}
)^{-1},\] here $\wp$ runs through all prime ideals of
$^{\ast}\mathcal{O}_{^{\ast}K}$.
\end{theorem}

\noindent \textbf{Proof:}
Let $M\in\mathbb{^{\ast}N}$ and let $\Omega_{^{\ast}K,M}$ be the
set of prime ideals with $N(\wp)\leq M$. So,
\[ \prod_{\wp\in\Omega_{^{\ast}K,M}}(1-(N(\wp))^{-s})^{-1} =
\prod_{\wp\in\Omega_{^{\ast}K,M}}( 1+ (N(\wp))^{-s} +
(N(\wp))^{-2s} +\ldots) = \sum_{\boldsymbol{a}\in\mathcal{M}}
\frac{1}{(N(\boldsymbol{a}))^{-s}},\] where $\mathcal{M}$ is the set
of ideals such that all its prime ideal factors are elements of
$\Omega_{^{\ast}K,M}$. By letting $\xi_{M}(s)=    \sum_{\boldsymbol{a}\in\mathcal{M}}
\frac{1}{(N(\boldsymbol{a}))^{-s}}$ one can proceed in an identical
manner to proposition \ref{L8}

\begin{flushright} \textbf{$\Box$} \end{flushright}

\section{Functional Equation}

Like the Riemann zeta function the Dedekind zeta function can
also be extended to a meromorphic function on $\mathbb{C}$ but
with a simple pole at $s=1$. In the nonstandard
setting a similar result is expected.

\begin{definition}[Partial Zeta Function] \label{partiald} Let
$\mathcal{I}\in Cl_{^{\ast}K}$ and define \[
\zeta_{^{\ast}K}(\mathcal{I},s) =
\sum_{\boldsymbol{a}\in\mathcal{I}}\frac{1}{(N(\boldsymbol{a})
)^{s}},\] where the sum is only over internal integral ideals.
Then \[ \zeta_{^{\ast}K}(s)=\sum_{\mathcal{I}\in Cl_{^{\ast}K}}
\zeta_{^{\ast}K}(\mathcal{I},s).\] \end{definition}

This expression can be simplified by the following lemma. For any fractional
ideal $\mathcal{b}$, $\mathcal{O}_{^{\ast}K}^{\times}$ acts on the set
$\boldsymbol{b}^{\times}= \boldsymbol{b}\setminus \{0\}$. The set of orbits is denoted
by $\boldsymbol{b}^{\times}\setminus\mathcal{O}_{^{\ast}K}^{\times}$

\begin{lemma}\label{debij}
Let $\boldsymbol{a}\in\Omega_{^{\ast}K}$ and $\mathcal{I}$ the class of
$\boldsymbol{a}^{-1}$. There is a bijection
\[\boldsymbol{b}^{\times}\setminus\mathcal{O}_{^{\ast}K}^{\times} \rightarrow \{
\boldsymbol{b}\in\mathcal{I}: \boldsymbol{b}\in\Omega_{^{\ast}K}\}, \qquad
\overline a \rightarrow \boldsymbol{b}=(a)\boldsymbol{a}^{-1}.\]
\end{lemma}

\noindent \textbf{Proof:}
Let $a\in\boldsymbol{a}^{\times}$, then
$(a)\boldsymbol{a}^{-1}$ is an integral ideal and lies in $\mathcal{I}$. Suppose
$(a)\boldsymbol{a}^{-1}=b\boldsymbol{a}^{-1}$ then $(a)=(b)$ so $ab^{-1}\in
\mathcal{O}_{^{\ast}K}$ and hence injectivity. The map is also surjective.
Indeed let $\boldsymbol{b}\in\mathcal{I}$ then
$\boldsymbol{b}=(a)\boldsymbol{a}^{-1}$ with $a\in\boldsymbol{a}\boldsymbol{b}
\subset \boldsymbol{a}$.

\begin{flushright} \textbf{$\Box$} \end{flushright}

In order to rewrite the partial zeta functions one uses the hyper Minkowski
space ($\boldsymbol{R}$) and the norm on it. In the special case of a principal
internal ideal $\boldsymbol{a}=(a)$ of $\mathcal{O}_{^{\ast}K}$. Let
$\omega_{1},\ldots, \omega_{n}$ be an integral basis of $\Omega_{^{\ast}K}$
then $a\omega_{1},\ldots ,a\omega_{n}$ is an integral basis of
$\boldsymbol{a}$. Let $A=(a_{ij})$ be the transition matrix $a\omega_{i}=\sum
a_{ij}\omega_{j}$.

Therefore for any $a\in\mathbb{^{\ast}K}^{\times}$
\[N((a))=|N_{\mathbb{^{\ast}K}\setminus\mathbb{^{\ast}Q}}(a)|=|\mathcal{N}(a)|\] and so
\[\zeta_{^{\ast}K}(\mathcal{I},s) =N(\boldsymbol{a})^{s} \sum_{ \overline a \in \boldsymbol{b}^{\times}\setminus
\mathcal{O}_{^{\ast}K}^{\times}} \frac{1}{|\mathcal{N}(\overline a)|^{s}}.\]
It has been shown in proposition \ref{completelattice}e that $\boldsymbol{a}$
forms a complete lattice in the hyper Minkowski space with volume
$\sqrt d_{\boldsymbol{a}}$
($d_{\boldsymbol{a}}=N(\boldsymbol{a})^{2}|d_{^{\ast}K}|.$

\subsection{A Higher Dimensional Hyper Gamma Function}

In order to define a higher dimensional hyper gamma function some analogues of
the one dimensional case are needed, for example the Haar measure and the
integration space. Let $\boldsymbol{^{\ast}R}$ be the hyper Minkoswki space
introduced in the above section. Recall $\boldsymbol{^{\ast}R}=\{z\in
^{\ast}K_{\mathbb{^{\ast}C}} :z=\overline z\}$ where $z=(z_{\tau})$ and
$\overline z=(\overline z_{\overline \tau})$. Now define analogues of
$\mathbb{^{\ast}R}\setminus\{0\}$ and $\mathbb{^{\ast}R}^{\times}_{+}$, \[
\boldsymbol{^{\ast}R}_{\pm}=\{ x\in\boldsymbol{^{\ast}R}: x=\overline x \}
\quad\text{ and } \quad \boldsymbol{^{\ast}R}_{+}^{\times}=\{ x\in
\boldsymbol{^{\ast}R}_{\pm} : x>0\}.\] The analogue of the upper
half plane can also be defined as $^{\ast}\boldsymbol{H}=
\boldsymbol{^{\ast}R}_{\pm} + i
\boldsymbol{^{\ast}R}_{+}^{\times}$. Two functions can be
defined, $| |: \boldsymbol{^{\ast}R}^{\times} \rightarrow
\boldsymbol{^{\ast}R}_{+}^{\times}$ and $\log :
\boldsymbol{^{\ast}R}^{\times}_{+} \rightarrow
\boldsymbol{^{\ast}R}_{\pm}$ which act on each component as
their hyperreal analogues. Using hyper complex exponentiation
$z^{p}$ can be defined for $z$ not a negative hyperreal or zero
by \[z^{p}=(z_{\tau}^{p_{\tau}}).\]

A Haar measure exists on
$\boldsymbol{^{\ast}R}^{\times}_{+}$. One is fixed by the following. Clearly
\[ \boldsymbol{^{\ast}R}^{\times}_{+}=
\prod_{\boldsymbol{p}} \boldsymbol{^{\ast}R}^{\times}_{+, \boldsymbol{p}} ,\]
where the product is over the hyper conjugation classes $\{\tau, \overline \tau
\}$ and $\boldsymbol{^{\ast}R}^{\times}_{+, \boldsymbol{p}}=
\mathbb{^{\ast}R}_{+}^{\times}$ when $\boldsymbol{p}$ is
hyperreal and $\boldsymbol{^{\ast}R}^{\times}_{+,
\boldsymbol{p}}=\{ (y,y):y\in \mathbb{^{\ast}R}_{+}^{\times}\}$
when hyper complex. Define the isomorphism
$\boldsymbol{^{\ast}R}^{\times}_{+, \boldsymbol{p}}\rightarrow
\mathbb{^{\ast}R}_{+}^{\times}$ by $y\mapsto y$ when
$\boldsymbol{p}$ is hyperreal and $y\mapsto y^{2}$ when
$\boldsymbol{p}$ is hyper complex. Together these result in an
isomorphism \[
\Phi:\boldsymbol{^{\ast}R}^{\times}_{+}\rightarrow
\prod_{\boldsymbol{p}}\mathbb{^{\ast}R}_{+}^{\times}.\] The
usual Haar measure on $\mathbb{^{\ast}R}_{+}^{\times}$, $dt/t$
can give a product measure on
$\boldsymbol{^{\ast}R}^{\times}_{+}$ which will be called the
canonical measure, denoted $dy/y$.

\begin{definition}\label{highg} For $s=(s_{\tau})\in
^{\ast}K_{\mathbb{^{\ast}C}}$ such that $^{\ast}\Re(s_{\tau})>0$
define the gamma function associated to $^{\ast}K$ by \[
\Gamma_{^{\ast}K}(s)=\int_{\boldsymbol{^{\ast}R}^{\times}_{+}}^{\ast}
\mathcal{N}(^{\ast}\exp(-y)y^{s})\frac{dy}{y}.\] \end{definition}

The study of this function reduces to to the hyper gamma function using the
isomorphism defined above. Indeed by the product decomposition \[ \Gamma_{^{\ast}K}(s)=
\prod_{\boldsymbol{p}}
\Gamma_{\boldsymbol{p}}(s_{\boldsymbol{p}}).\]

\begin{prop}\label{highgg}
For $\boldsymbol{p}$ hyperreal, $s_{\boldsymbol{p}}=s_{\tau}$ and
\[ \Gamma_{\boldsymbol{p}}(s_{p})
=^{\ast}\Gamma(s_{\boldsymbol{p}}).\] For $\boldsymbol{p}$ hyper
complex, $s_{\boldsymbol{p}}=(s_{\tau},\overline s_{\tau})$ and
\[ \Gamma_{\boldsymbol{p}}(s_{\boldsymbol{p}})
=2^{1-s_{\tau}-s_{\overline \tau}} \text{}
^{\ast}\Gamma(s_{\tau}+s_{\overline \tau})\]. \end{prop}

\noindent \textbf{Proof:}
For the hyperreal case,
$\Gamma_{\boldsymbol{p}}(s_{\boldsymbol{p}})
=\int_{\mathbb{^{\ast}R}_{+}^{\times}}^{\ast}
\mathcal{N}(^{\ast}\exp (-t)t^{s_{\tau}}) dt/t$ which is just the
usual gamma integral because the norm is simply the identity in
this case.

In the hyper complex case \begin{align}  \Gamma_{\boldsymbol{p}}(s_{\boldsymbol{p}})
&= \int_{\mathbb{^{\ast}R}_{+}^{\times}}^{\ast}
\mathcal{N}(^{\ast}\exp (-(t,t))(\sqrt t,\sqrt
t)^{(s_{\tau},s_{\overline \tau})}) dt/t, \notag \\
&=\int_{\mathbb{^{\ast}R}_{+}^{\times}}^{\ast}
\mathcal{N}(^{\ast}\exp (-2\sqrt t)\sqrt t ^{s_{\tau}+
s_{\overline \tau}}) dt/t,\notag\\
&= 2^{1-s_{\tau}-s_{\overline \tau}} \text{}
^{\ast}\Gamma(s_{\tau}+s_{\overline \tau}), \qquad
\text{$t\mapsto (t/2)^{2}$}.\notag
\end{align}

\begin{flushright} \textbf{$\Box$} \end{flushright}

This decomposition shows that this gamma integral converges for
$\boldsymbol{s}=(s_{\tau})$ with $^{\ast}\Re>0$ and admits a
Q-analytic continuation to all of $^{\ast}K_{\mathbb{^{\ast}C}}$
except for poles corresponding to those of the hyperreal gamma
function.

\subsection{Hyper Theta Functions}

A hyper theta function was introduced in the second chapter and
naturally there are many generalizations of these functions.

\begin{definition}\label{hight} For every complete lattice $L$ of
$^{\ast}\boldsymbol{R}$ define the theta series
for $z\in\boldsymbol{^{\ast}H}$ \[^{\ast}\Theta_{L}(z)
=\sum_{g\in L}^{\ast}\exp(\pi i \langle gz,z\rangle ).\]
\end{definition}

\begin{prop}\label{tconv}
The above theta series Q-converges absolutely and uniformly for
all $z\in ^{\ast}K_{\mathbb{^{\ast}C}}$ with $^{\ast}\Im
(z)>\delta$ for $\delta\in\mathbb{^{\ast}R}$ and $\delta>0$.
\end{prop}

\noindent \textbf{Proof:}
From the definition of a lattice $L=\mathbb{^{\ast}Z}r_{1} +
\ldots +\mathbb{^{\ast}Z}r_{n}$ where $(r_{1},\ldots, r_{n})$ is
aa basis of $L$ with $r_{i}$ in $\boldsymbol{R}$. Define a
partial theta function by \[  ^{\ast}\Theta_{L_{M}}(z)
=\sum_{g\in L_{M}}^{\ast}\exp(\pi i \langle gz,z\rangle ),\]
where $M\in\mathbb{^{\ast}M}$ and
$L_{M}=(\mathbb{^{\ast}Z}\setminus M\mathbb{^{\ast}Z})r_{1}
+\ldots + (\mathbb{^{\ast}Z}\setminus M\mathbb{^{\ast}Z})r_{n}.$
Then $ ^{\ast}\Theta_{L_{M}}$ Q-converges to $^{\ast}\Theta_{L}$
by transfer from the classical case. Using an identical method
of proof by transfer from proposition $11.2$ gives uniform
Q-convergence. Absolute Q-convergence follows from applying the
proof to the function \[^{\ast}\Theta_{L_{M},|.|}(z) =\sum_{g\in
L_{M}}|^{\ast}\exp(\pi i \langle gz,z\rangle )|.\]

\begin{flushright} \textbf{$\Box$} \end{flushright}

Hyper Schwartz functions can be defined on
$\boldsymbol{^{\ast}R}$ . Let $^{\ast}f\in
C^{\infty}(\boldsymbol{^{\ast}R})$ be the space of Q-smooth
functions $^{\ast}f:\boldsymbol{^{\ast}R}\rightarrow
\mathbb{^{\ast}C}$.

\begin{definition}\label{hihgs} Define the hyper Schwartz space
\[\mathcal{S} = \{ ^{\ast}f\in
C^{\infty}(\boldsymbol{^{\ast}R}): x\in \boldsymbol{^{\ast}R},
\lim_{x\rightarrow\infty} |x|^{m}
\frac{d^{k}\text{}^{\ast}f}{dx^{m}}=0, \forall
k,m\in\mathbb{^{\ast}N}\} \].\end{definition}

\begin{definition}\label{highfo} Let $^{\ast}f\in\mathcal{S}$ and
define its Fourier transform to be, \[
\widehat{f}(y)=\int^{\ast}_{\boldsymbol{^{\ast}R}^{\times}}
\text{} ^{\ast}f(x) ^{\ast}\exp(-2\pi i \langle x,y \rangle
dx.\] \end{definition}

\begin{prop}\label{selftran} The function
$^{\ast}h(x)=^{\ast}\exp(-\pi \langle x,x \rangle )$ is its own
Fourier transform. \end{prop}

\noindent \textbf{Proof:}
To perform the hyper integration an isometry can be used to
identify the euclidean vector space $\boldsymbol{^{\ast}R}$ with
$\mathbb{^{\ast}R}^{n}$. The Haar measure $dx$ becomes the
Lebesgue measure $dx_{1}\ldots dx_{n}$. As $^{\ast}h(x)$
decomposes as $h(x)=\prod_{i=1}^{n} \text{}^{\ast}\exp(-\pi
x_{i}^{2})$ the Fourier transform is then
$\widehat{h}(y)=\prod_{i=1}^{n}(^{\ast}\exp(-\pi x_{i}^{2}))$.
From the case $n=1$ (dealt with above) the result follows.

\begin{flushright} \textbf{$\Box$} \end{flushright}
 
\begin{prop}\label{lintran} Let $A$ be a linear transformation of
$\boldsymbol{^{\ast}R}$ and define the function $^{\ast}f_{A}(x)=
^{\ast}f(Ax)$. Then \[ \widehat{^{\ast}f}_{A}(y) = \frac{1}{|\det
A|} \widehat{^{\ast}f}(^{t}A^{-1}y),\] where $^{t}A$ is the
adjoint transformation of $A$. \end{prop}

\noindent \textbf{Proof:}
From the definition of the Fourier transform, \[ \widehat{f}_{A}
(y) =\int_{\boldsymbol{^{\ast}R}}^{\times} \text{} ^{\ast}f(Ax) ^{\ast}\exp(-2\pi
i \langle x,y\rangle) dx.\] Making a change of variable $x\mapsto
Ax$, \begin{align} \widehat{f}_{A}
(y) &=\int_{\boldsymbol{^{\ast}R}}^{\times} \text{} ^{\ast}f(x)
^{\ast}\exp(-2\pi i \langle A^{-1}x,y\rangle) |\det A|^{-1}dx,
\notag\\ &=|\det A|^{-1}\int_{\boldsymbol{^{\ast}R}}^{\times} \text{} ^{\ast}f(x)
^{\ast}\exp(-2\pi i \langle x,^{t}A^{-1}y\rangle) dx. \notag
\end{align}

\begin{flushright} \textbf{$\Box$} \end{flushright}
 
Using this proposition the Poisson summation formula can then be
proven.

\begin{theorem}\label{latdual}
Let $L$ be a complete lattice in $\boldsymbol{^{\ast}R}$ and
define its dual lattice \[ L'=\{ g' \in\boldsymbol{^{\ast}R}:
\langle g,g' \rangle\in\mathbb{^{\ast}Z} \forall g\in L \}.\]
Then for any hyper Schwartz function $^{\ast}f$: \[ \sum_{g\in L}
^{\ast}f(g) = \frac{1}{\operatorname{vol}(L)}\sum_{g' \in L'}
\widehat{^{\ast}f}(g').\]
\end{theorem}

\noindent \textbf{Proof:}
Consider the complete lattice $\mathbb{^{\ast}Z}^{n}$ then there
exists an invertible map, $A$, to $L$. Thus
$L=A\mathbb{^{\ast}Z}^{n}$ and $\operatorname{vol}=|\det A|$.
Clearly the lattice $\mathbb{^{\ast}Z}^{n}$ is self dual and for
all $\boldsymbol{n}  \in   \mathbb{^{\ast}Z}^{n}$ then
\begin{align} g'\in L &\Longleftrightarrow  ^{t}(A\boldsymbol{n}
)g' =^{t}\boldsymbol{n} ^{t}Ag'\in\mathbb{^{\ast}Z} ,\notag\\
&\Longleftrightarrow ^{t}Ag'\in\mathbb{^{\ast}Z}^{n},\notag\\
&\Longleftrightarrow g'\in ^{t}A^{-1}
\mathbb{^{\ast}Z}^{n}.\notag
\end{align}
Therefore $L'=A* \mathbb{^{\ast}Z}^{n}$ with $a*=^{t}A^{-1}$.
Then the statement to be proved becomes \[\sum_{\boldsymbol{n}
\in \mathbb{^{\ast}Z}^{n}}\text{}^{\ast}f_{A}(\boldsymbol{n}) =
\sum_{\boldsymbol{n} \in
\mathbb{^{\ast}Z}^{n}}\text{}^{\ast}\widehat{f}_{A}(\boldsymbol{n}).\]
Let
$^{\ast}g(x)=\sum_{\boldsymbol{^{\ast}k}\in\mathbb{^{\ast}Z}^{n}}
^{\ast}f_{A}(x+\boldsymbol{^{\ast}K})$. Via transfer this
function is absolutely and uniformly Q-convergent. It is also
clearly periodic as for all $\boldsymbol{n}\in \mathbb{^{\ast}Z}
^{n}$ \[^{\ast}g(x+\boldsymbol{n})=\text{} ^{\ast}g(x).\]
$^{\ast}g$ is also a hyper Schwartz function. Indeed as
$^{\ast}f$ is a hyper Schwartz function \[
|^{\ast}f(x+\boldsymbol{k})| ||\boldsymbol{k}||^{n+1} \leq C\]
for almost all $\boldsymbol{k}\in\mathbb{^{\ast}Z}^{n}$ and $x$
varying in a Q-compact domain.. Therefore $^{\ast}g(x)$ is
majorized by a constant multiple of the convergent series
$\sum_{\boldsymbol{k}\neq 0}
\frac{1}{||\boldsymbol{k}||^{n+1}}$. This argument can be
repeated for the partial derivatives of $^{\ast}f$ to show that
$^{\ast}g(x)$ is a $C^{\infty}$ function (where a $C^{\infty}$
function is defined in the Hyper Fourier Transform section of
chapter 2). Combining these observations suggests that it should
have some form of hyper Fourier expansion. One can just extend
the results of the one dimensional case to obtain an expansion
\[ ^{\ast}g(x)=\sum_{\boldsymbol{n}\in\mathbb{^{\ast}Z}^{n}}
a_{\boldsymbol{n}}^{\ast}\exp(2\pi i ^{t}\boldsymbol{n}x),\]
with $a_{\boldsymbol{n}}=\widehat{^{\ast}f}(\boldsymbol{n})$.
Therefore \begin{align}
\sum_{\boldsymbol{n}\in\mathbb{^{\ast}Z}^{n}}
\text{}^{\ast}f(\boldsymbol{n}) &=^{\ast}g(0),\notag\\
&=\sum_{\boldsymbol{n}\in\mathbb{^{\ast}Z}^{n}}
a_{\boldsymbol{n}},\notag\\
&=\sum_{\boldsymbol{n}\in\mathbb{^{\ast}Z}^{n}}
\widehat{^{\ast}f}(\boldsymbol{n}).\notag \end{align}

\begin{flushright} \textbf{$\Box$} \end{flushright}

The reason for the developing this Poisson summation formula is
to prove a transformation equation for the hyper theta function.

\begin{theorem}\label{highttr}
\[ ^{\ast}\Theta_{L}(-1/z)= \frac{\sqrt
\mathcal(z/i)}{\operatorname{vol}(L)}\text{}
^{\ast}\Theta_{L'}(z).\] \end{theorem}

\noindent \textbf{Proof:}
As both sides of the transformation are Q-holomorphic in $z$ it is
sufficient to check the identity for $z=iy$ with $y\in
\boldsymbol{^{\ast}R}_{+}^{\times}$. Let $t=y^{-1/2}$ then \[
z=\frac{i}{t^{2}} \qquad \text{and} \qquad -1/z=it^{2}.\] Then
\[ ^{\ast}\Theta_{L}(-1/z) = \sum_{g\in L}\text{}
^{\ast}\exp(-\pi \langle gt,gt \rangle).\] Let
$^{\ast}f(x)=^{\ast}\exp(-\pi \langle g,g \rangle)$ then
\[^{\ast}\Theta_{L}(-1/z)=\sum_{g\in L} \text{} ^{\ast}f(tg)
\qquad \text{and} \qquad ^{\ast}\Theta_{L'}(z)=\sum_{g'\in L'}
\text{} ^{\ast}f(t^{-1}g).\] Let $A$ be the self-adjoint
transformation of $\boldsymbol{R}$, $x\mapsto tx$, clearly this
has determinant $\mathcal{T}$. Then using propositions
\ref{lintran} and \ref{selftran} with the function $^{\ast}f(x)$
in combination with Poisson summation gives the result.

\begin{flushright} \textbf{$\Box$} \end{flushright}

\subsection{Integral Representation}

A theta series can be associated with the partial zeta function,
\[\theta(\boldsymbol{a},s)=\theta_{\boldsymbol{a}}(s/d_{\boldsymbol{a}}^{1/n})
=\sum_{a\in\boldsymbol{a}} ^{\ast}\exp(\pi i \langle
as/d_{\boldsymbol{a}}^{1/n}, a \rangle).\]

\begin{theorem}\label{whatt}
\[|d_{^{\ast}K}|^{s}\pi^{-ns}\Gamma_{^{\ast}K}(s)
\zeta_{^{\ast}K}(\mathcal{I},2s)=
\int_{\boldsymbol{R}_{+}^{\times}}^{\ast} \text{}^{\ast}g(y)
\mathcal{N}(y)^{s}\frac{dy}{y},\] with the series \[
^{\ast}g(y)=\sum_{a\in\mathcal{I}} \text{}
^{\ast}\exp(-\pi\langle
ay/d_{a}^{1/n},a \rangle).\]
\end{theorem}

\noindent \textbf{Proof:}
In the integral of the hyper higher dimensional gamma function
substitute $y\mapsto \pi |a|^{2}y/d_{\boldsymbol{a}}^{1/n}$. Then
\[|d_{^{\ast}K}|^{s}\pi^{-ns}\Gamma_{^{\ast}K}(s)
\frac{N(a)^{2s}}{|\mathcal{N}(a)^{2s}} =
\int_{\boldsymbol{R}_{+}^{\times}}^{\ast} \text{} ^{\ast}\exp(-\pi\langle
ay/d_{a}^{1/n}, a \rangle)
\mathcal{N}(y)^{s} \frac{dy}{y}.\] This can be summed over a
full system of representatives of $\boldsymbol{a}^{\times}
\setminus\mathcal{O}_{^{\ast}K}^{\times}$ and gives the result. The only
point which needed to be checked was the interchange of summation
and integration. This is just an extension of the result in one
dimension.

\begin{flushright} \textbf{$\Box$} \end{flushright}

The Euler factor at infinity can be defined as \[
^{\ast}Z_{\infty}(s)=
|d_{^{\ast}K}|^{s/2}\pi^{-ns/2}\Gamma_{^{\ast}K}(s/2),\] and
define \[^{\ast}Z(I,s)=\text{} ^{\ast}Z_{\infty}(s)
\zeta_{^{\ast}K}(\mathcal{I},s).\] The aim is to relate this to
the hyper theta function introduced above. The problem is that
the series $^{\ast}g$ is not over all $\boldsymbol{a}$ like the
hyper theta series. To overcome this problem a new measure is
introduced.

\subsection{Hyper Dirichlet Theorem}

Recall the $j$ map which embeds $^{\ast}K^{\times}$ into
$^{\ast}K_{\mathbb{^{\ast}R}}$. There is then the map from this
space to $[\prod_{\tau}\mathbb{^{\ast}R}]^{+}$ given by the
hyper logarithm acting on acting on the absolute value of each
component. This can be represented diagrammatically \[
^{\ast}K^{\times}\rightarrow^{j}
\text{}^{\ast}K_{\mathbb{^{\ast}R}} \rightarrow^{l}
\text{}[\prod_{\tau}\mathbb{^{\ast}R}]^{+} .\] Subgroups of these
groups can be considered.

In particular let $\boldsymbol{S}=\{x\in\boldsymbol{R}_{+}^{\times}:
\mathcal{N}(x)=1\}$ - this is called the norm one hypersurface.
Let $|.|:\boldsymbol{R}_{+}\rightarrow
\boldsymbol{R}_{+}^{\times}$ with $(x_{\tau}\rightarrow
(|x_{\tau}|)$. Then $\mathcal{O}_{^{\ast}K}^{\times}$ is
contained in this set. Every $y\in\boldsymbol{R}_{+}^{\times}$
can be written in the form \[ y=xt^{1/n}\] with
$x=y/\mathcal{y}^{1/n}$ and $t=\mathcal{N}(y)$. This gives a
decomposition \[\boldsymbol{R}_{+}^{\times}= \boldsymbol{S} \times
\mathbb{^{\ast}R}^{\times}_{+}.\]

By transfer there exists a unique Haar measure on
$\boldsymbol{S}$  such that the Haar measure $dy/y$ on
$\boldsymbol{R}_{+}^{\times}$ becomes the product measure \[
\frac{dy}{y}=d^{*}x \times \frac{dt}{t}.\]

The group $|\mathcal{O}_{^{\ast}K}^{\times}|^{2}$ acts quite
naturally on $\boldsymbol{S}$. Recall the logarithm map
$^{\ast}\log : \boldsymbol{R}_{+}^{\times} \rightarrow
\boldsymbol{R}_{\pm}$ given by $(x_{\tau})\mapsto (^{\ast} \log
x_{\tau}).$ Clearly this takes $\boldsymbol{S}$ to the trace zero
space $^{\ast}H=\{x\in\boldsymbol{R}_{\pm}: \operatorname{Tr}(x)=0$.

Let $^{\ast}\lambda$ be the composite map of the $j$ map, $|.|$ map and the
$^{\ast}\log$ map; from $^{\ast}K$ to $\boldsymbol{R}_{\pm}$. By restriction
$|\mathcal{O}_{^{\ast}K}^{\times}|$ lies in $^{\ast}H$ and be transfer

\begin{lemma}\label{dirlog}
Under the $^{\ast}\log$ map $|\mathcal{O}_{^{\ast}K}^{\times}|$ is taken to a
complete lattice $^{\ast}G$ in $^{\ast}H$.
\end{lemma}

A fundamental domain, $^{\ast}F$, for the above action can be taken to be the
preimage of an arbitrary fundamental mesh of the lattice $2^{\ast}G$.

\begin{prop}\label{himel}
The function $^{\ast}Z(\mathcal{I},2s)$ is the Mellin transform,
$^{\ast}Z(\mathcal{I},2s)(L(^{\ast}f,s)$, of the hyper function \[ ^{\ast}f(t)=
^{\ast}f_{F}(\boldsymbol{a},t)=\frac{1}{w}\int_{^{\ast}F}^{\ast}\text{}
^{\ast}\theta(\boldsymbol{a},ixt^{1/n}) d^{*}x,\] with
$w=$\#$\mu_{^{\ast}K}$ the number of roots of unity of
$^{\ast}K$. \end{prop}

\noindent \textbf{Proof:}
From the decomposition of $\boldsymbol{R}_{+}^{\times}$ and the definition of
$^{\ast}Z(\mathcal{I},2s)$ \[^{\ast}Z(\mathcal{I},2s)=
\int_{\mathbb{^{\ast}R}_{+}}^{\times} \int_{^{\ast}\boldsymbol{S}}^{\times}
\sum_{a\in\boldsymbol{a}}\text{} ^{\ast}\exp(-\pi \langle axt',a\rangle) d^{*}x
t^{s}\frac{dt}{t},\] with $t'=(t/d_{\boldsymbol{a}})^{1/n}$. By the definition
of the fundamental domain $^{\ast}F$, \[ \boldsymbol{S}=\bigcup_{\eta\in |
\mathcal{O}_{^{\ast}K}^{\times}|} \eta^{2}\text{}^{\ast}F.\] This is a
disjoint union. Consider a transformation of $\boldsymbol{^{\ast}S}$, $x\mapsto
\eta^{2}x$. This leaves the Haar measure invariant, by definition, and maps
$^{\ast}F$ to $\eta^{2}\text{}F$. Therefore \begin{align} \int_{^{\ast}\boldsymbol{S}}^{\times}
\sum_{a\in\boldsymbol{a}}\text{} ^{\ast}\exp(-\pi \langle axt',a\rangle) d^{*}x &=
\sum_{\eta\in |
\mathcal{O}_{^{\ast}K}^{\times}|} \int_{    \eta^{2}\text{}^{\ast}F}^{\ast}
\text{} ^{\ast}\exp(-\pi \langle axt',a\rangle) d^{*}x,\notag\\
&=\frac{1}{w}\int_{^{\ast}F}^{\times} \sum_{\epsilon\in
\mathcal{O}_{^{\ast}K}^{\times}} \text{} ^{\ast}\exp(-\pi \langle a
\epsilon xt',a\epsilon\rangle) d^{*}x,\notag\\
&=\frac{1}{w}\int_{^{\ast}F}^{\times} (^{\ast}\theta(\boldsymbol{a},
ixt^{1/n})-1)d^{*}x,\notag\\
&=^{\ast}f(t)-\text{}^{\ast}f(0).\notag\end{align}

The factor $1/w$ appears as the kernel of the map $
\mathcal{O}_{^{\ast}K}^{\times}\rightarrow |
\mathcal{O}_{^{\ast}K}^{\times}|.$ Indeed let $\gamma\in \mu_{^{\ast}K}$ and let
$\tau: ^{\ast}K\rightarrow \mathbb{^{\ast}C}$ be an embedding. Then
$^{\ast}\log|\tau(\gamma)|= \text{}\log 1=0$ and thus $\mu_{^{\ast}K} \subset
\ker (^{\ast}\lambda)$. For the converse let $\epsilon\in \ker (^{\ast}\lambda)$
which implies $^{\ast}\lambda(\epsilon)=l(j(\epsilon))=0$ Therefore for each
embedding $\tau: ^{\ast}K\rightarrow \mathbb{^{\ast}C}$, $|\tau(\epsilon)|=1$.
This then means that $j(\epsilon)=(\tau(\epsilon))$ lies in a bounded domain of
$\boldsymbol{^{\ast}R}$. Recall that the $j$ map of ideals leads to complete
lattices in $\boldsymbol{^{\ast}R}$ and in particular $j(\epsilon)$ is a point
of the lattice $j\mathcal{O}_{^{\ast}K}$. Hence $\ker(\lambda)$ can contain only
a hyperfinite number of elements. To finish this off the following lemma is
needed.

\begin{lemma}\label{subsub}
Let $^{\ast}H$ be a hyperfinite subgroup of $^{\ast}K^{\times}$. Then $^{\ast}H$
consists of roots of unity.
\end{lemma}

This follows directly from the transfer principle.

\begin{flushright} \textbf{$\Box$} \end{flushright}

The functional equation for the completed partial zeta function derives from the
hyper theta transformation formula. To achieve this a dual lattice is needed for
$\boldsymbol{a}$ and the volume of $^{\ast}F$ with respect to $d^{*}x$ is needed.

\begin{definition}\label{diff}
Let $\mathcal{^{\ast}D}^{-1}= \{ x\in\text{} ^{\ast}K:
\operatorname{Tr}(x\mathcal{O}_{^{\ast}K})\subset
\mathbb{^{\ast}Z}.$ Then define the different of
$^{\ast}K\setminus\mathbb{^{\ast}Q}$ to be the fractional inverse of
$\mathcal{^{\ast}D}^{-1}$
\end{definition}

\begin{lemma}\label{invo}
The lattice $L'$ which is dual to the lattice $L=\boldsymbol{a}$
is given by $*L'=(\boldsymbol{a}\mathcal{D})^{-1}$. Here the $*$
represents the involution $(x_{tau})\mapsto (\overline x_{\tau})$
on $\boldsymbol{^{\ast}R}$.
\end{lemma}

\noindent \textbf{Proof:}
Using the definition of $\langle,\rangle$ and of the dual lattice
\begin{align} *L'&= \{ *g\in\boldsymbol{^{\ast}R}: \langle
g,a\rangle \in \mathbb{^{\ast}Z} \forall a\in\boldsymbol{a} \},
\notag\\
&=\{ x\in\boldsymbol{^{\ast}R}:\operatorname{Tr}(x\boldsymbol{a})
\subset \mathbb{^{\ast}Z}\},\notag\\
&=\{ x\in\text{}^{\ast}K:\operatorname{Tr}(x\boldsymbol{a})
\subset \mathbb{^{\ast}Z}\}.\notag
\end{align}
The final condition in the
middle set ($\operatorname{Tr}(x\boldsymbol{a}) \subset
\mathbb{^{\ast}Z}$) implies $x\in \text{} ^{\ast}K$.
Using the definition of the inverse different $x\in \text{}
*L'\leftrightarrow \operatorname{Tr}(xa\mathcal{O}_{^{\ast}K})
\subset \mathbb{^{\ast}Z}$ for all $a\in\boldsymbol{a}$
$\leftrightarrow x\boldsymbol{a}\subset \mathcal{^{\ast}D}^{-1}
\leftrightarrow x\in (\boldsymbol{a}\mathcal{^{\ast}D})^{-1}$.

\begin{flushright} \textbf{$\Box$} \end{flushright}

\begin{definition}\label{regg}
Let $t=r_{1}+r_{2}-1$ and $\epsilon_{1},\ldots \epsilon_{t}$ be a
set of fundamental units of $^{\ast}K$. Define the regulator,
$^{\ast}R$, of $^{\ast}K$ to be the absolute value of
\[\det(\delta_{i}\text{}^{\ast}\log|\tau_{i}(\epsilon_{j})|)_{1\leq
i,j \leq r}.\] Here $\delta_{i}=1$ if $\tau_{i}$ is real and
$\delta_{j}=2$ if $\tau_{j}$ is complex. The embeddings are
written as $\tau_{1},\ldots,\tau_{r_{1}},\tau_{r_{1}+1},\ldots
,\tau_{r+1},$ $\overline \tau_{r_{1}+1},\ldots, \overline
\tau_{r+1}$. \end{definition}

\begin{lemma}\label{infreg}
Let $^{\ast}F$ be a fundamental domain of $\boldsymbol{S}$ then
with respect to $d^{*}x$
\[\operatorname{vol}(^{\ast}F)=2^{r-1}\text{} ^{\ast}R,\] where
$r$ is the number of infinite places.
\end{lemma}

\noindent \textbf{Proof:}
Recall the decomposition $\boldsymbol{R}_{+}^{\times}=
\boldsymbol{^{\ast}S} \times \mathbb{^{\ast}R}^{\times}_{+}$ This
isomorphism is given by $\alpha: (x,t)\mapsto xt^{1/n}.$ It
transforms the canonical measure into the product measure $d^{*}x
\times dt/t$. With respect to $dt/t$ the set $^{\ast}I=\{
t\in\mathbb{^{\ast}R}_{+}^{\times}:1\leq t \leq \text{} ^{\ast}
\exp(1) \}$ so $\operatorname{vol}(^{\ast}F)$ is also the volume
of $F\times I$ with respect to $d^{*}x \times dt/t$ and the volume
of $^{\ast}\alpha(^{\ast}F\times \text{}^{\ast}I)$ with respect
to $dy/y$. From above there are isomorphisms \[
\boldsymbol{R}_{+}^{\times} \rightarrow ^{^{\ast}\log}
\boldsymbol{R}_{\pm}\rightarrow^{\Phi}\prod_{\boldsymbol{p}}\mathbb{^{\ast}
R}\]. Let $\psi$ be the composite of these. Then $dy/y$ is
transformed into the Lebesgue measure on $\mathbb{^{\ast}R}^{r}$
as already mentioned above resulting in
\[\operatorname{vol}(^{\ast}F)=
\operatorname{vol}_{\mathbb{^{\ast}R}^{r}}(\psi\alpha(^{\ast} F
\times \text{} ^{\ast}I)).\]
Let $\boldsymbol{1}=(1,\ldots, 1)\in\boldsymbol{^{\ast}S}$. Then
\[ \Phi\alpha((\boldsymbol{1},t))=\boldsymbol{e}\text{}
^{\ast}\log (t^{1/n})=\frac{1}{n}\boldsymbol{e}\text{}
^{\ast}\log(t).\]
Here $\boldsymbol{e}=(\boldsymbol{e}_{\boldsymbol{p}_{1}}, \ldots
\boldsymbol{e}_{\boldsymbol{p}_{r}}\in\mathbb{^{\ast}R}^{r}$ with
$\boldsymbol{e}_{\boldsymbol{p}_{i}}=1$ or $2$ depending on
whether $\boldsymbol{e}_{\boldsymbol{p}_{i}}$ is real or complex
respectively. Using the definition of $^{\ast}F$ \[
\Phi\alpha(^{\ast}F\times \{1\})=2\text{}^{\ast}T.\]
Here $^{\ast}T$ is the fundamental mesh of the unit lattice $G$
in the trace zero space
$^{\ast}H=\{(x_{i})\in\mathbb{^{\ast}R}^{r}:\sum x_{i} =0\}$.
Then
\[\Phi\alpha(^{\ast}F\times\text{}^{\ast}I)=2^{\ast}T+
\text{}[0,1/n]^{\ast}\boldsymbol{e}.\] This is the parallelpiped
spanned by the vectors $2\boldsymbol{e}_{1},\ldots
2\boldsymbol{e}_{r-1},1/n \boldsymbol{e}$ with $\boldsymbol{e}_{1},\ldots
\boldsymbol{e}_{r-1}$ span the fundamental mesh. This volume is
$\frac{1}{n}2^{r-1}$ times the absolute value of the determinant
\[\det \left( \begin{array}{cccc} \boldsymbol{e}_{11}& \ldots &
\boldsymbol{e}_{r-1,1} & \boldsymbol{e}_{\boldsymbol{p}_{1}} \\
\vdots & & \vdots & \vdots \\ \boldsymbol{e}_{1r}& \ldots &
\boldsymbol{e}_{r-1,r} & \boldsymbol{e}_{\boldsymbol{p}_{r}}
\end{array}\right).\]
The usual matrix operations, which carry through by
transfer, can be used. In particular adding the first $r-1$ lines
to the last one makes all entries zero apart from the last one
which is $n-\sum\boldsymbol{e}_{\boldsymbol{p}_{i}}$. The matrix
above these zeros leads to the regulator when the determinant is
taken.

\begin{flushright} \textbf{$\Box$} \end{flushright}

\begin{prop}\label{nearly} The hyper function
$^{\ast}f_{^{\ast}F}(\boldsymbol{a}, 1/t)=t^{1/2}\text{}
^{\ast}f_{^{\ast}F^{-1}}((\boldsymbol{a}\mathcal{^{\ast}D})^{-1},
t)$ and \[ ^{\ast}f_{^{\ast}F}(\boldsymbol{a}, t)=
\frac{2^{r-1}}{w}R+O(^{\ast}\exp(-ct^{1/n}))\qquad \text{ for }
t\rightarrow\infty,c>0.\]
\end{prop}

\noindent \textbf{Proof:}
Let $L=\boldsymbol{a}$ be a lattice in $\boldsymbol{^{\ast}R}$
then by above this has volume $\operatorname{vol}(L)=
N(\boldsymbol{a})|d_{^{\ast}K}|^{1/2}$. The lattice dual to this
is given by $*L'=(\boldsymbol{a}\mathcal{^{\ast}D})^{-1}$.
$\langle *gz,*g\rangle = \langle gz,g \rangle $ and so
$^{\ast}\theta_{L'}(z)=^{\ast}\theta_{*L'}(z)$. Moreover
$d_{(\boldsymbol{a}\mathcal{^{\ast}D})^{-1}}=1/d_{\boldsymbol{a}}.$
The transformation $x\mapsto x^{-1}$ leaves $d^{*}x$ unchanged and
maps $^{\ast}F$ to the fundamental domain $^{\ast}F^{-1}$. Also
observe that $\mathcal{N}(x(td_{\boldsymbol{a}})^{1/n})=td_{
\boldsymbol{a}})$ for $x\in\text{}^{\ast}\boldsymbol{S}$. Using
these observations and the transformation formula for the hyper
theta function.
\begin{align} ^{\ast}f_{^{\ast}F}(\boldsymbol{a}, 1/t) &=
\frac{1}{w} \int_{^{\ast}F}^{\ast} \text{} ^{\ast}\theta_{
\boldsymbol{a}}(ix(td_{\boldsymbol{a}})^{-1/n})d^{*}x,\notag\\
&=\frac{1}{w} \int_{^{\ast}F^{-1}}^{\ast} \text{} ^{\ast}\theta_{
\boldsymbol{a}}(-(ix)^{-1}(td_{\boldsymbol{a}})^{-1/n})d^{*}x,\notag\\
&=
\frac{1}{w}\frac{(td_{\boldsymbol{a}})^{1/2}}{\operatorname{vol}
(\boldsymbol{a}}  \int_{^{\ast}F^{-1}}^{\ast} \text{}
^{\ast}\theta_{( \boldsymbol{a}\mathcal{^{\ast}D})^{-1}}(ix
(td_{\boldsymbol{a}})^{-1/n})d^{*}x,\notag\\
&=\frac{t^{1/2}}{w}\int_{^{\ast}F^{-1}}^{\ast} \text{}
^{\ast}\theta_{( \boldsymbol{a}\mathcal{^{\ast}D})^{-1}}(ix
(t/d_{(\boldsymbol{a}\mathcal{^{\ast}D)^{-1}}})^{-1/n})d^{*}x,\notag\\
&=t^{1/2}f_{^{\ast}F^{-1}}((\boldsymbol{a}\mathcal{^{\ast}D})^{-1},t).
\notag\end{align}

In order to prove the second part \[ ^{\ast}f_{^{\ast}F}(
\boldsymbol{a},t)= \frac{1}{w}\int_{^{\ast}F}^{\ast} d^{*}x +
\frac{1}{w}\int_{^{\ast}F}^{\ast}(^{\ast}\theta( \boldsymbol{a},
ixt^{1/n})-1) d^{*}x =\frac{ \operatorname{vol}(F)}{w}+ \text{}
^{\ast}r(t).\] The summands of the integrand are of the form
$^{\ast}\exp(-\pi\langle ax,a\rangle (t')^{1/n})$ with
$a\in\boldsymbol{a}\setminus\{0\}$ and $t'=t/d_{\boldsymbol{a}}.$
Also $x_{\tau}\geq \delta>0$ for all $\tau$ which results in
$\langle ax,a\rangle \geq \delta \langle a,a\rangle$ and so \[
^{\ast}r(t) \leq \frac{\operatorname{vol}(^{\ast}F)}{w}(
^{\ast}\theta_{\boldsymbol{a}}(i\delta t'^{1/n})-1).\] Let
$m=\min \{\langle a,a\rangle: a\in\boldsymbol{a}\setminus\{0\} \}$ and
$M=$\#$\{a\in\boldsymbol{a}:\langle a,a\rangle = m\}.$ Thus
\[^{\ast}\theta_{\boldsymbol{a}}(i\delta t'^{1/n})-1 = \text{}
^{\ast}\exp(-\pi\delta mt'^{1/n})(M+\sum){\langle a,a \rangle
>m} \text{} ^{\ast}\exp(-\pi \delta (\langle a,a \rangle -
m)t'^{1/n})) = O(^{\ast}\exp(-ct^{1/n}).\]

\begin{flushright} \textbf{$\Box$} \end{flushright}

\begin{theorem}\label{bigfunv} The function
$^{\ast}Z(\mathcal{I},s$ admits a Q-analytic continuation to
$\mathbb{^{\ast}C}\setminus\{0,1\}$ and satisfies the functional
equation \[ ^{\ast}Z(\mathcal{I}, s) =
\text{}^{\ast}Z(\mathcal{I}',1-s),\] where
$\mathcal{I}\mathcal{I}'=[\mathcal{^{\ast}D}]$, the ideal class
of $\mathcal{^{\ast}D}$. It has simple poles at $s=0$ and $s=1$
with residues $-2^{r}\text{}^{\ast}R/w$ and
$2^{r}\text{}^{\ast}R/w$ respectively.\end{theorem}

\noindent \textbf{Proof:}
Let
$^{\ast}f(t)=f_{^{\ast}F}(\boldsymbol{a},t)$
and $^{\ast}g(t)=f_{^{\ast}F^{-1}}((\boldsymbol{a}
\mathcal{^{\ast}D})^{-1},t)$ then the previous proposition
implies that \[ ^{\ast}f(t)=t^{1/2}\text{}^{\ast}g(t),\] and
$^{\ast}f(t)=a_{0}+O(^{\ast}\exp(-ct^{1/n})$, $^{\ast}g(t) =
a_{0} +O(^{\ast}\exp(-ct^{1/n}),$ with
$a_{0}=2^{r-1}\text{}^{\ast}R/w$. Using the section on the hyper
Mellin transform enables the Q-analytic continuation of the hyper
Mellin transforms of $^{\ast}f$ and $^{\ast}g$. Also the
functional equation is obtained \[ ^{\ast}L(^{\ast}f,s) = \text{}
^{\ast}L(^{\ast}g,1/2-s)\] with the simple poles of $^{\ast}L
(^{\ast}f,s)$ at $s=0$ and $s=1/2$ with residues $-a_{0}$ and
$a_{0}$ respectively. By the above proposition
$^{\ast}Z(\mathcal{I},s) = \text{} ^{\ast}L(^{\ast}f,s/2)$ and so
$^{\ast}Z(\mathcal{I},s)$ admits a Q-analytic continuation to
$\mathbb{^{\ast}C}/\{0,1\}$ with simple poles at $s=0$and $s=1$
with residues $-2^{r}\text{}^{\ast}R/w$ and
$2^{r}\text{}^{\ast}R/w$ respectively. Moreover it satisfies the
functional equation \[ \text{} ^{\ast}Z(\mathcal{I},s)=\text{}
^{\ast}L(^{\ast}f,s/2)=\text{} ^{\ast}L(^{\ast}g, (1-s)/2) =
\text{} ^{\ast}Z(\mathcal{I}',1-s).\]

\begin{flushright} \textbf{$\Box$} \end{flushright}

This theorem about the partial zeta functions then enables
similar results to be obtained for the completed zeta function
and the hyper Dedekind zeta function of $^{\ast}K$.

\begin{corollary}\label{ye}
$^{\ast}Z_{^{\ast}K}(s)=|d_{^{\ast}K}|^{s/2}\pi^{-ns/2}\text{}
^{\ast}\Gamma_{^{\ast}K}(s/2)\zeta_{^{\ast}K}(s)=
\sum_{\mathcal{I}}\text{} ^{\ast}Z(\mathcal{I},s)$ admits a Q-analytic
continuation to $\mathbb{^{\ast}C}/\{0,1\}$ and satisfies the
functional equation \[^{\ast}Z_{^{\ast}K}(s)=\text{}^{\ast}Z_{
^{\ast}K}(s).\] It has simple poles at $s=0$ and $s=1$ with
residues $-2^{r}h_{^{\ast}K}\text{}^{\ast}R/w$ and
$-2^{r}h_{^{\ast}K}\text{}^{\ast}R/w$ respectively.
\end{corollary}

Using $\zeta_{^{\ast}K}(s)= (
|d_{^{\ast}K}|^{s/2}\pi^{-ns/2}\text{}
^{\ast}\Gamma_{^{\ast}K}(s/2))^{-1}\text{}
^{\ast}Z_{^{\ast}K}(s)$ it is seen that the bracketed term is
only zero at $s=0$ so it cancels the pole of
$^{\ast}Z_{^{\ast}K}(s)$ at this point.

\begin{corollary}\label{yo}
$\zeta_{^{\ast}K}(s)$ has a Q-analytic continuation to
$\mathbb{^{\ast}C}/\{1\}$. It has a simple pole at $s=1$ with
residue \[ \frac{2^{r_{1}}(2\pi)^{r_{2}}}{w|d_{^{\ast}K}|^{1/2}}
h_{^{\ast}K}\text{}^{\ast}R.\] It also satisfies the functional
equation \[\zeta_{^{\ast}K}(1-s)=\text{}^{\ast}A(s)
\zeta_{^{\ast}K}(s).\] Here \[^{\ast}A(s)=
2^{n}\times (2\pi)^{-ns} |d_{^{\ast}K}|^{1/2 -s}(^{\ast}\cos(\pi
s/2))^{r_{1}+r_{2}}(^{\ast}\sin(\pi s/2)^{r_{2}}
(^{\ast}\Gamma(s))^{n}.\]
\end{corollary}
    
\chapter{Nonstandard Interpretation of $p$-adic Interpolation}\label{chadic}

\section{Mahler's Theorem}

The classical problem of interpolation is well known.
For example, interpolating the factorials, $ n!$
$(n\in\mathbb{N})$, is solved by noting that \[
\int_{0}^{\infty} \exp(-x) x^{n}dx=n!. \] Hence finding a
continuous function $f(s)$ taking the value of $n!$ for $s=n$
results in \[ f(s)=\int_{0}^{\infty}\exp(-x)x^{s}dx
=\Gamma(s+1),\] where $\Gamma(s)$ is the familiar gamma
function.

Fix a finite prime $p$ then the $p$-adic analogue (of
interpolation) has no direct solution since $|n!|_{p}\rightarrow
0$ as $|n|\rightarrow\infty$, where $|.|$ is the standard
archimedean metric on $\mathbb{R}$. In many ways it is more
natural to consider \[\prod_{k<n, (k,p)=1}k, \] as a $p$-adic
factorial.

$p$-adic interpolation is not one fixed method. It consists of
many ideas with the aim of obtaining a $p$-adic object from a
real object. There are many ways of achieving this aim.
\begin{itemize}

\item Mahler interpolation uses the topological property of
taking a function on a dense subset of $\mathbb{Z}_{p}$ (for
example $\mathbb{Z}$ or $\mathbb{N}$) and using the natural
extension to give a function on $\mathbb{Z}_{p}$.

\item As an extension of the previous idea the function can be
$p$-adically analytically continued to a larger domain using
ideas of Washington. [Example: $n^{s}$ for $n\in\
1+p\mathbb{Z}_{p}$.]

\item $p$-adic integration as used for the
$p$-adic version of L-functions, the Riemann zeta function and
Eisenstein series among others.

\item Interpolation via a twist of the original function.
[Example: $p$-adic Hurwitz zeta function.]

\item Replacing a function in a variable with a $p$-adic variable
where this makes sense. [Example: $p$-adic L-functions.]
\end{itemize}

In many cases interpolation is a combination of the above ideas
and more techniques as well. The natural question is why
interpolate? In many ways this question is specific to the actual
problem being considered. Though a clear answer is to learn more
about a function from a $p$-adic perspective by creating a
$p$-adic function. Often this $p$-adic function is a simpler
function to work with compared to the original function. This
enables a problem to be tackled from various aspects. By
studying interpolation from a $p$-adic perspective one hopes to
obtain new interpretations of this standard problem.

In more detail the problem of $p$-adic Mahler interpolation is
derived as follows. Suppose there is a sequence $\{ c_{k}
\}_{k=1}^{\infty}$ with $c_{k}\in\mathbb{Q}_{p}$ (where $p$ is a
fixed finite prime). This sequence can be encoded as a function
$g:\mathbb{N}\rightarrow\mathbb{Q}_{p}$ with $g(n)=c_{n}$. Then
when does there exist a continuous function
$f:\mathbb{Z}_{p}\rightarrow\mathbb{Q}_{p}$ such that
$f(n)=g(n)$ $\forall n\in\mathbb{N}$? Since $\mathbb{Z}_{p}$ is
compact, $f$ must be uniformly continuous and bounded
(\cite{M}). This can be stated as, \[ (\forall
m\in\mathbb{N})(\exists n\in\mathbb{N})(\forall
x,y\in\mathbb{Z}_{p})(|x-y|_{p}\leq p^{-n} \Rightarrow
|f(x)-f(y)|_{p}\leq p^{-m}).\] In particular this is true
$\forall x,y\in \mathbb{N}$ as
$\mathbb{N}\subset\mathbb{Z}_{p}$. Then as $f(n)=g(n)$ $\forall
n\in\mathbb{N}$, \[ (\forall m\in\mathbb{N})(\exists
n\in\mathbb{N})(\forall x,y\in\mathbb{N})(|x-y|_{p}\leq p^{-n}
\Rightarrow |g(x)-g(y)|_{p}\leq p^{-m}).\] Thus necessary
conditions are that $g$ is uniformly continuous and bounded
$p$-adically.

Conversely suppose that $g$ is uniformly continuous and bounded.
Then a continuous function
$f:\mathbb{Z}_{p}\rightarrow\mathbb{Q}_{p}$ can be constructed.
For $x\in\mathbb{Z}_{p}$ let $\{x_{i}\}$ be a sequence of
integers tending to $x$ and define
\[ f(x)=\lim_{i\rightarrow\infty} f(x_{i}).\]
A simple exercise shows that this function is well-defined
because of the uniform continuity of $g$.

This can be stated in the following theorem.

\begin{theorem}\cite[112--113]{M} \label{ma1}
Let $f:\mathbb{N}\rightarrow \mathbb{Q}_{p}$ be uniformly
continuous on $\mathbb{N}$. Then there exists a unique function
$F:\mathbb{Z}_{p}\rightarrow \mathbb{Q}_{p}$ which is uniformly
continuous and bounded on $\mathbb{Z}_{p}$, and
$F(x)=f(x)$ if $x\in\mathbb{N}$.
\end{theorem}

In a certain sense $p$-adic interpolation is trivial because when
a function is defined on a dense subset of $\mathbb{Z}_{p}$ (for
example $\mathbb{Z}$ or $\mathbb{N}$) it has a continuous
extension to all of $\mathbb{Z}_{p}$.

Returning to the $p$-adic gamma function. Let $a_{n}=\prod_{k<n,
(k,p)=1}k$ then it is easily checked (essentially Wilson's
theorem) that $a_{n+p^{s}}\equiv -a_{n} (\operatorname{mod}\text{ $p^{s}$})$. So by
making a slight sign adjustment in $a_{n}$ they give a uniformly
continuous function on $\mathbb{N}$ and can be $p$-adically
interpolated using the above theorem to give the $p$-adic gamma
function ($p\neq 2$): \[ \Gamma_{p}(n)=(-1)^{n}\prod_{k<n,
(k,p)=1}k, \] and define $\Gamma_{p}(0)=1.$ Thus $\Gamma(x)$ is
a uniformly continuous function for $x\in\mathbb{Z}_{p}$. This
construction is due to Morita (\cite{Mo}).

Mahler's theorem can be interpreted as a $p$-adic analogue of
a classical theorem of Weierstrass which states that a
continuous function on a closed interval can be uniformly
approximated by polynomials.

\begin{theorem}[Mahler: \cite{M}, chapters 9 and 10] \label{ma2}
Suppose $f:\mathbb{Z}_{p}\rightarrow\mathbb{Q}_{p}$ is continuous
and let
\[a_{n}(f)=\sum_{n\in\mathbb{N}}(-1)^{n-k}\binom{n}{k}f(k).\]
Then $|a_{n}(f)|_{p}\rightarrow 0$ as $|n|\rightarrow\infty$ and
the series \[ \sum_{n\in\mathbb{N}}\binom{x}{n}a_{n}(f) \]
converges uniformly in $\mathbb{Z}_{p}$. Moreover,
\begin{displaymath}
f(x)=\sum_{n\in\mathbb{N}}\binom{x}{n}a_{n}(f).
\end{displaymath} \end{theorem}

The aim of this chapter is to view this interpolation in terms of
the nonstandard shadow map. It begins with some definitions in
nonstandard mathematics. Then it proceeds to demonstrate that the
natural extension of a continuous function to its hyper function
can have an interpretation (with some extra conditions) via the
$p$-adic shadow map to classical Mahler interpolation.

\section{Interpolation Series}

In order to investigate Mahler interpolation from a nonstandard perspective it
seems reasonable to consider functions of the form
$f:\mathbb{N}\rightarrow\mathbb{Q}_{p}$. As has been previously described
functions can be extended their hyper functions
$^{\ast}f:\mathbb{^{\ast}N}\rightarrow\mathbb{^{\ast}Q}_{p}$. To investigate
these functions further the necessary constraint of $p$-adic Q-uniform
continuity will be considered in order to write the corresponding hyper
function as a series and for interpolation.

\begin{definition}\label{binom}
\begin{itemize}
\item For $n,k\in\mathbb{N}$ the classical definition of the
binomial symbol gives
 \[\binom{n}{k}=\left\{ \begin{array}{ll}
\frac{n!}{k!(n-k)!} & \mbox{$0\leq k\leq n$,} \\
0 & \mbox{otherwise.}
\end{array}
\right. \]
\item This can be extended to all $x\in\mathbb{R}$ and
$x\in\mathbb{Z}_{p}$ by \[\binom{x}{k}=\left\{ \begin{array}{ll}
x(x-1)(x-2)\ldots (x-k+1)/k! & \mbox{$k\geq 1$,} \\ 1 &
\mbox{$k=0$.} \end{array} \right. \] \end{itemize}
\end{definition}

Naturally in the latter expression for $x\in\mathbb{N}$ this
agrees with the former definition. It takes integer values for
$x\in\mathbb{Z}$ and has properties such as $| \binom{x}{n}
|_{p}\leq 1 \forall x\in\mathbb{Z}_{p}$ and for $x\in\mathbb{R}$
(with  $x>0$) $\binom{-x}{k}=(-1)^{k}\binom{x+k-1}{k}$. By the
transfer principle, the above definition of the binomial symbol
extends to $\mathbb{^{\ast}R}$. By using the identity for
$\binom{-x}{k}$ one only needs to consider the transfer principle
applied to the classical definition of the binomial symbol to
obtain a definition for $x\in\mathbb{^{\ast}Z}$. Also the
following result is true.

\begin{lemma}\label{les1}
$|\binom{x}{k}|_{p}\leq 1 \forall x\in\mathbb{^{\ast}Z}$.
\end{lemma}

\noindent In order to develop interpolation series some basic
results are needed about the binomial symbol.

\begin{lemma}[Binomial Inversion Formula] \label{l1}
Let $\{a_{k}\}$ be a set of values in $\mathbb{^{\ast}Q}_{p}$
then \[ b_{n}=\sum_{k\in\mathbb{^{\ast}N}}\binom{n}{k}a_{k}
\Longleftrightarrow
a_{n}=\sum_{k\in\mathbb{^{\ast}N}}\binom{n}{k}(-1)^{n-k}b_{k}.\]
\end{lemma} \begin{corollary}[Orthogonality] \label{l2} \[
\sum_{k\in\mathbb{^{\ast}N}}(-1)^{k}\binom{n}{k}\binom{k}{m}=\left\{
\begin{array}{ll} (-1)^{m} & \mbox{if $n=m$,} \\ 0 &
\mbox{otherwise.} \end{array} \right. \] \end{corollary} Both of
these results are simple exercises which are easy to verify
directly by the same proofs as the classical cases.

\begin{definition}\label{coef}
For $n\in\mathbb{^{\ast}N}$ and
$^{\ast}f:\mathbb{^{\ast}N}\rightarrow\mathbb{^{\ast}Q}_{p}$
define
\[a_{n}(^{\ast}f)=\sum_{k\in\mathbb{^{\ast}N}}\binom{n}{k}(-1)^{n-k}\text{}^{\ast}f(k).\]
\end{definition}

This hypersum is well defined since the binomial symbol is defined
in the previous section and $(-1)^{n-k}$ extends to
$\mathbb{^{\ast}Z}$ by the transfer principle. Also the sum is
hyperfinite because $\binom{n}{k}=0$ for $k>n$. In the case when
f is a standard function then the coefficient is the same as that
found in any text on $p$-adic Mahler interpolation.

\begin{prop} \label{pro1}
Let $^{\ast}f:\mathbb{^{\ast}N}\rightarrow\mathbb{^{\ast}Q}_{p}$
be a function. Then for all
$n\in\mathbb{^{\ast}N}$,
\[
^{\ast}f(n)=\sum_{k\in\mathbb{^{\ast}N}}\binom{n}{k}a_{k}(^{\ast}f).\]
\end{prop}

\noindent\textbf{Proof:} This result follows from the binomial
inversion formula. In the formula in lemma 7 let
$a_{n}=a_{n}(f)$ and $b_{k}=\text{}^{\ast}f(k)$. The result
follows.

\begin{flushright} \textbf{$\Box$}\end{flushright}

\section{Interpretation of Interpolation}

\begin{theorem} \label{thint1} Let $f:\mathbb{N}\rightarrow
\mathbb{Q}_{p}$ be a uniformly continuous function, with respect
to the $p$-adic metric, on $\mathbb{N}$ and let
$^{\ast}f:\mathbb{^{\ast}N}\rightarrow \mathbb{^{\ast}Q}_{p}$ be
the extension to its hyper function. Then
$\operatorname{sh}_{p}(^{\ast}f):\mathbb{Z}_{p}\rightarrow\mathbb{Q}_{p}$
is the the unique $p$-adic function obtained by Mahler
interpolation. \end{theorem}

\noindent This gives a natural interpretation of $p$-adic
interpolation in a nonstandard setting.

\noindent\textbf{Proof:}
As discussed in the previous section $^{\ast}f$
can be written as a hyper finite sum for all
$n\in\mathbb{^{\ast}N}$, \[
^{\ast}f(n)=\sum_{k\in\mathbb{^{\ast}N}}\binom{n}{k}
a_{k}(^{\ast}f). \]

This function can be extended to all hyper integers by using
this sum and noting that for $x>0$, \[
\binom{-x}{k}=(-1)^{k}\binom{x+k-1}{k},\] and defining for all
$n\in\mathbb{^{\ast}Z}$ with $n<0$,
\[^{\ast}f(n)=\sum_{k\in\mathbb{^{\ast}N}}(-1)^{k}\binom{-n+k-1}{k}
a_{k}(^{\ast}f). \] To check that $^{\ast}f$ is Q-convergent it
is required to show the nonstandard case of convergence in the
$p$-adic metric. Classically,

\begin{theorem} \label{t1}
For $c_{n}\in\mathbb{Q}_{p}$, a $p$-adic series
\[\sum_{n\in\mathbb{N}}c_{n},\] converges if and only if
$|c_{n}|_{p}\rightarrow 0$ as $|n|\rightarrow\infty$.
\end{theorem}

This result extends by transfer to $\mathbb{^{\ast}Q}_{p}$ since
it is $^{\ast}-$complete with respect to the extended $p$-adic
metric.

\begin{prop}\label{pro2} Let $^{\ast}f$ be as in theorem
\ref{thint1} and $a_{n}(^{\ast}f)$ as defined in definition
\ref{coef} then $|a_{n}(^{\ast}f(v)|_{p}\rightarrow 0$ as
$|v|_{\eta}\rightarrow \infty$. \end{prop}

\noindent\textbf{Proof:}
This is a slightly adapted proof given in \cite[61--63]{RM} and
is based on a proof of Bojanic. As $^{\ast}f$ is uniformly
Q-continuous, \[(\forall s\in\mathbb{^{\ast}N})(\exists
t\in\mathbb{^{\ast}N})(\forall
x,y\in\mathbb{^{\ast}N})(|x-y|_{p}\leq p^{-t} \Rightarrow
|\text{}^{\ast}f(x)-\text{}^{\ast}f(y)|_{p}\leq p^{-s}).\] Let
$x=k\in\mathbb{^{\ast}N}$ and $y=k+p^{t}$ then $|\text{}
^{\ast}f(k+p^{t})-\text{} ^{\ast}f(k)|_{p}\leq p^{-s}$.
$\mathbb{^{\ast}N}$ is bounded in the $p$-adic metric and
$^{\ast}f$ is uniformly continuous on $\mathbb{^{\ast}N}$ it is
also bounded there. Let
$p^{u}=\max_{n\in\mathbb{^{\ast}N}}|\text{}^{\ast}f(n)|_{p}$.
One can replace $^{\ast}f$ by $^{\ast}g=p^{u}\text{}^{\ast}f$ so
that $|\text{}^{\ast}g(n)|_{p}\leq 1$ for all
$n\in\mathbb{^{\ast}N}$. So without loss of generality assume
that $|^{\ast}f(n)|_{p}\leq 1$ for all $n\in\mathbb{^{\ast}N}$.

From the definition of $a_{n}(^{\ast}f)$, $|\binom{n}{k}|_{p}\leq
1$ and properties of the metric,

\begin{align}
|a_{n}(^{\ast}f)|_{p} &\leq | \sum_{k\in\mathbb{^{\ast}N}}
(-1)^{n-k} \binom{n}{k}\text{}^{\ast}f(k)|_{p},
\notag \\
&\leq
\max_{k\in\mathbb{^{\ast}N}}|\text{}^{\ast}f(k)|_{p},\notag\\
&\leq 1.\notag
\end{align}

To continue several lemmas are needed to enable a more precise
bound to be put on values of $a_{n}(^{\ast}f)$.

\begin{definition} \label{diff}Let $^{\ast}f$ be as above and
$n\in\mathbb{^{\ast}N}$ then define the difference operators:
\[ D^{n}(f(x))=\sum_{k=0}^{n}\binom{n}{k}(-1)^{n-k}f(x+k). \]
\end{definition}

\begin{lemma} \label{diff1}For $m\in\mathbb{^{\ast}N}$,
\[ D^{n}(f(x))=\sum_{j=0}^{m}\binom{m}{j}D^{n+j}(f(x-m)).\]
\end{lemma}

\noindent\textbf{Proof of lemma:} Using the definition of $D$,
\begin{align}
\sum_{j=0}^{m}\binom{m}{j}D^{j}(f(x-m)) &=
\sum_{j=0}^{m}\binom{m}{j}\sum_{k=0}^{j}\binom{j}{k}(-1)^{j-k}f(x-m+k),\notag\\
&=\sum_{k=0}^{m}(-1)^{k}f(x-m+k)\sum_{j=0}^{m}\binom{m}{j}\binom{j}{k}(-1)^{j}.\notag
\end{align}
Using corollary 8 the inner sum is zero unless $k=m$ and in this
case it is equal to $(-1)^{m}$. Then by applying $D^{n}$ to each
side gives the result.
\begin{flushright} \textbf{$\Box$}\end{flushright}

\begin{lemma}\label{bound} Let
$^{\ast}f:\mathbb{^{\ast}N}\rightarrow\mathbb{^{\ast}Q}_{p}$ be
any function and define,
\[
a_{n}(^{\ast}f)=\sum_{k=0}^{n}(-1)^{n-k}\binom{n}{k}\text{}^{\ast}f(k).\]
Then,
\[\sum_{j=0}^{m}\binom{m}{j}a_{n+j}(^{\ast}f)=\sum_{k=0}^{n}(-1)^{n-k}
\binom{n}{k}\text{}^{\ast}f(k+m).\] \end{lemma}

\noindent\textbf{Proof of lemma:} (Note that this definition of
$a_{n}$ agrees with the previous but is a little more general as
there are no conditions on the function.) For $m=0$ it is the
definition of $a_{n}$. By the previous lemma, \[
D^{n}(^{\ast}f(m))=\sum_{j=0}^{m}\binom{m}{j}D^{n+j}(^{\ast}f(0)).\]
However from the definition of the difference operators,
\[D^{n}(^{\ast}f(m))=\sum_{k=0}^{n}\binom{n}{k}(-1)^{n-k}\text{}^{\ast}f(m+k).\]
Then letting $m=0$ in the previous equation gives,
\[D^{n}(^{\ast}f(0))=a_{n}(^{\ast}f).\] Putting these three
equations together gives the result. \begin{flushright}
\textbf{$\Box$}\end{flushright}

This lemma is now used
in conjunction with the formula for $a_{n}$ to give,
\[a_{n+p^{t}}(^{\ast}f)=-\sum_{j=1}^{p^{t}-1}\binom{p^{t}}{j}a_{n+j}(^{\ast}f)
+\sum_{k=0}^{n}(-1)^{n-k}\binom{n}{k}(^{\ast}f(k+p^{t}) -\text{}
^{\ast}f(k)).\] From elementary number theory $p$ divides
$\binom{p^{t}}{j}$ for each $j$ in the first sum. Also using the
estimates on the uniform continuity of $^{\ast}f$ at the
beginning of the proof one finds,

\begin{align}
|a_{n+p^{t}}(^{\ast}f)|_{p}&\leq \max_{1\leq j < p^{t}} \{
p^{-1}|a_{n+j}(^{\ast}f)|_{p},p^{-s}\},\notag\\ &\leq p^{-1}
\text{  for $n\geq p^{t}$},\notag \end{align}

since $|a_{n}(^{\ast}f)|_{p}\leq 1$. The above argument can be
repeated by $n\rightarrow n+p^{t}$ in the penultimate inequality
to get, \[ |a_{n}(^{\ast}f)|_{p}\leq p^{-2} \text{  for $n\geq
2p^{t}$} .\] Repeating this argument $(s-1)$ times gives \[
|a_{n}(^{\ast}f)|_{p}\leq p^{-s} \text{  for $n\geq sp^{t}$} .\]
Therefore, $|a_{n}(^{\ast}f)|_{p}\rightarrow 0$ as
$|n|_{\eta}\rightarrow \infty$. \begin{flushright}
\textbf{$\Box$}\end{flushright}

This proposition gives that $|a_{n}(^{\ast}f)|_{p}\rightarrow 0$
as $|n|\rightarrow \infty$. Using theorem \ref{t1} the sum for
the function $^{\ast}f(n)$ converges for all
$n\in\mathbb{^{\ast}Z}$.

\begin{prop}\label{lim} $\forall n\in\mathbb{^{\ast}Z}$
$^{\ast}f(n)\in\mathbb{^{\ast}Q}_{p}^{\lim_{p}}$.
\end{prop}
\noindent\textbf{Proof of proposition:} Returning to the original
function $f$ the binomial inversion formula can be used in the
classical setting to give the classical analogue of proposition
\ref{pro1}.
\[f(n)=\sum_{k\in\mathbb{N}}\binom{n}{k} a_{k}(f),\]
where $a_{k}(f)$ are defined in theorem 2. Similarly it can
be extended to give function on $\mathbb{Z}$. From the definition
of $a_{k}(^{\ast}f)$ it is seen that
$a_{k}(f)=a_{k}(^{\ast}f)$ since
$^{\ast}f(n)=f(n)$ $\forall n\in\mathbb{N}$ (from the definition
of extending a function to its hyper function). To show the
result it is required to demonstrate that $|^{\ast}f(n)|_{p}$ is
limited for all $n\in\mathbb{^{\ast}Z}$. Let
$f_{N}(n)=\sum_{n\leq N}\binom{n}{k} a_{k}(f)$ ($N\in\mathbb{N}$)
and also $^{\ast}f_{N}(n)=\sum_{n\leq N}\binom{n}{k}
a_{k}(^{\ast}f)$ ($N\in\mathbb{^{\ast}N}$). Classically
$\{f_{N}\}$ converges absolutely to $f$ (\cite[chapter 5]{RM}).
By the transfer principle $\{^{\ast}f_{N}\}$ Q-converges
absolutely to $^{\ast}f$. (Note that for $N\in\mathbb{N}$,
$^{\ast}f_{N}=f_{N}$ and $^{\ast}f_{N}$ is the extension of
$f_{N}$ to nonstandard $N$.) Hence,
\begin{equation}
(\forall u\in\mathbb{^{\ast}R}^{+})(\exists
M\in\mathbb{^{\ast}N})(\forall N\geq M)(|^{\ast}f-
\text{}^{\ast}f_{N}|_{p}<u).\label{eq1} \end{equation}

A bound now needs to be put on $^{\ast}f_{N}$. Classically the
$p$-adic valuation satisfies for all $x,y\in\mathbb{Q}_{p}$
$|x+y|_{p}\leq \max \{|x|_{p}, |y|_{p} \}$. By induction, for all
$n\in\mathbb{N}$ and $x_{i}\in\mathbb{Q}_{p}$ ($1\leq i\leq n$):
\[ |\sum_{1\leq i\leq n} x_{i}|_{p} \leq \max_{1\leq i\leq n}
|x_{i}|_{p}. \]
By the transfer principle this result extends to all
$n\in\mathbb{^{\ast}N}$ and $x_{i}\in\mathbb{^{\ast}Q}_{p}$. This
is now applied below.

\begin{align}
|^{\ast}f_{N}(n)|_{p} &=
|\sum_{k\leq N}\binom{n}{k} a_{k}(^{\ast}f)|_{p},\notag\\
&\leq \max_{k\leq N}
\{|\binom{n}{k} a_{k}(^{\ast}f)|_{p}\},\notag\\
&\leq \max_{k\leq N} \{|\binom{n}{k}|_{p}
|a_{k}(^{\ast}f)|_{p}\},\notag\\
&\leq \max_{k\leq N} \{|a_{k}(^{\ast}f)|_{p}\}\text{ (lemma
6).}\label{dave}
\end{align}

The following lemma examines the nonstandard values of
$a_{n}(^{\ast}f)$. It looks at the nonstandard part of a
$p$-adically Q-convergent sequence and is the analogue of the
real case found in chapter 5 of \cite{Go}.

\begin{lemma} \label{lem21}
Let $\{c_{n}\}$ be a convergent of
$p$-adic numbers in $\mathbb{Q}_{p}$. This sequence naturally
extends to a $p$-adically Q-convergent sequence $\{c_{n}\}$ in
$\mathbb{^{\ast}Q}_{p}$ for all $n\in\mathbb{^{\ast}N}$. (See
for example chapter 5 of \cite{Gv}.) Then for all
$n\in\mathbb{^{\ast}N}-\mathbb{N}$,
$c_{n}\in\mathbb{^{\ast}Q}_{p}^{\inf_{p}}$. \end{lemma}

\noindent\textbf{Proof:} Suppose there exists an
$N\in\mathbb{^{\ast}N}-\mathbb{N}$ and an
$\epsilon(N)\in\mathbb{R}^{+}$ such that
$|c_{N}|_{p}>\epsilon(N)$.

As the sequence is convergent, by definition
$|c_{n}|_{p}\rightarrow 0$ as $n\in\mathbb{N}\rightarrow\infty$:
\[ (\forall \delta\in\mathbb{R})(\exists
m(\delta)\in\mathbb{N})(\forall v>m)(|c_{v}|_{p}<\delta).\] In
particular choosing $\delta=\epsilon(N)\in\mathbb{R}^{+}$ gives
that for all $n\in\mathbb{N}$, with $n>m(\epsilon(N))$,
$|c_{n}|_{p}<\epsilon(N)$.

Using the transfer principle on the above logical statement gives
\[ (\forall \delta\in\mathbb{^{\ast}R})(\exists
m(\delta)\in\mathbb{^{\ast}N})(\forall
v>m)(|c_{v}|_{p}<\delta).\]
Again choose $\delta=\epsilon(N)\in\mathbb{R}^{+}$ but this this
time note the transfer principle enables the deduction that for
all $n\in\mathbb{^{\ast}N}$, with $n>m(\epsilon(N))$,
$|c_{n}|_{p}<\epsilon(N)$. This is a contradiction as
$m(\epsilon(N))\in\mathbb{N}$.

\begin{flushright} \textbf{$\Box$}\end{flushright}

This lemma is applied to the sequence $a_{n}(^{\ast}f)$ to
deduce that $a_{n}(^{\ast}f)\in\mathbb{^{\ast}Q}^{\inf_{p}}_{p}$.

From proposition \ref{pro2}, $|a_{n}(^{\ast}f)|_{p}\rightarrow 0$
as $n\rightarrow\infty$. The previous lemma shows that for
nonstandard $n$,
$a_{n}(^{\ast}f)\in\mathbb{^{\ast}Q}_{p}^{\inf_{p}}$. As already
stated $a_{n}(^{\ast}f)=a_{n}(f)$ $\forall n\in\mathbb{N}$ so
$a_{n}(^{\ast}f)\in\mathbb{^{\ast}Q}_{p}^{\inf_{p}}$ (for
$n\in\mathbb{^{\ast}N}-\mathbb{N}$) or
$a_{n}(^{\ast}f)\in\mathbb{Q}_{p}$ (for $n\in\mathbb{N}$). In
both cases $a_{n}(^{\ast}f)\in\mathbb{^{\ast}Q}_{p}^{\lim_{p}}$.
Therefore using equation \ref{dave}
\begin{equation}
|^{\ast}f_{N}(n)|_{p}\leq p^{C},\label{eq2}
\end{equation}
here $p^{C}=\max_{n\in\mathbb{^{\ast}N}} \{
|a_{n}(^{\ast}f)|_{p}\}$.

Finally for some $N$ to be chosen below:
\begin{align}
|^{\ast}f|_{p} &=|(^{\ast}f - \text{}^{\ast}f_{N}) +
^{\ast}f_{N}|_{p},\notag\\
&\leq \max \{ | ^{\ast}f- \text{}^{\ast}f_{N}|_{p}, |
^{\ast}f_{N}|_{p} \} \notag. \end{align} Equation \ref{eq1}
enables $N$ to be chosen such that
$|^{\ast}f-\text{}^{\ast}f_{N}|_{p}<p^{C}$. Thus
$|^{\ast}f|_{p}<p^{C}$ and therefore
$^{\ast}f\in\mathbb{^{\ast}Q}_{p}^{\lim_{p}}$.

\begin{flushright}\textbf{$\Box$}\end{flushright}

Applying the definition of the $p$-adic shadow map to the
function $^{\ast}f$ and using the previous proposition gives a
domain consisting of $\mathbb{Z}_{p}$. To
determine the exact nature of the function the values
$\operatorname{sh}_{p}(^{\ast}f)$ need to be investigated.

Let $z\in\mathbb{^{\ast}Z}$ then $\operatorname{sh}_{p}$ can be
taken of $^{\ast}f(z)$ by the previous proposition.

\[ \operatorname{sh}_{p}(^{\ast}f(z)) =\operatorname{sh}_{p}\left(
\sum_{n\in\mathbb{^{\ast}N}}
\binom{z}{n} a_{k}(^{\ast}f) \right).\]

$\operatorname{sh}_{p}$ is a ring homomorphism and as the sum is
absolutely Q-convergent,

\begin{equation} \label{eq3}
\operatorname{sh}_{p}(^{\ast}f(z))
=\operatorname{sh}_{p} \left( \sum_{n\in\mathbb{N}}
\binom{z}{n}a_{n}(^{\ast}f) \right) +\operatorname{sh}_{p}
\left( \sum_{n\in\mathbb{^{\ast}N}-\mathbb{N}}
\binom{z}{n}a_{n}(^{\ast}f) \right).
\end{equation}

The second term vanishes. Indeed, in the classical case
(\cite[corollary 4.1.2]{Gv}) consider a convergent infinite
series $\sum_{n\in\mathbb{N}} b_{n}$ ($b_{n}\in\mathbb{Q}_{p}$)
then \[ |\sum_{n\in\mathbb{N}}b_{n}|_{p}\leq
\max_{n\in\mathbb{N}}\{|b_{n}|_{p}\}.\]
By the transfer principle  the following lemma will hold.

\begin{lemma} \label{lem22}
Let $\sum_{n\in\mathbb{^{\ast}N}} b_{n}$
($b_{n}\in\mathbb{^{\ast}Q}_{p}$) be a Q-convergent infinite series then
\[|\sum_{n\in\mathbb{^{\ast}N}} b_{n}|_{p}\leq
\max_{n\in\mathbb{^{\ast}N}}\{|b_{n}|_{p}\}.\]
\end{lemma}

This lemma is applied to the infinite sum $b_{n}=0$
($n\in\mathbb{N}$) and $b_{n}=\binom{z}{n}a_{n}(^{\ast}f)$
($n\in\mathbb{^{\ast}N}-\mathbb{N}$). It is Q-convergent by
proposition \ref{pro2} and theorem \ref{t1}. Then
\[|\sum_{n\in\mathbb{^{\ast}N}}
b_{n}|_{p}\leq \max_{n\in\mathbb{^{\ast}N}}\{|b_{n}|_{p}\}\leq
\max_{n\in\mathbb{^{\ast}N}-\mathbb{N}} \{|b_{n}|_{p}\},\] since
$b_{n}=0$ for $n\in\mathbb{N}$. However, \[
|b_{n}|_{p}=|\binom{z}{n}|_{p}|a_{n}(^{\ast}f)|_{p}\leq
|a_{n}(^{\ast}f)|_{p},\] because $|\binom{z}{n}|_{p}|\leq 1$.
Proposition \ref{lim} implies
$|b_{n}|_{p}\in\mathbb{^{\ast}Q}_{p}^{\inf_{p}}$ for all
$n\in\mathbb{^{\ast}N}-\mathbb{N}$. Therefore \[
\sum_{n\in\mathbb{^{\ast}N}-\mathbb{N}}
\binom{z}{n}a_{n}(^{\ast}f)\in\mathbb{^{\ast}Q}_{p}^{\inf_{p}},\]
and the the second term does vanish in equation (\ref{eq3})
leaving
\begin{equation}\label{eq4}
\operatorname{sh}_{p}(^{\ast}f(z))
=\operatorname{sh}_{p} \left( \sum_{n\in\mathbb{N}}
\binom{z}{n}a_{n}(^{\ast}f) \right).
\end{equation}
As already stated in theorem \ref{props} $\operatorname{sh}_{p}$
is a ring homomorphism. By using induction the
following is true for all $N\in\mathbb{N}$ and
$b_{n}\in\mathbb{^{\ast}Q}_{p}^{\lim_{p}}$,
\[\operatorname{sh}_{p}\left( \sum_{n\leq N} b_{n} \right)=
\sum_{n\leq N}\operatorname{sh}_{p}\left( b_{n} \right).\] Using
the transfer principle this statement becomes true for all
$N\in\mathbb{^{\ast}N}$. In particular taking $N$ nonstandard
implies that \[\operatorname{sh}_{p}\left( \sum_{n\in\mathbb{N}}
b_{n} \right)= \sum_{n\in\mathbb{N}}\operatorname{sh}_{p}\left(
b_{n} \right).\] This is true because one could define $b_{n}=0$
for $n\in\mathbb{^{\ast}N}-\mathbb{N}$. Or alternatively note
that $\mathbb{N}\subset \{n\in\mathbb{^{\ast}N}:n\leq N\}$. Thus
equation (\ref{eq4}) becomes
\begin{equation}\label{eq5}
\operatorname{sh}_{p}(^{\ast}f(z))
= \sum_{n\in\mathbb{N}}
\operatorname{sh}_{p} \left(\binom{z}{n}a_{n}(^{\ast}f) \right).
\end{equation}
Again using the ring homomorphism property,
\[\operatorname{sh}_{p} \left(\binom{z}{n}a_{n}(^{\ast}f)
\right)= \operatorname{sh}_{p} \left(\binom{z}{n}
\right) \operatorname{sh}_{p} \left( a_{n}(^{\ast}f) \right).\]
In the proof of proposition 17 it was showed that
$a_{n}(^{\ast}f)=a_{n}(f)$ for $n\in\mathbb{N}$. Therefore
$a_{n}(^{\ast}f)\in\mathbb{^{\ast}Q}_{p}$ and
$\operatorname{sh}_{p} \left( a_{n}(^{\ast}f) \right)=
a_{n}(^{\ast}f).$ This results in
 
\begin{equation}\label{eq6}
\operatorname{sh}_{p}(^{\ast}f(z))
= \sum_{n\in\mathbb{N}}
\operatorname{sh}_{p} \left(\binom{z}{n}\right) a_{n}(^{\ast}f) .
\end{equation}
 
Finally for all $z\in\mathbb{^{\ast}Z}$ and $n\in\mathbb{N}$,
$\binom{z}{n}=\frac{z(z-1)(z-2)\ldots (z-n+1)}{n!}$. The
numerator has a finite number of terms and the denominator is
standard. Therefore with a final application of the ring
homomorphism property:
\begin{align}
\operatorname{sh}_{p} \left(\binom{z}{n}
\right)&=\frac{\operatorname{sh}_{p}(z)(\operatorname{sh}_{p}(z)-1)
(\operatorname{sh}_{p}(z)-2)\ldots
(\operatorname{sh}_{p}(z)-n+1)}{n!},\notag\\
&=\binom{\operatorname{sh}_{p}(z)}{n}\notag.
\end{align}
Putting all this together gives

\[\operatorname{sh}_{p}(^{\ast}f(z)) = \sum_{n\in\mathbb{N}}
\binom{\operatorname{sh}_{p}(z)}{n}a_{n}(f).\]

It has already been proved that the domain of this function is
$\mathbb{Z}_{p}$ and the values lie in $\mathbb{Q}_{p}$.
Further this is identical to the function that would have been
obtained via Mahler $p$-adic interpolation (\cite[chapter
4(2.3)]{R}):

\[ f(x)= \sum_{n\in\mathbb{N}}
\binom{x}{n}a_{n}(f).\]

The identification $x=  \operatorname{sh}_{p}(z)$ is made and the
theorem is proven.

\begin{flushright} \textbf{$\Box$}\end{flushright}

Graphically the above process can be displayed as follows:
\begin{center}
$f:\mathbb{N}\rightarrow\mathbb{Q}_{p},$ \\
$\downarrow$\\
(hyper extension),\\
$\downarrow$\\
$^{\ast}f:\mathbb{^{\ast}N}\rightarrow\mathbb{^{\ast}Q}_{p},$\\
$\downarrow$ \\
(extend to $\mathbb{^{\ast}Z}$),\\
$\downarrow$\\
$^{\ast}f:\mathbb{^{\ast}Z}\rightarrow\mathbb{^{\ast}Q}_{p},$\\
$\downarrow$ \\
($p$-adic shadow map),\\
$\downarrow$\\
$f:\mathbb{Z}_{p}\rightarrow\mathbb{Q}_{p}$,\\
($p$-adic interpolated function).
\end{center}

\section{Extension}

In the previous section I have obtained a nonstandard interpretation of
Mahler's theorem. This work strengthens that result by showing
that all continuous $p$-adic functions
$f:\mathbb{Z}_{p}\rightarrow\mathbb{Q}_{p}$ can be viewed as
functions
$^{\ast}g:\mathbb{^{\ast}N}\rightarrow\mathbb{^{\ast}Q}^{\lim_{p}}$.

\begin{theorem}\label{stronger} Let $f:\mathbb{N}\rightarrow
\mathbb{Q}_{p}$ be a uniformly continuous function, with respect
to the $p$-adic metric, on $\mathbb{N}$. Then there exists a hyper function
$^{\ast}g:\mathbb{^{\ast}N}\rightarrow \mathbb{^{\ast}Q}^{\lim_{p}}$ such that
$\operatorname{sh}_{p}(^{\ast}g):\mathbb{Z}_{p}\rightarrow\mathbb{Q}_{p}$ is
the the unique $p$-adic function obtained by Mahler interpolation.\end{theorem}

\subsection{$p$-adic Functions}

\begin{definition}\label{puni} Let the set of
$p$-adically uniformly continuous functions from
$\mathbb{Z}_{p}$ to $\mathbb{Q}_{p}$ be denoted by
$\mathcal{C}_{p}(\mathbb{Z}_{p},\mathbb{Q}_{p})$.
\end{definition}

By restriction of a function in $\mathcal{C}_{p}(\mathbb{Z}_{p},
\mathbb{Q}_{p})$ to $\mathbb{N}$ gives a map to
$\mathcal{C}_{p}(\mathbb{N},\mathbb{Q}_{p})$. By Mahler's
theorem this map is an isomorphism.

Consider the nonstandard space consisting of $p$-adically
uniformly Q-continuous functions
$\mathcal{C}_{p}(\mathbb{^{\ast}N},\mathbb{^{\ast}Q}^{\lim_{p}})$.
There exists a map from this space to
$\mathcal{C}_{p}(\mathbb{Z}_{p},\mathbb{Q}_{p})$, the $p$-adic
shadow map ($\operatorname{sh}_{p}$) acting on functions. The
definition of this map has been given above. It is
defined for all $^{\ast}g:\mathbb{^{\ast}N}\rightarrow
\mathbb{^{\ast}Q}^{\lim_{p}}$ because $\operatorname{sh}_{p}
(\mathbb{^{\ast}N})=\mathbb{Z}_{p}$ and $\operatorname{sh}_{p}
(\mathbb{^{\ast}Q}^{\lim_{p}})=\mathbb{Q}_{p}$ and are also
surjective maps. This suggests that the $p$-adic shadow map
acting on functions should also be surjective. By the definition
of this map the image of $\operatorname{sh}_{p}(
\mathcal{C}_{p}(\mathbb{^{\ast}N}, \mathbb{^{\ast}Q}^{\lim_{p}})$
is contained in $\mathcal{C}_{p}(\mathbb{Z}_{p},
\mathbb{Q}_{p})$.

By Mahler's theorem any $f\in\mathcal{C}_{p}(\mathbb{Z}_{p},
\mathbb{Q}_{p})$ can be written as \[f(x)=\sum_{n\in\mathbb{N}}
a_{n}(f)\binom{x}{n},  \forall x\in\mathbb{Z}_{p}.\]
Here $a_{n}(f)=\sum_{k=0}^{n}(-1)^{n-k}\binom{n}{k}f(k)$ and
$|a_{n}|_{p}\rightarrow 0$ as $n\rightarrow\infty$. This theorem
extends to the nonstandard space
$\mathcal{C}_{p}(\mathbb{^{\ast}Z}_{p}, \mathbb{^{\ast}Q}_{p})$
by the previous section. In particular as
$\mathcal{C}_{p}(\mathbb{^{\ast}N}, \mathbb{^{\ast}Q}^{\lim_{p}})
\subset  \mathcal{C}_{p}(\mathbb{^{\ast}Z}_{p},
\mathbb{^{\ast}Q}_{p})$ any
$^{\ast}g\in\mathcal{C}_{p}(\mathbb{^{\ast}N},
\mathbb{^{\ast}Q}^{\lim_{p}})$ can be written as
\[^{\ast}g(m)=\sum_{n\in\mathbb{^{\ast}N}}b_{n}(^{\ast}g)\binom{m}{n},
 \forall m\in\mathbb{^{\ast}N}.\] Here
$b_{n}(^{\ast}g)=\sum_{k=0}^{n}
(-1)^{n-k}\binom{n}{k}^{\ast}g(k)$ and $|b_{n}|_{p}\rightarrow 0$
as $n\rightarrow\infty$. The $p$-adic shadow map acts on this
hyperfinite sum since for any $m\in\mathbb{^{\ast}N}$,
$\binom{m}{n}=0$ for $n>m$. Properties of the $p$-adic shadow map
are used from above.

\begin{equation}\label{eq1} \operatorname{sh}_{p}(^{\ast}g(m))=
\sum_{n\in\mathbb{^{\ast}N}}\operatorname{sh}_{p}(b_{n}(^{\ast}g))
\operatorname{sh}_{p}(\binom{m}{n}).\end{equation}

However from the definition of the the image of the $p$-adic
shadow map

\begin{equation}\label{eq2} f(t)=\operatorname{sh}_{p}(^{\ast}g(m))=
\sum_{n\in\mathbb{N}}a_{n}(f)\binom{\operatorname{sh}_{p}(m)}{n},
\end{equation}
Here $t=\operatorname{sh}_{p}(m)$. By equating these two
equations at values of $t=m\in\mathbb{N}$ it is easily shown that
$a_{n}(f)=\operatorname{sh}_{p}(b_{n}(^{\ast}g))$ for all
$n\in\mathbb{N}$ and also
$f(t)=\operatorname{sh}_{p}(^{\ast}g(m))$ for all
$m=t\in\mathbb{N}$.

In order for the two equalities to hold in general and by the
previous paragraph $\operatorname{sh}_{p}(b_{n}(^{\ast}g))=0$ for
all $n\in\mathbb{^{\ast}N}\setminus\mathbb{N}$. This implies that
$b_{n}(^{\ast}g)\in\mathbb{^{\ast}Q}^{\inf_{p}}$ for these values
of $n$.

One question which can be asked is what is the kernel of the
shadow map. $\operatorname{ker}(\operatorname{sh}_{p})=\{
^{\ast}g\in\mathcal{C}_{p}(\mathbb{^{\ast}N},
\mathbb{^{\ast}Q}^{\lim_{p}}):
\operatorname{sh}_{p}(^{\ast}g)=0\}.$ By the previous paragraphs
this kernel is equal to
$\mathcal{C}_{p}(\mathbb{^{\ast}N},
\mathbb{^{\ast}Q}^{\inf_{p}})$.

For the surjectivity consider an $f\in\mathcal{C}_{p}(
\mathbb{N},\mathbb{Q}_{p})$. Then this can be written as
\[ f(n)=\sum_{k\in\mathbb{N}}a_{k}(f)\binom{n}{k}.\]
In particular for each $n\in\mathbb{N}$, $f(n)\in\mathbb{Q}_{p}$
and let $q_{n}=f(n)$. Then by the series expansion for $p$-adic
numbers
\[q_{n}=p^{r_{n}}q_{n_{r{n}}}+p^{r_{n}+1}q_{n_{r{n}+1}}+\ldots
.\]
Let \[ [q_{n}]_{m}=p^{r_{n}}q_{n_{r{n}}}+\ldots+
p^{r_{n}+m}q_{n_{r{n}}+m}\in\mathbb{Q}.\]
Then $[q_{n}]_{m}\rightarrow q_{n}$ as $n\rightarrow \infty$ and
$|[q_{n}]_{m}|_{p}=|q_{n}|_{p}$ $\forall m\in\mathbb{N}$. Now
define a function $f_{m}:\mathbb{N}\rightarrow\mathbb{Q}$ by
$f_{m}(n)=[q_{n}]_{m}$. Since $|[q_{n}]_{m}|_{p}=|q_{n}|_{p}$
$f_{m}\in\mathcal{C}_{p}(\mathbb{N},\mathbb{Q})$. Using the
classical norm on $\mathcal{C}_{p}(
\mathbb{N},\mathbb{Q}_{p})$ with
$|f|_{p}=\sup_{n\in\mathbb{N}}|a_{n}(f)|_{p}$ then
$f_{m}\rightarrow f$ as $m\rightarrow \infty$. Now consider the
nonstandard extension of these functions $f_{m}$ to
$m\in\mathbb{^{\ast}N}$. These are going to be functions of
the form
$^{\ast}f_{m}:\mathbb{^{\ast}N}\rightarrow\mathbb{^{\ast}Q}\in
\mathcal{C}_{p}(\mathbb{^{\ast}N},\mathbb{^{\ast}Q})$. Using the
previous work on Mahler interpolation it can be
shown that in fact
$^{\ast}f_{m}\in\mathcal{C}_{p}(\mathbb{^{\ast}N},\mathbb{^{\ast}Q}^{\lim_{p}}).$
Moreover the classical norm on functions can also be extended to
the nonstandard functions.

For $^{\ast}g\in\mathcal{C}_{p}(\mathbb{^{\ast}N},
\mathbb{^{\ast}Q}^{\lim_{p}})\setminus\mathcal{C}_{p}(\mathbb{^{\ast}N},
\mathbb{^{\ast}Q}^{\inf_{p}})$ define
$|^{\ast}g|_{p}=\sup_{n\in\mathbb{N}}|a_{n}(^{\ast}g)|_{p}$. This
is well defined because it has already been shown that for
$n\in\mathbb{^{\ast}N}\setminus\mathbb{N}$ that
$a_{n}(^{\ast}g)\in\mathbb{^{\ast}Q}^{\inf_{p}}$. Since $^{\ast}g
\notin\mathcal{C}_{p}(\mathbb{^{\ast}N},
\mathbb{^{\ast}Q}^{\inf_{p}})$ there exists at least one
$a_{n}(^{\ast}g)$ (for some $n\in\mathbb{N}$) with
$a_{n}(^{\ast}g)\in\mathbb{^{\ast}Q}^{\lim_{p}}
\setminus\mathbb{^{\ast}Q}^{\inf}_{p}$. Thus
\[\sup_{n\in\mathbb{^{\ast}N}}|a_{n}(^{\ast}g)|_{p}=
\sup_{n\in\mathbb{N}}|a_{n}(^{\ast}g)|_{p}.\] It can be easily
shown that this norm satisfies the triangle inequality.

For each $^{\ast}g\in\mathcal{C}_{p}(\mathbb{^{\ast}N},
\mathbb{^{\ast}Q}^{\lim_{p}})\setminus\mathcal{C}_{p}(\mathbb{^{\ast}N},
\mathbb{^{\ast}Q}^{\inf_{p}})$ define the monad of $^{\ast}g$ to
be
\[\mu(^{\ast}g)=\{^{\ast}h\in\mathcal{C}_{p} (\mathbb{^{\ast}N},
\mathbb{^{\ast}Q}^{\lim_{p}}): ^{\ast}h\simeq_{p} ^{\ast}g\} \cup
\{ j\in\mathcal{C}_{p} (\mathbb{N},
\mathbb{Q}_{p}): j\simeq_{p} g\}.\]
Here $\simeq_{p}$ is the notion of two elements being
infinitesimally close with respect to the norm. That is
$a\simeq_{p} b$ iff $|a-b|_{p}=p^{-N}$ where
$N\in\mathbb{^{\ast}N}\setminus\mathbb{N}$. It is clear that this
is an equivalence relation. Moreover

\begin{lemma} \label{mumu}The only element of $\mu(^{\ast}g)$ which is an
element of $\mathcal{C}_{p} (\mathbb{N},
\mathbb{Q}_{p})$ is $\operatorname{sh}_{p}(^{\ast}g)$.
\end{lemma}
\noindent \textbf{Proof:}
Dealing with uniqueness first. Suppose there are two distinct
elements ($h_{1}$ and $h_{2}$) of $\mathcal{C}_{p} (\mathbb{N},
\mathbb{Q}_{p})$ in $\mu(^{\ast}g)$. Then by the transitive
property of $\simeq_{p}$, $h_{1}\simeq_{p} h_{2}$ which is a
contradiction since both functions are standard. Hence
$h_{1}=h_{2}$.

The existence requires examination of the shadow map of
$^{\ast}g$. Using the above work
\[|^{\ast}g-\operatorname{sh}_{p}(^{\ast}g)|_{p}
=\sup_{n\in\mathbb{N}} |a_{n}(^{\ast}g)-\operatorname{sh}_{p}
(a_{n}(^{\ast}g))|_{p}.\] By the definition of the $p$-adic
shadow map each of the above values is infinitesimal and hence so
is the value of the norm.

\begin{flushright} \textbf{$\Box$} \end{flushright}

Similar definitions hold for elements of $\mathcal{C}_{p}(\mathbb{^{\ast}N},
\mathbb{^{\ast}Q}^{\inf_{p}})$ with the shadow map of any element
being zero and the monad of any element being $\mathcal{C}_{p}(\mathbb{^{\ast}N},
\mathbb{^{\ast}Q}^{\inf_{p}})$ and the $0$ map in $\mathcal{C}_{p} (\mathbb{N},
\mathbb{Q}_{p})$.

Returning to any $f\in\mathcal{C}_{p} (\mathbb{N},
\mathbb{Q}_{p})$ there exists a sequence of functions
$^{\ast}f_{m}:\mathbb{^{\ast}N}\rightarrow\mathbb{^{\ast}Q}^{\lim_{p}}$
converging to $f$. In particular for
$m\in\mathbb{^{\ast}N}\setminus\mathbb{N}$, $f\simeq_{p}\text{}
^{\ast}f_{m}$. Hence $f=\operatorname{sh}_{p}(^{\ast}f_{m})$ and
the $p$-adic shadow map on functions is surjective.

In conclusion $\operatorname{sh}_{p}: \mathcal{C}_{p}(
\mathbb{^{\ast}N},\mathbb{^{\ast}Q}^{\lim_{p}})
\rightarrow \mathcal{C}_{p}(\mathbb{N},\mathbb{Q}_{p})$ is
surjective and this map provides an isomorphism

\[  \mathcal{C}_{p}(\mathbb{N},\mathbb{Q}_{p})
\cong
\mathcal{C}_{p}(
\mathbb{^{\ast}N},\mathbb{^{\ast}Q}^{\lim_{p}})\setminus \mathcal{C}_{p}(
          \mathbb{^{\ast}N},\mathbb{^{\ast}Q}^{\inf_{p}}).\]

As a final remark to this section this result extends my
nonstandard work on Mahler interpolation. This can be displayed
graphically
\begin{center}
$f:\mathbb{N}\rightarrow\mathbb{Q}_{p},$ \\
$\downarrow$\\
(hyper extension),\\
$\downarrow$\\
$^{\ast}f:\mathbb{^{\ast}N}\rightarrow\mathbb{^{\ast}Q}^{\lim_{p}},$\\
$\downarrow$ \\  ($p$-adic shadow map),\\  $\downarrow$\\
$f:\mathbb{Z}_{p}\rightarrow\mathbb{Q}_{p}$,\\ ($p$-adic
interpolated function). \end{center}
The reason for looking at Mahler interpolation in this way is two
fold. Firstly a nonstandard perspective can often provide a
different view point and in this case could also simplify the
work since it deals with functions on $\mathbb{^{\ast}N}$ which
are potentially simpler than those on $\mathbb{Z}_{p}$. Secondly
the possibility of working with more than one prime at the
same time would provide different approach to $p$-adic analysis.
The development of this begins in the next chapter.
    
\chapter{Applications and Double Interpolation}\label{chapp}

The previous chapter showed how one form of $p$-adic
interpolation, Mahler interpolation, could be viewed as the
$p$-adic shadow map of a certain hyper function. In the
introduction it was stated that there are a variety of methods
which can be used for $p$-adic interpolation. Some of these are
now developed from a nonstandard perspective culminating in
interpolation with respect to two or more primes.

\section{Morita Gamma Function}
The gamma function has been studied since the 1700s with its
importance realised by Euler and Gauss. It occurs at one solution
to the problem of finding a continuous function of real or
complex variable that agrees with the factorial function at the
integers.

\begin{definition}[Gamma Function]\label{euler}
For $z\in\mathbb{C}$ with $\Re(z)>0$ \[\Gamma(z)= \int_{0}^{
\infty} t^{z-1}\exp(-t) dt,\] and functional equation for
$z\notin -\mathbb{N}$ \[\Gamma(z+1)=z\Gamma(z).\]
\end{definition}

From the $p$-adic perspective the aim is to find a $p$-adically
continuous function extending the factorial function. One
solution was found by Morita which extended the restricted
factorial \[ n!_{p}:=\prod_{1\leq j <n, p\nmid n}j.\] The
interpolated function is the Morita gamma function.

\begin{definition}[Morita Gamma Function]\label{mori}
This is the $p$-adically continuous function \[\Gamma_{p}:
\mathbb{Z}_{p}\rightarrow\mathbb{Z}_{p}\] that extends
$f(n):= (-1)^{n}n!_{p}$ for $n\geq 2$. \end{definition}

This is quite an attractive function to try and view from a
nonstandard perspective. Indeed, define for all
$n\in\mathbb{^{\ast}N}$ with $n\geq2$ \[ ^{\ast}\Gamma_{p} :
\mathbb{^{\ast}N}\rightarrow \mathbb{^{\ast}Z}, \qquad
^{\ast}\Gamma_{p}(n)= ( -1)^{n}\prod_{1\leq j < n, p\nmid n}
j,\] with $^{\ast} \Gamma(0) =-\text{} ^{\ast} \Gamma(1) = 1$.
By applying Wilson's theorem in a nonstandard setting this hyper
function is $p$-adically Q-continuous. It also satisfies the
functional equation \[ ^{\ast}\Gamma_{p}(n+1)=\text{}
^{\ast}h_{p}(n)\text{} ^{\ast}\Gamma_{p}(n),\] where \[
^{\ast}h_{p}(n)=\left\{ \begin{array}{ll} -n & \mbox{$p\nmid
n$,}\\ -1 & \mbox{$p|n$.} \end{array} \right.\] As the function
lies in $\mathbb{^{\ast}Z}$ the $p$-adic shadow map can be taken
leading to \[ \operatorname{sh}_{p}(^{\ast}\Gamma_{p}):
\mathbb{Z}_{p}\rightarrow\mathbb{Z}_{p},\] and it satisfies
$\operatorname{sh}_{p}(^{\ast}\Gamma_{p}(n))= \Gamma_{p} (n)$
for all $n\in\mathbb{N}$. By the properties of the shadow map
$\operatorname{sh}_{p}(^{\ast}\Gamma_{p})$ is $p$-adically
continuous. Thus as these two $p$-adically continuous functions
agree on a dense subset of $\mathbb{Z}_{p}$ they must be
equal. Alternatively one can take this construction as the
definition of the Morita gamma function. So one sees that the set
$\mathbb{^{\ast}N}\setminus\mathbb{N}$ determines the non natural
$p$-adic values of the Morita gamma function.

\begin{lemma}\label{shga}
\[\operatorname{sh}_{p}(^{\ast}\Gamma_{p})=\Gamma_{p}.\]
\end{lemma}

It can easily be checked that the $p$-adic shadow map preserves
all the properties of $\operatorname{sh}_{p}(^{\ast}\Gamma_{p})$.
For example taking the $p$-adic shadow map of $^{\ast}h_{p}(x)$
leads to $h_{p}(x)$ where \[
h_{p}(x)=\left\{ \begin{array}{ll} -x &
\mbox{$x\in\mathbb{Z}_{p}^{\times}$,}\\ -1 & \mbox{$x \in p
\mathbb{Z}_{p}$,} \end{array} \right. \] and $\Gamma_{p}(x+1)=
h_{p}(x)\Gamma(x)$ for all $x\in\mathbb{Z}_{p}$. In many ways it
is a lot easier to work with a function in the integers or hyper
integers than in $\mathbb{^{\ast}Z}_{p}$.

\section{The Kubota-Leopoldt Zeta Function}

For many years the value of the Riemann zeta function at integer
points has attracted many mathematicians. In particular the value
at negative integers lead to work in the $p$-adic area and the
discovery of the $p$-adic Riemann zeta function of Kubota and
Leopoldt in 1964 (\cite{Kub-L}).

The value of $\zeta_{\mathbb{Q}}(s)$ at
negative integers was originally derived by Euler using divergent
series. These values are related to the Bernoulli numbers
($B_{n}$) which are given by
\[\frac{t}{\exp(t)-1}=\sum_{k\in\mathbb{N}}B_{k}\frac{t^{k}}{k!}.\]
\begin{lemma}\label{bber} For all $k\in\mathbb{N}$,
\[\zeta_{\mathbb{Q}}(-k)=(-1)^{k}\frac{B_{k+1}}{k+1}.\] In
particular as $B_{k}=0$ for $k>1$ and odd, ($k\in\mathbb{N}$)
\[\zeta_{\mathbb{Q}}(-1-k)=-\frac{B_{k+2}}{k+2}.\]
\end{lemma}

The idea of $p$-adic interpolation is obtain a $p$-adic function
which has similar properties to the original real or complex
function. The aim is to learn more about the original function by
recasting the problem in the $p$-adic world which can often be
simpler. In the case of the Riemann zeta function the
interpolation is based on its values at negative integers. In
order for the resulting $p$-adic function to be $p$-adically
continuous the Riemann zeta function is modified by removing
the $p$-Euler factor. Let \[
\zeta_{\mathbb{Q},p}(s)=(1-p^{-s})\zeta_{\mathbb{Q}}(s),\] then
for all $k\in\mathbb{N}$
\[\zeta_{\mathbb{Q},p}(-k)=-(1-p^{k})\frac{B_{k+1}}{k+1}.\] It
should be noted than the factor $(-1)^{k}$ is not needed because
when $k=0$ the $p$-Euler factor is equal to zero and so the
modified zeta function is also zero.

There are now several ways to develop the $p$-adic Riemann zeta
function; Kubota-Leopoldt, Iwasawa theory, and Mazur measures and
$p$-adic integration. Each $p$-adic zeta function has $p-1$
branches. Initially define for $n\in\mathbb{N}\setminus \{0\}$
\[\zeta_{p}(1-n)=(1-p^{n-1})\zeta_{\mathbb{Q},p}(1-n)=
-(1-p^{n-1})\frac{B_{n}}{n}.\] The $p$-adic zeta function
interpolates these numbers.

\begin{definition} \label{zeta}Let $p\geq 5$ be a prime and fix
$s_{0}\in \{0,1,\ldots,p-2\}$. Then define the $p$-adic zeta
function by
\begin{align}
\zeta_{p,s_{0}}:\mathbb{Z}_{p}&\rightarrow
\mathbb{Q}_{p},\notag\\
s& \mapsto \lim_{t_{\alpha}}\zeta_{p}(
1-(s_{0}+(p-1)t_{\alpha})).\notag
\end{align}
Here $s$ is a $p$-adic integer with $\{t_{\alpha}\}_{\alpha\geq
1}$ a sequence of natural numbers which $p$-adically converges to
$s$. The case $s_{0}=0=s$ is excluded and dealt with below.
\end{definition}

It should be remarked that the $p$-adic zeta function
$\zeta_{p,s_{0}}(s)$ interpolates the zeta function
$\zeta_{\mathbb{Q},p}$ at negative integer values $s$ by
\[\zeta_{p,s_{0}}(s)=\zeta_{p}(1-n),\] where $n\equiv s
(\operatorname{mod} p-1)$ and $n=s_{0}+s(p-1)$. The continuity
properties are deduced from the Kummer congruences.
 
\begin{theorem} \label{kummer}
Suppose $i\in\mathbb{N}$, $m\geq 2$ and $p$ a prime with
$(p-1)\nmid j$. if $i\equiv j (\operatorname{mod} p^{n}(p-1))$
then \[ (1-p^{i-1})\frac{B_{i}}{i}\equiv
(1-p^{j-1})\frac{B_{j}}{j} (\operatorname{mod}p^{N+1}).\]
\end{theorem}

For some fixed $s_{0}$ let $s,t\in\mathbb{N}\setminus \{0\}$ with
$s\equiv t (\operatorname{mod} p^{N})$ then for some
$k\in\mathbb{N}$ $s=t+kp^{N}$. Now let $i=s_{0}+s(p-1)$ and
$j=s_{0}+t(p-1)=s_{0}+s(p-1)+kp^{N}(p-1)$ then
$i\equiv j\operatorname{mod} p^{N}(p-1)$ and by the Kummer
congruences \[  (1-p^{i-1})\frac{B_{i}}{i}\equiv
(1-p^{j-1})\frac{B_{j}}{j} (\operatorname{mod}p^{N+1}).\]
Hence,
\[\zeta_{p,s_{0}}(s)\equiv \zeta_{p,s_{0}}(t) (\operatorname{mod}
p^{N+1}),\]
and $\zeta_{p,s_{0}}$ is uniformly continuous on $\mathbb{Z}_{p}$
since $\mathbb{N}\setminus \{0\}$ is dense in $\mathbb{Z}_{p}$. This also
proves the uniqueness of the $p$-adic zeta function by the
interpolation property.

One should note that in the case when $s_{0}\in \{ 1,3,\ldots
p-2\}$ (the odd congruence classes) gives the zero function since
for such $s_{0}$, $B_{s_{0}+(p-1)k}=0$. This only leaves the even
congruence classes and the case $s_{0}=0$. The latter case will
enable the zeta function to be defined for $p=2$ and $p=3$ since
the only congruence class is $s_{0}=0$.

Let $s_{0}=0$ then for $s\neq 0$ definition \ref{zeta} can be
applied. For $s=0$
\[ \zeta_{p,0}:\mathbb{Z}_{p}\rightarrow \mathbb{Q}_{p},\quad
\zeta_{p,0}(0)=\lim_{t_{\alpha}\rightarrow
0}(\zeta_{p}(1-(p-1)t_{\alpha})=\zeta_{p}(1).\] This can be
viewed as a pole of the $p$-adic zeta function as in the case of
the Riemann zeta function. This can be seen more clearly in terms
of $p$-adic integration.

\subsection{A Nonstandard Construction}\label{621}
\begin{definition}\label{hypbb}
\begin{enumerate}

$\zeta_{\mathbb{^{\ast}Q}}(s)=\sum_{n\in\mathbb{^{\ast}N} \setminus
\{0\}}n^{-s}.$
\item The hyper Bernoulli numbers are defined by
$\frac{t}{^{\ast}\exp(t)-1}=\sum_{k\in\mathbb{^{\ast}N}}
\text{}^{\ast}B_{k}t^{k}/k!.$
\end{enumerate}
\end{definition}

\begin{lemma}\label{hypzb}
For $k\in\mathbb{^{\ast}N}$ and $k>1$
\[\zeta_{\mathbb{^{\ast}Q}}(1-k)=-\frac{^{\ast}B_{k}}{k}.\]
\end{lemma}

The proof of this result can derived in an identical manner to
the classical case. In a similar manner to the classical
construction of the previous section the $p$-Euler factor can be
removed for some fixed standard prime, $p$. Define
\[\zeta_{\mathbb{^{\ast}Q},p}(s)=(1-p^{-s})\zeta_{\mathbb{^{\ast}Q}},\]
and
\[\zeta_{\mathbb{^{\ast}Q},p}(1-k)=-(1-p^{k-1})\frac{^{\ast}B_{k}}{k}.\]

Now define a nonstandard function.
\begin{align}
^{\ast}f: \mathbb{^{\ast}N}&\rightarrow
\mathbb{^{\ast}Q},\notag\\
k&\mapsto
-(1-p^{k})\frac{^{\ast}B_{k+1}}{k+1}=\zeta_{\mathbb{^{\ast}Q},p}(-1-k).\notag
\end{align}

\noindent This function can be used for the $p$-adic continuation
of the Riemann zeta function as in the standard case. Indeed let
$\sigma_{0}\in \{-1,0,1,2,\ldots,p-3\}$ and
$\sigma\in\mathbb{^{\ast}N}$. Let $p\geq 5$ and $\sigma_{0}\neq
-1$ \begin{align} ^{\ast}f_{\sigma_{0}}:
\mathbb{^{\ast}N}&\rightarrow \mathbb{^{\ast}Q},\notag\\
\sigma&\mapsto -(1-p^{\sigma_{0}+\sigma(p-1)})
\frac{^{\ast}B_{\sigma_{0}+\sigma(p-1)+1}}{\sigma_{0}+\sigma(p-1)+1}
=\zeta_{\mathbb{^{\ast}Q},p}(-1-\sigma_{0}-\sigma(p-1)).\notag
\end{align}

For $\sigma_{0}$ even the value of the Bernoulli number is zero
for all $\sigma\in\mathbb{^{\ast}N}$. So only odd $\sigma_{0}$
needs to be considered. This function is $p$-adically
continuous for a fixed $\sigma_{0}$. Indeed let
$\sigma\equiv \tau (\operatorname{mod} p^{N})$,
$\tau=\sigma+kp^{N}$ (for some $k,N\in\mathbb{^{\ast}N}$.)
The Kummer congruences also carry through in the nonstandard
setting and imply $^{\ast}f_{\sigma_{0}}(\sigma) \equiv$$^{\ast}f_{\sigma_{0}}( \tau)
(\operatorname{mod}p^{N+1}).$ The Kummer congruences also imply
$|^{\ast}B_{k}/k|_{p}\leq 1$ for $p-1\nmid k$ and so
$^{\ast}f_{\sigma_{0}}(\sigma)\in\mathbb{^{\ast}Q}^{\lim_{p}}$
for all $\sigma\in\mathbb{^{\ast}N}$. This enables the $p$-adic
shadow map to be taken. Using work in Goldblatt (\cite{Go},
chapter 18)
\[\operatorname{sh}_{p}(\mathbb{^{\ast}N})=\mathbb{Z}_{p}.\]
Thus \[\operatorname{sh}_{p}(^{\ast}f_{\sigma_{0}}(\sigma))=
f_{p,\sigma_{0}}(\operatorname{sh}_{p}(\sigma)): \mathbb{Z}_{p}
\rightarrow \mathbb{Q}_{p},\] with
\[f_{p,\sigma_{0}}(n)=-(1-p^{\sigma_{0}+n(p-1)})
\frac{B_{\sigma_{0}+n(p-1)+1}}{\sigma_{0}+n(p-1)+1}, \forall
n\in\mathbb{N}.\]

The $p$-adic shadow map preserves continuity so
$f_{p,\sigma_{0}}$ is $p$-adically continuous on
$\mathbb{Z}_{p}$. Moreover \[f_{p,\sigma_{0}}(n)=
\zeta_{p,\sigma_{0}+1}(n), \forall n\in\mathbb{N}.\]

As these continuous functions agree on a dense set in
$\mathbb{Z}_{p}$
\[f_{p,\sigma_{0}}=\zeta_{p,s_{0}}\]
where $s_{0}=\sigma_{0}+1$.

This leaves the case $\sigma_{0}=-1$. For $\sigma\neq 0$ the same
definition of the nonstandard function above can be used to
define $^{\ast}f_{-1}$ for all primes. So let $r=n-1$
then for all $r\in\mathbb{^{\ast}N}$, $^{\ast}f_{-1}(r)
=-(1-p^{-1+(r+1)(p-1)})^{\ast}B_{(r+1)(p-1)}{(r+1)(p-1)}.$ From
classical results transferred to the nonstandard setting this
function is $p$-adically continuous and
$|^{\ast}f_{p,-1}(r)|_{p}\leq 1$ for all $r\in\mathbb{^{\ast}N}$
thus the shadow map can be taken and a function is obtained

\[\operatorname{sh}_{p}(^{\ast}f_{-1}(\sigma))=
f_{p,-1}(\operatorname{sh}_{p}(\sigma)): \mathbb{Z}_{p}
\rightarrow \mathbb{Q}_{p},\]
with \[f_{p,-1}(n)=-(1-p^{-1+n(p-1)})
\frac{B_{n(p-1)}}{n(p-1)}, \forall
n\in\mathbb{N}\setminus \{0\}.\]

Then $f_{p,-1}(n)=\zeta_{p,0}(n)$ for all $n\in\mathbb{N}\setminus
\{0\}.$ As these functions agree are on a dense set in
$\mathbb{Z}_{p}$ and are continuous then they are equal, proving the following
theorem.
\begin{theorem}\label{klzetash21}
For a fixed $\sigma_{0}\in\{-1,1,3,\ldots p-3\}$ \[\operatorname{sh}_{p}
((^{\ast}f_{\sigma_{0}}(\sigma)) = \zeta_{p,\sigma_{0}+1}( \operatorname{sh}_{p}
(\sigma)).\] \end{theorem}

\section{Double Interpolation}

The ideas and techniques behind $p$-adic interpolation have been
discussed in various parts throughout this work. The previous
sections show that some forms of $p$-adic interpolation can be
viewed from a nonstandard perspective. This interpretation lends
itself to extending interpolation. In particular the possibility
of interpolating a function or set of numbers with certain
properties with respect to more than one prime. Initially the
case of interpolating with respect to two distinct primes $p$
and $q$ will be considered. The range of functions which can be
interpolated is reduced because of conditions required with
respect to two primes rather than a single prime.

In the standard world double interpolation does not seem
possible. The main reason being that an interpolated function
would have to have a domain consisting of $\mathbb{Z}_{p}$ and
$\mathbb{Z}_{q}$ and possibly a range of values in
$\mathbb{Q}_{p}$ and $\mathbb{Q}_{q}$. It appears very difficult
to achieve this because the fields $\mathbb{Q}_{p}$ and
$\mathbb{Q}_{q}$ are not isomorphic. In a nonstandard setting
these problems are removed since one can just consider the spaces
$\mathbb{^{\ast}N}$ (or $\mathbb{^{\ast}Z}$) and
$\mathbb{^{\ast}Q}$. This is the first restriction on functions
which can be interpolated, one has functions of the form
$f:\mathbb{N}\rightarrow \mathbb{Q}$ and not with values in
$\mathbb{Q}_{p}$ or $\mathbb{Q}_{q}$ like standard interpolation.
The next condition relates to the continuity. For any hope of
double interpolation $f$ must be uniformly continuous and bounded
when considered in the $p$-adic and $q$-adic valuation. This
ensures any interpolated function is a continuous extension of
$f$. Given these conditions one constructs the nonstandard
extension of $f$, $^{\ast}f:\mathbb{^{\ast}N}\rightarrow
\mathbb{^{\ast}Q}$. Since the original function was assumed to be
bounded in both valuations the actual function is \[ ^{\ast}f :
\mathbb{^{\ast}N}\rightarrow \mathbb{^{\ast}Q}^{\lim_{p}} \cap
\mathbb{^{\ast}Q}^{\lim_{q}}.\] This also enables the shadow maps
to be taken leading to a $p$-adic function $f_{p}: \mathbb{Z}_{p}
\rightarrow \mathbb{Q}_{p}$ and a $q$-adic function $f_{q}:
\mathbb{Z}_{q} \rightarrow \mathbb{Q}_{q}$ These correspond to
the single interpolation of $f$. This can de displayed
graphically

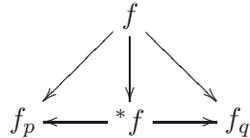
\begin{displaymath}
\xymatrix{ & f \ar[d] \ar[dl] \ar[dr] & \\
f_{p} & ^{\ast}f \ar[l] \ar[r] & f_{q} }
\end{displaymath}

Why may double interpolation be of any use? The main reason is
the same as for single prime interpolation which is often
specific to the problem being considered. In general the aim is
to learn more about a function from the $p$-adic function which
share some special values (though sometimes they are twisted).
Often the $p$-adic function is easier to examine and so enables
the original problem from various aspects. By using double
interpolation one hopes that another aspect is added to the
methods of looking at a problem.

\section{Riemann Zeta Function}

In this section attempts are made to double interpolate the
Riemann zeta function. Ideally one would like to choose the set
\[ \{ (1-p^{m})(1-q^{m})\zeta_{\mathbb{^{\ast}Q}}(-m)
\}_{m\in\mathbb{^{\ast}N}} .\] This seems a sensible choice based
on the nonstandard work on the Kubota-Leopoldt zeta function, by
 removing the $p$ and $q$ Euler factors. This also is symmetric
with respect to the prime factors. In the single prime case
the continuity was deduced from the Kummer congruences. These can
be extended slightly, in a symmetrical way to take account of the
extra Euler factor. Recall the Kummer congruences

\begin{theorem}\label{kumkum}
If $(p-1)\nmid i$ and $i\equiv j\pmod{p^{n}(p-1)}$ then \[
(1-p^{i-1})\frac{B_{i}}{i}\equiv (1-p^{j-1})\frac{B_{j}}{j}\pmod{
p^{n+1}}.\]
\end{theorem}

As $q$ is distinct from $p$ and still assume the conditions in
the theorem

$p\nmid (q^{i-1}-1)$ and so by Euler's theorem \[
q^{i-1}-1 \equiv q^{j-1}-1 \pmod{p^{n+1}}.\] The one case when
this is not true is when $(p-1)|(i-1)$ and the congruence holds
trivially. Combining this observation with the Kummer congruences
gives

\begin{corollary}\label{kumkumkum}
If $(q,p)=1$, $(p-1)\nmid i$

and $i\equiv j\pmod{p^{n}(p-1)}$ then \[
(1-q^{i-1})(1-p^{i-1})\frac{B_{i}}{i}\equiv
(1-q^{j-1})(1-p^{j-1})\frac{B_{j}}{j}\pmod{ p^{n+1}}.\]
\end{corollary}

This can be made symmetrical by including further conditions on
$i$ and $j$.

\begin{corollary}\label{kkm}
If $(q,p)=1$, $(p-1)\nmid i$,

$(q-1)\nmid i$

, $i\equiv j\pmod{q^{n}(q-1)}$ and $i\equiv j\pmod{p^{n}(p-1)}$
then \[ (1-q^{i-1})(1-p^{i-1})\frac{B_{i}}{i}\equiv
(1-q^{j-1})(1-p^{j-1})\frac{B_{j}}{j}\pmod{ p^{n+1}},\] and \[
(1-q^{i-1})(1-p^{i-1})\frac{B_{i}}{i}\equiv
(1-q^{j-1})(1-p^{j-1})\frac{B_{j}}{j}\pmod{ q^{n+1}},\]
\end{corollary}

As in the single prime case branches are considered. Firstly
define

\begin{align}
^{\ast}f_{p,q}: \mathbb{^{\ast}N}&\rightarrow\mathbb{^{\ast}Q},
\notag\\ k&\mapsto -(1-p^{n-1})(1-q^{n-1})\frac{^{\ast}B_{n}}{n}.
\notag \end{align}

\noindent This function is not Q-continuous with respect to
either prime. Let $p,q>5$ and without loss of generality assume
$p>q$ and let $\sigma_{0}\in \{ -1,0,1,2,\ldots (p-1)(q-1)-2 \}$.
Due to the symmetry only the $p$ case will be considered. Let
$\sigma_{0}\notin \{-1,p-1, 2(p-1),\ldots (q-2)(p-1)\}$ then
define
\begin{align}
^{\ast}f_{p,q,\sigma_{0}}: \mathbb{^{\ast}N}&\rightarrow
\mathbb{^{\ast}Q},\notag\\ \sigma&\mapsto
-(1-p^{\sigma_{0}+\sigma(p-1)(q-1)})(1-q^{\sigma_{0}+
\sigma(p-1)(q-1)})
\frac{^{\ast}B_{\sigma_{0}+\sigma(p-1)(q-1)+1}}{\sigma_{0}
+\sigma(p-1)(q-1)+1}.\notag
\end{align}

For fixed $\sigma_{0}$ as above this function is $p$-adically
Q-continuous. The extended Kummer congruences then give this
result. The Kummer congruences also give $|^{\ast}B_{k}/k|_{p}
\leq 1$ for $(p-1)\nmid k$ which shows the function actually lies
in $\mathbb{^{\ast}Q}^{\lim_{p}}$. The shadow map leads to a
$p$-adic function.

For $\sigma_{0}\in \{ -1,0,1,2,\ldots (p-1)(q-1)-2 \}$ with
$\sigma\neq 0$ the same definition of the nonstandard function
can be given. The case $\sigma=0$ corresponds to a pole.

So with the double interpolated zeta function defined one also
has new $p$-adic functions defined. The next section tries to
gain information about these functions by using $p$-adic
measures. As a final note on this section it is clear that this
method generalises to a finite set of primes for interpolation.

\subsection{Katz's Theorem}
Katz's theorem on $p$-adic measures follows directly from the
work on $p$-adic interpolation by Mahler. A proof can be
found in chapter 3 of \cite{Hi}.

\begin{definition}
A $\mathbb{Q}_{p}$-linear map
$\phi:\mathcal{C}_{p}(\mathbb{Z}_{p},
\mathbb{Q}_{p}) \rightarrow\mathbb{Q}_{p}$ (where $
\mathcal{C}_{p}(\mathbb{Z}_{p}, \mathbb{Q}_{p})$ is defined in
\ref{puni}) is called a bounded $p$-adic measure if there exists
a constant $B\geq 0$ such that $|\phi(f)|_{p}\leq B|f|_{p}$ for
all $f\in\mathcal{C}_{p}(\mathbb{Z}_{p}, \mathbb{Q}_{p})$.
\end{definition}

\begin{theorem}\label{k1}
Consider a bounded sequence of numbers $\{b_{n}\}$ in
$\mathbb{Q}_{p}$. Let $f:\mathbb{Z}_{p}\rightarrow\mathbb{Q}_{p}$
be a $p$-adically continuous function and $a_{n}(f)$ the
Mahler coefficients. Then a uniquely bounded $p$-adic measure
$\phi$ can be defined by \[ \int_{\mathbb{Z}_{p}}f
d\phi=\sum_{n\in\mathbb{N}}b_{n}a_{n}(f).\]
Moreover the measure satisfies \[
\int_{\mathbb{Z}_{p}}\binom{x}{n}d\phi = b_{n},\] for all
$n\in\mathbb{N}$. All bounded measures on $\mathbb{Z}_{p}$ are
obtained in this way. \end{theorem}

\begin{corollary}\label{k2}
Each bounded $p$-adic measure having values in $\mathbb{Q}_{p}$
is uniquely determined by its values at all monomials $x^{m}$
($m\in\mathbb{N}$).
\end{corollary}

One of the applications of this theorem is to develop the
$p$-adic measure of the Riemann zeta function and again \cite{Hi}
is the reference.

\begin{theorem}\label{k3}
Let $\zeta_{\mathbb{Q}}$ be the Riemann zeta function. Also let
$a\in\mathbb{N}$ with $a\geq 2$ and $(a,p)=1$. Then for all
$m\in\mathbb{N}$
\[ (1-a^{m+1})\zeta_{\mathbb{Q}}(-m)=\left(
t\frac{d}{dt}\right) ^{m} \Psi(t)\mid_{t=1}.\]
where $\Psi(t)=(1-t^{a})^{-1}\sum_{b=1}^{a}\xi(b)t^{b}$ and

\begin{align}
\xi:\mathbb{Z}&\rightarrow \mathbb{Z}\notag, \\
n&\mapsto \left\{ \begin{array}{ll} 1 & \mbox{$a\nmid n$,}\\ 1-a
& \mbox{$a \mid n$.} \end{array} \right. \end{align}
\end{theorem}

\begin{corollary}\label{k4}
Let $a\in\mathbb{N}$ with $a\geq 2$ and $(a,p)=1$. Then there
exists a unique $p$-adic measure $\zeta_{a}$ on
$\mathbb{Z}_{p}$ such that  \[
\int_{\mathbb{Z}_{p}}x^{m}d\zeta_{a}=(1-a^{m+1})\zeta_{\mathbb{Q}}(-m),\]
for all $m\in\mathbb{N}$.
\end{corollary}

The central result in this section is generalizing theorem
\ref{k3} and corollary \ref{k4} with the following results.

\begin{theorem}\label{bc1}
Let $a\in\mathbb{N}$ with $a\geq 2$ and $(a,p)=1$. Also let
$r\in\mathbb{N}$, $(r,p)=1$. Then for all $m\in\mathbb{N}$
\[ (1-a^{m+1})r^{m}\zeta_{\mathbb{Q}}(-m)=\left( t\frac{d}{dt}
\right) ^{m}\Psi_{r}(t)\mid_{t=1}.\]
Here $\Psi_{r}(t)=(1-t^{ra})^{-1}\sum_{b=1}^{a}\xi_{r}(br)t^{br}$
and

\begin{align}
\xi_{r}:\mathbb{Z}&\rightarrow \mathbb{Z}\notag, \\
n&\mapsto \left\{ \begin{array}{ll}
0& \mbox{$r\nmid n$,}\\
 1 & \mbox{$r\mid n$, $ra\nmid n$,}\\ 1-a
& \mbox{$r\mid n$, $ra \mid n$.} \end{array} \right. \end{align}
\end{theorem}

From this theorem \ref{k3} is a special case of the above when
$r=1$. To begin the proof the following short lemma is needed.

\begin{lemma}\label{bc2}
\[\sum_{b=0}^{a-1}\xi_{r}(br)=\sum_{b=1}^{a}\xi_{r}(br)=0.\]
\end{lemma}

\noindent \textbf{Proof:}
From the definition, $\xi_{r}(b)\neq 0$ only for the multiples of
$r$. Hence the sum will consist of only $a$ terms.
\[\sum_{b=1}^{a}\xi_{r}(rb)=\xi_{r}(r)+\xi_{r}(2r)+\ldots +
\xi_{r}((a-1)r) +\xi_{r}(ar).\] Only the final term is divisible
by $ar$ thus the other $a-1$ terms have value 1.
\[\sum_{b=1}^{a}\xi_{r}(br)= 1+1 +\ldots+ 1 +1 +(a-1)=0.\] The
same method can be used to show the second equality as $ra\mid
0$.

\noindent \textbf{Proof of theorem \ref{bc1}:}
Consider the polynomial
\[ \Phi_{r}(t)=\frac{t^{r}+1}{t^{r}-1} -
a\frac{t^{ar}+1}{t^{ar}-1}.\]
Using classical results the cotangent function satisfies
\[\pi \cot(\pi z)
=\frac{1}{z}-2\sum_{k\in\mathbb{N}}\zeta_{\mathbb{Q}}(2k)z^{2k-1}.\]
Letting $e(z)=\exp(2\pi iz)$ then
\begin{equation} \label{phi}
\Phi_{r}(e(z))=-(i\pi)^{-1}\sum_{k=1}^{\infty}2(1-a^{2k})\zeta_{\mathbb{Q}}(2k)
z^{2k-1}r^{2k-1}. \end{equation} The expression for $\Phi_{r}(t)$
can be rewritten using the $\xi_{r}$ function.

\begin{align}
\Phi_{r}(t)& =
\frac{(t^{r}+1)\left(\frac{t^{ar}-1}{t^{r}-1}\right) -
a(t^{ar}+1)}{t^{ar}-1}\notag,\\
&=\frac{(t^{r}+1)(1+t^{r}+t^{2r}+\ldots+t^{r(a-1)})-at^{ar}-a}{t^{ar}-1},\notag\\
&=\frac{
(1-a)+t^{ar}(1-a)+2(t^{r}+\ldots+t^{r(a-1)})}{t^{ar}-1},\notag\\
&=\frac{ (1-a)(t^{ar}+1) +2((t^{r}+\ldots+t^{r(a-1)})-(\xi_{r}(0)
+\xi_{r}(1)+\ldots+\xi_{r}((a-1)r)))}{t^{ar}-1},\notag\\
&=\frac{ (1-a)(t^{ar}-1)+2(
(t^{r}-1)+(t^{2r}-1)+\ldots+(t^{ar}-1))}{t^{ar}-1},\notag\\
&=\frac{(1-a)(1+t^{r}+t^{2r}+\ldots+t^{r(a-1)})+2\sum_{b=1}^{a-1}(1+t^{r}
+t^{2r}+\ldots+t^{r(b-1)})}{1+t^{r}+\ldots+t^{r(a-1)}},\notag\\
&=-(1-a)+2\frac{\sum_{b=1}^{a}\xi_{r}(br)(1+t^{r}+\ldots+t^{r(b-1)})}{1+t^{r}
+\ldots+t^{r(a-1)}},\notag \end{align}

Now let
\[
\Psi_{r}(t)=-\frac{\sum_{b=1}^{a}\xi_{r}(br)(1+t^{r}+\ldots+t^{r(b-1)})}{1+t^{r}
+\ldots+t^{r(a-1)}}=\frac{
\sum_{b=1}^{a}\xi_{r}(br)t^{rb}}{1-t^{ar}}.\] The last equality
follows from lemma \ref{bc2}.

\begin{equation}\label{phi1}
-2\Psi_{r}(t)=\Phi_{r}(t)-(a-1).
\end{equation}

The proof continues in two separate parts depending on whether
$m$ is odd or even.

For the even case use equation \ref{phi} and let
$v\in\mathbb{N}$ and $v\geq 1$.

\begin{align}
\left(\frac{d}{dz}\right) ^{2v}\Psi_{r}(e(z))&= -\frac{1}{2}
\left(\frac{d}{dz}\right) ^{2v}\Phi_{r}(e(z))\notag,\\
&=(2i\pi)^{-1}\sum_{k=1}^{\infty}2(1-a^{2k})\zeta_{\mathbb{Q}}(2k)r^{2k-1}(2k-1)
\ldots (2k-2v)z^{2k-2v-1}.\notag \end{align}

Therefore
\[\left(\frac{d}{dz}\right) ^{2v}\Psi_{r}(e(z))\mid_{z=0}=0.\]
By letting $t=\exp(2\pi iz)$ gives
\[\left(t\frac{d}{dt}\right) ^{2v}\Psi_{r}(t)\mid_{t=1}=0.\]
Since $\zeta_{\mathbb{Q}}(s)$ has trivial zeros at $s\in
-2\mathbb{N} (s\neq 0)$ the theorem is proved for even $m$.

For the odd case use equation \ref{phi} and let
$v\in\mathbb{N}$ and $v\geq 1$.

\begin{align}
-2\left(\frac{d}{dz}\right) ^{2v-1}\Psi_{r}(e(z))&=
\left(\frac{d}{dz}\right) ^{2v-1}\Phi_{r}(e(z))\notag,\\
&=-(i\pi)^{-1}\sum_{k=1}^{\infty}2(1-a^{2k})\zeta_{\mathbb{Q}}(2k)r^{2k-1}(2k-1)
\ldots (2k-2v)z^{2k-2v}.\notag
\end{align}

Therefore
\[-2\left(\frac{d}{dz}\right) ^{2v-1}\Psi_{r}(e(z))\mid_{z=0}=
-(i\pi)^{-1}(2v-1)!2(1-a^{2v})\zeta_{\mathbb{Q}}(2v)r^{2v-1}.\]

Using the functional equation for $\zeta_{\mathbb{Q}}$
\[ (1-a^{2v})\zeta_{\mathbb{Q}}(1-2v)=(2\pi
i)^{-2v}(2v-1)!2(1-a^{2v})\zeta_{\mathbb{Q}}(2v),\] and
substituting $t=\exp(2\pi iz)$ gives

\[\left( t\frac{d}{dt}\right)^{2v-1}\Psi_{r}(t)\mid_{t=1}
=(1-a^{2v})r^{2v-1}\zeta_{\mathbb{Q}}(-(2v-1)).\] The theorem is
then proved.
\begin{flushright} \textbf{$\Box$} \end{flushright}

An application of this theorem in combination with theorem
\ref{k1} leads to a slightly modified $p$-adic measure,
comparable to corollary \ref{k4} for the Riemann zeta function.

\begin{lemma}\label{bc3}
Let $a\in\mathbb{N}$ with $a\geq 2$ and $(a,p)=1$. Then there
exists a unique $p$-adic measure $\zeta_{a,p}$ on
$\mathbb{Z}_{p}$ having values in $\mathbb{Z}_{p}$ such that for
all $m\in\mathbb{N}$
\[ \int_{\mathbb{Z}_{p}}x^{m}d\zeta_{a,p,q}=
(1-a^{m+1})(1-q^{m})\zeta_{\mathbb{Q}}(-m).\]
\end{lemma}

\noindent \textbf{Proof:}
In the proof of theorem \ref{k3} Katz showed that
\[
(1-a^{m+1})\zeta_{\mathbb{Q}}(-m)=\left(t\frac{d}{dt}\right)^{m}\Psi(t)\mid_{t=1}.\]
Let $\binom{x}{n}=\sum_{k=0}^{n}c_{n,k}x^{k}$ with
$c_{n,k}\in\mathbb{Q}$. (These numbers are the Stirling numbers
of the second kind divided by $n!$.) The uniqueness of such a
modified $p$-adic Riemann zeta function measure is guaranteed by
corollary \ref{k2}. The sequence of numbers $\{
(1-a^{m+1})(1-q^{m})\zeta_{\mathbb{Q}}(-m)\}$ is bounded in
$\mathbb{Q}$ because the classical case by Katz gives that
$|(1-a^{m+1})\zeta_{\mathbb{Q}}(-m)|_{p}\leq1$ and as
$(1-q^{m})\in\mathbb{Z}$, $|(1-q^{m})|_{p}\leq 1$. In order to
show the existence of the measure it is required to show that $
\int_{\mathbb{Z}_{p}}\binom{x}{n}d\zeta_{a,p,q}$ is bounded.

\begin{align}
|\int_{\mathbb{Z}_{p}}\binom{x}{n}d\zeta_{a,p,q}|_{p}&=|\sum_{m=0}^{n}c_{n,m}
(1-a^{m+1})(1-q^{m})\zeta_{\mathbb{Q}}(-m)|_{p},\notag\\
&=|\left(\sum_{m=0}^{n}c_{n,m}
(1-a^{m+1})\zeta_{\mathbb{Q}}(-m)\right) -\left(
\sum_{m=0}^{n}c_{n,m}
(1-a^{m+1})q^{m}\zeta_{\mathbb{Q}}(-m)\right)|_{p},\notag\\ &
\leq \max \left( | \sum_{m=0}^{n} c_{n,m} (1-a^{m+1})
\zeta_{\mathbb{Q}}(-m)|_{p}, | \sum_{m=0}^{n}c_{n,m}
(1-a^{m+1})q^{m}\zeta_{\mathbb{Q}}(-m)|_{p} \right),\notag\\
&\leq\max\left(1, |\sum_{m=0}^{n}c_{n,m}
(1-a^{m+1})q^{m}\zeta_{\mathbb{Q}}(-m)|_{p}\right).\notag
\end{align} This last inequality follows from the classical
result in theorem \ref{k3}. So the proof now only requires the
last term to be bounded. Let $r\in\mathbb{N}$ and $(p,r)=1$.

\begin{align}
\sum_{m=0}^{n}c_{n,m}(1-a^{m+1})r^{m}\zeta_{\mathbb{Q}}(-m)
&=\sum_{m=0}^{n}c_{n,m}\left(
t\frac{d}{dt}\right) ^{m}\Psi_{r}(t)\mid_{t=1},\notag\\
&=\binom{t\frac{d}{dt}}{n}\Psi_{r}(t)\mid_{t=1},\notag\\
&:=\delta_{n}\Psi_{r}(t)\mid_{t=1}.\notag
\end{align}

Properties of the differential operator $\delta_{n}$ are known
via the classical results of Katz. The proofs can be found in
\cite{Hi}.

\begin{lemma}\label{k5}
\[ \delta_{n}=\frac{t^{n}}{n!}\frac{d^{n}}{dt^{n}}.\]
\end{lemma}

\begin{lemma}\label{k6}
Let
$R'=\{P(t)/Q(t):P(t),Q(t)\in\mathbb{Z}_{p}[t],|Q(1)|_{p}=1\}.$
Then $R'$ is a ring and stable under the action of $\delta_{n}$
for all $n\in\mathbb{N}$.
\end{lemma}

From the definition of $\Psi_{r}(t)$ it is apparent that
$\Psi_{r}(t)\in R'$ and by the previous lemma
$\delta_{n}\Psi_{r}(t)\in R'$. Therefore
\[ \delta_{n}\Psi_{r}=\frac{P_{r}}{Q_{r}},\] for
$P_{r}(t), Q_{r}(t)\in\mathbb{Z}_{p}[t]$ with $|Q_{r}(1)|_{p}=1$.
Thus $P_{r}(1)\in\mathbb{Z}_{p}$ and for all $r,n\in\mathbb{N}$
\[|\delta_{n}\Psi_{r}(t)|_{t=1}|_{p}=\frac{|P_{r}(1)|_{p}}{|Q_{r}(1)|_{p}}
\leq 1.\]
Hence in the case $r=q$
\[|\int_{\mathbb{Z}_{p}}\binom{x}{n}d\zeta_{a,p,q}|_{p}\leq 1.\]

\begin{flushright} \textbf{$\Box$} \end{flushright}

It should be noted that for the prime $q$ the conditions in the
corollary also prove that a $q$-adic measure $\zeta_{a,q,p}$
exists.

\begin{equation}\label{pqzeta} \int_{\mathbb{Z}_{q}}x^{m}d\zeta_{a,q,p}=
(1-a^{m+1})(1-p^{m}) \zeta_{\mathbb{Q}}(-m).\end{equation}

\subsubsection{Single Interpolation}
Without loss of generality let the following analysis be
completed with the prime $p$ (as the $q$ case is identical). The
reason that the new $p$-adic measure has been constructed above
is to try and find measures which corresponds to the single
interpolation of the double interpolated set of numbers.

\[\int_{\mathbb{Z}_{p}^{\times}}x^{m}d\zeta_{a}=
(1-a^{m+1})(1-p^{m}) \zeta_{\mathbb{Q}}(-m). \]

\begin{lemma}\label{bbc}
\[\int_{\mathbb{Z}_{p}^{\times}}x^{m}d\zeta_{a,p,q}=
(1-a^{m+1})(1-p^{m})(1-q^{m}) \zeta_{\mathbb{Q}}(-m). \]
\end{lemma}

In fact a slightly more general version will be proven.

\begin{theorem}\label{bbc1}
Let $\phi$ be a locally constant function on $\mathbb{Z}_{p}$
then \[\int_{\mathbb{Z}_{p}^{\times}}\phi(x)
x^{m}d\zeta_{a,p,q}= \sum_{n>1} \left((\phi(n)-
a^{m+1}\phi(na))n^{m} -
(\phi(nq)-a^{m+1}\phi(naq))n^{m}q^{m} \right).\]
\end{theorem}

\noindent \textbf{Proof:}
As $\phi$ is locally constant it is going to be constant on the
classes modulo $p^{k}$ for some $k\in\mathbb{N}$. Let $F\in R'$
and define $[\phi](F)$ by the Fourier inversion formula
\[ [\phi](F)(t)=\frac{1}{p^{k}}\sum_{b \pmod{p^{k}}} \phi(b)
\sum_{\zeta^{p^{k}}=1}\zeta^{-b}F(\zeta t).\] Then by \cite{Ka},
page 85, $[\phi](F)\in R'$ and  \[ \int_{\mathbb{Z}_{p}} f(x)
d\mu_{F} = \int_{\mathbb{Z}_{p}}f(x)d\mu_{[\phi]F},\] where for
$F\in R'$, $d\mu_{F}$ is the measure associated to $F$ via the
work above. In particular it has been shown above that associated
to the function $\Psi_{r}(t)\in R'$ is the measure $\zeta_{a,r}$
satisfying \[ \int_{\mathbb{Z}_{p}} x^{m}d\zeta_{a,r} =
(1-a^{m+1})r^{m}\zeta_{\mathbb{Q}}(-m).\]
\begin{align}
\int_{\mathbb{Z}_{p}}\phi(x)x^{m}d\zeta_{a,r} &=
\int_{\mathbb{Z}_{p}}\phi(x)x^{m}d\mu_{\Psi_{r}},\notag\\
&= \int_{\mathbb{Z}_{p}} x^{m} d\mu_{[\phi]\Psi_{r}}, \notag\\
=\left( t\frac{d}{dt} \right) ^{m}([\phi]\Psi_{r})|_{t=1}.\notag
\end{align}
Recall that $\Psi_{r}(t)=(\sum_{n=1}^{r}\xi_{r}(nr)t^{rn})/
(1-t^{ar})=\sum_{n=1}^{ar}\xi_{r}(n)t^{n})/
(1-t^{ar})$ which can be rewritten as \[\Psi_{r}(t)=
\frac{\sum_{n=1}^{arp^{k}}\xi_{r}(n)t^{n}}{1-t^{arp^{k}}}.\]
Then by simple calculation and using the definition of $\Psi_{r}$

\begin{align} [\phi]\Psi_{r}(t) &= \frac{\sum_{n=1}^{arp^{k}}
\phi(n)\xi_{r}(n)t^{n}}{1-t^{arp^{k}}},\notag\\
&=\sum_{n\geq 1}\phi(n)\xi_{r}t^{n},\notag\\
&=\sum_{n\geq 1}(\phi(nr)t^{nr}-a\phi(nar^{2})t^{nar}
 ) .\notag\\ \end{align}
 
\noindent Then  \[ \left(
t\frac{d}{dt}\right)^{m}([\phi]\Psi_{r})|_{t=1} = \sum_{n>1}
(\phi(nr)-a^{m+1}\phi(nar))n^{m}r^{m}.\]

Recall that \[ \int_{\mathbb{Z}_{p}}x^{m}d\zeta_{a,p,q}=
\left(t\frac{d}{dt}\right)^{m}(\Psi_{1}(t)-\Psi_{q}(t))|_{t=1},\]
then
\begin{align}
 \int_{\mathbb{Z}_{p}}\phi(x)
x^{m}d\zeta_{a,p,q} &=
\left(t\frac{d}{dt}\right)^{m}([\phi]\Psi_{1}(t)-
[\phi]\Psi_{q}(t))|_{t=1},\notag\\
&=\sum_{n>1} \left((\phi(n)-
a^{m+1}\phi(na))n^{m} -
(\phi(nq)-a^{m+1}\phi(naq))n^{m}q^{m} \right).\notag
\end{align}

\begin{flushright} \textbf{$\Box$} \end{flushright}

In particular choosing $\phi=\chi_{\mathbb{Z}_{p}^{\times}}$ (the
characteristic function of $\mathbb{Z}_{p}^{\times}$). Then as
$(a,p)=(q,p)=1$, $\phi(n)=\phi(naq)=\phi(nq)=\phi(na)$
and the above theorem reduces to
\begin{align}
 \int_{\mathbb{Z}_{p}}\chi_{\mathbb{Z}_{p}^{\times}}
x^{m}d\zeta_{a,p,q} &= \int_{\mathbb{Z}_{p}^{\times}}
x^{m}d\zeta_{a,p,q} ,\notag\\
&= (1-p^{m})(1-q^{m})(1-a^{m+1})\zeta_{\mathbb{Q}}(-m). \notag
\end{align}

\begin{definition}\label{pqqq}
Let $p,q$ be finite distinct primes in $\mathbb{Q}$ and for all
$k\in\mathbb{N}\setminus\{1\}$

\[
\zeta_{p,q}(1-k)=(1-p^{k-1})(1-q^{k-1})\zeta_{\mathbb{Q}}(1-k).\]
\end{definition}

By the corollary above

\[\zeta_{p,q}(1-k)=\frac{1}{1-a^{k}}\int_{\mathbb{Z}_{p}^{\times}}
x^{k-1} d\zeta_{a,p,q}.\]

It should be noted that the expression on the right does not
depend on $a$. Suppose $b\neq a$ with $(b,p)=(b,q)=1$ and $b\geq
2$. Then by the definition

\[\frac{1}{1-b^{k}}\int_{\mathbb{Z}_{p}^{\times}} x^{k-1}
d\zeta_{b,p,q}= \frac{1}{1-a^{k}}\int_{\mathbb{Z}_{p}^{\times}} x^{k-1}
d\zeta_{a,p,q}, \] since both equal
$(1-p^{k-1})(1-q^{k-1})\zeta_{\mathbb{Q}}(1-k)$.

\begin{definition}\label{pqzeta1}
Fix $s_{0}\in\mathbb\{0,1,2,\ldots,p-2\}.$ For
$s\in\mathbb{Z}_{p}$ ($s\neq 0$ if $s_{0}=0$) define
\[
\zeta_{p,q,s_{0}}(s)=\frac{1}{1-a^{s_{0}+(p-1)s}}\int_{\mathbb{Z}_{p}^{\times}}
x^{s_{0}+(p-1)s-1}d\zeta_{a,p,q}.\]
\end{definition}

If $k\in\mathbb{N}$ congruent to $s_{0}$ (mod $p-1$)
($k=s_{0}+(p-1)k_{0}$) then $\zeta_{p,q}(1-k)=\zeta_{p,q,s_{0}}
(k_{0}).$ In the case $s=0$ when $s_{0}=0$ the denominator
vanished. Using the above equality
$\zeta_{p,q}(1-k)=\zeta_{p,q,s_{0}}(k_{0})$ it follows that the
excluded case corresponds to $\zeta_{p,q}(1)$ and so has a pole
at $s=1$.

\begin{lemma}\label{bc5}
For $p,q$ fixed and $s_{0}$ as above the function
$\zeta_{p,q,s_{0}}$ is a continuous function of $s$ that does not
depend upon the choice of $a\in\mathbb{N}$ with $a\geq 2$ and
$(a,p)=(a,q)=1$.
\end{lemma}

\noindent \textbf{Proof:} The continuity is clear from standard
$p$-adic analysis. The independence of $a$ is established as
follows. Let $b\in\mathbb{N}$ with $b\geq 2$ and $(b,p)=(b,q)=1)$
then the two functions

\[\frac{1}{1-a^{s_{0}+(p-1)s}}\int_{\mathbb{Z}_{p}^{\times}}
x^{s_{0}+(p-1)s-1}d\zeta_{a,p,q} \] and
\[\frac{1}{1-b^{s_{0}+(p-1)s}}\int_{\mathbb{Z}_{p}^{\times}}
x^{s_{0}+(p-1)s-1}d\zeta_{b,p,q}\]

\noindent agree whenever $s_{0}+(p-1)s=k$ is an integer greater
than zero since in both cases the value of the function is
$(1-p^{k-1})(1-q^{k-1})\zeta_{\mathbb{Q}}(1-k)$. So both
functions agree on the set of non-negative integers which is
dense in $\mathbb{Z}_{p}$ and hence are equal.

\begin{flushright} \textbf{$\Box$} \end{flushright}

This is of course the function obtained by the shadow map of the double
interpolation function.

\subsubsection{Questions}

Initial attempts at proving theorem \ref{bc3} tried to use the
proof of theorem \ref{k1} directly. Instead of introducing the
new power series $\Psi_{r}(t)$ the power series $\Psi_{1}(t)$ was
used as in the work of Katz but introducing a different
differential operator. So in trying to prove the boundedness of
$\int_{\mathbb{Z}_{p}}\binom{x}{n}d\zeta_{a,r}$ the integral was
rewritten

\begin{align}
\int_{\mathbb{Z}_{p}}\binom{x}{n}d\zeta_{a,r}&=\sum_{m=0}^{n}
c_{n,m}r^{m}(1-a^{m+1})\zeta_{\mathbb{Q}}(-m),\notag\\
&=\sum_{m=0}^{n}c_{n,m}\left( rt\frac{d}{dt} \right) ^{m}
\Psi_{1}(t)|_{t=1},\notag\\
&:=\Delta_{n}\Psi_{1}(t)|_{t=1}.\notag
\end{align}

By the proof of lemma \ref{bc3} it is known that
$|\Delta_{n}\Psi_{1}(t)|_{t=1}|_{p}\leq 1$.

\noindent\textbf{Questions:}
\begin{itemize}
\item Does there exist a proof of lemma \ref{bc3} more $p$-adic
in nature? That is, following directly from the proof of theorem
\ref{k1} instead of having to work via power series evaluated at
$t=1$.
\item Does there exist a closed form for $\Delta_{n}$ for all
$n\in\mathbb{N}$?
\item Is $R'$ stable under $\Delta_{n}$ for all
$n\in\mathbb{N}$?
\end{itemize}

\subsubsection{The Measure on Sets}
There are several ways to view a $p$-adic measure. It has been
seen above that one way to define a measure is via a bounded
sequence. Two other methods include a definition by power
series in $\mathbb{Z}_{p}[[t]]$ and a definition by the open sets
of $\mathbb{Z}_{p}$. The power series can be easily deduced from
the bounded sequence method (for the details see section 3.5 of
\cite{Hi}). Of potentially more interest is the definition via
the open sets. For example corollary \ref{k4} establishes the
existence of a $p$-adic measure for the Riemann zeta function.
Work by Mazur establishes the existence of this measure from the
definition of open sets via the Bernoulli distributions.

Recall that in the $p$-adic topology all sets of the form
$a+p^{N}\mathbb{Z}_{p}$ are both open and closed. An integration
theory can be developed over such sets, see chapter 2 of
\cite{K}. An equivalent definition of a $p$-adic measure can be
established via sets and distributions.

\begin{definition}
Suppose $X$ is a compact-open subset of $\mathbb{Q}_{p}$.  A
$p$-adic distribution, $\mu$, on $X$ is an additive map from the
set of compact-open sets in $X$ to $\mathbb{Q}_{p}$. Thus, if
$U\subset X$ is a disjoint union of compact-open sets
$U_{1},\ldots,U_{n}$, then \[ \mu(U)=\mu(U_{1})+\cdots +
\mu(U_{n}).\]
\end{definition}

\begin{definition}
A distribution $\mu$ of $X$ is called a $measure$ if there is a
constant $B$ such that \[|\mu(U)|_{p}\leq B,\] for all
compact-open $U\subset X$.
\end{definition}

$\mathbb{Q}_{p}$ has a basis of open sets consisting of all
sets of the form $a+p^{N}\mathbb{Z}_{p}$ for $a\in\mathbb{Q}_{p}$
and $N\in\mathbb{Z}$.  Thus any open set of $\mathbb{Q}_{p}$ is a
union of open subsets of this type. Hence any distribution, and
measure, is determined by the values on such sets, see chapter 7
of \cite{RM} or chapter 2 of \cite{K} for information on this and
p-adic integration.

The aim of this section is to find a method which establishes the
action of a $p$-adic measure on open sets of $\mathbb{Z}_{p}$
directly from a measure defined by a bounded sequence. This
action is unique from classical results on $p$-adic measures. I
have such a method but it is entirely numerical in nature and
although gives a solution it does so only one set at a time
rather than give a general form of the measure.

From above the open sets of $\mathbb{Z}_{p}$ are of
the form $b+p^{n}\mathbb{Z}_{p}$ with $n\in\mathbb{N}$ and $0\leq
b< p^{n}$. Define a characteristic function of an open set.

\begin{definition}\label{characteristic} Let
$b+p^{n}\mathbb{Z}_{p}$be an open set of $\mathbb{Z}_{p}$ then
its characteristic function is given by
\[\chi_{b+p^{n}\mathbb{Z}_{p}}(x)= \left\{ \begin{array}{ll} 1 &
\mbox{$x\in b+\mathbb{Z}_{p}$,}\\ 0 & \mbox{$x\notin
b+\mathbb{Z}_{p}$.} \end{array} \right.\]
\end{definition}

Then by the definition of the $p$-adic integral with respect to a
$p$-adic measure $\mu$

\[\mu(b+p^{n}\mathbb{Z}_{p})=\int_{b+p^{n}\mathbb{Z}_{p}} d\mu =
\int_{\mathbb{Z}_{p}}\chi_{b+p^{n}\mathbb{Z}_{p}}(x)d\mu.\]

By defining a $p$-adic measure by an the integral of $x^{k}$ for
$k\in\mathbb{N}$ enables the integral of $p$-adic polynomials and
power series using classical results. Indeed let
$f(x)=\sum_{k\in\mathbb{N}}a_{k}(f)\binom{x}{k}$ be a continuous
function on $\mathbb{Z}_{p}$ in terms of its Mahler series. Let
$\{d_{k}\}$ be a bounded sequence used to define a $p$-adic
measure $\mu$ with $\int_{\mathbb{Z}_{p}}x^{k}d\mu=d_{k}$ for
all $k\in\mathbb{N}$. As in the first
section let $\binom{x}{k}=\sum_{m=0}^{k}c_{k,m}x^{m}$.

\begin{align}
\int_{\mathbb{Z}_{p}}f(x)d\mu=
&\int_{\mathbb{Z}_{p}}
\sum_{k\in\mathbb{N}}a_{k}(f)\binom{x}{k}d\mu,\notag\\
&=\sum_{k\in\mathbb{N}}a_{k}(f)\left( \int_{\mathbb{Z}_{p}}
\binom{x}{k}d\mu\right),\notag\\
&=\sum_{k\in\mathbb{N}}a_{k}(f)d_{k}.\label{conv}
\end{align}

The middle equality follows from classical results due to the
Mahler series converging (section 5, \cite{R}). Here $d_{k}$ is
also a bounded sequence by theorem \ref{k1}.
\[d_{k}=\int_{\mathbb{Z}_{p}}\binom{x}{k}d\mu=\sum_{m=0}^{k}c_{k,m}d_{k}.\]
It should be noted that equation \ref{conv} converges since by
Mahler's theorem $|a_{k}(f)|_{p}\rightarrow 0$ as
$k\rightarrow\infty$ and $d_{k}$ is bounded.

Returning to the characteristic function it is continuous and has
a Mahler series
\[\chi_{b+p^{n}\mathbb{Z}_{p}}(x)=\sum_{k\in\mathbb{N}}a_{k}(b,n)\binom{x}{k}.\]
Therefore
\begin{equation}\label{calc}
\mu(b+p^{n}\mathbb{Z}_{p})=\int_{\mathbb{Z}_{p}}\chi_{b+p^{n}\mathbb{Z}_{p}}
(x)d\mu =\sum_{k\in\mathbb{N}}a_{k}(b,n)d_{k}.\end{equation}

This expression enables $\mu(b+p^{n}\mathbb{Z}_{p})$ to be
numerically calculated for all open sets of $\mathbb{Z}_{p}$. By
some initial calculations it appears that in general extensive
calculations are needed to calculate the sequence $d_{k}$ and
then the coefficients $a_{k}(b,n)$. Unless either one or both of
these sets of numbers take a closed form or particularly simple
form then this method does not appear greatly practical.

As an example suppose one only knew that the Riemann zeta
function measure was defined by
\[\int_{\mathbb{Z}_{p}}x^{k}d\zeta_{a}=(1-a^{k+1})\zeta_{\mathbb{Q}}(-k).\]
Could it then be deduced that using the method outlined above
results in
\[\zeta_{a}(b+p^{n}\mathbb{Z}_{p})=\frac{1}{a}\left[\frac{ab}{p^{n}}\right]
+ \frac{(1/a)-1}{2}?\]
Here $[.]$ is the integral part function. This now leaves a
search for another method which, like the closed form for the
Riemann zeta function measure, may specific to a problem rather
than the general approach considered above.

The reason why such a process may be necessary is that it may

be interesting to see if a closed form of
$\zeta_{a,p,q}(b+p^{n}\mathbb{Z}_{p})$ and related measures
exists and in particular if they involves any Bernoulli
polynomials.

\section{A Double Morita Gamma Function}

An extension of the Morita gamma function is the double Morita
gamma function which is a $p$-adically and $q$-adically
Q-continuous function in $\mathbb{^{\ast}Z}$. A naive attempt
would be to define for odd distinct primes \[
\Gamma_{p,q}(n)=\prod_{1\leq j<n, p\nmid j, q\nmid j}j \quad
(j\geq 2).\] The problem is that this is not continuous in either
valuation. The reason is due to the generalized Wilson's theorem
not being applicable. Focussing on the prime $p$; by removing
elements divisible by $q$ one does not always have an inverse for
elements in $(\mathbb{Z}\setminus p^{r}\mathbb{Z})^{\times}$ so the
product of elements is not a unit in this group. So the search is
for the double Morita gamma function of the form (for $n\geq 2$)
\[ \Gamma_{p,q}(n)=C(n)\prod_{1\leq j <n, A(p,q,j)}j.\] Here
$A(p,q,j)$ is a set of conditions on $j$ depending on $p$ and
$q$. (As an example the naive attempt had $A(p,q,j)=\{p\nmid j,
q\nmid j\}$.) $A(p,q,j)$ must contain the conditions $p\nmid j$
and $q\nmid j$ otherwise there can be no form of continuity.

Initial attempts at finding an appropriate $A(p,q,j)$ revolve
around at trying to get some inverse residues for each $j$. So fix a
representation of $\mathbb{Z}/n\mathbb{Z}$ to be $\{ 1,2, \ldots, n\}$. To
explain, as $p,q\nmid j$ each such $j$ is invertible in $\mathbb{Z}/
p^{r}\mathbb{Z}$ ($r\in\mathbb{^{\ast}N}$) and $\mathbb{Z}/ q^{s}\mathbb{Z}$
($s\in\mathbb{^{\ast}N}$). The problem is that the inverse in $\mathbb{Z}/
p^{r}\mathbb{Z}$ may be divisible by $q$ and so will not be in the product
(similarly for the $p$-divisible case). So really one wants to include only
those $j$ such that for each $r,s\in\mathbb{^{\ast}N}$, the inverse residue of
$j$ in $\mathbb{Z}/ p^{r}\mathbb{Z}$ is not divisible by $q$ and in
$\mathbb{Z}/ q^{s}\mathbb{Z}$ is not divisible by $p$. Let
$S(p,q)\subset\mathbb{^{\ast}N}$ be the set of such $j$. Then for all $r$ and
$s$ \[ \prod_{1\leq j<p^{r}, j\in S(p,q)}j\equiv 1 \pmod{p^{r}} \quad \text{
and } \quad \prod_{1\leq j<q^{s}, j\in S(p,q)}j\equiv 1 \pmod{q^{s}}.\] The
immediate problem regards the elements in $S(p,q)$. Clearly $1\in S(p,q)$ but
are there any other? It is also known that if there is another element then
there are an infinite number of elements because of residue inverses for $j$ in
all the rings $\mathbb{Z}/ p^{r}\mathbb{Z}$ and $\mathbb{Z}/ q^{s}\mathbb{Z}$.
If it could be shown that $S(p,q)\equiv \{1\}$ for all $p$ and
$q$ then the double Morita gamma function (and all finite sets
of primes version) is trivial as the minimum conditions have
been imposed to gain continuity. If the sets are not trivial
then this would be a pleasant solution to an elementary number
theory problem but would not lead to a solution of the double
Morita gamma function directly because the resulting products
may not be continuous.

\section{Triviality}
\begin{theorem}\label{ivanko}
For all distinct primes $p$ and $q$, $S(p,q)\equiv \{1\}$.
\end{theorem}
\begin{corollary}
The double Morita gamma function is identical to 1.
\end{corollary}

\noindent\textbf{Proof:}
Let $j$ be an integer greater than $1$ such that for every
positive integer $r$ and $s$ the remainder of the inverse of $j
\pmod{p^{r}}$ is prime to $p$ and $r$ and the remainder of the
inverse of $j\pmod{q^{s}}$ is prime to $q$. Clearly $p\nmid j$
and $q\nmid j$. Then
\[ 1/j = a_{0}+a_{1}p+a_{2}p^{2}+\ldots \in \mathbb{Z}_{p},\]
with $a_{0}\in\{ 1,2,\ldots ,p-1\}$ and for all $i>0$,
$a_{0}\in\{ 0,1,2,\ldots, p-1\}$.

It is well known that the sequence $a_{0},a_{1},\ldots$ is
periodic so there exists an $n\in\mathbb{N}$ such that for all
$i\geq 0$, $a_{i+n}=a_{i}$.

Define \[ A_{mn}=\sum_{r=0}^{mn-1} a_{r}p^{r}.\]
Then by periodicity, \[ A_{mn}=(1+p^{n}+p^{2n}+ \ldots
+p^{n(m-1)})(a_{0}+a_{1}p+\ldots +a_{n-1}p^{n-1}),\] and
\[1/j \equiv A_{mn}\pmod{p^{mn}}.\]

Let $b\equiv p^{n}\pmod{q}$, then
\[ (1+p^{n}+\ldots + p^{n(m-1)})\equiv 1+b+b^{2}+ \ldots +
b^{m-1} \pmod{q}.\]

Suppose now that $b=1$ then in the case $m=q$, \[ q| 1+p^{n}+
\ldots p^{n(m-1)},\] and so is $A_{qn}$, a contradiction.

In the remaining cases, $b>1$, so for $m=q-1$, \[1+b+b^{2} +
\ldots + b^{q-1} = \frac{b^{q-1}-1}{b-1},\] which is divisible by
$q$ by Euler's theorem and hence so is $A_{n(q-1)}$, a
contradiction.

\begin{flushright} \textbf{$\Box$} \end{flushright}

\section{The Set $S(p,q)$}

Before proving the result in theorem \ref{ivanko} I looked at the
problem from a numerical perspective and I present some results
below which although elementary I have not seen in the
literature.

As an initial attempt to tackle the problem in theorem
\ref{ivanko} consider the integer 2. The inverses of 2 in
$\mathbb{Z}/ p^{r}\mathbb{Z}$ and $\mathbb{Z}/ q^{s}\mathbb{Z}$
are $(p^{r}+1)/2$ and $(q^{s}+1)/2$ respectively. The conditions
for 2 not to lie in $S(p,q)$ are $q\mid (p^{r}+1)/2$ and $p\mid
(q^{s}+1)/2$ for some $r,s\in\mathbb{^{\ast}N}$. As $p$ and $q$
are odd it follows that the conditions become $p^{r}\equiv
-1\pmod{q}$ and $q^{s}\equiv -1 \pmod{p}$ with $r<q$ and $s<p$.
So for any odd primes $p$ and $q$ if one of the two congruences
is soluble then $2\notin S(p,q)$ - this is a necessary
condition. Immediately one knows this is true when one of the
primes has a primitive root of the other. This leaves the case
when neither prime has the other as a primitive root.

To begin with choose a
primitive root $g$ for $p$. Then there exists integers
$1\geq j(q),j(-1)<p$ such that $g^{j(q)}\equiv q \pmod{p}$ and
$g^{j(-1)}\equiv -1\pmod{p}$. Putting these together \[
g^{j(q)s}\equiv g^{j(-1)}\pmod{p},\] \[j(q)s\equiv
j(-1)\pmod{p-1}.\]
This final equation is just a linear congruence which can be
solved for $s$ $\iff$ gcd$(p-1,j(q))\mid j(-1)$. In fact this can
not always be solved to give the actual value of $s$. Indeed let
the order of $q$ in $p$ be $2u$ (the order has to be even for a
solution to the necessary condition). Then by basic congruences
$q^{u}\equiv -1 \pmod{p}$ giving $s=u$. An example of where it
does not work is given by $p=5$ and $q=27$. Then $5^{27}\equiv 1
\pmod{109}$.

Considering just the classes modulo $p^{r}$; if the necessary
condition on $(p^{r}+1)/2$ is not satisfied then further
conditions can be investigated. For example the inverse of
$(p^{r}+1)/2$ in $\mathbb{Z}/p^{s}\mathbb{Z}$ also must not be
divisible by $q$ for $s>r$ (the case $s\leq r$ is already
covered by the original conditions since $(p^{r}+1)/2\equiv
(p^{s}+1)/2\pmod{p^{s}}$). The problem of finding these inverses
is given by the following lemma.

\begin{lemma}  \label{lem}
Let $s>r$ and let $n\in\mathbb{N}$ be maximal in satisfying
$nr<s$ (that is $(n+1)r\geq s$). Then the inverse, $x$, of
$(p^{r}+1)/2$ in $(\mathbb{Z}/ p^{s}\mathbb{Z})^{\times}$ is
given by \[ x=2((-1)^{n}p^{nr}+(-1)^{n-1}p^{(n-1)r}+\ldots
-p^{r}+1)+(1-(-1)^{n})(p^{s}/2),\] (if the leading term is
negative ($n$ is odd) then  $p^{s}$ is added to put the inverse
in the range $0<x<p^{s}$). \end{lemma}

\noindent\textbf{Proof:}
As $(p^{r}+1)/2\in (\mathbb{Z}p^{s}\mathbb{Z})^{\times}$ it
has an inverse $x$ satisfying $0<x<p^{s}$ and \[\frac{p^{r}+1}{2}
x\equiv 1 \pmod{p^{s}}.\] Let $x=2x_{1}$ then
\[(p^{r}+1)x_{1}\equiv 1 \pmod{p^{s}}.\] Now let $x_{1}=x_{2}+1$
then \[ p^{r}x_{2}+x_{2}+p^{r}\equiv 0 \pmod{p^{s}},\] which
implies the existence of a $k\in\mathbb{N}$ such that
\begin{equation}\label{zero} p^{r}x_{2}+x_{2}+p^{r} =kp^{s}.\end{equation} Viewing this
modulo $p^{r}$ shows that $x_{2}$ can be written as
$x_{2}=p^{r}x_{3}$
satisfying \begin{equation}\label{one} p^{r}x_{3}+x_{3}+1=kp^{s-r}.\end{equation}
Reducing this equation modulo $p^{r}$ leads to
$x_{3}=-1+x_{4}p^{r}$ satisfying \begin{equation}\label{two}
 p^{r}x_{4}+x_{4}-1=kp^{s-2r}.\end{equation} This can be reduced
modulo $p^{r}$ again to give $x_{4}=1+x_{5}p^{r}$ satisfying \begin{equation}
\label{three} p^{r}x_{5}+x_{5}+1=kp^{s-3r}.\end{equation} This is
of the same form as equation \ref{one}. In total $n$ reductions
modulo $p^{r}$ can be taken reducing to the equation \begin{equation}
\label{four}
p^{r}x_{n+2}+x_{n+2}+(-1)^{n-1}=kp^{s-nr}.\end{equation} and the
solution is currently in the form \[x=2(1-p^{r}+p^{2r} +\ldots +
(-1)^{mr}p^{mr} +\ldots +x^{n}p^{nr}).\] As $s-nr<r$ the equation
can only be reduced modulo $p^{s-nr}$. This gives a solution of
the form $x_{n}=(-1)^{nr}+x_{n+1}p^{s-nr}.$ Therefore the
solution is
\begin{align}
x &= 2(1-p^{r}+p^{2r} +\ldots +(-1)^{mr}p^{mr} +\ldots
+(-1)^{nr}p^{nr}+(-1)^{nr}x_{n+1}p^{s},\notag\\
&\equiv 2(1-p^{r}+p^{2r} +\ldots +(-1)^{mr}p^{mr} +\ldots
+(-1)^{nr}p^{nr} \pmod{p^{s}}.\notag
\end{align}
Induction can be trivially used to show that for any
$n\in\mathbb{N}$, \[|\sum_{m=0}^{n}(-1)^{mr}p^{mr}|<p^{nr}.\]
As $p$ is an odd prime $2p^{nr}<p^{nr+1}\leq p^{sr}$ and hence
$p^{nr}<p^{sr}/2$. Therefore $0<|x|<p{s}$ and it is clear that
for even $n$, $x>0$ and for odd $n$, $x<0$ so the solution is
given by $x+p^{s}$.

\begin{flushright} \textbf{$\Box$} \end{flushright}

The general method of solving
a linear congruence of the form $uv\equiv 1 \pmod{k}$ is by
Euclid's algorithm. One of the implications from the above lemma
is that if one tried to calculate an inverse of $(p^{r}+1)/2$ in
$(\mathbb{Z}/ p^{s}\mathbb{Z})^{\times}$ for specific $r$ and $p$
using Euclid's algorithm then one knows the parity of the number
of steps in the algorithm before the algorithm is even carried
out, it simply is the parity of $n$.

One now has more necessary conditions for $2\neq S(p,q)$. One
looks at the $q$ integrality of the inverse found in the lemma
(and similarly for the $q$ case and $p$-integrality). Of course
these inverses have inverses in $(\mathbb{Z}\setminus
p^{v}\mathbb{Z})^{\times}$ and so on. (For example more lemmas
can be given to find these inverses such as the inverse of
$2-2-p^{2}$ is $(p^{2}-2p-3)/4 +p^{2} +(p+1)/2 \pmod{p^{3}}$. Of
course this process can be carried out for any integer not just
2. Hopefully one can now see how a chain of elements of $S(p,q)$
can be found given an element (not 1) in $S(p,q)$.

An example can be given for $p=3$. One then proceeds to find the
necessary inverses in $(\mathbb{Z}/ 3^{s}\mathbb{Z})^{\times}$.
Hence \[ 2\rightarrow \{2\}_{3}\rightarrow \{5\}_{9}\rightarrow
\{11,14\}_{27} \rightarrow \{29,41,59,65\}_{81}\rightarrow\{
83,86,122,146,173,176,191,221\}_{243}\rightarrow\ldots.\] One
immediately sees that the primes $q=5,7,11,13,17,29,41,43,59, 61
,73,173,183,191$ already occur in the inverses as a divisor thus
implying that $2\notin S(3,q)$ for those $q$ above.

The above lemma \ref{lem} can be extended to other residue
classes.

\begin{corollary}  \label{lema}
Let $l|(p^{r}+1)$. Also let $s>r$ and $n\in\mathbb{N}$ be
maximal in satisfying $nr<s$ (that is $(n+1)r\geq s$). Then the
inverse, $x$, of $(p^{r}+1)/l$ in $(\mathbb{Z}\setminus
p^{s}\mathbb{Z})^{\times}$ is given by \[
x=l((-1)^{n}p^{nr}+(-1)^{n-1}p^{(n-1)r}+\ldots
-p^{r}+1)+(1-(-1)^{n})(p^{s}/2),\] (if the leading term is
negative ($n$ is odd) then  $p^{s}$ is added to put the inverse
in the range $0<x<p^{s}$). \end{corollary}

\noindent\textbf{Proof:} Let $x=n(x_{1}+1)$ then the condition
for finding the inverse becomes \[ px_{1}+x_{1}+p=kp^{s},\] which
is the same as equation \ref{zero}.

\begin{flushright} \textbf{$\Box$} \end{flushright}

The most general form of the lemma is given by the following.

\begin{corollary}\label{lemb}
Let $m,r,t,v,n,s\in\mathbb{N}$ with $v|(mp^{r}+t)$. Let $s>r$ and
$n$ be maximal in satisfying $nr<s$. Then the inverse, $x$, of
$(mp^{r}+t)/v$ in $(\mathbb{Z}\setminus p^{s}\mathbb{Z})^{\times}$ is
given by \[ x= v( t_{s}-t_{s}^{2}mp^{r}+\ldots
(-1)^{l}t_{S}^{l+1}(mp^{r})^{l}+\ldots
+(-1)^{n}t_{s}^{n+1}(mp^{r})^{n})+(1-(-1)^{n})(p^{s}/2),\]
where $t_{s}\in \{1,2,\ldots p^{s}-1\}$ satisfies $tt_{s}\equiv 1
\pmod{p^{s}}$.
\end{corollary}

The proof is identical in nature to the above. The only problem
with this corollary is the almost circular argument used for
finding the inverse of t in $(\mathbb{Z}\setminus
p^{s}\mathbb{Z})^{\times}$. In the lemma and first corollary this
inverse is simple to find and so the result is useful. The
advantage is that for fixed $s$ and $t$ the inverse $t_{s}$ can
be found using Euclid's algorithm and the second corollary
provides a useful result for finding other inverses having only
used Euclid's algorithm once.

One can even investigate ways of
finding $t_{s}$ simply or more generally the inverse of a $j\in
(\mathbb{Z}/ n\mathbb{Z})^{\times}$. The inverse of $j$ is going
to be of the form $(yn+1)/j$. For example if $n\equiv -1
\pmod{j}$ then $y=1$ or if $n\equiv 1 \pmod{j}$ then $y=j-1$. The
value of $y$ depends on $n\pmod{j}$. Indeed one needs $yn\equiv
-1 \pmod{j}$.

For the more general integer $j\neq 1,2$ it is more difficult
because in general an inverse residue to $j$ cannot just be
written down like it was done for the case $j=2$ above,
different cases modulo $j$ have to be examined.

\section{A Universal $p$-adic Function}
The Riemann zeta function provided the first example of double
interpolation and in fact of interpolation with respect to a
finite set of primes. This next example shows that interpolation
can take place with respect to all primes. For the Riemann zeta
function this was not possible because removing all the Euler
factors would result in something trivial.

\subsection{Translated Ideals}
The open sets of $\mathbb{Z}_{p}$ are an essential part of
$p$-adic analysis (for example in some forms of $p$-adic
measures and integration). These open sets are translated ideals
of $\mathbb{Z}_{p}$; $p^{n}\mathbb{Z}_{p}$ for
$n\in\mathbb{N}$. Let
\[\mathcal{S}_{p}=\{u+p^{v}\mathbb{Z}_{p}:v\in\mathbb{N}\cup\{\infty\},0\leq
u<p^{v}\} .\] (The point at infinity accounts for the trivial
ideal.) To find a nonstandard version of these ideals in
$\mathbb{^{\ast}Z}$ one studies its translated ideals. The
ideals of $\mathbb{^{\ast}Z}$ are of the form
$m\mathbb{^{\ast}Z}$ for $m\in\mathbb{^{\ast}N}$. In analogy
with $\mathcal{S}_{p}$ define
\[
\mathcal{S}=\{a+m\mathbb{^{\ast}Z}:m\in\mathbb{^{\ast}N},0\leq a
<m\}.\] There is a natural connection between these two sets and
that is via the $p$-adic shadow map restricted to
$\mathbb{^{\ast}Z}$ (in this case it is the same as the map
described by \cite{Gv}, chapter 16). It is defined by
\begin{align}
\operatorname{sh}_{p}: \mathbb{^{\ast}Z}&\rightarrow
\mathbb{Z}_{p}\notag , \\
x&\mapsto \langle x \operatorname{mod p}, x \operatorname{mod
p^{2}}, \ldots \rangle.\notag
\end{align}
It is a homomorphism and surjective, along with other properties.

\begin{lemma}\label{lem1}
Let $m\mathbb{^{\ast}Z}$ be an ideal of $\mathbb{^{\ast}Z}$ then
$\operatorname{sh}_{p}(m\mathbb{^{\ast}Z})$ is an ideal of
$\mathbb{Z}_{p}$.
\end{lemma}
\noindent \textbf{Proof:} Using properties of the $p$-adic
shadow map $\mathbb{^{\ast}Z}^{\lim_{p}}=\mathbb{^{\ast}Z}$
which means $\operatorname{sh}_{p}$ is defined for all
$\mathbb{^{\ast}Z}$. Therefore it is defined on all sets of
$\mathbb{^{\ast}Z}$. Now the basic definition of an ideal is
used. Let
$I_{m}=\operatorname{sh}_{p}(m\mathbb{^{\ast}Z})\subset
\mathbb{Z}_{p}$.

Firstly, $\forall r,s\in m\mathbb{^{\ast}Z}$, $r+s\in
m\mathbb{^{\ast}Z}$. Let $a=\operatorname{sh}_{p}(r)\in I_{m}$
and $b=\operatorname{sh}_{p}(s)\in I_{m}$. Also let $t=r+s$ and
$c=\operatorname{sh}_{p}(t)\in I_{m}$. Then using the
homomorphism property of the $p$-adic shadow map,
\[a+b=\operatorname{sh}_{p}(r)+\operatorname{sh}_{p}(s)
=\operatorname{sh}_{p}(r+s)=c\in I_{m}\].

Secondly, $\forall r\in m\mathbb{^{\ast}Z}$, $\forall w\in
m\mathbb{^{\ast}Z}$, $wr\in m\mathbb{^{\ast}Z}$. Let $a=\operatorname{sh}_{p}(r)\in I_{m}$
and $l=\operatorname{sh}_{p}(w)\in \mathbb{Z}_{p}$. Also let
$t=wr$ and $h=\operatorname{sh}_{p}(t)\in I_{m}$. Then using
the homomorphism property of the $p$-adic shadow map,
\[la=\operatorname{sh}_{p}(l)\operatorname{sh}_{p}(r)
=\operatorname{sh}_{p}(lr)=h\in I_{m}\].

Therefore $I_{m}$ is an ideal of $\mathbb{Z}_{p}$ for all
$m\in\mathbb{^{\ast}N}$. So for some $n(m)\in\mathbb{N}$,
\[\operatorname{sh}_{p}(m\mathbb{^{\ast}Z})=p^{n(m)}\mathbb{Z}_{p}.\]

\begin{flushright} \textbf{$\Box$} \end{flushright}

\begin{corollary} \label{cora}
Let $a+m\mathbb{^{\ast}Z}\in\mathcal{S}$ then
$\operatorname{sh}_{p}(a+m\mathbb{^{\ast}Z})\in\mathcal{S}_{p}.$
\end{corollary}
\noindent \textbf{Proof:} As $a\in\mathbb{^{\ast}Z}$,
$\operatorname{sh}_{p}(a)\in\mathbb{Z}_{p}$. So one can write \[
\operatorname{sh}_{p}(a)=\sum_{k\in\mathbb{N}}a_{k}p^{k}.\]
Then
\begin{align} \operatorname{sh}_{p}(a+m\mathbb{^{\ast}Z}) &=
\operatorname{sh}_{p}(a)
+\operatorname{sh}_{p}(m\mathbb{^{\ast}Z}),\notag\\
&=\operatorname{sh}_{p}(a) + p^{n(m)}\mathbb{Z}_{p},\notag\\
&=(a_{0}+a_{1}p+\ldots
+a_{n(m)-1}p^{n(m)-1})+p^{n(m)}(a_{n(m)}+\ldots)
+p^{n(m)}\mathbb{Z}_{p},\notag\\
&=(a_{0}+a_{1}p+\ldots
+a_{n(m)-1}p^{n(m)-1})+p^{n(m)}\mathbb{Z}_{p},\notag\\
&=b+p^{n(m)}\mathbb{Z}_{p}\in\mathcal{S}_{p},\notag
\end{align}
where $b=a_{0}+a_{1}p+\ldots
+a_{n(m)-1}p^{n(m)-1}\in\mathbb{N}$, in particular $0\leq
b<p^{n(m)}.$

\begin{flushright} \textbf{$\Box$} \end{flushright}

\begin{lemma} \label{lem2}Let $m\mathbb{^{\ast}Z}$ be an ideal of
$\mathbb{^{\ast}Z}$. Let $m=p^{r}e$ with $r\in\mathbb{^{\ast}N}$,
$e\in\mathbb{^{\ast}N}$ and $(e,p)=1$. Then,
\[ \operatorname{sh}_{p}(m\mathbb{^{\ast}Z})=\left\{
\begin{array}{ll}
p^{r}\mathbb{Z}_{p} & \mbox{$m\in\mathbb{^{\ast}Z}$,  $r\in\mathbb{N}$,} \\
p^{\infty}\mathbb{Z}_{p}=0 & \mbox{$m\in\mathbb{^{\ast}N}$,
$r\in\mathbb{^{\ast}N}\setminus \mathbb{N}$.}
\end{array} \right. \]
\end{lemma}
\noindent \textbf{Proof:} Let
$\operatorname{sh}_{p}(m\mathbb{^{\ast}Z}) =
p^{n(m)}\mathbb{Z}_{p}.$ Firstly suppose
$r\in\mathbb{^{\ast}N}\setminus\mathbb{N}$. Then $\forall
r\in\mathbb{N}$, $p^{r}|m$. Then for $x\in m\mathbb{^{\ast}Z}$,
\[\operatorname{sh}_{p}:x\mapsto \langle 0 \operatorname{mod p},
0 \operatorname{mod p^{2}},\ldots \rangle = 0.\] Thus
$n(m)=\infty$.

In the second case suppose $r\in\mathbb{N}$. This is then split
into two further cases. Suppose $r=0$ then $(m,p)=1$ and there
exists solutions to equations of the form $km+p^{t}s=1$ with
$t\in\mathbb{N}\setminus{0}$ and $s,m\in\mathbb{Z}$. This implies
there exist solutions to $m\equiv u \operatorname{mod p^{t}}$
with $(1\leq u<p^{t})$. So for $x\in m\mathbb{^{\ast}Z}$ it can
be written as $x=my$ with $y\in\mathbb{^{\ast}Z}$. Under the
shadow map $y$ maps to an element in $\mathbb{Z}_{p}$. However
as $m\equiv v \operatorname{mod p^{t}}$ $(v\neq 0)$ it follows
that $\operatorname{sh}_{p}(m)\in\mathbb{Z}_{p}^{\times}$. Hence
$\operatorname{sh}_{p}(m\mathbb{^{\ast}Z})\in\mathbb{Z}_{p}$.

In the second subcase suppose ($r>0$). Then $x\equiv 0
\operatorname{mod p^{s}}$ for $0<s\leq r$. So for
$x\in m\mathbb{^{\ast}Z}$,
\[\operatorname{sh}_{p}(x)=\langle 0\operatorname{mod p},\ldots,
0\operatorname{mod p^{r}}, x \operatorname{mod
p^{r+1}},\ldots\rangle.\]
Therefore $\operatorname{sh}_{p}(x)\in p^{r}\mathbb{Z}_{p}.$

\begin{flushright} \textbf{$\Box$} \end{flushright}

\begin{corollary}\label{corb}
Consider $a+m\mathbb{^{\ast}Z}\in\mathcal{S}$ then
\[\operatorname{sh}_{p}(a+m\mathbb{^{\ast}Z})=\left\{
\begin{array}{ll}
b+p^{\operatorname{ord}_{p}(m)}\mathbb{Z}_{p} & \mbox{
$\operatorname{ord}_{p}(m)\in\mathbb{N}$,} \\
\operatorname{sh}_{p}(a) &
\mbox{$\operatorname{ord}_{p}(m)\in\mathbb{^{\ast}N}\setminus
\mathbb{N}$.} \end{array} \right. \] Here let
$\operatorname{sh}_{p}(a)=\sum_{k\in\mathbb{N}}a_{k}p^{k}$ then
$b=\sum_{k=0}^{\operatorname{ord}_{p}(m)-1} a_{k}p^{k}$.
\end{corollary}

\subsection{Characteristic Functions}
Choose a $N\in\mathbb{^{\ast}N}\setminus\mathbb{N}$ and set $\mathcal{P}=
\prod_{p\leq N,\text{ $p$ prime }}p.$ Then this number in some
ways acts as a generic finite rational prime. For example it has
already been seen above that $\operatorname{sh}_{p}
(\mathcal{P}^{n}\mathbb{^{\ast}Z})=p^{n}\mathbb{Z}_{p}$ for all
finite primes $p$ and $n\in\mathbb{N}$. By looking at translated
ideals $a +\mathcal{P}^{n}\mathbb{^{\ast}Z}$ under the shadow
maps translated ideals in $\mathbb{Z}_{p}$ are obtained. These
ideas can be extended to look at characteristic functions. For
example set $\phi_{\mathbb{^{\ast}Z}}:\mathbb{^{\ast}Q}
\rightarrow \{ 0,1\}$ to be the characteristic function of
$\mathbb{^{\ast}Z}$. Taking the $p$-adic shadow map of this function leads
to $\phi_{\mathbb{Z}_{p}}$. Similarly defining a characteristic
function on the ideal $\mathcal{P}^{n}\mathbb{^{\ast}Z}$ leads to
$\phi_{p^{n}\mathbb{Z}_{p}}$ and in the same manner for
translated ideals.

\subsection{The Universal Function}
Classically the $p$-adic interpolation of the function $n^{s}$ is
performed to give a $p$-adic continuous function $n\mapsto n^{s}$
with $n\in 1+p\mathbb{Z}_{p}$ and $s\in\mathbb{Z}$. This is
proved (in \cite{Gv}, 127--133) using the binomial theorem.
This function can be looked at in a nonstandard
setting to obtain a nonstandard function for a fixed
$n\in 1+p\mathbb{^{\ast}Z}$
\begin{align}
f:\mathbb{^{\ast}Z}&\rightarrow
\mathbb{^{\ast}Q}^{\lim_{p}},\notag\\ s&\mapsto n^{s}.\notag
\end{align}

This function is given explicitly by the $p$-adic convergent sum
\[n^{s}=(1+(n-1))^{s}=\sum_{k\in\mathbb{^{\ast}N}}\binom{s}{k}(n-1)^{k}.\]
By taking the $p$-adic shadow map one obtains the classical
$p$-adic function.

At this point one wonders if double interpolation can take place.
That is interpolation with respect to two finite distinct primes
($p$ and $q$). In the previous work this was done. An $n\equiv 1
\operatorname{mod pq}$ was fixed and the same binomial expansion
gave a function which was $p$-adically and $q$-adically
convergent. Moreover taking the respective shadow maps lead to
the classical functions.

The next step is to consider interpolation with respect to all
finite primes. This has to be carried out in a nonstandard space
so one has the relevant congruence. For interpolation of
$n^{s}$ with respect to each prime ($p$) it is required that
\[n\equiv 1 \operatorname{mod p}.\]
This means one needs to consider
$\mathcal{P}\in\mathbb{^{\ast}N}\setminus\mathbb{N}$ with
$\mathcal{P}=\prod_{q \text{ prime},q\leq M}q$,
($M\in\mathbb{^{\ast}N}\setminus\mathbb{N}$). (It will become
clear that one could consider any nonstandard $M$ as all that
matters is that $\mathcal{P}$ is divisible by all standard
primes.)

The 'local' information can
then be gathered into a 'global' equivalence \[ n\equiv 1
\operatorname{mod \mathcal{P}}.\] Equivalently \[n\in
1+\mathcal{P}\mathbb{^{\ast}Z}.\] For such an $n$ define a
function on $\mathbb{^{\ast}Z}$ with values in
$\mathbb{^{\ast}Q}$ by
\[g_{n}(s)=n^{s}=(1+(n-1))^{s}=\sum_{k\in\mathbb{^{\ast}N}}\binom{s}{k}(n-1)^{k}.\]

\begin{lemma}\label{lem3}
The function $g_{n}$ is $p$-adically uniformly continuous and
convergent with respect to every finite prime $p$ in
$\mathbb{Q}$. \end{lemma}

\noindent \textbf{Proof:}
\textbf{Continuity:} It is the case that $n\equiv 1
\operatorname{mod \mathcal{P}}$ or equivalently there exists
$t\in\mathbb{^{\ast}Z}$ such that $n=1+\mathcal{P}t$. Using the
binomial theorem on this last expression gives
$n^{\mathcal{P}^{m}}\equiv 1 \operatorname{mod
\mathcal{P}^{m+1}}.$ Therefore $n^{k+p^{m}}\equiv
n^{k}\operatorname{mod \mathcal{P}^{m+1}}$. Continuity is then
established since the last expression implies $n^{k+p^{m}}\equiv
n^{k}\operatorname{mod p^{m+1}}$.

\textbf{Convergence:} It has to be shown that the series
$b_{k}=\binom{s}{k}(n-1)^{k}$ is a null sequence in each $p$-adic
norm. Since $|\binom{s}{k}|_{p}\leq 1$ for all
$s,k\in\mathbb{^{\ast}Z}$ and $p$ prime one needs only consider
the coefficients $a_{k}=(n-1)^{k}$. For a prime $p$,
\begin{align}
|a_{k}|_{p}&=|n-1|_{p}^{k},\notag\\
&=|\mathcal{P}t|_{p}^{k} \mbox{  (for some
$t\in\mathbb{^{\ast}Z}$),}\notag\\
&=|\mathcal{P}/p|_{p}^{k}|pt|_{p}^{k},\notag\\
&\leq p^{-k}.\notag
\end{align}
Thus $a_{k}$ is a null sequence with respect to each prime.

\begin{flushright} \textbf{$\Box$} \end{flushright}

\begin{lemma}\label{lem4}
For all $s\in\mathbb{^{\ast}Z}$,
$g_{n}(s)\in\mathbb{^{\ast}Q}^{\lim_{p}}$.
\end{lemma}
\noindent \textbf{Proof:} As $n\in
1+\mathcal{P}\mathbb{^{\ast}Z}$, $|n|_{p}\leq 1$ for all $p$.
Raising to the power $s\in\mathbb{^{\ast}Z}$, $|n|_{p}^{s}\leq
1^{s}=1$. So
$n^{s}\in\mathbb{^{\ast}Q}^{\lim_{p}}$ for all $p$.
\begin{flushright} \textbf{$\Box$} \end{flushright}

This last lemma enables the $p$-adic shadow map to be taken for
each prime. This is done by my work on nonstandard interpolation.

\begin{align}
\operatorname{sh}_{p}(g_{n}(s))&=\operatorname{sh}_{p}(n)^{\operatorname{sh}_{p}(s)},\notag\\
&=\sum_{k\in\mathbb{N}}\binom{\operatorname{sh}_{p}(s)}{k}
(\operatorname{sh}_{p}(n)-1)^{k}.\notag
\end{align}
This is the classical $p$-adic function because the function is
defined on all of $\operatorname{sh}_{p}(s)\in\mathbb{Z}_{p}$.
Moreover using the results of section 3, the resulting function
is defined on all of $1+p\mathbb{Z}_{p}$ because
$\operatorname{sh}_{p}(1+\mathcal{P}\mathbb{^{\ast}Z})=1+p\mathbb{Z}_{p}$.
At first glance this is quite surprising because in the
nonstandard set $(1+\mathcal{P}\mathbb{^{\ast}Z})$ the only
standard integer, in fact only standard number, is 1. Yet the
respective shadow maps produce the required integers
$\operatorname{mod} p$ and $p$-adic integers to produce the
ideals $1+p\mathbb{Z}_{p}$.

Returning to the choice of $\mathcal{P}$. In the above proof the
only point which matters is that $\mathcal{P}=e \times\prod_{q
\text{ prime},q\leq M }q$ with $(e,p)=1$ for all finite
primes. The fact that $e$ is $p$-integral for all p enables the
shadow map to map $\mathcal{P}\mathbb{^{\ast}Z}$ to
$p\mathbb{Z}_{p}$ rather than $p^{r}\mathbb{Z}_{p}$ if some
power of $p$ divided $e$. As was seen in section 3 all the
shadow map 'cares' about is the $p$-part and with $e$
$p$-integral there is no effect on the final translated ideal.
Similarly with the choice of $M$ it does not affect the shadow
map.

In conclusion the nonstandard function $g_{n}$ acts as a
source of the $p$-adic function $n^{s}$ for all finite prime
$p$.

\section{Double Hurwitz Zeta Function}
This problem is more
difficult for the Hurwitz zeta function because interpolation is based on
a twisted Hurwitz zeta function. The actual numbers interpolated are
explicitly dependent on the prime $p$ in the $p$-adic
interpolation. So to find the correct numbers for double
interpolation they have to depend explicitly on $p$ and $q$
which can be achieved by two potential methods.
\begin{itemize}
\item Define an analogue of the Teichm\"{u}ller character which
explicitly depends on $p$ and $q$.
\item View the Teichm\"{u}ller character as a Dirichlet
character.
\end{itemize}

The second approach potentially leads onto the ideas of a double $L$-function
but are not considered here.

\subsection{Nonstandard Teichm\"{u}ller Character}
Before trying to find a double version it would be prudent to
find a nonstandard version. Let
$^{\ast}\omega_{p}:(\mathbb{^{\ast}Z}/p\mathbb{^{\ast}Z})^{\times}
\rightarrow\mathbb{^{\ast}N}$ be defined for standard prime $p$
by
\begin{itemize} \item $^{\ast}\omega_{p}(n)\equiv n\pmod{p}$,
\item $^{\ast}\omega_{p}(n)^{p-1}\simeq_{p} 1.$ \end{itemize}

The first condition implies that there exists
$t_{n}\in\mathbb{^{\ast}N}$ such that $^{\ast}\omega_{p}(n) = n
+t_{n}p$. Now the second condition can be used to determine the
value(s) of $t_{n}$ or it will show that such a function does not
exist. Using the binomial theorem,

\begin{align}
^{\ast}\omega_{p}(n)^{p-1} &= (n+t_{n}p)^{p-1}\notag,\\
&=\sum_{r=0}^{p-1}\binom{p-1}{r}(t_{n}p)^{r}n^{p-1-r},\notag\\
&\simeq_{p} 1.\notag
\end{align}

A version of Hensel's lemma in $\mathbb{^{\ast}Z}$ cannot be used
to obtain a value of $t_{n}$ because of lack of completeness.
Instead one looks at properties of the surjective $p$-adic shadow
map ($\operatorname{sh}_{p}:\mathbb{^{\ast}N}\rightarrow
\mathbb{Z}_{p}$). Then for each $x\in\mathbb{Z}_{p}$ there exists
a $n_{x}\in\mathbb{^{\ast}N}$ such that $\operatorname{sh}_{p}(
n_{x})=\omega_{p}(x)$. Moreover by the properties of the shadow
map $n_{x}\equiv x \pmod{p}$ and $n_{x}^{p-1}\simeq_{p} 1$.
Therefore it satisfies the defining properties of the
Teichm\"{u}ller character. Using these observations one can
define a (non-unique) nonstandard Teichm\"{u}ller character
$^{\ast}\omega_{p}: \mathbb{^{\ast}N}\rightarrow \mathbb{^{\ast}
N}$ such that the shadow map is the standard Teichm\"{u}ller
character. In the case of $p|n$ define $^{\ast}\omega_{p}(n)=0$.
There is no unique nonstandard character because for any two such
nonstandard characters the value taken by both for a given
$n\in\mathbb{^{\ast}N}$ lie in the same monad.

\subsection{Double Teichm\"{u}ller Character}
In searching for an initial interpretation of this function the
following conditions should be included in the definition. Let
the double Teichm\"{u}ller character initially be a function
defined \[ ^{\ast}\omega_{p,q}:(\mathbb{^{\ast}Z}/
pq\mathbb{^{\ast}Z})^{\times}\rightarrow \mathbb{^{\ast}N}.\] The
obvious extension to $\mathbb{^{\ast}N}$ is given by
$^{\ast}\omega_{p,q}(n)=0$ for $p|n$ or $q|n$ and
$^{\ast}\omega_{p,q}(n) =\text{}^{\ast}\omega_{p,q}(n \pmod{pq})$
otherwise. From the definition of the Teichm\"{u}ller character
in its natural $p$-adic setting $^{\ast}\omega_{p,q}$ should
satisfy for all $n\in (\mathbb{^{\ast}Z}/
pq\mathbb{^{\ast}Z})^{\times}$
\begin{itemize}
\item $^{\ast}\omega_{p,q}(n)\equiv n \pmod{p}$,
\item $^{\ast}\omega_{p,q}(n)\equiv n \pmod{q}$,
\item $^{\ast}\omega_{p,q}(n)^{p-1} \simeq_{p} 1$,
\item $^{\ast}\omega_{p,q}(n)^{q-1} \simeq_{q} 1$.
\end{itemize}
Suppose such a function exists then the shadow map with respect
to one of the primes will lead to a function which agrees
with the respective Teichm\"{u}ller character but is defined on a
slightly smaller set due to the $q$ part in the original
definition.

A function satisfying the first two conditions is trivial to
find, just set $\omega_{p,q}(n)=n+kpq$ for some hyper natural
number $k$. The last two conditions are equivalent to
 $\omega_{p,q}(n)^{p-1}=1 +u_{p,n}p^{N_{p,n}}$ and
$\omega_{p,q}(n)^{q-1}=1 +u_{q,n}q^{N_{q,n}}$ where
$N_{p,n},N_{q,n} \in \mathbb{^{\ast}N}\setminus\mathbb{N}$ and
$u_{p,n},u_{q,n}\in\mathbb{^{\ast}N}$. In order to apply the
methods of finding a nonstandard character as in the previous
section one needs to determine whether or not for each
$n\in\mathbb{^{\ast}N}$ the following subset of
$\mathbb{^{\ast}N}$ is empty
\[ \mu_{p}(\omega_{p}(n))\cap \mu_{q}(\omega_{q}(n)).\] In the
case that it is non-empty then a (non-)unique double nonstandard
character could be defined by setting $^{\ast}\omega_{p,q}(n)$ to
be equal to an element in the above subset. All the properties
desired are then satisfied.

In solving this problem one could consider it in a slightly more
general setting. Given any two distinct primes $p$ and $q$ let
$x_{p}\in\mathbb{Q}_{p}$ and let $y_{q}\in\mathbb{Q}_{q}$. The
problem is to determine whether ot not the following set is empty
\[\mathcal{S}(x_{p},y_{q})=\mu_{p}(x_{p})\cap \mu_{q}(y_{q}).\]
The first case to consider is when $x_{p},y_{q}\in\mathbb{N}$
and are distinct. Elements $N\in\mathcal{S}(x_{p},y_{q})$ can be
written in the form $N=x_{p}+k_{p}p^{M_{p}}$ and
$N=y_{q}+k_{q}q^{M_{q}}$ for some
$k_{p},k_{q}\in\mathbb{^{\ast}N}$ and
$M_{p},M_{q}\in\mathbb{^{\ast}N}\setminus\mathbb{N}$. Without loss of
generality suppose $x_{p}>y_{q}$ and fix nonstandard values of
$M_{p}$ and $M_{q}$. Then the problem reduces to solving the
linear equation \[t=x_{p}-y_{q}=k_{p}p^{M_{p}}-k_{q}q^{M_{q}}.\]
This is soluble because $(p,q)=1 | t$.

This process can be applied to the non-natural elements of
$\mathbb{Z}_{p}$ and $\mathbb{Z}_{q}$. Indeed let
$x_{p}\in\mathbb{Z}_{p}\setminus\mathbb{N}$ and
$y_{q}\in\mathbb{Z}_{q}\setminus\mathbb{N}$. Then let $U\in\mu_{p}
(x_{p})$ and $V\in\mu_{q}(y_{q}$, clearly
$U,V\in\mathbb{^{\ast}N}\setminus\mathbb{N}$. If there exists some
$N\in\mathcal{S}(x_{p},y_{q})$ then it can be written in the form
$N=S+k_{p}p^{M_{p}}$ and $N=T+k_{q}q^{M_{q}}$ for
some $k_{p},k_{q}\in\mathbb{^{\ast}N}$ and
$M_{p},M_{q}\in\mathbb{^{\ast}N}\setminus\mathbb{N}$. As above a linear
equation results and as $S-T\in\mathbb{Z}$ one has basically the
same problem as above and is soluble as $(p,q)=1$.

This can then be used to find a double version of the $p$-adic
function defined for $x\in\mathbb{Z}_{p}^{\times}$ $<x>_{p}:=x /
\omega_{p}(x)\in 1+\mathbb{Z}_{p}$. Then for
$n\in\mathbb{^{\ast}}$ with $p,q\nmid n$ \[ <n>_{p,q} = n/
\text{} ^{\ast}\omega_{p,q}(n).\] THen taking the shadow map for
a given $n$ leads to the $p$-adic or $q$-adic $<.>$.

\subsection{Hurwitz Zeta Function}
Pick a nonstandard double Teichm\"{u}ller character $^{\ast}
\omega_{p,q}$ such that $^{\ast} \omega_{p,q}(1)=1$. The problem
is that the nonstandard version cannot be viewed as a hyper
Dirichlet character as it is not a multiplicative homomorphism.
This is due to the condition on it not being an exact $p-1$ root
of unity just infinitely close to being one (similarly in the $q$
case). It is the same for the single nonstandard Teichm\"{u}ller
character. So instead of considering a Dirichlet $L$-series
consider the following $L$-function \[ ^{\ast}L(s,\text{} ^{\ast}
\omega_{p,q})=\sum_{n\in\mathbb{^{\ast}N}}\frac{^{\ast}\omega_{p,q}
(n)}{n^{s}}.\]
Using the comparison test this Q-converges absolutely for at
least $^{\ast}\Re(s)>1$. Using the hyper Hurwitz zeta function the
above can be written as \[ ^{\ast}L(s,\text{} ^{\ast}
\omega_{p,q})=(pq)^{-s}\sum_{r=1}^{pq}\text{}
^{\ast}\omega_{p,q}(r) \zeta_{\mathbb{^{\ast}Q}}(s,r/pq).\]
Q-analytic continuation of this L-function is then established
via the Q-analytic continuation of the hyper Hurwitz zeta
function. Therefore for $n\in\mathbb{^{\ast}N}$ \[
^{\ast}L(1-n,\text{}
^{\ast}\omega_{p,q})=-\frac{(pq)^{n-1}}{n}\sum_
{r=1}^{pq}\text{}^{\ast}\omega_{p,q}(r)\text{} ^{\ast}
B_{n}(r/pq),\] where $^{\ast} B_{n}(x)$ is the n-th hyper
Bernoulli polynomial. Thus \[ ^{\ast}L(1-n,\text{}
^{\ast}\omega_{p,q})=-\frac{(pq)^{n-1}}{n}\sum_
{r=1}^{pq}\text{}^{\ast}\omega_{p,q}(r)\sum_{k\in\mathbb{
^{\ast}N}} \binom{n}{k}(r/pq)^{n-k}\text{}^{\ast}B_{k}.\]

Of particular importance is the Hurwitz zeta function and
finding a double interpolation of it. As a first attempt consider
the following function for $n\in\mathbb{^{\ast}N}\setminus\{0\}$,
$b,F\in\mathbb{N}\setminus\{0\}$, $b<F$ and $p,q|F$ \[
^{\ast}H_{p,q}(1-n,b,F)= -\frac{1}{n}
\frac{1}{F}\text{}^{\ast}<b>_{p,q}^{n}\sum_{k\in\mathbb{ ^{\ast}
N}} \binom{n}{k}(F/b)^{k}\text{}^{\ast}B_{k}.\] Since the sum is
hyperfinite the actual values of this function lie in
$\mathbb{^{\ast}Q}$. Moreover for all $n$ the sum lies in
$\mathbb{^{\ast}Q}^{\lim_{p}}\cap \mathbb{^{\ast}Q}^{\lim_{q}}$
by the von Staudt-Clausen theorem. By the properties developed
above the shadow maps can be taken of $<.>_{p,q}$ for $p,q\nmid
b$. Since $1/F\in\mathbb{Q}$ the shadow map is trivial and so one
is just left to deal with $1/n$. For the shadow map to be defined
on this term it is required that $|n|_{p}$ and $|n|_{q}$ are not
infinitesimal. (Under the shadow maps this last requirement
corresponds to the pole of the $p$-adic or $q$-adic Hurwitz zeta
function.) This suggests defining the $p-q$-adic Hurwitz zeta
function for $n\in -\mathbb{^{\ast}N}$ as \[ ^{\ast}\zeta_{p,q}
(n,b,F) = -\frac{1}{1-n}
\frac{1}{F}\text{}^{\ast}<b>_{p,q}^{1-n}\sum_{k\in\mathbb{
^{\ast} N}} \binom{1-n}{k}(F/b)^{k}\text{}^{\ast}B_{k}.\]

\chapter{The Work of Shai Haran}\label{chshai}

\section{Overview}

The dictionary between arithmetic and geometry is a fascinating
area of mathematics having been developed by many great
mathematicians from Kummer and Kronecker to Artin and Weil. The
analogy begins with $\mathbb{Z}$ in arithmetic and with $k[x]$
(the ring of polynomials in one variable over a field $k$) in
geometry. By choosing a prime $p$ of $\mathbb{Z}$ the $p$-adic
integers $\mathbb{Z}_{p}=\lim_{\leftarrow}\mathbb{Z}/ p^{n}$
and field of fractions $\mathbb{Q}_{p}=\mathbb{Z}_{p}[1/p]$ are
obtained. Similarly for a prime $f$ of $k[x]$ the geometric
analogues are $k_{f}[[f]]=\lim_{\leftarrow}k[x]/f^{n}$ and the
field of Laurent series $k_{f}((f))=k_{f}[[f]][1/f]$. This
dictionary extends a lot further but there are two anomalies,
that is two constructions in the geometric picture which have no
"obvious" analogue in the arithmetic picture.

\begin{enumerate}
\item To produce theorems in geometry the change is made from
affine to projective geometry. To the affine line the point at
infinity is added which corresponds to the ring $k[[1/x]]$ and
its field of fractions $k((1/x))$. The analogue of $\infty$ for
$\mathbb{Q}$ is the real prime (denoted by $\eta$). The
associated field is $\mathbb{Q}_{\eta}=\mathbb{R}$ but there is
no analogue $\mathbb{Z}_{\eta}$ of $k[[1/x]]$. By following the
definition for finite primes $\mathbb{Z}_{p}=\{x\in\mathbb{Q}_{p}
:|x|_{p}\leq 1\}$, $\mathbb{Z}_{\eta}\text{}'=' [-1,1]$ but this
is not closed under addition.

\item The second problem refers to tensor products.
By taking the product of given geometrical objects a new
geometrical object is obtained. For example the affine plane
($\mathbb{A}^{2}$) is the product of the affine line
($\mathbb{A}^{1}$) with itself. This corresponds to the tensor
product of two polynomial rings $k[x_{1}]$ and $k[x_{2}]$, in
the category of $k$-algebras, and is equal to the ring of
polynomials in two variables, $k[x_{1},x_{2}]$. Trying to find
the corresponding arithmetical surface it is found that in the
category of commutative rings $\mathbb{Z}\otimes
\mathbb{Z}=\mathbb{Z}$ (the surface reduces to the diagonal). In
fact in any geometry based on rings the arithmetical surface
reduces to the diagonal since
$\operatorname{Spec}(\mathbb{Z}) \otimes
\operatorname{Spec}(\mathbb{Z}) =
\operatorname{Spec}(\mathbb{Z})$. So does there
exist a category in which absolute Decartes powers
$\operatorname{Spec} \mathbb{Z} \ldots \times
\operatorname{Spec}\mathbb{Z}$ do not reduce to the diagonal?
\end{enumerate}

Much of Shai Haran's work has been trying to resolve these two
issues. The bulk of this chapter will review Haran's attempts to
find an interpretation of the first point. Work in the direction
of resolving the second point is still very much in its infancy
and so only a small amount of space will be dedicated to the
theory of non-additive geometry and a new language between
arithmetic and geometry.

\section{The Real Integers}

Almost all Haran's work in the search of the real integers is
found in his book Mysteries of the Real Prime which in his own
words is "very condense and hard to read". Only very recently
have some lecture notes become available which expand on parts
of his book. So the references for this work are the chapters
1--9 of his book \cite{Ha1} and chapters 0--4 and 7 of these
lecture notes.

\begin{itemize}
\item \textbf{Aims:} For the $p$-adic integers there are three
important inverse limit expressions which rely upon reduction
modulo $p$. \[ \mathbb{P}^{1} = \lim_{\leftarrow_{N}}
\mathbb{P}^{1} (\mathbb{Z}/p^{N}\mathbb{Z}),\qquad
\mathbb{Z}_{p}=\lim_{\leftarrow_{N}}\mathbb{Z}/p^{N} \mathbb{Z},
\qquad
\mathbb{Z}_{p}^{\ast}=\lim_{\leftarrow_{N}}(\mathbb{Z}/p^{N}
\mathbb{Z})^{\ast}.\] The aim is to understand the real analogue
of these not by reduction of points but by "reduction" of certain
complex valued functions.

\item \textbf{Techniques:} Initial studies of Markov chains
enable the above $p$-adic limits to be viewed as boundaries of
certain trees. The $q$-world is used to construct chains which
generalise the $p$-adic chains but in the limit $q\rightarrow 0$
the $p$-adic chains are recovered. Moreover in the limit
$q\rightarrow 1$ a real chain is obtained. By interpreting the
real chain as an analogue of the $p$-adic chain certain
deductions are made regarding properties of the real prime.

\item \textbf{Philosophy:} The historical tendency is to develop
results in the $p$-adic setting from the corresponding ones in
the real setting. This is based on the human perception of
macroscopic 'reality' as many aspects of everyday life are based,
to a good approximation, on the real numbers and measuring
distances using the absolute value. Given early mathematics,
engineering, physics,$\ldots$ were influenced by the physical
world it is not surprising that the real results far outnumber
those of the $p$ adic numbers in the Platonic mathematical world.
Often results in the $p$-adic setting are simpler and more
natural than those in the real setting so by studying the reals
from the $p$-adic perspective it is hoped to find new results.

\end{itemize}

\subsection{Markov Chains}
\begin{definition} A Markov chain consists of a state space, a
set $X$, and transition probabilities, a function $P_{X}:X\times
X\rightarrow [0,1]$ satisfying $P_{X}(x,X)=1$ for all $x\in X$.
\end{definition}

The source of $p$-adic chains for his work is derived from graph
theory and to begin with trees. Indeed let $X$ be a tree with
root $x_{0}\in X$. For $n\in\mathbb{N}$ let $X_{n}=\{ x\in X:
d(x_{0},x)=n\}$ then $X=\sqcup_{n\in\mathbb{N}} X_{n}$, a
disjoint union. The boundary, $\delta X$, is defined as the
inverse limit of the sets $X_{n}$. The important result relating
to this is the following theorem.

\begin{theorem} Let $\mathcal{M}_{1}(\delta X)$ be the set of
probability measures on the boundary. Then there is a one-to-one
correspondence between $\tau\in\mathcal{M}_{1}(\delta X)$ and
Markov chains on $X$. \end{theorem}

In the constructions used in the proof it is shown that there
exists a special probability measure, the harmonic measure. From
this probability measure it can be shown that it leads to a
probability measure on $X_{n}$ and an associated Hilbert space
$H_{n}=l_{2}(X_{n},\tau_{n})$. The key observation is
that there is a unitary embedding $H_{n}\hookrightarrow
H_{n+1}$ and an orthogonal projection from $H_{n+1}$ onto the
subspace $H_{n}$. Similarly on the boundary there is the Hilbert
space $H=l_{2}(\delta X,\tau)$, unitary embedding $H_{n}
\hookrightarrow H$ and orthogonal projection $H$ onto $H_{n}$.
These Hilbert spaces (and in particular the orthogonal bases) are
probably the most important aspect of the work.

In order to deal with the real and $q$-chains the above theory
has to be extended to chains which are not trees. One of the
characterizations of a chain which is a tree is that the boundary
is totally disconnected. Therefore non-tree chains also include
ones which have continuous boundaries. The theory is based around
harmonic functions. Given a chain with state space $X=\sqcup
X_{n},X_{0}=\{x_{0}\}$ and transition probabilities $P_{X}$ then
$P_{X}$ can be regarded as a matrix over $X\times X$ with entry
0 if two points are not connected. Further $P_{X}$ can be
regarded as an operator on $l_{\infty}(X)$ by
$P_{X}f(x)=\sum_{x'\in X} P_{X}(x, x') f(x')$ .
\begin{definition} A function $f:X\rightarrow [0,\infty )$ is
called harmonic if $P_{X}f=f$ and $f(x_{0})=1$.
\end{definition}

The collection of all harmonic
functions is a convex set, denote this by
$\operatorname{Harm}(X)$. This decomposes
in the standard way of a convex set: $\operatorname{Harm}(X)=
\operatorname{Harm}(X)_{\text{non-ext}} \cup
\operatorname{Harm}(X)_{\text{ext}}$ where
$\operatorname{Harm}(X)_{\text{non-ext}} = \{ \lambda_{0}f_{0} +
\lambda_{1}f_{1}: f_{0},f_{1}\in\operatorname{Harm}(X),
\lambda_{0}, \lambda_{1}>0, \lambda_{0}+\lambda_{1}=1\}$ and
$\operatorname{Harm}(X)_{\text{ext}}$ be those functions which
are not non-extreme. This leads to a decomposition of the
boundary, $\delta X= \delta X_{\text{ext}} \cup \delta
X_{\text{non-ext}}$. Indeed the boundary is the compactification
of $X$ with respect to the Martin metric which itself is derived
from the Martin kernel and the Green kernel. Here the Green
kernel is an operator on $X\times X$ given by $G(x,y)=
\sum_{m\in\mathbb{N}} (P_{X})^{m}(x,y)$ and the Martin kernel is
given by $K(x,y)=G(x,y)/G(x_{0},y)$.

\begin{theorem} For the general Markov chain there is a
one-to-one correspondence between the harmonic functions on $X$
and the probability measures on $\delta X_{\text{ext}}$.
\end{theorem}

In particular the constant function is clearly harmonic and the
corresponding unique measure is called the harmonic measure.

\subsubsection{A $p$-adic Beta Chain}

The actual chain of considerable use is the non-symmetric
$p$-adic $\beta$ chain. The next stage is calculations. Given the
chain the harmonic measure, boundary, Hilbert spaces, $\ldots$
can all be calculated. They are lengthy calculations and so do
not appear in an explicit way in his book though really to most
readers they can just be accepted. To others they form an
extensive set of exercises which would probably double the length
of his book.

The symmetric $\beta$ chain on
$\mathbb{P}^{1}(\mathbb{Q}_{p})/ \mathbb{Z}_{p}^{\ast}$ is very
complicated so a chain is considered on the tree $\mathbb{P}^{1}
(\mathbb{Q}_{p})/ \mathbb{Z}_{p}^{\ast}\ltimes \mathbb{Z}_{p}$
where $\ltimes$ is a semidirect product.

In summary this chain has the state space $X_{n}=\{ (i,j)\in
\mathbb{N}\times \mathbb{N} : i+j=n\}$. The state space can be
identified with $\mathbb{N}\times\mathbb{N}$ by the following
pararmetrization \[  X_{n} \ni(i,j)\mapsto
(1:p^{n-j})=(1:p^{i})\in \mathbb{P}^{1}(\mathbb{Z}/p^{n})/
(\mathbb{Z}/p^{n})^{\ast} \ltimes (\mathbb{Z}/p^{n}).\] The
explicit identification with $\mathbb{N}\times\mathbb{N}$ can
be found in chapter 4 of \cite{Ha1}. The transition probabilities
are given for $\alpha,\beta >0$ by \[ P_{X}((i,j),(u,v)) =
\left\{ \begin{array}{ll}
\frac{1-p^{-\beta}}{1-p^{-\beta-\alpha}} & \textrm{ if
$(i,j)=(0,0)$ and $(u,v)=(0,1)$,} \\ \frac{(1-p^{-\alpha})
p^{-\beta}}{1-p^{-\alpha-\beta}} & \textrm{ if $(i,j)=(0,0)$ and
$(u,v)=(1,0)$,} \\ p^{-\beta} & \textrm{ $j=v=0$, $i\geq 1$ and
$u=i+1$,} \\ 1-p^{-\beta} & \textrm{$j=0$, $v=1$, $i\geq 0$ and
$u=i$,} \\ 1 &\textrm{$i\geq 0$, $u=i$, $j\geq 1$ and $v=j+1$,}
\\ 0 & \textrm{otherwise.} \end{array} \right.\]
The boundary is
$\mathbb{P}^{1}(\mathbb{Z}_{p})/\mathbb{Z}_{p}^{\ast} \ltimes
\mathbb{Z}_{p}\cong p^{\mathbb{N}} \cup \{0\}$ with harmonic
measure the gamma measure. This is defined as $\tau^{\beta}_{
\mathbb{Z}_{p}}= \phi_{\mathbb{Z}_{p}}(x)|x|_{p}^{\beta}
d^{\ast}x /\zeta_{p}(\beta)$ where $\zeta_{p}$ is the local
component of the completed Riemann zeta function.

The key calculation is that of the orthogonal bases. On the
boundary the basis for the associated Hilbert space, $H=
\bigoplus_{m\in\mathbb{N}} \mathbb{C}\phi_{p,m} $, is called the
$p$-adic Jacobi basis ($\{\phi_{p,m}\}$). On the finite
dimensional Hilbert spaces the basis of $H_{n}=\bigoplus_{0\leq
m\leq N}\mathbb{C}\phi_{p,N,m}$ is given by the $p$-adic Hahn
basis ($\{\phi_{p,N,m}\}$). The most important observation is
reduction $\bmod{p^{N}}$. This reduction takes place between the
two bases using integration against the Martin kernel since the
embeddings and projections mentioned above can be given in terms
of the Martin kernel. In particular the projection, or what can
be referred to as reduction $\bmod{p^{N}}$ from $H$ to $H_{N}$ is
given by \[ K_{p,N}\phi_{p,m}=\left\{ \begin{array}{ll}
\phi_{p,N,m} & \textrm{ $0\leq m\leq N$,}\\ 0 &\textrm{$N<m$.}
\end{array}\right. \] It is this interpretation which is sought
initially in the $q$ case and then the real case to give reduction
of certain complex polynomials $\bmod{\eta^{N}}$.

The initial problem encountered by Haran is finding the "correct"
real chain as there are many candidates ranging from expanding a
real number with some fixed base to the use of continued
fractions. Haran uses the $q$-chains to justify his choice of
real chain as the correct one.

A $q$-chain is a $q$-interpolation between a $p$-adic chain and
the real chain. The $q$ world is extensive in interpolating
between the $p$-adic and real numbers.  In this world the
continuum $\mathbb{Q^{\ast}}_{\eta}/\mathbb{Z^{\ast}}_{\eta}=
\mathbb{R}^{+}$ is approximated by $q^{\mathbb{Z}}$ which
resembles the $p$-adic $\mathbb{Q^{\ast}}_{p}/\mathbb{Z^{\ast}}_{p}=
p^{\mathbb{Z}}$. Recent work can be found in the papers
\cite{Ki1}-\cite{Ki9}, \cite{KKuW}, \cite{KuW}, \cite{KuWY} and
\cite{KuOW}.  Chapter 6 of \cite{Ha1} provides a concise
introduction to the $q$-world.

\subsubsection{The $q$-Beta Chain}

Let the state space be $X_{(q)}=\mathbb{N}\times \mathbb{N}$ and
define the transition probabilities for $\alpha,\beta >0$ as
\[ P_{X_{q}}((i,j),(u,v)) = \left\{ \begin{array}{ll}
\frac{(1-q^{\beta+j})}{(1-q^{\alpha+\beta+i+j)}} &
\textrm{$u=i$ and $v=j+1$,}\\
\frac{(1-q^{\alpha+i)}q^{\beta+j}}{(1-q^{\alpha+\beta+i+j)}} &
\textrm{$u=i+1$ and $v=j$}
\end{array} \right.\]

Clearly this is not a tree. As it is a $q$-chain (interpolating
between the $p$-adic $\beta$ chain and the real $beta$ chain) it
has the special property that taking the $p$-adic limit
($q=p^{-N}$, $\alpha:=\alpha /N$, $\beta:=\beta/ N$,
$N\rightarrow\infty$) converges to the $p$-adic $\beta$ chain
introduced above. Moreover in the real limit ($q:=q^{2/N}$,
$\alpha:=\alpha/2$, $\beta:=\beta/2$, $N\rightarrow\infty$) a
new chain is obtained and this is the real $\beta$ chain.

\subsubsection{The Real Beta Chain}
The state space remains $\mathbb{N}\times\mathbb{N}$,
$X_{\eta}=\sqcup_{n\in\mathbb{N}}X_{\eta(n)}$ with
$X_{\eta(n)}=\{ (i,j):i+j=n\}$ and the transition probabilities
become \[ P_{X_{\eta}}((i,j),(u,v)) = \left\{ \begin{array}{ll}
\frac{\beta + 2j}{\alpha +\beta +2(i+j)} & \textrm{$u=i$
and $v=j+1$,}\\
\frac{\alpha +2i}{\alpha + \beta +2(i+j)} &
\textrm{$u=i+1$ and $v=j$} \end{array} \right.\]

The relevant calculations can be carried out and the resulting
boundary is $\delta X_{\eta}=\mathbb{P}^{1}(\mathbb{R})/\{ \pm
1\} =[0,\infty]$ and the harmonic measure on the boundary is the
real beta measure with finite approximation
$\tau^{\alpha,\beta}_{\eta{n}}(i,j)=(n!/(i!j!))(\zeta_{\eta}(
\alpha+2i,\beta+2j)/\zeta_{\eta}(\alpha,\beta)$ ($n=i+j$ and
$\zeta_{\eta}( , )$ is the beta function).  The beta measure is
defined as the product of two gamma measures on the real plane
which is then projected onto $\mathbb{P}^{1}(\mathbb{R})$. Using
these calculations the Hilbert spaces can be examined but the
most important aspect is to find the relations (the ladder
structure) between them analogous to the structure discussed
above on trees.

Let the finite Hilbert spaces be denoted by $H_{\eta(n)}^{\alpha,
\beta}$ and on the boundary by $H_{\eta}^{\alpha,\beta}$. Recall
in the $p$-adic case the orthogonal bases were found by using the
relations between the spaces. In this case integration against
the Martin kernel does not lead to embeddings or projections. The
theory relies on difference operators, which in his work are
developed in relation to the $q$-beta chain and follow in the
real case by taking the real limit. The difference operator
$D_{n}:H_{\eta(n)}^{\alpha,\beta}\rightarrow
H_{\eta(n-1)}^{\alpha +2, \beta +2}$ and its adjoint $D^{+}_{n}:
H_{\eta(n-1)}^{\alpha +2, \beta +2} \rightarrow
H_{\eta(n)}^{\alpha,\beta}$ are defined as \begin{align} D_{n}
\phi(i,j) &= \left(\frac{\alpha+\beta}{2}+n\right)
(\phi(i,j+1)-\phi(i+1,j)),\notag\\   D_{n}^{+} \phi(i,j)&=
\left( \frac{\alpha+\beta}{2}+n\right)^{-1} (j(\alpha/2
+i)\phi(i,j-1) - i(\beta/2 +j)\phi(i-1,j)).\notag\end{align}
These operators satisfy the Heisenberg relation \[ D_{n}D^{+}_{n}
- D^{+}_{n-1}D_{n-1} = \left( \frac{\alpha +\beta}{2}\right)
\operatorname(id)_{H_{\eta(n-1)}}^{\alpha+2,\beta+2}.\] 
This shows that a constant multiple of the identity operator is
obtained for the difference in going down and up this "ladder".
Moreover the constant function is characterised by the equation
$D_{n}\mathbf{1}=0$ which means that $D_{n}$ can be regarded as
an annihilation operator, $D^{\ast}$ a creation operation and
the constant function as the vacuum. Using these observations
the orthogonal basis $\varphi_{\eta(n),m}^{\alpha,\beta} =
(-1)^{m}/m!
(D^{+})^{m}\mathbf{1}_{H_{\eta(n-m)}^{\alpha+2m,\beta+2m}}$, the
real Hahn basis. This process is repeated for the boundary with
difference operators, satisfying the Heisenberg relation,
creating the ladder structure. Ultimately they lead to the real
Jacobi basis $\varphi_{\eta,m}^{\alpha,\beta}= (-1)^{m}/m!
(D^{+})^{m}\mathbf{1}_{H_{\eta}^{\alpha+2m,\beta+2m}}$. The link
between the Hilbert space on the boundary and those on the
finite sets is also defined via the Martin kernel as \[
K^{\alpha,\beta}_{\eta(n)}\phi(i,j)=\int_{\delta X}
K((i,j),x)\phi(x)\tau^{\alpha,\beta}_{\eta}(x),\] where $K$ is
the Martin kernel on the boundary. Then \[ K^{\alpha,\beta}_{
\eta(n)} \varphi^{\alpha,\beta}_{\eta, m} = \left\{\begin{array}{ll}
 \varphi_{\eta(n),m}^{\alpha,\beta} & \textrm{if $0\leq m
\leq n$,} \\ 0 & \textrm{if $m>n$.}   \end{array}\right.\] This
relates the real Jacobi and Hahn basis. The interpretation is
"reduction $\bmod \eta^{n}$", analogous to the $p$-adic case. So
this gives an interpretation to the inverse limit for
$\mathbb{P}^{1}(\mathbb{R}$ in terms of complex polynomials.

Hopefully the basic method can be seen in how an interpretation
of complex polynomials can be given. In order to "view" the
inverse limits for the other two problems suitable chains with
boundary $\mathbb{Z}_{p}$ and $\mathbb{Z}_{p}^{\ast}$ are
needed. In fact both chains can be deduced from the beta chain.

\subsubsection{Gamma Chain}
For the beta chain taking the limit $\alpha\rightarrow\infty$
leads to the gamma chain. Unfortunately in the real case this
reduces to the unit shift walking from the origin to $0\in\delta
X_{\eta}$. However taking the appropriate normalized limit
$\alpha\rightarrow \infty$ of the above operators $D$ and $D^{+}$
and of the real Jacobi polynomials leads to the Laguerre
polynomials. This is an orthogonal basis for the boundary space
$H^{\beta}_{\eta}=L^{2}(\mathbb{R}/\{\pm 1\},
\tau^{\beta}_{\mathbb{Z}_{\eta}})$ where
$\tau_{\mathbb{Z}_{\eta}}^{\beta}$ is the real gamma measure
defined by $\tau_{\mathbb{Z}_{\eta}}^{\beta}/(x)=\exp(-\pi x^{2})
|x|^{\beta}_{\eta}\pi^{\beta/2}d^{\ast}x/\Gamma(\beta/2)$. This
can be also be derived from the $q$-beta chain by taking
appropriate limits. Haran does not make clear the interpretation
of this in terms of the real limit of the additive problem
initially stated. A problem is that he does state that the real
Laguerre basis should be considered as the interpretation in the
introduction but this is not developed in the book at any stage.
I believe part of the problem is the lack of a real chain, there
are no non-trivial finite Hilbert spaces and hence no reduction
type operator. Further what is the real chain with boundary
$\mathbb{Z}_{\eta}$? These thoughts were confirmed when reading
the appendix in the lecture notes where as a problem he asks
what is the real gamma chain and what is the finite real
Laguerre basis?

\subsubsection{The Real Units}
In the $p$-adic case a non-trivial chain can be obtained with
boundary $\mathbb{Z}_{p}^{\ast}$ by taking the limit $\alpha,
\beta\rightarrow\infty$ of the $p$-adic beta chain. The harmonic
measure on the boundary is just the multiplicative Haar measure
on $\mathbb{Z}_{p}^{\ast}$ normalized by
$d^{\ast}(\mathbb{Z}_{p}^{\ast})=1$. Since a real beta chain
exists taking limits $\alpha,\beta\rightarrow\infty$ should lead
to a similar chain. Indeed it does with transition probabilities
of $1/2$ on the same state space as for the real beta chain. On
this chain the Hilbert spaces exist with difference operators
which lead to orthogonal bases. In this case the basis on the
boundary is given by the real Hermite polynomials \[
\varphi_{\eta,m}(w)=\exp(\pi w^{2})\frac{(-1)^{m}}{m!} \left(
\frac{\delta}{\delta w}\right)^{m}\exp(-\pi w^{2}).\] On the
finite layers the basis is essentially given by the elementary
symmetric functions \[
\varphi_{\eta(n),m}=\frac{(2\pi)^{m/2}}{(N(N-1)\ldots
(N-m+1))^{1/2}} \sigma_{N,m}(1,\ldots, 1,-1\ldots, -1),\] where
$\sigma_{N,m}$ is the $m$-th elementary symmetric function of $N$
variables. Again it is implied that the interpretation for
$\mathbb{Z}_{\eta}^{\times}$ is the Hermite polynomials while for
the approximations to $(\mathbb{Z}/\eta^{n}\mathbb{Z})^{\ast}$
are given by the symmetric functions above. What is not made
clear is the process of reduction $\bmod{\eta^{n}}$ and I think
this is another problem in need of resolution.

\subsubsection{Higher Dimensional Theory}
The theory developed above has a natural extension to global
chains. By considering restricted direct products of chains one
can obtain chains over the finite adeles. So far there is little
to comment on in this area because all the properties of these
global chains follow from the local chains. Of more interest is
the theory of higher dimensions. In the one dimensional theory
all the $q$ chains of interest have been constructed from the
products and semi-direct products of the basic chain. This has
state space $\mathbb{N}$ and transition probabilities given by \[
P_{\mathbb{N}(q)}^{\beta}(i,j)=\left\{ \begin{array}{ll}
q^{\beta+i} & \textrm{$i=j$,} \\ 1-q^{\beta+i} &
\textrm{$i+1=j$,}\\ 0 &\textrm{otherwise.} \end{array} \right.\]
For example the $q$-gamma chain is constructed as the product of
the basic chain $P^{\beta}_{\mathbb{N}(q)}$ and the unit shift
chain $P^{\infty}_{\mathbb{N}(q)}$. The $q$-beta chains are also
constructed from the non-stop (semi-direct) products of just two
basic chains. The next step is to consider non-stop (semi-direct)
products of $r>2$ basic chains. This construction leads
to higher dimensional beta chains.

\subsection{Remarks}
In many areas of mathematics there is often a deeper theory
underlying a set of results. Following Tate's thesis and the
subsequent Langland's program arithmetic problems are translated
into representation theory of adelic groups. This takes place in
an elegant way in which although the computations are done in
representation theory the flavour is algebraic geometry. At this
stage the usual consequence is to study the local objects,
whatever they are, and glue the results to the adeles (the
method of Langlands). What Haran proposes is the alternative of
whether the local objects can be put on an equal footing.

Using this approach Haran tries to attack the Riemann hypothesis
by using the interpolating $q$-objects. (See the next section for
full details.) The rest of his work searches for the
"correct" local object which can capture representation theories
of all local fields - a non-trivial question. Where does one
begin this search? The fragments of information Haran uses are
special functions which are only glimpses of representation
theory. Haran's process of using chains is one way of building
these interpolating functions. The functions can be studied
without this process and one can check properties and limits to
show that they are matrix elements in representations of
$GL_{n}(F)$ for all local fields. An advantage of Haran's process
is that in defining some of the objects above he is able to give
an interpretation  of $\bmod{\eta^{n}}$. Future work will show
whether or not these interpretations are the correct ones.

\section{Applications of The Quantum World}
As mentioned in the previous section one of the main areas of
Haran's work has been in attacking the Riemann hypothesis using
$q$-interpolating objects.

\begin{itemize}
\item \textbf{Aims:} To prove the Riemann hypothesis.
\item \textbf{Techniques:} The work is based on $q$-interpolating
objects. The particular objects of interest are the explicit sums
of Weil and the Reisz potential.
\end{itemize}

\subsection{The Riemann Zeta Function}
The completed zeta function is defined via the product of the
local factors to give a global zeta function which is related to
certain operators acting on the adeles, full details can
be found in the opening of chapter 13 of \cite{Ha1}.  \[
\zeta_{\mathbb{A}}(s)=\prod_{p\geq \eta} \zeta_{p}(s), \qquad
\textrm{$\Re(s)>1$.}\] Here the product is over all primes
including the real prime. The local factors are defined as
$\zeta_{p}(s)=(1-p^{-s})^{-1}$ for the finite primes and
$\zeta_{\eta}(s)=\pi^{-s/2}\Gamma(s/2)$. Although the functional
equation is well known Haran views it from the perspective of the
Heisenberg group to obtain \[
\zeta_{\mathbb{A}}(1/2+s)=\zeta_{\mathbb{A}}(1/2-s),\] with
simple poles at $\pm 1/2$ and residues $\pm 1$. There are
numerous reformulations of the Riemann hypothesis and the main
one Haran develops are the explicit sums.

For any $f\in C_{c}^{\infty}(\mathbb{R}^{+})$ there exists a
Mellin transform of it
$\hat{f}(s)=\int_{0}^{\infty}f(a)a^{s}d^{\ast}a$. Then using the
residue theorem and properties of the completed zeta function \[
\sum_{\zeta_{\mathbb{A}}(1/2 +s)=0}\hat{f}(s)- \hat{f}(1/2) -
\hat{f}(-1/2)=-\sum_{p\geq \eta}\mathcal{W}_{p}(f),\] where
$\mathcal{W}_{p}(f)$ is the Weil distribution defined by \[
\mathcal{W}_{p}(f)=\frac{1}{2\pi i}\int_{-i\infty}^{i\infty}
\hat{f}(s)d\log \frac{\zeta_{p}(1/2 +s)}{\zeta_{p}(1/2-s)}.\] For
the finite primes there is a closed expression for the Weil
distribution \[ \mathcal{W}_{p}(f)=\log (p) \sum_{n\neq 0}
p^{-|n| /2}f(p^{n}).\] So essentially the sum of $\hat{f}(s)$
over the zeros of $\zeta_{\mathbb{A}}(1/2 +s)$ is "equal" to the
weighted sum of $f$ over the prime powers. Moreover for finite
$p$, $f\geq 0$ implies $\mathcal{W}_{p}(f)\geq 0$. The
contribution of the real prime cannot be ignored. There are
several ways to write the real Weil distribution in a finite
closed form, for example \[ \mathcal{W}_{\eta}(f)= \int_{0}^{
\infty} \frac{f(a)-f(1)}{1-\min(a,a^{-1})^{2}} \min(a,a^{-1})
^{1/2} d^{\ast}a +(\gamma+\pi/2 +\log(8\pi))f(1),\] where
$\gamma$ is the Euler constant. This formula is obtained by Haran
by introducing $q$-Weil distributions and taking the real limit.
The $p$-adic limit leads to the $p$-Weil distribution as
expected. The importance of the Weil distributions is in relation
to the Riemann hypothesis. The Mellin transform takes
multiplicative convolutions of functions to multiplication of
functions. Applying this to the explicit sums  \[
\sum_{\zeta_{\mathbb{A}}(1/2 +s)=0}\hat{f}(s)\overline{ \hat{f}
(-\overline{s})} -2\Re( \hat{f}(1/2)\overline{
\hat{f}(-1/2)}=-\sum_{p\geq \eta}\mathcal{W}_{p}(f \circ
f^{\ast}),\] where $\circ$ is the multiplicative convolution and
$^{\ast}$ is the adjoint function. Without loss of generality it
can be assumed that $\hat{f}(\pm1/2)=0$. Then the Riemann
hypothesis is $\zeta_{\mathbb{A}}(1/2 +s)=0$ implies
$s=-\overline{s}$ which is equivalent to \[
\sum_{\zeta_{\mathbb{A}} (1/2 +s)=0}\hat{f}(s)\overline{ \hat{f}
(-\overline{s})} \geq 0,\qquad \textrm{for all $f$.}\] Therefore
the Riemann hypothesis is equivalent to the positivity of \[
-\sum_{p\geq \eta}\mathcal{W}_{p}(f \circ  f^{\ast})\geq 0.\]

\subsection{Riesz Potentials}
Haran's next step is finding a connection between the Weil
distributions and the Riesz potentials.

Let $\mathcal{F}_{p}$ be the Fourier transform
then the Riesz potential $R^{s}_{p}$ can be formulated for all
primes acting on functions on $\mathbb{Q}_{p}$ by
\[R^{s}_{p}\varphi \mathcal{F}
|x|_{p}^{-s}\mathcal{F}^{-1}\varphi.\]

\begin{theorem}[Local Formula] \begin{align} \mathcal{W}_{p}(f)
&= ( \mathcal{F}_{p} \log (|x|_{p}^{-1}) \mathcal{F}_{p}^{-1}
|x|_{p}^{-1/2}f(|x|_{p}))(1)\notag,\\ &=\frac{\delta}{\delta s}
\bigg| _{s=0}(R_{p}^{s}|x|_{p}^{-1/2}f(|x|_{p}))(1).\notag
\end{align} \end{theorem} The main advantage of this
theorem is that it enables the Weil distribution to be written
in the form of the trace of an operator in at least one way. One
operator to consider is the operator on $L_{2}(\mathbb{Q}_{p}$
given by $R_{p}^{s}\pi(f(|.|_{p}))\phi_{\mathbb{Z}_{p}}$ where
$\pi$ is the unitary action of the multiplicative group
($\pi(a)f(x)=|a|_{p}^{-1/2}f(a^{-1}x)$) and
$\phi_{\mathbb{Z}_{p}}$ is the operator of multiplication by the
characteristic function of $\mathbb{Z}_{p}$. Then \[
\mathcal{W}_{p}(f)=\frac{\delta}{\delta s}
\bigg| _{s=0} \frac{1}{\zeta_{p}(s)}\operatorname{tr}
(R_{p}^{s}\pi(f(|x|_{p}))\phi_{\mathbb{Z}_{p}} | L_{2}
(\mathbb{Q}_{p})).\]  The Taylor expansion is given by \[
\frac{1}{\zeta_{p}(s)}\operatorname{tr}
(R_{p}^{s}\pi(f(|x|_{p}))\phi_{\mathbb{Z}_{p}} | L_{2}
(\mathbb{Q}_{p}))=f(1)+\mathcal{W}_{p}(f)s+O(s^{2}), \qquad
\textrm{as $s\rightarrow\infty$.}\] There are other versions of
writing the Weil distributions locally including work by
Connes using cut off operators \cite{Co1}.

Indeed, let $B_{c}(x)=1$ if $|x|_{p}\leq c$ and $0$ otherwise.
This can be viewed as a projection on $L_{2}(\mathbb{Q}_{p})$.
The dual projection is given by $\hat{B_{c}}=\mathcal{F}_{p}
B_{c} \mathcal{F}_{p}^{-1}$. Further let $\pi^{0}(f(|x|_{p}))
\varphi(x)=\int d^{0}a f(|a|_{p})|a|_{p}^{1/2}\varphi(ax)$ where
the normalized multiplicative measure is given by $d^{0}a=
n_{p}da/|a|_{p}$ with $n_{p}=\log(p)/(1-p^{-1})$ for finite
primes and $1/2$ for the real prime.

\begin{theorem}
For all $N$ and as $c\rightarrow\infty$ \[ \operatorname{tr}(
\hat{B_{c}}B_{c}\pi^{0}(f(|x|_{p})))=(2\log(c))f(1)+\mathcal{W}
_{p}(f)+O(c^{-N}).\]
\end{theorem}

The next step is to provide global formulae. These provide
reformulations of the Riemann hypothesis. Let
$\mathcal{R}_{\mathbb{A}}=\sum_{p\geq \eta} \mathcal{R}_{p}$
where $\mathcal{R}_{p}=\delta/\delta s |_{s=0} R^{s}_{p}$. A
function $f\in C_{c}^{\infty}(\mathbb{R}^{+})$ a function on
$\mathbb{A^{\ast}}$ can be associated to it
$\tilde{f}=\phi_{\hat{\mathbb{Z}}^{\ast}} \otimes f$, further let
$\tilde{f}_{q}(x)=\tilde{f}(qx)$ for $q\in\mathbb{Q}^{\ast}$ and
$x\in\mathbb{A}^{\ast}$.

\begin{theorem}[Global Formula] \[ \sum_{p\geq\eta}
\mathcal{W}_{p}(f) = \hat{f}(1/2) +\hat{f}(-1/2) - \sum_{\zeta_{
\mathbb{A}}(1/2 +s)=0} \hat{f}(s) = \sum_{q\in \mathbb{Q}^{\ast}}
\mathcal{R}_{\mathbb{A}}(x^{-1/2}f)^{~}(q).\]
\end{theorem}

There also exist global trace formulae whose positivity implies
the Riemann hypothesis. In a simplistic one way of stating this
is \[ \textrm{Riemann Hypothesis} \leftrightarrow \sum_{q\in
\mathbb{Q}^{\ast}} \operatorname{CT}_{s=0} \operatorname{tr}
("\otimes_{p\geq \eta}R^{s}_{p} "
\pi(\hat{f})\pi(\hat(f))^{\ast} \pi(q))\geq 0.\]
Here $\operatorname{CT}_{s=0}$ is the constant term in the Taylor
series about $s=0$.

This is very similar to the work of Connes, in fact he has
rewritten the global trace formulas using the tools of
non-commutative geometry.

In Haran's work he claims that his work involving global trace
formula provides a proof of the Riemann hypothesis for the
function field case but this is not the case as he conceded
himself. Whether or not this method can be used to prove the
Riemann hypothesis in either the function field case or the
number field case is yet to be resolved.

\section{Non-Additive Geometry}
This section gives a very narrow introduction of what little
theory exists regarding the second problem of the
arithmetic-geometry dictionary - the arithmetic surface. In many
ways this is closely related to the first problem. Indeed in
trying to consider $\mathbb{Z}_{\eta}$ as the interval $[-1,1]$
the problem was it not being closed under addition. So why not
abandon addition? This process has led to the consideration of
the "field with one element".  A detailed overview can be found
in \cite{Du}

\subsection{The Field Of One Element}
The conception of the "field of one element" (denoted
$\mathbb{F}_{1}$) was in group theory in the work of Jacques
Tits. Without going into detail the "the field of one element"
was thought of/defined  as $G(\mathbb{F}_{1})=W$ where $G$ is a
Chevalley group scheme over $\mathbb{Z}$ and $W$ is its Weyl
group. On the face of it the idea is absurd since fields by
definition must have at least two elements. However there is
plenty of evidence that something like it exists, essentially
using the ideas of Tits. With finite fields, $\mathbb{F}_{q}$
($q$ a power of a prime), there are numerous formula for
counting structures on projective spaces over $\mathbb{F}_{q}$.
In these formula taking $q=1$ gives results for finite sets. As
a very basic example the number of maximal flags in an
$n$-dimensional vector space over $\mathbb{F}_{q}$ is the
$q$-factorial $[n]!=[1][2] \ldots [n]$, where
$[n]=(q^{n}-1)/(q-1)$. Taking $q\rightarrow 1$ gives the value
$n!$, which is the number of ways to order a set with $n$
elements. This analogy extends a lot deeper and is considered so
powerful that the interpretation of finite sets is as projective
spaces over the "field with one element", $\mathbb{F}_{1}$. With
this interpretation and "existence" of a "field with one
element" further interpretations and definitions are made. This
was, and perhaps still is, a set of suggestions and results based
on these in order to find the "correct" approach to dealing with
such an entity.

One must not forget that the main aim from a number theory point
of view is to be able to use geometrical methods for solving
number theoretical problems. Central to this is to view
$\operatorname{Spec}\mathbb{Z}$ as a curve over some field. In
this case some of the well known conjectures of arithmetic
(Riemann hypothesis, ABC,$\ldots$) become easy theorems in the
geometric analogue of a curve $C$ over a finite field since the
surface $C\times C$ can be formed.

This speculation began several decades after Tits' work with
Manin suggesting possible zeta functions over $\mathbb{F}_{1}$
(\cite{Ma1}). The possible geometrical theorems associated to
$\mathbb{F}_{1}$ continued with a category of varieties over
$\mathbb{F}_{1}$ in \cite{So}, derivations over
$\mathbb{F}_{1}$ in \cite{KuOW} and more recently work by Deitmar
on schemes over $\mathbb{F}_{1}$, the related cohomology and zeta
functions.

Part of Soul\"{e}'s work seems to have been inspired by Manin as
he gives a definition of zeta functions over $\mathbb{F}_{1}$.
Indeed a starting point is to view $\mathbb{F}_{1}$ as the base
field of $\operatorname{Spec} \mathbb{Z}$. This means that any
variety over $\mathbb{F}_{1}$ must have a base change over
$\mathbb{Z}$, which is an algebraic variety over $\mathbb{Z}$.
This base change is found in several of the works above.  In
order to define a variety over $\mathbb{F}_{1}$ it has been
suggested, in unpublished work by Kapranov and Smirnov, that
$\mathbb{F}_{1}$ should have an extension $\mathbb{F}_{1^{n}}$
given by adjoining roots of unity suggesting that \[
\mathbb{F}_{1^{n}} \otimes_{\mathbb{F}_{1}}
\mathbb{Z}=\mathbb{Z}[T]/(T^{n}-1),\] which is $R_{n}$, the ring
of functions on the affine group scheme of $n$-th roots of unity.

\begin{definition} A variety over $\mathbb{F}_{1}$ is a covariant
functor $X$ from the category $\mathcal{R}$ (with objects $R_{n}$
and their finite tensor products) to finite sets. For
$R\in\operatorname{Ob}\mathcal{R}$ there is the natural inclusion
$X(R)\subset X_{\mathbb{Z}}(R)$ where $X_{\mathbb{Z}}$ is a
variety over $\mathbb{Z}$. There is also a further condition
(universal property) which relates to functors. \end{definition}
Unfortunately this category of varieties over
$\operatorname{Spec}\mathbb{F}_{1}$ does not contain
$\operatorname{Spec}\mathbb{Z}$.

From Weil's work a local zeta function can be defined for a
scheme of finite type over $\mathbb{Z}$, $X$. \[
Z_{X}(p,T)=\exp\left( \sum_{n=1}^{\infty}\frac{T^{n}}{n}
\#X(\mathbb{F}_{p^{n}})\right),\] where $p$ is a prime number.
Soul\'{e} considered the condition that there exists a
polynomial $N(x)$ with integer coefficients for every prime $p$
and $N\in\mathbb{N}$ with
$\#X(\mathbb{F}_{p^{n}})=N(p^{n})$. He then defines
\[\zeta_{X|\mathbb{F}_{1}}(s) =\lim_{p\rightarrow 1} \frac{
Z_{X}(p, p^{-s})^{-1}}{(p-1)^{N(1)}} =s^{a_{0}}(s-1)^{a_{1}}
\cdots (s-n)^{a_{n}}.\] The second equality follows by setting
$N(x)=a_{0}+\cdots +a_{n}x^{n}$. For example
$\zeta_{\operatorname{Spec}\mathbb{F}_{1}}(s)=s.$

Zeta functions over arbitrary $\mathbb{F}_{1}$-schemes
were first defined in Deitmar's work on schemes. He shows that the condition on requiring
$\#X(\mathbb{F}_{q})=N(q)$ can be made slightly weaker, by the
following theorem, which enables the zeta functions to be
defined.

\begin{theorem} Let $X$ be a $\mathbb{Z}$-scheme defined over
$\mathbb{F}_{1}$. Then there exists $e\in\mathbb{N}$ and
$N(x)\in\mathbb{Z}[x]$ such that for every prime power $q$ \[
(q-1,e)=1 \Rightarrow \#X_{\mathbb{Z}}(\mathbb{F}_{q}) = N(q). \]
The condition determines the zeta polynomial of $X$, $N$,
uniquely. \end{theorem}

The zeta function of an arbitrary $\mathbb{F}_{1}$-scheme, $X$,
can then be defined via the zeta polynomial ($N_{X}(x)= a_{0}
+\cdots +a_{n}x^{n}$) as \[ \zeta_{X|\mathbb{F}_{1}}(s)=
s^{a_{0}} \cdots (s-n)^{a_{n}}.\] Explicit calculations of zeta
functions are found in \cite{K} who also defines the Euler
characteristic as \[ \#X(\mathbb{F}_{1})=\sum_{k=0}^{n} a_{k}.\]

\subsection{Haran}

An important aspect arising from \cite{KuOW} is the building the
theory of $\mathbb{F}_{1}$ on multiplication only - losing
addition. Supposing that $\operatorname{Spec}\mathbb{Z}$ is a
"curve" over $\mathbb{F}_{1}$ then derivations should exist. If
the concept of $\mathbb{Z}$-linearity is removed and one only
considers the Leibniz rule then the result is derivations. This
idea of losing additivity/$\mathbb{Z}$-linearity is applied to
other $\mathbb{F}_{1}$ objects in the form of a (forgetful)
functor from objects over $\mathbb{Z}$. This idea is used by
Deitmar to define schemes and is a starting point in Haran's
work.

Haran's input in the study of the "field with one
element" area is different from the
other work in this area.  He uses his work on the real prime as
guidance in particular the concept that by abandoning addition
the real integers become a real object.

The initial seed for the work comes from the work of \cite{KuOW}
as they interpret vector spaces over $\mathbb{F}_{1}$ as finite
pointed sets. This combined with a hint of the real prime leads
to the definition of a category $\mathbb{F}$ which has objects
finite sets and arrows as partial bijections. He chooses to use
this as a model for the "field with one element". The actual
connection with the previous work in this area is not made clear.
Clearly there is a form of connection because the objects can
also be taken to be vector spaces over $\mathbb{F}_{1}$ with
certain maps related to the real integers. So unlike the previous
work he has taken the two problems, initially stated about the
arithmetic-geometric dictionary, and combined them to search for
some answers/interpretations. Again this is an interpretation as
the theory in his paper can be taken for a category $F$ with two
symmetric monoidal structures $\oplus$, $\otimes$, the unit
element ($[0]$) for $\otimes$ is the initial and final object of
$F$, $\otimes$ is distributive over $\oplus$ and it respects $X
\otimes [0]=[0]$ ($X\in |F|$).

From this definition various structures are imposed on
$\mathbb{F}$. One of the most important is the $\mathbb{F}$-ring.

\begin{definition} A $\mathbb{F}$-ring is a category $A$ with
objects $|\mathbb{F}|$ and arrows $A_{Y,X}=\operatorname{Hom}_{A}
(X,Y)$ containing the arrows of $\mathbb{F}$. This means there is
a faithful functor $\mathbb{F}\rightarrow A$ which is the
identity on objects. \end{definition}

As examples of these  he shows that there is a
functor from commutative rings to $\mathbb{F}$-rings. The
$\mathbb{F}$-ring associated to a commutative ring $R$ has
morphisms which are matrices with values in $R$. The most
important example is that relating to the real integers. Let
$X\in |\mathbb{F}|$ and let $\mathbb{R}\cdot X$ denote the real
vector space with inner product having $X$ as an orthonormal
basis. Then for $a=(a_{x})\in\mathbb{R}\cdot X$ there is a norm
$|a|_{\eta}=(\sum_{x\in X}|a_{x}|^{2})^{1/2}$. There is also the
related operator norm for a real linear map
$f\in\operatorname{Hom}(\mathbb{R}\cdot X, \mathbb{R}\cdot Y)$,
$|f|_{\eta}=\sup_{|a|_{\eta}\leq 1} |f(a)|_{\eta}$. Then the
$\mathbb{F}$-ring, $\mathcal{O}_{\mathbb{R}}$, is defined to have
morphisms \[ (\mathcal{O}_{\mathbb{R}})_{Y,X}=\{ f\in
\operatorname{Hom}_{\mathbb{R}}(\mathbb{R}\cdot X,
\mathbb{R}\cdot Y), |f|_{\eta}\leq 1\}.\] Haran also defines
the residue field of $\mathcal{O}_{\mathbb{R}}$ to be the
$\mathbb{F}$-ring of partial isometries ($\mathbb{F}_{\eta}$ ) \[
(\mathbb{F}_{\eta})_{Y,X}=\{f:V\xrightarrow{\sim} W:V\subset
\mathbb{R}.X, W\subset\mathbb{R}.y, \text{real sub-vector
spaces and $f$ is a real linear isometry}\}    .\] These can be
generalised to any number field, $k$, and $\eta:k\rightarrow
\mathbb{C}$ a real or complex prime. In the real case he simply
mentions that these constructions are analogous to the $p$-adic
integers. In future work I think these $\mathbb{F}$-rings may
play an important role.

The rest of the paper sets about developing the language of
geometry using $\mathbb{F}$-rings. This naturally starts with the
definition of modules, submodules and ideals of an
$\mathbb{F}$-ring. An $A$-module of a $\mathbb{F}$-ring $A$ is
basically a collection of sets with maps which are compatible
with $A$. An $A$-submodule $M'$ of $M$ is a collection of subsets
of $M$ which are closed under the maps which are used to define
$M$. In the special case of $M=A$ then the $A$-submodule is
called an ideal. A clear example following on from the
$\mathbb{F}$-rings are the modules of a commutative ring. Under
the functor $\mathbb{F}$ the modules become modules of an
$\mathbb{F}$-ring. The main emphasis is placed on ideals in order
to develop prime ideals, $\operatorname{Spec}$ and schemes.

The ideals which he uses are called $H$(omogeneous)-ideals.
A normal ideal $\mathbf{a}$ is a collection of subsets $\{
\mathbf{a}_{Y,X} \subset A_{Y,X}\}$ which are closed under the
functors $\otimes$, $\oplus$ and under composition $\circ$ of
maps. The $H$-ideals are a subset of ideals with the property
that an ideal $\bold{a}$ is generated by $\bold{a}_{[1] ,[1]}$
where $[n]=\{0,1,\ldots n\}$. An application of Zorn's lemma
gives
\begin{theorem} Every $\mathbb{F}$-ring contains a maximal
(proper) $H$-ideal. \end{theorem}
Naturally there is a notion of a prime $H$-ideal. A $H$-ideal
$\bold{p}\in A_{[1],[1]}$ is called prime if $A_{[1],[1]}\
\bold{p}$ is multiplicative closed: $f,g\in A_{[1],[1]}\ \bold{p}
\rightarrow f.g\notin \bold{p}$. The set of prime ideals of $A$
is denoted by $\operatorname{Spec}(A)$.

As an example he takes the $\mathbb{F}$-ring of real "integers",
$\mathcal{O}_{\mathbb{R}}$ and states that
$\bold{m}_{\eta}=\{x\in \mathbb{R}:|x|_{\eta}<1\}$ is the unique
maximal$H$-ideal of $\mathcal{O}_{\mathbb{R}}$. With these ideas
in place he can then define a Zariski topology on
$\operatorname{Spec}A$. The closed sets are given for a set
$\mathcal{U}\subset A_{[1],[1]}$, $V_{A}(\mathcal{U})= \{
\bold{p}\in\operatorname{Spec}A:\bold{p}\supseteq \mathcal{U} \}
$ and for a $H$-ideal, $\bold{a}$, generated by $\mathcal{U}$,
$V_{A}(\mathcal{U})=V_{A}(\bold{a})$. The closed sets of the
Zariski topology on $\operatorname{Spec}A$ are $\{ V_{A}(
\bold{a}): \bold{a}\in H-\text{id}(A)\}$.

The theory of localization of a $\mathbb{F}$-ring is almost
identical to that of commutative rings due the multiplicative
theory. This enables $\mathbb{F}$-ringed spaces to be defined
and a category of Zariski-$\mathbb{F}$-schemes. So far this work
is in setting the stage for actual applications and progress
towards finding the arithmetical surface. Based on the other work
I would also expect zeta functions to be developed for the
schemes introduced by Haran.

\section{A Nonstandard Beginning}

The q-world is vast and varied in nature but in this work some
of the aspects which Shai Haran chose to examine are looked at
through a nonstandard view point. 

\subsection{$q$-Integers}

\begin{definition}
Let q be an indeterminate which can be considered in the real
field. For any integer $s$ define the $q$-integer $[
s]_{q}=\frac{1-q^{s}}{1-q}$. This definition in fact holds for
any $s\in\mathbb{C}$, the non-symmetric $q$-numbers.
\end{definition}
From this the nonstandard non-symmetric $q$-numbers can be
defined.

\begin{definition} Let $q\in\mathbb{^{\ast}R}$ and
$s\in\mathbb{^{\ast}C}$ then define $^{\ast}[
s]_{q}=\frac{1-q^{s}}{1-q}$.
\end{definition}
This is a function
$^{\ast}[.]_{(.)}:\mathbb{^{\ast}C}\times\mathbb{^{\ast}R}\rightarrow
\mathbb{^{\ast}C}$. Of particular interest are the values
of the function within $\mu_{\eta}(1)$.

Consider $q\in\mu_{\eta}(1)$ then for some fixed infinitesimal
$\delta$, $q=1+\delta$ and take $s$ to lie in
$\mathbb{^{\ast}C}^{\lim_{\eta}}$ so the $\eta$-shadow map can be
taken. In particular for $q$ as above and
$s\in\mathbb{^{\ast}C}^{\lim_{\eta}}$, $^{\ast}[s]_{1+\delta}\in
\mathbb{^{\ast}C}^{\lim_{\eta}}$.

 \begin{align}
^{\ast}[q]_{1+\delta}&=
\frac{1-(1+\delta)^{s}}{1-(1+\delta)},\notag\\
&=\frac{1-\text{}^{\ast}\exp(s^{\ast}\log(1+\delta))}{-\delta}
,\notag\\
&=\frac{s^{\ast}\log(1+\delta)+(1/2)s^{2}\text{}^{\ast}
\log^{2}(1+\delta)+\ldots }{\delta},\notag\\
&=s+\delta(s^{2}-s/2)+\delta^{2}(s/3-s^{2}+s^{3})+\ldots +
\delta^{n}S_{n}(s)+\ldots.\notag \end{align}

Here for all $n\in\mathbb{^{\ast}N}$, $S_{n}(s)$ is a monic
hyperpolynomial of degree $n+1$. This sum is absolutely
Q-convergent (by the properties of the functions). From the
convergence properties of the standard function it can be
deduced that $|\delta^{n}S_{n}(s)|\simeq_{\eta}0$ for all
$n\in\mathbb{^{\ast}N}\setminus\mathbb{N}$ and
$s\in\mathbb{^{\ast}C}^{\lim_{\eta}}$.

\begin{align}
\operatorname{sh}_{\eta}(\frac{1-(1+\delta)^{s}}{1-(1+\delta)}
)&=
\operatorname{sh}_{\eta}(\sum_{n\in\mathbb{^{\ast}N}}\delta^{n}S_{n}(s))
,\notag\\
&=\sum_{n\in\mathbb{^{\ast}N}}\operatorname{sh}_{\eta}(\delta^{n}S_{n}(s))
,\notag\\ &= s.\notag \end{align} Since the nonstandard terms are
zero by above and the standard terms, apart from the leading one,
are zero since the delta part is infinitesimal while $S_{n}(s)$
is standard leading to an infinitesimal term. Therefore the only
term remaining is the first one and \[
\operatorname{sh}_{\eta}(\frac{1-q^{s}}{1-q})=s,\] for
$s\in\mathbb{^{\ast}C}^{\lim_{\eta}}$ and $q\in\mu_{\eta}(1)$.
This corresponds to the standard case \[ \lim_{q\rightarrow
1}[s]_{q}=s.\]

There also exists other special values. For example let $q$ be an
infinitesimal of the form $r^{-1/\delta}$ where
$r\in\mathbb{^{\ast}R}$ and $\delta\in\mu_{\eta}(0)$
(an infinitesimal) Let $s\in\mathbb{^{\ast}C}^{\lim_{\eta}}$ and
define $t=s\delta\in\mu_{\eta}(0)$ (also an infinitesimal). Then
by an analogous method as above
\[
\operatorname{sh}_{\eta}(^{\ast}[s\delta]_{r^{-1/\delta}})=1-
\operatorname{sh}_{\eta}(r)^{-\operatorname{sh}_{\eta}(s)}.\]
This example shows that the function is not Q-continuous. This is
due to terms of the form $\epsilon^{\delta}$, where both are
infinitesimals. This makes the nonstandard function quite
interesting because given two values of the function given by
infinitesimals $(s_{1},q_{1})$ and $(s_{2},q_{2})$ it does not
follow that $^{\ast}[s_{1}]_{q_{1}}$ and $^{\ast}[s_{2}]_{q_{2}}$
are infinitely close.

Naturally for other basic functions, such as the non-symmetric
$q$-numbers, special values exist and the real shadow map of
these corresponds to the standard limits.

\subsection{$q$-Zeta}
One definition for the $q$ zeta function is given for
$s\in\mathbb{C}$ and $q\in\mathbb{R}^{+}$ by
\[ \zeta_{q}(s)=\prod_{n\in\mathbb{N}}(1-q^{s+n})^{-1}.\]
One important limit of this function is when $q=p^{-N}$ ($p$
prime) and $s:=s/N$ with $N\rightarrow\infty$ then
\[\zeta_{p^{-N}}(s/N)\rightarrow_{N\rightarrow\infty}\zeta_{p}(
s)=(1-p^{-s})^{-1}.\]

To consider this from a nonstandard
perspective one considers a nonstandard function
$^{\ast}\zeta_{q}(s)=\prod_{n\in\mathbb{N}} (1-q^{s+n})^{-1}$
with nonstandard arguments. To find an analogy to the limit
consider an infinitesimal $\delta\in\mu_{\eta}(0)$ and a prime
$p$ setting $q=p^{-1/\delta}$ and for
$t\in\mathbb{^{\ast}C}^{\lim_{\eta}}$ set $s=\delta t$. Then
\[ ^{\ast}\zeta_{p^{-1/\delta}}(s)=\prod_{n\in\mathbb{N}}
(1-p^{-t}p^{-n/\delta})^{-1}.\] Using classical analysis this
product is absolutely Q-convergent because clearly
$\sum_{n\in\mathbb{N}}p^{-n/\delta}$ is absolutely Q-convergent.

\[\operatorname{sh}_{\eta}(^{\ast}\zeta_{p^{-1/\delta}}(s))=
(1-p^{-\operatorname(sh)_{\eta}(s)})^{-1}\operatorname{sh}_{\eta}
(\prod_{n\in\mathbb{N}\setminus\{0\}} (1-p^{-t}p^{-n/\delta})^{-1}.\]
From the basic properties of the real shadow map and for all
$N\in\mathbb{^{\ast}N}$ with
$a_{n}\in\mathbb{^{\ast}C}^{\lim_{\eta}}$ \[
\operatorname{sh}_{\eta}(\prod_{n=1}^{N}
a_{n})=\prod_{n=1}^{N} \operatorname{sh}_{\eta}(a_{n}).\] In
particular choosing $a_{n}=(1-p^{-t}p^{-n/\delta})^{-1}$ for
$n\in\mathbb{N}$ and $a_{n}=1$ for
$n\in\mathbb{^{\ast}N}\setminus\mathbb{N}$ then \[
\operatorname{sh}_{\eta}(^{\ast}\zeta_{p^{-1/\delta}}(s)) = (1-p^
{-\operatorname{sh}_{\eta}(s)})^{-1}.\]

\subsection{Chains}
One of the main approaches in investigating the real prime is via
Markov chains. To interpolate between the real gamma chains and
$p$-adic gamma chains $q$-chains are used.

\subsubsection{$q$-Gamma Chain}
\begin{definition}
Let $(X^{0}_{q},\mathcal{P}^{\beta}_{X^{0}_{q}})$ be the
$q$-gamma chain where the state space is \[
X^{0}_{q}=\coprod_{n\in\mathbb{N}} X^{0}_{q(n)},\qquad
X^{0}_{q(n)}=\{(n,j)|0\leq j\leq n \}, \] and the transition
probabilities are given by
\[\mathbb{P}^{\beta}_{X^{0}_{q}}((i,j),(u,v))= \left\{
\begin{array}{ll} q^{\beta+j} & \textrm{$u=i+1$, $v=j$},\\
1-q^{\beta+j} & \textrm{$u=i+1$, $v=j+1$},\\ 0 &
\textrm{otherwise}. \end{array}\right.\]
\end{definition}
There are two approaches to considering this in the nonstandard
setting. The methods both use shadow maps to obtain the standard
objects.

\subsubsection{Real Shadow Map}
Consider the following extension to the $q$-gamma chain given by
$(^{\ast}X^{0}_{q},\mathbb{^{\ast}P}^{\beta}_{^{\ast}X^{0}_{q}})$
with state space
$^{\ast}X^{0}_{q}=\coprod_{n\in\mathbb{^{\ast}N}}
^{\ast}X^{0}_{q(n)}$ ($^{\ast}X^{0}_{q(n)}=\{(n,j)|0\leq j\leq n
\}$) and transition probabilities as above but defined for all
nonstandard natural numbers and $q\in\text{} ^{\ast}(0,1)$ and
$\beta\in\mathbb{^{\ast}R}^{+}$. Special values of this chain
exist for a prime $p$ and let $q+p^{-1/\delta}$ for some
$\delta\in\mu_{\eta}$ and consider the chain
$(^{\ast}X^{0}_{p^{-1/\delta}},
\mathbb{^{\ast}P}^{\beta\delta}_{^{\ast}X^{0}_{p^{-1/\delta}}})$
with $\beta\in\mathbb{^{\ast}R}^{+}$. The shadow map can be taken
for all $\beta\in (\mathbb{^{\ast}R}^{+})^{\lim_{\eta}}$ leading
to transition probabilities
$\mathbb{^{\ast}P}^{\operatorname{sh}_{p}(\beta)}_{^{\ast}X^{0}_{p}}$.
The state space is $X^{0}_{p}$. This is the $p$-adic gamma chain.

As a note consider the case when $p$ is not necessarily prime but
$q=r^{-1/\delta}$ for some $r\in
(\mathbb{^{\ast}R}^{+})^{\lim_{\eta}}$. The above holds as there
was no use made of r being prime. Taking the shadow map results
in a chain
$\mathbb{^{\ast}P}^{\operatorname{sh}_{r}(\beta)}_{^{\ast}X^{0}_{p}}$.
In particular for $r\in\mathbb{N}$ then the resulting chain under
the shadow map has a harmonic measure consisting of harmonic
measures relating to the primes divisors of $r$. (This follows
from the basic results on $g$-adic numbers ($g\in\mathbb{N}$)
given by Mahler. These $r$-gamma chains are almost "handmade"
rather than resulting from some natural structure as in the
$p$-adic gamma chains (\cite{Ha1}, chapter 4).

If instead one considers $n\in\mathbb{N}$ and considers a tree
\[Y^{0}_{n}=\coprod_{N\in\mathbb{N}} Y^{0}_{n(N)}, \qquad
Y^{0}_{n(N)}=\mathbb{Z}/n^{N}\mathbb{Z},\qquad
Y^{0}_{n(0)}=\{0\}.\]
One can think of this chain as a $n$-adic expansion with an edge
$y\in Y^{0}_{n(N)}$ going to its $n$ pre-images in
$Y^{0}_{n(N+1)}$. The boundary is given by
\[ \partial Y^{0}_{n}=\lim_{\leftarrow_{N}}Y^{0}_{n(N)}=\mathbb{Z}_{n}.\]
Here $\mathbb{Z}_{n}$ are the $n$-adic integers as described in
Mahler. (\cite{M}, chapter 5). These numbers have a
decomposition as a direct sum of the $p$-adic integers with
$p|n$. On this boundary there is a harmonic measure
$\tau^{\beta}_{\mathbb{Z}_{n}}$ which can be decomposed as a
direct sum of $\tau^{\beta}_{\mathbb{Z}_{p}}$ for $p|n$. This
measure induces a unique chain on this tree. Moreover the measure
is invariant under the $\mathbb{Z}_{n}^{\ast}$ action and a
quotient chain is obtained. The structure of the chain is
determined by the number of distinct prime factors of $n$. For
example when $n$ is prime the chain is the $p$-adic gamma chain
as in Haran. Explicit examples can be easily given. I have yet to
find any use for these chains in the same way there is little use
for general $n$-adic numbers compared to $p$-adic numbers. These
chains are different from the ones arising from the shadow map,
hence the term "handmade" for the shadow map chains.

\subsubsection{$p$-adic Shadow Maps}
The second approach to the gamma chains relies on the other
shadow maps of $\mathbb{^{\ast}Q}$. Choose a nonstandard integer
$N\in\mathbb{^{\ast}N}\setminus\mathbb{N}$ and consider the set
$\mathcal{P}=\{ p\in\mathbb{^{\ast}N}: \text{ $p$ prime, } p\leq
N\}$ (a hyperfinite set of primes) and the product of primes \[
P=\prod_{p\in\mathcal{P}} p.\] This is a hyperfinite product.

\begin{lemma}
For all $p\in\mathcal{P}$ there exists a
non-unique $u\in\mathbb{^{\ast}Q}$ such that for some non-unique
$k_{p},r_{p}\in\mathbb{^{\ast}N}$, $u=p+
k_{p}p^{r_{p}}$ which in particular is equivalent for standard
primes to  $\operatorname{sh}_{p}(u)=p$.   \end{lemma}

\noindent\textbf{Proof:}
Fix a set of $\{r_{p}\}$ simply by choosing a
$r_{p}\in\mathbb{^{\ast}N}\setminus\mathbb{N}$ for each
$p\in\mathcal{P}$. Then it is required to solve the system of
equations $\{u\equiv p (\operatorname{mod }p^{r_{p}} )\}$. Since
$(p^{r_{p}},q^{r_{q}})=1$ for all distinct primes the Chinese
remainder theorem in a nonstandard setting can be used to obtain
a solution $u$. Any two solutions of $u$ for a particular set of
$\{r_{p}\}$ differ by a multiple of
$\prod_{p\in\mathcal{P}}p^{r_{p}}$. In particular for all
standard primes $u\simeq_{p}p$ and by taking the $p$-adic shadow
map $\operatorname{sh}_{p}(u)=p$.
\begin{flushright} \textbf{$\Box$} \end{flushright}

The non-uniqueness needs to be considered and perhaps some form
of equivalence relation could be used but this will be dealt with
later. For the moment suppose some set $\{r_{p}\}$ has been fixed
along with a solution $u$. For such a $u$ a hyper chain will be
constructed based on the $p$-gamma chain. Take a state space
\[ X^{0}_{u}=\coprod_{N\in\mathbb{N}}X^{0}_{u(N)},\quad
X^{0}_{u(N)}=\{(N,j):0\leq j\leq N\}.\]
For the transition probabilities let
\[\mathbb{^{\ast}P}^{\beta}_{X^{0}_{u}}((i,j),(c,d))= \left\{
\begin{array}{ll} u^{-\beta} & \textrm{$c=i+1$, $d=j=0$},\\
1-u^{-\beta} & \textrm{$c=i+1$, $d=1$, $j=0$},\\
1 & \textrm{$c=i+1$, $d=j+1$, $j\geq 1$, $d\geq 2$},\\
0 &\textrm{otherwise}. \end{array}\right.\]
In order for the shadow maps to be taken $\beta$ has to be an
element of $\mathbb{^{\ast}C}^{\lim_{p}}$ for all standard $p$
and if the real shadow map is to be taken then also for $p$ the
real prime. As an initial investigation just let
$\beta\in\mathbb{N}$. For these values of $\beta$ the $p$-adic
shadow map acting on the chain leads to the $p$-adic gamma chain
with value $\beta$ for all finite $p$, since
$\operatorname{sh}_{p}(u^{-\beta})=p^{-\beta}$. Moreover for
these values of $\beta$, $\operatorname{sh}_{\eta}(u^{-\beta})=0$
and a chain is obtained which is the real gamma chain (the unit
shift chain to the right). So the initial nonstandard chain leads
to all the $p$-gamma chains (for all standard primes and the
real prime) by an application of the respective shadow map.

The problem is, of course, the restrictive nature of the values
of $\beta$. Suppose the set of values for $\beta$ is enlarged to
contain $\mathbb{Q}$ then the in general the $p$-adic shadow maps
cannot be taken because $u^{-\beta}$ is irrational. Suppose
instead it is taken to be $\mathbb{^{\ast}N}$ then one obtains a
chain based on a $p$-adic number which is a "handmade" chain for
a $p$-adic number instead of just real numbers considered above.
Even so their use appears limited meaning the "maximal" set for
beta is $\mathbb{N}$. Can a topology be put on chains in order to
look at sequences of chains?

The next option to consider is the transition probabilities as
above but with the enlargement of the state space
\[ ^{\ast}X^{0}_{u}=\coprod_{N\in\mathbb{^{\ast}N}}\text{}
^{\ast}X^{0}_{u(N)}, \quad ^{\ast}X^{0}_{u(N)}=\{(N,j):0\leq
j\leq N\}.\]
My problem with this state space is taking the shadow map. For
all nonstandard $N\in\mathbb{^{\ast}N}\setminus\mathbb{N}$, \[
\mathbb{N}\subset \operatorname{sh}_{p}(^{\ast}X^{0}_{u(N)})
\subset \mathbb{Z}_{p}.\] The standard chain which results from
the $p$-adic shadow map contains $X^{0}_{p}$ but also subsets of
$\mathbb{Z}_{p}$. At the moment I am unsure of how this "fits in"
with this chain.

\end{document}